\documentclass[12pt]{report}
\usepackage{subfiles}
\usepackage{amssymb,amsmath,amsthm,amstext,amscd, scalefnt}
\usepackage[pdftex]{graphicx}
\usepackage{subfig}
\usepackage{url}
\usepackage{float}
\usepackage{tikz}
\usetikzlibrary{calc,positioning,shapes,arrows, positioning,svg.path,plotmarks,arrows}
\usepackage{pdflscape}
\usepackage{multirow}
\usepackage{enumerate}
\usepackage{xstring}
\usepackage{ifthen}
\usepackage[nomessages]{fp}
\usepackage{forest} 
\usepackage{placeins}
\usepackage[nospace,noadjust]{cite}
\usepackage{mathtools}
\usepackage{xr} 
\usepackage{xr-hyper} 
\usepackage[hidelinks=true]{hyperref}
\usepackage[numbers]{natbib} 

 \def\biblio{\bibliographystyle{apa}\bibliography{\main/caref}}  

\usepackage{geometry}
\renewcommand{\baselinestretch}{1.3} 

\DeclarePairedDelimiter{\ceil}{\lceil}{\rceil}
\DeclarePairedDelimiter{\floor}{\lfloor}{\rfloor}

\tikzset{
  solid node/.style={circle,draw,inner sep=1.2,fill=black},
  hollow node/.style={circle,draw,inner sep=1.2},
  left label/.style={above left,midway},
  right label/.style={above right,midway}
}

\def\centerarc[#1](#2)(#3:#4:#5)
    { \draw[#1] ($(#2)+({#5*cos(#3)},{#5*sin(#3)})$) arc (#3:#4:#5); }
    
\newcommand{\tvec}[2]{ 
	\left( \begin{smallmatrix} {#1} \\ {#2} \end{smallmatrix} \right)
}
\newcommand{\tvecal}[2]{ 
	\left( \begin{smallmatrix*}[r] {#1} \\ {#2} \end{smallmatrix*} \right)
}
\newcommand{\CII}[2]{ 
\left(\begin{smallmatrix} #1 \\ #2 \end{smallmatrix} \right)
}
\newcommand{\thvec}[3]{ 
	\left( \begin{smallmatrix} {#1} \\ {#2} \\ {#3} \end{smallmatrix} \right)
}
\newcommand{\thvecal}[3]{ 
	\left( \begin{smallmatrix*}[r] {#1} \\ {#2} \\ {#3} \end{smallmatrix*} \right)
}

\numberwithin{equation}{section}

\newtheorem{thm}{Theorem}[section]
\newtheorem{cor}[thm]{Corollary}
\newtheorem{lem}[thm]{Lemma}

\newtheorem{prop}[thm]{Proposition}
\newtheorem{cnj}[thm]{Conjecture}

\theoremstyle{definition}
\theoremstyle{definition}
\newtheorem{rmkdfn}[thm]{Remark/Definition}
\newtheorem{dfn}[thm]{Definition}
\theoremstyle{definition}
\newtheorem{ex}[thm]{Example}
\theoremstyle{definition}
\newtheorem{rmk}[thm]{Remark}
\theoremstyle{definition}
\newtheorem{ntn}[thm]{Notation}
\theoremstyle{definition}
\newtheorem{alg}[thm]{Algorithm}
\theoremstyle{definition}
\newtheorem{apdxprf}[thm]{Proof}
\theoremstyle{definition}
\newtheorem{qn}[thm]{Question}
\newtheorem{case}{Case}

\newcommand{\A}{\mathcal{A}}

\newcommand{\C}{\mathbb{C}}
\newcommand{\N}{\mathbb{N}}
\newcommand{\K}{\mathbb{K}}
\newcommand{\Z}{\mathbb{Z}}
\newcommand{\F}{\mathcal{F}}

\newcommand{\la}{\langle}
\newcommand{\ra}{\rangle}
\newcommand{\circled}[2][]{%
  \tikz[baseline=(char.base)]{%
    \node[shape = circle, draw, inner sep = 1pt]
    (char) {\phantom{\ifblank{#1}{#2}{#1}}};%
    \node at (char.center) {\makebox[0pt][c]{#2}};}}
\newcommand{\circtx}[1]{\circled[00]{#1}}
\newcommand\pref[1]{(\ref{#1})} 

\newcommand\scalemath[2]{\scalebox{#1}{\mbox{\ensuremath{\displaystyle #2}}}}

\newcommand{\wt}[1]{\widetilde{#1}}
\newcommand{\und}[1]{\underline{#1}}
\newcommand{\mmut}[2]{#1 _{[{#2}]} } 
\newcommand{\pth}[1]{[\und{#1}]}
\newcommand{\pths}[2]{[\und{#1}, \und{#2}]}
\newcommand{\var}{\operatorname{var}}
\newcommand{\cl}{\operatorname{cl}}
\newcommand{\sd}{\operatorname{sd}}
\newcommand{\dsd}{\deg.\operatorname{sd}}
\newcommand{\dv}{\text{\bf{d}}} 
\newcommand{\dcl}{\dv.\cl} 
\newcommand{\dslof}[2]{\dcl \left( {#2} \right) |_#1} 
\newcommand{\dslp}[2]{\dcl {#1} |_#2} 
\newcommand{\dslpm}[3]{\dcl_{#1} {#2} |_#3} 
\newcommand{\sfr}[1]{\displaystyle{#1^\circ}} 

\newcommand{\esseq}{\equiv}
\newcommand{\cluster}[1]{(#1)}
\newcommand{\mrep}[1]{ m({#1}) } 

\newcommand{\prpth}[3]{#2 *_{#3} #1} 
\newcommand{\Prpth}[3]{[\prpth{#1}{#2}{#3}]} 
\newcommand{\cfg}[1]{\kappa({#1})}
\newcommand{\ssmat}[3]{  
	\StrLen{#1}[\alength]
	\FPeval{\ralength}{clip(\alength-1)}
	\StrLen{#2}[\blength] 
	\FPeval{\rblength}{clip(\blength-1)}
	\StrLen{#3}[\clength]
	\FPeval{\rclength}{clip(\clength-1)}
	\left( \begin{smallmatrix}
	        0 & \IfBeginWith{#1}{-}{ \StrRight{#1}{\ralength} }{-#1} &  \IfBeginWith{#3}{-}{ \StrRight{#3}{\rclength} }{-#3}\\
	        {#1}& 0 &\IfBeginWith{#2}{-}{ \StrRight{#2}{\rblength} }{-#2}\\
	        {#3}&{#2} & 0 \end{smallmatrix} \right)
}

\newenvironment{xsmallmatrix}[1] 
  {\renewcommand\thickspace{\kern#1}\smallmatrix}
  {\endsmallmatrix}
        
\def\Arrow#1{\raisebox{-.5\height}{\scalebox{1.6}{  $\xrightarrow{\mbox{ \tiny{#1} } }$  }}}
\def\Image#1{\raisebox{-.5\height}{{#1}}}

\newcommand{\striL}[2]{ 
\Image{
\begin{tikzpicture}[node distance=1.6cm, auto, font=\footnotesize]                                          
  \node (A1) []{#1};                                                                                                                                                                                              
  \node (A3) [below right of = A1]{$*$};   
  \node (A2) [above right of  = A3]{#2};        
  \draw[-latex] (A1) to node {$2$} (A2);      
  \draw[dotted, -latex] (A3) to node [] {$W$} (A1); 
  \draw[dotted, latex-] (A3) to node [swap] {$w$} (A2);   
\end{tikzpicture}}
}

\newcommand{\striR}[2]{
\Image{
\begin{tikzpicture}[node distance=1.6cm, auto, font=\footnotesize]                                          
  \node (A1) []{#1};                                                                                                                                                                                              
  \node (A3) [below right of = A1]{$*$};   
  \node (A2) [above right of  = A3]{#2};        
  \draw[-latex] (A1) to node {2} (A2);      
  \draw[dotted, -latex] (A3) to node [] {$w$} (A1); 
  \draw[dotted, latex-] (A3) to node [swap] {$W$} (A2);   
\end{tikzpicture}}
}

\graphicspath{{images}{../images/}}

\begin{document}
\def\biblio{}

\pagenumbering{roman} 

\author{Thomas Booker-Price}
\title{Applications of Graded Methods to Cluster Variables in Arbitrary Types}
\date{2017}


\begin{titlepage}
	\centering
	{\scshape\Large \phantom{.} \par}
	\vspace{1.5cm}
{\huge {Applications of Graded Methods to Cluster Variables in Arbitrary Types}\par}
	\vspace{1.5cm}
     {\Large {Thomas Booker-Price}\par}	
	\vspace{7cm}
	{\Large	A thesis submitted for the degree of \\ Doctor of Philosophy \par}	
	\vspace{1cm}
	{\Large Department of Mathematics and Statistics \\ Lancaster University \\ June 2017 \par}
	\vspace{1cm}	
	\vfill
\end{titlepage}

{
\begin{center}
{\large {Applications of Graded Methods to Cluster Variables in Arbitrary Types}\par}
	\vspace{0.3cm}
     { {Thomas Booker-Price}\par}	
	\vspace{0.3cm}     
	{A thesis submitted for the degree of Doctor of Philosophy \par}	
	{ June 2017 \par}
	\vspace{1cm}	

	{ \bf Abstract}
\end{center}
This thesis is concerned with studying the properties of gradings on several examples of cluster algebras, primarily of infinite type. We start by considering two classes of finite type cluster algebras: those of type $B_n$ and $C_n$. We give the number of cluster variables of each occurring degree and verify that the grading is balanced. These results complete a classification in \cite{G15} for coefficient-free finite type cluster algebras.

We then consider gradings on cluster algebras generated by $3 \times 3$ skew-symmetric matrices. We show that the mutation-cyclic matrices give rise to gradings in which all occurring degrees are positive and have only finitely many associated cluster variables (excepting one particular case). For the mutation-acyclic matrices, we prove that all occurring degrees have infinitely many variables and give a direct proof that the gradings are balanced.

We provide a condition for a graded cluster algebra generated by a quiver to have infinitely many degrees, based on the presence of a subquiver in its mutation class. We use this to study the gradings on cluster algebras that are (quantum) coordinate rings of matrices and Grassmannians and show that they contain cluster variables of all degrees in $\N$.

Next we consider the finite list (given in \cite{FSTFinite}) of mutation-finite quivers that do not correspond to triangulations of marked surfaces. We show that $\A(X_7)$ has a grading in which there are only two degrees, with infinitely many cluster variables in both. We also show that the gradings arising from $\wt{E}_6$, $\wt{E}_7$ and $\wt{E}_8$ have infinitely many variables in certain degrees.

Finally, we study gradings arising from triangulations of marked bordered $2$-dimensional surfaces (see \cite{FST}). We adapt a definition from \cite{Muller} to define the space of \emph{valuation functions} on such a surface and prove combinatorially that this space is isomorphic to the space of gradings on the associated cluster algebra. We illustrate this theory by applying it to a family of examples, namely, the annulus with $n+m$ marked points. We show that the standard grading is of mixed type, with finitely many variables in some degrees and infinitely many in the others. We also give an alternative grading in which all degrees have infinitely many cluster variables.

}
\newpage

\section*{Acknowledgements}

Firstly, I would like to thank my supervisor, Jan Grabowski, for his generous support and guidance throughout the course of my studies. I am extremely grateful for all the time he gave up on the many useful discussions we have had regarding my research, as well as for his knack for finding pertinent results in the literature that have helped me make progress. I thank him also for his constant willingness to help provide interesting research topics, which has made my PhD experience all the more enjoyable.

I also thank my examiners, Anna Felikson and Mark MacDonald, for making many useful suggestions that have improved my thesis.

There are many people I am grateful to have met during my time at Lancaster. I thank the other PhD students I have met, many of whom I have become good friends with, for being such a down-to-earth group who have made the social life here a great experience. I am also grateful for the non-mathematical friends that I have made here (and for their reminders that it is inappropriate for a group of PhD students to repeatedly make lame mathematical jokes in their presence). Additionally, I thank the staff in the Maths and Stats department, who help to make it an exceptionally friendly and positive place.

I owe an enormous debt of gratitude to my parents, Moshe and Lesley, who have always believed in me and supported me in everything I have done. I also thank the rest of my family, including my sister Maddy, for all their love and support.

Last but not least, I owe a special thanks to Faye for her love, and all the encouragement she has given me (especially towards the end of my studies). I thank her also for providing an open problem far deeper than any I could encounter during my PhD: how to consistently make a cup of tea that isn't ``vile". I am grateful to her for all the chocolate we shared, the places we visited together and the wonderful times we have had---and for always putting up with my rambling during all of the above. 

\newpage

\section*{Declaration}

This thesis is my own work, and has not been submitted in substantially the same form for the award of a higher degree elsewhere.

\tableofcontents

\biblio

\chapter{Introduction}
\pagenumbering{arabic}

Cluster algebras are commutative, unital subalgebras of $\mathcal{F}:=\mathbb{Q}(X_1, \dots, X_n)$ generated by a certain iterative process which endows them with a very particular combinatorial structure.
They were introduced by Fomin and Zelevinsky in the series of papers \cite{FZCAI,FZCAII,FZBCAIII,FZCAIV} (the third with co-author Berenstein) in an attempt to capture certain combinatorial patterns that had been observed in a number of algebras associated to some particularly nice classes of geometric objects. The authors of these founding papers realised that, though the combinatorial patterns appear complex, they arise from relatively simple rules, iterated over and over again. The original goal of cluster algebras was to provide an algebraic framework for the study of dual canonical bases, and the related notion of total positivity, in semisimple Lie groups and their quantum analogues. Since their introduction, however, they (and their quantum analogues) have appeared in many areas of mathematics, such as algebraic and symplectic geometry, noncommutative algebra, mathematical physics and others.

The iterative process that produces the generators of a cluster algebra starts with only a small finite subset $\{x_1,\dots,x_n\} \subset \mathcal{F}$ of all of the generators we will eventually obtain. These generators are known as \emph{cluster variables} and the finite subset we start with is called the \emph{initial cluster}. From this initial cluster, more clusters are produced, each replacing one old variable with a new one, in a process called \emph{mutation}. How the mutation of the initial cluster is carried out is controlled by an $n \times n$ skew-symmetrisable matrix called the \emph{initial exchange matrix}, whose entries in fact determine the entire cluster algebra. When the initial exchange matrix is used to mutate the initial cluster, it too mutates, so that we obtain a new exchange matrix along with a new cluster (together called a \emph{seed}) which in turn can be mutated again. This process of seed mutation is iterated repeatedly, and the union of all clusters thus obtained is the set of generators for the algebra. A consequence of this method of construction is that we often end up with many more generators than we require. Indeed, there are finitely generated algebras (two important examples being 
coordinate rings of matrices and Grassmannians)
that have cluster algebra structures with infinitely many generators. On the other hand, the relations between these generators are always of a very restricted and relatively simple form.

While the exchange matrix dictates exactly how mutation is carried out, the new cluster variable $z_k'$ that is obtained from a given cluster $\{z_1,\dots,z_n\}$ via a single \emph{mutation in direction $k$} always satisfies an \emph{exchange relation} of the form
\[ z_k z_k'= M_1 + M_2,\]
where $M_1$ and $M_2$ are monomials in the variables $\{z_1,\dots,z_n\} -\{z_k\}$ with no common divisors. By definition, any two cluster variables (from any two clusters) can be obtained from each other via some sequence of mutations. Thus, any cluster variable may be expressed as a rational function in the variables of any given cluster. One of the fundamental behaviours that a cluster algebra exhibits is the \emph{Laurent phenomenon} (\cite[Theorem 3.1]{FZCAI}), which says that any cluster variable is in fact a Laurent polynomial in the variables of any given cluster. This property is remarkable, as mutating a cluster to obtain a new cluster variable requires dividing by a Laurent polynomial that in general will have many terms in the numerator. The fact that such a cancellation must always occur is therefore highly surprising. Another property, which arises from the form of exchange relations above, is that a cluster variable can always be represented as a \emph{subtraction-free} rational function in the variables of any given cluster. This reflects the connection between cluster algebras and the theory of total positivity.

Cluster algebras are either of \emph{finite type} or \emph{infinite type}, corresponding to whether the set of cluster variables is finite or infinite. That is, a cluster algebra is of finite type if the iterative process of mutation exhausts all possible clusters in a finite number of steps, and of infinite type otherwise. Cluster algebras of finite type are well-understood and have been fully classified. This classification is an incarnation of the well-known Cartan-Killing classification of  simple Lie algebras:
\begin{thm}[{\cite[Theorem 1.8]{FZCAII}}]
Let $\A$ be a cluster algebra. The following are equivalent.
\begin{enumerate}[(i)]
\item $\A$ is of finite type.
\item For every seed $(\und{x}, B)$ in $\A$, the entries of the matrix $B=(b_{ij})$ satisfy ${|b_{ij}b_{ji}| \leq 3}$ for all $i,j$.
\item $B$ is mutation equivalent to a matrix $B'$ whose Cartan counterpart is a Cartan matrix of finite type.
\end{enumerate}
\end{thm}

On the other hand, while cluster algebras of infinite type make up the vast majority, much less is understood about them. One of the reasons for this has been the lack of good tools to control the infiniteness present in them. This thesis is concerned with studying one of the first of these tools, namely, gradings on cluster algebras.

The introduction of \cite{G15} contains a brief historical overview of the use of cluster algebra gradings in the literature.
The first notion of a grading for cluster algebras was that of a $\Z^n$-grading introduced by Fomin and Zelevinsky (\cite[Proposition 6.1]{FZCAIV}), where $n$ is the \emph{rank} of the cluster algebra (the number of cluster variables in each cluster). The notion of a graded quantum seed was defined by Berenstein and Zelevinsky \cite[Definition 6.5]{BZ}, which gives rise to a module grading but not an algebra grading. Another notion of grading is due to Gekhtman et al.~(\cite[Section 5.2]{GSV}) via the language of toric actions. The definition given by Gekhtman et al.~is equivalent to the one we will eventually settle on, but more appropriate for the geometric perspective of the setting in which it arises.

The definition of grading in the form we are interested in was introduced by Grabowski and Launois in \cite{GL}, where the authors use it to prove that quantum Grassmannians $\K_q[Gr(k,n)]$ admit a quantum cluster algebra structure. This definition, which is of a $\Z$-grading, is generalised to $\Z^k$-gradings in \cite{G15}, which also studies and fully classifies gradings for finite type cluster algebras without coefficients (the classification for types $B_n$ and $C_n$ due to work that appears in Chapter \ref{chp:finite_type} of this thesis).

In general, when attempting to endow an algebra with a grading structure, we assign degrees to the generators and check that the relations are homogeneous, and for cluster algebra gradings we essentially do the same. However, instead of assigning degrees to all our generators at once, we initially only assign degrees to the generators in our initial cluster, according to a certain grading condition. The degrees of all the other cluster variables are then determined recursively from those of the initial cluster, and it turns out that the exchange relations mean all cluster variables will be homogeneous.

Prior to this thesis, almost nothing was known about the properties of gradings on infinite type cluster algebras, and almost any question one could ask was open. The goal of the thesis is to study several classes of examples of graded cluster algebras (including the matrix and Grassmannian cases mentioned above) and discover what structure the gradings have. In particular, we are interested in describing the cluster variables in terms of their degrees. More precisely, given a $Z^k$-graded cluster algebra, we aim to answer the following questions:
\begin{itemize}
\item Is the grading balanced, that is, is there a bijection between the variables of degree $d$ and those of degree $-d$? (We only ask this when applicable; sometimes we end up with $\N$-gradings.)
\item Do there exist cluster variables of infinitely many different degrees, or only finitely many?
\item Of the degrees that occur, how are the variables distributed among them; are there
\begin{enumerate}
\item finitely many variables of each occurring degree,
\item infinitely many variables of each occurring degree, or,
\item a mixture, with some degrees associated with infinitely many variables and others only finitely many?
\end{enumerate}
\end{itemize}

These questions are posed in a loosely increasing order of difficulty, though the real determinant of difficulty for any such question is the rank of the cluster algebra we are asking it about. For rank two cluster algebras, there are no nontrivial gradings; increasing the rank by just one, we immediately find gradings with complex structures and a rich variety of behaviours. In most cases, we are able to answer the questions above for rank three cluster algebras generated by skew-symmetric matrices. Even in this low-rank setting, though, some questions require significant effort to answer (and, for one particular case, some remain open). As the rank continues to grow, the complexity of cluster algebras and their graded structures increases extremely rapidly in a way that is difficult to control. For higher rank cases, therefore, our ambitions must be much more modest, even when the examples we look at have a great deal of additional structure.  An exception to this occurs in cluster algebras associated to triangulated marked surfaces (see \cite{FST} for background). Here, the associated structure allows us to describe cluster variables geometrically, and this proves such a powerful tool when applied to gradings that results which would otherwise be extremely difficult to obtain can essentially be read off from the associated geometric combinatorial objects.

It was anticipated that heavy use of representation theoretic techniques would be required in order to make progress on the questions above, but it turns out that a surprising amount can be said using direct, combinatorial methods (though some of the combinatorial facts we use do themselves rely on underlying representation theory). In this thesis, we have persisted with this elementary approach and attempted to see how far it can be pushed in obtaining useful results. 

An additional tool we use at times is that of computer-aided calculation. While only a few results strictly rely on these calculations (all contained in Chapter \ref{chp:finite_mutation_type}), in practice, several other results and conjectures (mainly in Chapter \ref{chp:matrix_grassmannian}) would not have been possible to obtain otherwise. The code developed while carrying out this research, along with accompanying documentation, is available along with the electronic version of the thesis.

\section{Thesis overview}
\interfootnotelinepenalty=10000
This thesis is organised as follows. In Chapter \ref{chp:notation_preliminaries} we summarise the basic theory of graded cluster algebras and their associated structures, and collect notation and definitions we will need throughout. In particular, we say what it means for two exchange matrices to be equivalent (which we term \emph{essentially equivalent}), from the perspective of mutation directions.

In Chapter \ref{chp:finite_type} we study the gradings associated to finite type cluster algebras generated by matrices whose Cartan counterparts are of types $B_n$ and $C_n$. We do so using a formula (\cite[Corollary 4.2]{G15}) that arises from a known bijection (\cite[Theorem 1.9]{FZCAI}) between the cluster variables in a cluster algebra of finite type and the almost positive roots in the root system of the corresponding type in the Cartan-Killing classification. The results we obtain complete the classification of finite type cluster algebras without coefficients, which is presented in \cite{G15}.

In Chapter \ref{chp:3_by_3} we turn to the first nontrivial infinite type case: gradings on rank $3$ skew-symmetric matrices. Our goal is to classify as much of this entire family as we can. We start by noting that, up to mutation and essential equivalence, every matrix is of the form $\ssmat{-a}{-b}{c}$ for some $a,b,c \in \N_0$. Similarly, we note that up to integer scaling there is only one grading vector for such a matrix, and that the entries of the grading vector are contained in the matrix. We give a conjectured classification of grading behaviours in Table \ref{tab:classification3x3} and turn to working through the table in the remainder of the chapter, starting by showing that the division between finitely many and infinitely many degrees is as conjectured. Next we consider the only mixed type case, $\ssmat{-2}{-1}{1}$, and prove that its grading behaves as such. In the next section we give an algorithm (similar to one in \cite{ABBS}) which determines whether a $3 \times 3$ skew-symmetric matrix is mutation-cyclic\footnote{A matrix is mutation-cyclic if all the matrices in its mutation class correspond to quivers that are \emph{cyclic} (i.e.~consist of an oriented cycle). In general, a skew-symmetric matrix is \emph{mutation-acyclic} if there is a matrix in its mutation class corresponding to an \emph{acyclic} quiver (one without any oriented cycles).}
or not and gives a minimal representative of the matrix in its mutation class. We prove that, aside from one particular case, mutation-cyclic matrices give rise to gradings with finitely many cluster variables in each degree. This is achieved by showing that the growth of degrees under mutation is too fast to allow for more than finitely many variables per degree. The penultimate section of this chapter considers mutation-acyclic
  matrices that are mutation-infinite, which require the most work to classify. We prove that the corresponding graded cluster algebras all have infinitely many variables in each occurring degree. A crucial result we use to do this is \cite[Corollary 3]{CalderoKeller}, which concerns connectedness of subgraphs of an important object associated with a cluster algebra: its \emph{exchange graph}. Another result which is helpful here is the classification of exchange graphs for $3 \times 3$ acyclic matrices, whose structures are identified in \cite[Theorem 10.3]{War}. Finally, in the last section, we consider the so-called \emph{singular cyclic} case: the mutation-cyclic matrices of the form $\ssmat{-a}{-a}{2}$, $a \geq 3$, that do not have the same fast growth behaviour as other mutation-cyclic matrices. We show that infinitely many of the occurring degrees contain infinitely many variables, and conjecture that in fact all degrees contain infinitely many variables. We briefly indicate why proving this is difficult and suggest some possible approaches.

Chapter \ref{chp:growing_subquivers} collects some results that provide a condition for identifying when a graded cluster algebra has infinitely many degrees in general, some of which we will apply in Chapter \ref{chp:matrix_grassmannian}. These involve inferring information about a \emph{degree quiver} (a diagrammatic way of representing a skew-symmetric matrix along with a compatible grading vector) from information about its subquivers. As such, they are some of our first results of this kind, and this line of research remains largely undeveloped (although generalising these results is likely an approachable goal). In the final section we discuss this in more detail and introduce a way to sum degree quivers together to give another degree quiver. This gives one potential way of viewing gradings on cluster algebras as being made up of smaller gradings. We define the notion of an \emph{irreducible} degree quiver and give an example of a restricted setting in which using sums of degree quivers is useful for inferring results from the summands.

As mentioned, the coordinate rings $\mathcal{O}(M(k,l))$ and $\mathcal{O}(Gr(k,k+l))$ of matrices and Grassmannians, along with their quantum counterparts, are important examples of well-understood algebras that turn out to have cluster algebra structures. In Chapter \ref{chp:matrix_grassmannian}, we consider gradings on these cluster algebras (the specific grading chosen being the one defined in \cite{GL}). These are the first cluster algebras we encounter that have \emph{frozen vertices} (or \emph{coefficients}), vertices of the exchange quiver (equivalently, variables present in every cluster) which we are never allowed to mutate. These frozen vertices are needed in order to allow the quantum versions of the graded cluster algebras to have nontrivial gradings, as quantum cluster algebras without coefficients do not admit gradings. This is also the first high-rank setting we encounter. 
Our goals here are to decide whether the cluster algebras, which are $\N$-graded, have variables of infinitely many different degrees, and, if so, whether this includes every positive degree. Both of these properties are suggested by the properties of $\mathcal{O}(M(k,l))$ and $\mathcal{O}(Gr(k,k+l))$ considered as graded algebras in the traditional sense. We show that this is also the case for the cluster algebra gradings, using some results from Chapter \ref{chp:growing_subquivers}. Mostly, this reduces to finding appropriate subquivers present in the degree quiver mutation class. While this is simple in principle, the high rank means such subquivers are very sparsely distributed amongst the mutation class, and this is one of the situations in which we have benefited from computer-aided calculation.

For $n \geq 3$, it is known (\cite[Theorem 6.1]{FSTFinite}) that $n \times n$ skew-symmetric matrices with a finite mutation class are either adjacency matrices of triangulations of marked bordered $2$-dimensional surfaces, or one of a finite list of eleven exceptional matrices. In Chapter \ref{chp:finite_mutation_type} we investigate what can be said about gradings on cluster algebras generated by some of these exceptional matrices. We use computer-aided calculation to find the size of the degree quiver class and the set of occurring degrees corresponding to several of the matrices. For one, $X_7$, we prove the grading gives rise to infinitely many variables in each degree (which would usually be very difficult for a cluster algebra of such high rank, but is possible here due to the large amount of symmetry in the quiver for $X_7$). For $\wt{E}_6$, $\wt{E}_7$ and $\wt{E}_8$, we show that certain of the occurring degrees contain infinitely many cluster variables (and we conjecture that the other degrees have only one or a finite number of variables).

Having considered the mutation-finite matrices that are not adjacency matrices of triangulated surfaces, we turn our attention in Chapter \ref{chp:surface_type} to the class of matrices that are. The process by which a cluster algebra arises from a bordered marked surface $(\bf{S},\bf{M})$ is described in \cite{FST}. The authors define the \emph{arc complex} of $(\bf{S},\bf{M})$ as the simplicial complex whose ground set is the set of all arcs between marked points in $\bf{M}$ and whose maximal simplices are the \emph{ideal triangulations} (collections of non-intersecting arcs that fully triangulate the surface). They show that the arc complex of a surface is isomorphic to the \emph{cluster complex} of the associated cluster algebra, and that the dual graph of the arc complex is isomorphic to the exchange graph. In particular, this means the cluster variables are in bijection with the arcs and seeds are in bijection with the ideal triangulations. In Section 3.5 (see also Section 7.5) of {\cite{Muller}}, Muller gives a notion of an \emph{endpoint grading} which assigns degrees to arcs based on values of the marked points at their endpoints. This notion does not match our definition of grading, but we make an adjustment to it to obtain one that is compatible with gradings in our sense. We define the space of \emph{valuation functions}---assignments of values to marked points according to certain compatibility conditions---for a marked surface and prove (using only the combinatorics of the surface theory) that it is isomorphic to the space of gradings on the associated cluster algebra. This allows us to translate grading questions into questions about degrees of arcs. As we noted above, this lets us capture a great deal of information about the grading structure with much less effort than would be expected. Lastly, we indicate some promising possibilities for further research involving this theory.

\biblio

\chapter{Notation and preliminaries} \label{chp:notation_preliminaries}

In this chapter we summarise basic definitions, results and notation which we will need throughout. We fix $n, m \in \N$ throughout the chapter.

\section{Definition of a cluster algebra}

\begin{dfn}
Let $X$ be a set of cardinality $n$ and let $l:\{1, \dots , n\} \to X$ be a bijection. We will call the pair $(X,l)$ a \emph{labelled set}, where $l$ is the \emph{labelling} of $X$. For $i \in \{1, \dots , n\}$, we will call $l(i)$ the element of $X$ \emph{labelled $i$}. We will often refer to $X$ by itself as a labelled set.
\label{dfn:labelled_set}
\end{dfn}
We may write any labelled set $(X,l)$ in the form $X=\{x_{l(1)},\dots,x_{l(n)}\}$.
This way we can think of a labelled set as a tuple when needed, and it will usually be convenient to do so as the labelled sets we will be working with will be \emph{clusters}, where the order of the elements is important for \emph{seed mutation}. (These terms will be defined below.) Hereafter, if a labelled set $X$ is written down without reference to $l$, it is to be assumed that it has been written in this form. If it is written as a tuple, this is understood as shorthand for the above.

\begin{ntn}
Let $\sigma \in S_n$ be a permutation of $\{1,\dots,n\}$ and let $(X,l)$ be a labelled set. We will write $\sigma(X)$ to mean the labelled set $( X, l \circ \sigma )$.
\end{ntn}

A cluster algebra has several associated structures which we now define before giving the definition of a cluster algebra itself. Most of the definitions in the remainder of this section (excluding those associated with \emph{mutation paths} and \emph{essential equivalence}) follow \cite{GSV}. An original source can be found in \cite{FZCAI,FZCAII,FZBCAIII,FZCAIV}, which are the seminal papers introducing cluster algebras.

\begin{dfn}
The \emph{coefficient group} $\mathcal{P}$ 
 is the free multiplicative abelian group of rank $m$ with generators $y_1,\dots,y_m$.

Let $\K$ denote the field of fractions of the group ring $\Z\mathcal{P}$. The \emph{ambient field} $\F$ is the field of rational functions in the $n$ indeterminates $X_1 \dots, X_n$ with coefficients in $\K$, that is, $\K(X_1,\dots,X_n)$. 
\end{dfn}
 
\begin{dfn}
Given an $(n+m) \times n$ matrix $\wt{B}$, its \emph{principal part}, $B$, is the $n \times n$ submatrix of $\wt{B}$ obtained by deleting rows $n+1$ to $n+m$. If $B$ is skew-symmetrisable (i.e.\ can be multiplied by an invertible diagonal matrix to obtain a skew-symmetric matrix) then we say that $\wt{B}$ also has this property. 
\end{dfn}

\begin{ntn}
We retain the notation $b_{ij}$ (rather than using $\wt{b}_{ij}$) to refer to the $(i,j)$ entry of $\wt{B}$. We will also use commas between indices in some situations, as in $b_{i,j}$, to aid readability.
\end{ntn}

\begin{dfn}[Cluster, extended cluster, seed]
The notions of \emph{cluster} and \emph{seed} are defined iteratively. Though they require the notion of \emph{mutation}, we will briefly suspend this definition as doing so allows for a neater presentation.

An \emph{initial cluster} in $\F$ is a labelled set $\und{x}_0$ of cardinality $n$ that is a transcendence base of $\F$ over $\K$. 
(Recall that this means $\und{x}_0$ is an algebraically independent set such that the field extension $\F / \K(\und{x}_0)$ is algebraic, that is, such that the elements of $\F$ are all algebraic over $\K(\und{x}_0)$.)
Given a fixed initial cluster $\und{x}_0$, a \emph{cluster} (in the corresponding cluster algebra) is either $\und{x}_0$ itself or a labelled set $\und{x}$ that satisfies the same requirements as an initial cluster and that can be obtained from $\und{x}_0$ by \emph{mutation}. 
The elements of a cluster are called \emph{cluster variables} and the elements of the initial cluster are called the \emph{initial cluster variables}. 

Given a cluster $\und{x}=\{x_1,\dots,x_n\}$, define the \emph{extended cluster} $\wt{\und{x}}$ to be the labelled set $\{x_1,\dots,x_n,y_1,\dots,y_m\}$, where the $y_i$ are the generators of $\Z\mathcal{P}$ defined above. (We extend the labelling $l$ on $\und{x}$ to $\wt{\und{x}}$ by setting $l(n+i)=y_i$, so we may refer to $y_i$ as $x_{l(n+i)}$.) The variables $y_{1},\dots, y_{m}$ are called the \emph{coefficient variables}. They are described as \emph{stable} or \emph{frozen} variables, whereas the cluster variables $x_1,\dots,x_n$ are called \emph{mutable}.

An \emph{initial seed} in $\F$ is a pair $\Sigma_0 = (\wt{\und{x}_0},\wt{B}_0)$, where $\wt{\und{x}_0}$ is an extended cluster and $\wt{B}_0$ is an $(n+m) \times n$ skew-symmetrisable integer matrix. Given a fixed initial seed $\Sigma_0$, a \emph{seed} (in the corresponding cluster algebra) is either $\Sigma_0$ itself or a pair $\Sigma = (\wt{\und{x}},\wt{B})$ that satisfies the same requirements as an initial seed and that can be obtained from $\Sigma_0$ by \emph{seed mutation}. Given a seed $\Sigma = (\wt{\und{x}},\wt{B})$, the matrix $\wt{B}$ is called an \emph{extended exchange matrix} and its principal part $B$ is called an \emph{exchange matrix}. The matrix in the initial seed is called the \emph{initial (extended) exchange matrix}.
\end{dfn}

All extended clusters $\und{\wt{x}}$ in a cluster algebra share the same stable variables and so we will usually replace the extended cluster $\und{\wt{x}}$ by the cluster $\und{x}$ in the associated notation. Similarly, when the context is clear, we may simply refer to an extended cluster as a cluster. We keep in mind when needed (such as when carrying out \emph{mutation}, as we now define) that the stable variables are present.

\begin{dfn}[Cluster mutation]
Let $\Sigma = ({\und{x}},\wt{B})$ be a seed and let $l$ be the labelling of the cluster $\und{x}=\{x_{l(1)},\dots,x_{l(n)}\}$. The \emph{mutation of $\und{x}$} in direction $k \in \{1,\dots,n\}$ is the labelled set 
\[\mu_{\Sigma,k}(\und{x}) = \big(\und{x} \backslash \{ x_{l(k)}\} \cup \{\mu_{\Sigma,k}(x_{l(k)})\} ,  l'\big),\]
where $\mu_{\Sigma,k}(x_{l(k)})$ is defined by the \emph{exchange relation}
\begin{equation}
 \mu_{\Sigma,k}(x_{l(k)}) x_{l(k)} = 
 \prod_{\substack{1 \leq i \leq n+m \\b_{i,k}>0}} x_{l(i)}^{b_{i,k}} + \prod_{\substack{1 \leq i \leq n+m \\ b_{i,k}<0}} x_{l(i)}^{-b_{i,k}},
 \label{eqn:cluster_mutation}
\end{equation}
with the empty product being taken to be equal to 1, and where
\[
 l'(i) = \begin{cases}
        \mu_{\Sigma,k}(x_{l(k)}) & \text{if } i=k, \\
        l(i)  & \text{otherwise}.
        \end{cases}
\]
\end{dfn}

\begin{dfn}
The indices $(n+1),\dots (n+m)$ in Equation \ref{eqn:cluster_mutation} correspond to the stable variables of $\wt{\und{x}}$. The corresponding parts of the monomials formed from these are called \emph{coefficients}. More precisely, define
\begin{align}
 p_k^+ &= \prod_{\substack{1 \leq i \leq m \\b_{(n+i),k}>0}} x_{l(n+i)}^{b_{(n+i),k}} = \prod_{\substack{1 \leq i \leq m \\b_{(n+i),k}>0}} y_{i}^{b_{(n+i),k}} \\
\intertext{and} p_k^- &= \prod_{\substack{1 \leq i \leq m \\b_{(n+i),k}<0}} x_{l(n+i)}^{-b_{(n+i),k}} = \prod_{\substack{1 \leq i \leq m \\b_{n+i,k}<0}} y_{i}^{-b_{(n+i),k}}
 \label{eqn:coefficients}
\end{align}
and call these respectively the \emph{positive and negative coefficient} associated to the exchange relation \pref{eqn:cluster_mutation}.
\end{dfn}

\begin{dfn}[Matrix mutation]
\label{dnf:matrix_mutation}
Given an $(n+m) \times m$ skew-symmetrisable integer matrix $\wt{B}$, we define the \emph{mutation of $\wt{B}$} in  direction $k \in \{1,\dots,n\}$ to be the ${(n+m) \times m}$ matrix $\wt{B'} = \mu_k(\wt{B})$ given by
\begin{equation}
 b'_{ij}= 	\begin{cases}
 -b_{ij} &\text{if } i=k \text{ or } j=k,\\
 b_{ij} +  [b_{ik}]_+b_{kj} + b_{ik}[-b_{kj}]_+ &\text{otherwise,}
 \end{cases}
\end{equation}
where we define $[-]_+ $ by
\begin{equation}
[x]_+=\begin{cases}
    x &\text{if } x \geq 0,\\
    0 & \text{otherwise}.
  \end{cases}
\end{equation}  
\end{dfn}

It is fundamental to the theory of cluster algebras that $\mu_k(B)$ is again skew-symmetrisable and that  $\mu_k (\mu_k(B))=B$. This is straightforward to show. It is also easy to see that $\mu_k(-B)=-\mu_k(B)$.

\begin{dfn}[Seed mutation]
Given a seed $\Sigma = ({\und{x}},\wt{B})$, the \emph{mutation of $\Sigma$} in direction $k \in \{1,\dots,n\}$ is the seed $\Sigma' = \mu_k(\Sigma)=(\mu_{\Sigma,k}(\und{x}),\mu_k(\wt{B}))$.
\end{dfn}

It is straightforward to show that $\mu_k(\mu_k(\Sigma)) =\Sigma$.

We will often wish to make several mutations in a row, so we make the following definitions. 

\begin{dfn}
A \emph{mutation path} is a list $\pth{p}=[p_t,\dots,p_1]$ with each $p_i \in \{1,\dots,n\}$, written from right to left. Its \emph{length}, $l(\pth{p})$, is $t$. A \emph{subpath} of $\pth{p}$ is a mutation path $[\und{p}']$ of the form $[p_j,\dots, p_i]$ for some $t\geq j \geq i \geq 1$, and a \emph{rooted subpath} of $\pth{p}$ is a subpath of $\pth{p}$ of the form $[p_j,\dots, p_1]$. The \emph{reverse path} of $\pth{p}$ is the path $[\und{p}^{-1}]=[p_1,\dots,p_t]$. The \emph{empty path} is the path $[\,]$ consisting of no mutation directions.
\end{dfn}

\begin{rmk}
Our notation for mutation paths is not the standard notation found in the literature, which uses $\mu_{k_n} \circ \cdots \circ \mu_{k_1}$ for $[k_n, \dots, k_1]$. We will need to write down mutation paths frequently and some will be quite long. We use our notation to reduce clutter and improve readability.
\end{rmk}

\begin{ntn}
The concatenation of two paths $\pth{p}$, $\pth{q}$ will be denoted by $[\und{p},\und{q}]$. We will use the notation $(p_t,\dots,p_1)^r$ to mean that the list $p_t,\dots,p_1$ is to be repeated $r$ times. We use $(p_t,\dots,p_1)_r$ to denote the first $r$ terms of the infinitely repeated list $(p_t,\dots,p_1)^{\infty}$. 
For example, the path $[(1,2,3)^2,(9,8,7)_4]$ is equal to $[1,2,3,1,2,3,7,9,8,7]$.
\end{ntn}

\begin{ntn}
\label{ntn:sd_cl_B_p}
Let $\wt{B}$ be an $(n+m) \times n$ matrix, $[\und{p}]=[p_t,\dots,p_1]$ a mutation path and ${\und{x}}=\{x_1, \dots, x_{n+m}\}$ a cluster. We denote by $\mmut{\wt{B}}{\und{p}}$ the matrix obtained by applying $\pth{p}$ to $\wt{B}$. (That is, we first mutate $\wt{B}$ in direction $p_1$, then the resulting matrix in direction $p_2$, etc.) We denote by $\var_{\wt{B}}[\und{p}]({\und{x}})$ the newest cluster variable obtained after applying $[\und{p}]$ to the seed $({\und{x}}, \wt{B} )$.
If $\wt{B}$ is clear from the context, we may omit it. 
Similarly, we denote, the seed and cluster obtained after applying $\pth{p}$ by $\sd_{\wt{B}} \pth{p} ({\und{x}})$ and $\cl_{\wt{B}}\pth{p}({\und{x}})$, respectively, though if ${\und{x}}$ is the initial cluster we will often omit it from the notation. Alternatively we may use $\pth{p} ({\und{x}},\wt{B}) $ for $\sd_{\wt{B}} \pth{p} ({\und{x}})$.
\end{ntn}

\begin{dfn}
Two clusters, matrices or seeds are \emph{mutation equivalent} if there exists a mutation path $\pth{p}$ under which the latter respective object is obtained from the former. The \emph{mutation class} of a cluster, matrix or seed is the set consisting of all clusters, matrices or seeds (respectively) that are mutation equivalent to the original object. A matrix is \emph{mutation-infinite} if its mutation class is infinite and \emph{mutation-finite} otherwise.
\end{dfn}

We are now ready to give the definition of a cluster algebra.

\begin{dfn}[Cluster algebra]
Let $\wt{\und{x}_0} = \{x_1,\dots,x_n,y_1,\dots,y_m\}$ be an extended cluster, $\wt{B}_0$ be an $(n+m) \times n$ extended exchange matrix and let $\Sigma_0=\big( \wt{\und{x}_0},\wt{B}_0 \big)$.
Let $\mathbb{A}$ be a unital subring of $\Z \mathcal{P}$ which contains all the coefficients $p_k^{\pm}$ for all $1 \leq k \leq n$ for every seed that is mutation equivalent to $\Sigma_0$.
The \emph{cluster algebra} $\mathcal{A}(\Sigma_0)$ over $\mathbb{A}$ is the $\mathbb{A}$-subalgebra of $\F$ generated by all cluster variables (the ones that occur as an entry in the cluster of some seed that is mutation equivalent to $\Sigma_0$). The ring $\mathbb{A}$ is called the \emph{ground ring} of $\mathcal{A}$ and the number $n$ of mutable variables is the \emph{rank} of $\mathcal{A}$.
A cluster algebra is said to be of \emph{finite type} its set of cluster variables is finite and of \emph{infinite type} otherwise.
\end{dfn}

For cluster algebras of \emph{geometric type} we take $\mathbb{A}$ to be $\Z[X_{n+1}, \dots, X_{n+m}]$. As we will see in Proposition \ref{prop:geomtype_coeffts} below, this choice of $\mathbb{A}$ does indeed contain all the required coefficients to be a valid one.

\begin{rmk}
An alternative way of defining the cluster algebra $\mathcal{A}$ above is as the $\mathbb{A}$-subalgebra of $\F$ generated by the set $ \{ \var_{\wt{B}_0}\pth{p} \mid \pth{p} \text{ is a mutation path} \}$.
\end{rmk}

\emph{A priori}, a cluster variable is just a rational function in the ambient field $\F$. However, it turns out that cluster variables possess the following remarkable property.

\begin{thm}[{\cite[Theorem 3.1]{FZCAI}}]
\label{thm:laurent_phenomenon}
Any cluster variable is a Laurent polynomial in the variables from the initial (or any other) cluster with coefficients in $\Z \mathcal{P}$.
\end{thm}

This is known as the Laurent phenomenon. In the case of cluster algebras of geometric type, we have an even stronger result:

\begin{prop}[{\cite[Proposition 3.20]{GSV}}]
\label{prop:geomtype_coeffts}
In a cluster algebra of geometric type, any cluster variable is a Laurent polynomial in the initial (or any other) cluster, whose coefficients are polynomials in the stable variables.
\end{prop}

\begin{rmk} \label{rmk:mut_is_rat_fn}
Let $\pth{p} = \pths{p_2}{p_1}$ be a mutation path and let $\und{x}=\cl_{\wt{A}} \pth{p_1}$ and $\wt{B}=\mmut{\wt{A}}{\und{p_1}}$. Then $\var_{\wt{A}} \pth{p}=\var_{\wt{B}} \pth{p_2}(\und{x})$. We can think of $\var_{\wt{B}} \pth{p_2}$ as a rational function $f$ in $n$ variables with coefficients in $\Z \mathcal{P}$. Then $\var_{\wt{A}} \pth{p} = f(\und{x})$. (Here we think of $\und{x}$ as a tuple, so that $f(\und{x})$ makes sense.)
\end{rmk} %

A cluster algebra often contains multiple seeds whose clusters are equal up to  relabelling the respective sets and which produce the same set of exchange relations when mutated in all possible directions. Such seeds are the same in the sense that respective mutation is the same up to permutation of direction. 

\begin{dfn}
Let $\wt{A}$ and $\wt{B}$ be $(n+m) \times n$ skew-symmetrisable integer matrices. 
We will say $\wt{A}$ and $\wt{B}$ are \emph{essentially equivalent (by $\sigma$)}
 if there exists a permutation $\sigma \in S_n$ such that
\begin{equation}
\label{eqn:esseq}
\var_{\wt{A}}[k]\big(\{X_1,\dots,X_n\}\big) = \var_{\wt{B}} [\sigma^{-1}(k)]\big( \sigma\{X_1,\dots,X_n\}\big) 
\end{equation} 
 for all $k \in \{1, \dots, n\}$.
We then write $\wt{B} \esseq \sigma (\wt{A})$. Similarly, we say that two seeds $\big( \und{x},\wt{A}\big)$ and $\big( \und{y},\wt{B}\big)$ are essentially equivalent (by $\sigma$) if there exists $\sigma \in S_n$ such that $\wt{B} \esseq \sigma(\wt{A})$ and $\und{y} = \sigma(\und{x})$. We may then write  $\big( \und{y},\wt{B}\big) \esseq \sigma \big( \und{x},\wt{A}\big)$.
\end{dfn}

\begin{ex} \label{ex:esseq}
Let $a,b,c \in \N$. The seed 
$\left( \{x_1,x_2,x_3\} , A=\ssmat{-a}{-b}{c} \right)$
is essentially equivalent to \sloppy
$\left( \{x_1,x_3,x_2\} ,B=\ssmat{-c}{-b}{a} \right)$ 
by the permutation $\sigma = (2~3)$. Explicitly we have
\begin{align*}
\scalemath{0.975}{ \var_A[1](\{x_1,x_2,x_3\}) = \frac{x_2^a+x_3^c}{x_1} = \var_B[1](\{x_1,x_3,x_2\}) =\var_B[\sigma^{-1}(1)](\sigma\{x_1,x_2,x_3\}), }\\
\scalemath{0.975}{ \var_A[2](\{x_1,x_2,x_3\}) = \frac{x_1^a+x_3^b}{x_2} = \var_B[3](\{x_1,x_3,x_2\}) =\var_B[\sigma^{-1}(2)](\sigma\{x_1,x_2,x_3\}),}\\
\scalemath{0.975}{ \var_A[3](\{x_1,x_2,x_3\}) = \frac{x_1^c+x_2^b}{x_3} = \var_B[2](\{x_1,x_3,x_2\}) =\var_B[\sigma^{-1}(3)](\sigma\{x_1,x_2,x_3\}).}
\end{align*}
\end{ex}

\begin{rmk}
We have defined essential equivalence in terms of mutation directions. In the literature, such matrices and seeds are known as \emph{equivalent}, which is defined from a slightly different perspective: saying \pref{eqn:esseq} holds for all $k$ is the same as saying $\wt{A}$ and $\wt{B}$ differ (up to sign) by a simultaneous permutation of rows and columns. 

More precisely, let $\sigma$ be a permutation of $n$ objects. Define $ N_{\sigma}$ to be the function that takes a matrix $\wt{A}$ whose principal part is $n\times n$ and simultaneously permutes its rows and columns by $\sigma$.
We have
\begin{equation} \label{eqn:esseq_vs_equivalent}
\wt{A} \esseq \sigma(\wt{B}) \iff \wt{A}= \pm N_{\sigma}(\wt{B}).
\end{equation}

\sloppy We favour our perspective since it preserves the intuition that mutating ${\big( \sigma (\und{x}),N_{\sigma}(\wt{A}) \big)}$ in one direction is the same as mutating $\big( \und{x},\wt{A}\big)$ in another.
\end{rmk}

Let us check that essential equivalence of matrices is transitive; doing so is useful as it lets us write down an explicit permutation involved.

\begin{lem} \label{lem:esseq_transitivity}
Let $\wt{A}$, $\wt{B}$ and $\wt{C}$ be $(n+m) \times n$ skew-symmetrisable matrices. Suppose $\wt{A}$ is essentially equivalent to $\wt{B}$ with $\wt{B} \esseq \sigma(\wt{A})$ and $\wt{B}$ is essentially equivalent to $\wt{C}$ with $\wt{C} \esseq \rho(\wt{B})$, for some permutations $\rho$ and $\sigma$. Then $\wt{A}$ is essentially equivalent to $\wt{C}$ and we have $\wt{C} \esseq (\sigma  \rho) (\wt{A})$.
\end{lem}

\begin{proof}
Since $\wt{C} \esseq \rho(\wt{B})$ and $\wt{B} \esseq \sigma(\wt{A})$, we have
\begin{align*}
 && \var_{\wt{B}}[k]\big(\{X_1,\dots,X_n\}\big) &= \var_{\wt{C}} [\rho^{-1}(k)]\big( \rho \{X_1,\dots,X_n\}\big) & \forall k\\
\implies && \var_{\wt{B}}[k]\big( \sigma \{X_1,\dots,X_n\}\big) &= \var_{\wt{C}} [\rho^{-1}(k)]\big( \sigma (\rho \{X_1,\dots,X_n\}) \big) & \forall k \\
\implies && \var_{\wt{B}}[\sigma^{-1}(k)]\big( \sigma \{X_1,\dots,X_n\}\big) &= \var_{\wt{C}} [\rho^{-1}(\sigma^{-1}(k))]\big( \sigma (\rho \{X_1,\dots,X_n\}) \big) &  \forall k \\
\overset{ }{\implies} && \var_{\wt{A}}[k]\big(\{X_1,\dots,X_n\}\big) &= \var_{\wt{C}} [(\sigma \rho)^{-1}(k))]\big( (\sigma\rho) \{X_1,\dots,X_n\} \big)  &  \forall k
\end{align*} 
which shows that $\wt{A}$ is essentially equivalent to $\wt{C}$ with $\wt{C} \esseq (\sigma \rho)(\wt{A})$.
\end{proof}

\begin{ntn}
In the situation above, where we have $\wt{B} \esseq \sigma(\wt{A})$ and $\wt{C} \esseq \rho(\wt{B})$, we may write $\wt{C} \esseq \rho \big( \sigma(\wt{A} )\big)$.
\end{ntn}

\begin{rmk}
By Equation \ref{eqn:esseq_vs_equivalent}, if $\wt{C} \esseq \rho \big( \sigma(\wt{A} ) \big)$ then $\wt{C} = N_{\rho} \circ N_{\sigma}(\wt{A})$. Note that $N_{\rho} \circ N_{\sigma}$ is equal to $N_{\sigma \circ \rho}$ (not  $N_{\rho  \circ \sigma}$), as required by Lemma \ref{lem:esseq_transitivity}.
\end{rmk}
	
\begin{lem}
\label{lem:esseq}
Suppose the seeds $\big(\und{x},\wt{A}\big)$ and $\big(\und{y},\wt{B}\big)$ are essentially equivalent by $\sigma$. Let $\pth{p}$ be a mutation path.  We have the following.
\begin{enumerate}[(a)] 
\item  $\wt{A}_{\pth{p}}$ is essentially equivalent to $\wt{B}_{[\sigma^{-1}(\und{p})]}$ with $\wt{B}_{[\sigma^{-1}(\und{p})]} = \sigma(\wt{A}_{\pth{p}})$ (that is, by the same permutation $\sigma$ as above).
\item The seeds $\big(\cl_{\wt{A}}\pth{p}(\und{x}),\wt{A}_{\pth{p}} \big)$ and $ \big( \cl_{\wt{B}}[\sigma^{-1}(\und{p})](\und{y}), \wt{B}_{[\sigma^{-1}(\und{p})]} \big)$ are essentially equivalent, again by $\sigma$.
\end{enumerate}
\end{lem}

\begin{proof}
Part $(a)$ is easy to see when $\pth{p}$ has length one. The full result then follows by induction. For $(b)$, when $\pth{p}$ has length one, the result is immediate using part $(a)$. Otherwise, suppose $\pth{p}$ has length $r>1$ and assume the result is true for the rooted subpath of length $r-1$, that is, the path $[p_{r-1},\dots,p_1]=:[\und{p}']$. We have
\begin{equation}
\cl_{\wt{A}}\pth{p}(\und{x}) = \cl_{\wt{A}_{ [\und{p}'] }}[p_r] \big( \cl_{\wt{A}}[\und{p}'](\und{x}) \big)
 =: \cl_{\wt{A}_{[\und{p}']}}[p_r] \big( \und{x}' \big)
\end{equation}
and
\begin{equation}
\cl_{\wt{B}}[\sigma^{-1}(\und{p})](\und{y}) = \cl_{\wt{B}_{[\sigma^{-1}(\und{p}')]}}[\sigma^{-1}(p_r)] \big(\cl_{\wt{B}}[\sigma^{-1}(\und{p}')](\und{y}) \big)
=: \cl_{\wt{B}_{[\sigma^{-1}( \und{p}')]}}[\sigma^{-1}(p_r)] \big( \und{y}' \big).
\end{equation}
Now we wish to show that 
$\big( \cl_{\wt{A}_{[\und{p}']}}[p_r] \big( \und{x}' \big) ,\wt{ A}_{\pth{p}}\big)$ and $\big( \cl_{\wt{B}_{[\sigma^{-1}(\und{p}')]}}[p_r] \big( \und{y}' \big), \wt{B}_{[\sigma^{-1}(\und{p})]} \big)$ are essentially equivalent by $\sigma$. But this is true, since $(b)$ is established for paths of length one, since we have by assumption that  
$\big(  \und{x}' , \wt{A}_{[\und{p}']}\big)$ and $\big( \und{y}', \wt{B}_{[\sigma^{-1}(\und{p}')]}\big)$ are essentially equivalent by $\sigma$, and since $(\wt{A}_{[\und{p}']})_{[p_r]} =\wt{A}_{\pth{p}}$ and $(\wt{B}_{[\sigma^{-1}(\und{p}')]})_{[\sigma^{-1}(p_r)]} = \wt{B}_{[\sigma^{-1}(\und{p})]}$. Thus $(b)$ follows by induction.
\end{proof}

Lemma \pref{lem:esseq} immediately implies the following formula, which is useful for relating mutation paths involving essentially equivalent seeds.

\begin{cor}
\label{cor:esseq_var}
Suppose $\wt{A}$ and $\wt{B} \esseq \sigma(\wt{A})$ are essentially equivalent matrices. Then
\begin{equation}
 \var_{\wt{A}} \pth{p}(\und{x}) = \var_{\sigma(\wt{A})}[\sigma^{-1}(\und{p})](\sigma(\und{x})).
 \label{eqn:esseq_var}
\end{equation}
\end{cor}

Notice that for any skew-symmetrisable matrix $\wt{B}$, we have that $\wt{B}$ and $-\wt{B}$ are \emph{trivially essentially equivalent} (i.e.\ essentially equivalent by the identity permutation). Thus, by \pref{eqn:esseq}, reversing the sign of a matrix does nothing as far variables are concerned. Therefore we may at times be liberal with the signs of exchange matrices in such a context.

\begin{prop}
\label{prop:esseq_isom_graded_algebras}
If $\wt{A}$ and $\wt{B}$ are essentally equivalent matrices then the initial seeds $(\und{x},\wt{A})$ and $(\und{x},\wt{B})$ generate isomorphic algebras.
\end{prop}

\begin{proof}
Write $\und{x} = (x_1, \dots, x_n)$ as usual. Let $z$ be a cluster variable in $\mathcal{A}(\und{x},\wt{A})$. Then $z= \var_{\wt{A}} \pth{p}(\und{x})$ for some path $\pth{p}$. Let $f_{\wt{A},\pth{p}}$ be the rational function corresponding to this path with respect to $\wt{A}$. Let $\sigma \in S_n$ be a permutation such that $B \esseq \sigma(\wt{A})$. By Corollary \ref{cor:esseq_var}, $f_{\wt{A},\pth{p}}$ is then the rational function $f_{\wt{B},[\sigma^{-1}(\und{p})]} \circ \sigma$, where $f_{\wt{B},[\sigma^{-1}(\und{p})]}$ is the function corresponding to the path $[\sigma^{-1} (\und{p})]$ in $\A(\und{x},\wt{B})$.
Therefore any cluster variable $z$ in $\A(\und{x},\wt{A})$ corresponds to the cluster variable $y$ in $\A(\und{x},\wt{B})$ that is obtained by replacing every occurrence of initial variable $x_j$ in $z$ by $x_{\sigma(j)}$. Since the cluster variables generate the algebras, it follows that they are isomorphic.
\end{proof}

\begin{dfn}
Let $\wt{B}$ be an $(n+m) \times n$ skew-symmetrisable matrix. Mutation directions $i$ and $j$ are called \emph{essentially equivalent} with respect to $\wt{B}$ if mutation of $\left(\{x_1,\dots,x_n\} ,\wt{B}\right)$ in direction $i$ yields the same cluster variable as that for mutation in direction $j$, but with $x_j$ in the latter swapped with $x_i$. 
\end{dfn}

\begin{ex}
Let $B=\left( \begin{smallmatrix}0&a&-c\\-a& 0&c\\c&-c&0\end{smallmatrix} \right)$.
Then directions 1 and 2 are essentially equivalent with respect to $B$ since mutation in directions 1 and 2 produces the elements $(x_2^a+x_3^c)/{x_1}$ and $(x_1^a+x_3^c)/{x_2}$, respectively.
\label{ex:esseq_directions}
\end{ex}

\begin{lem}
Let $\sigma \in S_n$ be a permutation which sends each mutation direction $i$ to an essentially equivalent direction $j$ with respect to the matrix $\wt{B}$. Then
\begin{equation}
  \var \pth{p}_{\wt{B}}(\sigma(\und{x})) = \var[\sigma(\und{p})]_{\wt{B}}(\und{x}).
\end{equation}
\end{lem}

\begin{proof}
This is immediate from the definition.
\end{proof}

\begin{ex}
With $B$ as in Example \ref{ex:esseq_directions}, we have $\var_B [1](x_2,x_1,x_3) = \\ (x_1^a+x_3^c)/{x_2}= \var_B[2](x_1,x_2,x_3)$.
\end{ex}

Often we will want to repeatedly apply a mutation path that leaves the mutated degree seed essentially equivalent to the initial one, but not trivially so. In this case we will need to permute the path along the way to ensure we are mutating the new seed in an analogous way. The following gives a notation for this.

\begin{ntn}
Let $\pth{p}$ be a mutation path and $\sigma$ a permutation. We define
\[ \Prpth{\und{p}}{k}{\sigma}:=[\sigma^{-(k-1)}(\und{p}),\dots,\sigma^{-1}(\und{p}),\und{p}]. \]
\end{ntn}

Finally, there are three further important structures associated to a cluster algebra.

\begin{dfn}
The \emph{exchange tree} of a cluster algebra $\A$ is the $n$-regular tree $\mathbb{T}_n$ with each vertex $t_i$ associated to a seed $\Sigma(t_i)$ and with two vertices $t_i$ and $t_j$ connected by an edge labeled $k$ if and only if $\Sigma(t_i) = \mu_k(\Sigma(t_j))$. We usually denote the vertex associated to the initial seed by $t_0$.
\end{dfn}

\begin{dfn}
The \emph{exchange graph} of $\A$ is the exchange tree modulo essential equivalence of seeds. As it no longer makes sense to label edges by mutation direction, we may instead choose to label them by the relevant exchange relation, that is, the right hand side of Equation \ref{eqn:cluster_mutation}.
\end{dfn}

\begin{dfn}
The cluster complex of $\A$, denoted $\Delta(\A)$, is the simplicial complex whose ground set is the set of all cluster variables and whose maximal simplices are are the clusters.
\end{dfn}

\section{Graded cluster algebras}

There are several notions of gradings for cluster algebras in the literature. Our definition of gradings will follow \cite{G15}.

\begin{dfn}
Let $\A$ be a cluster algebra of rank $n$. A \emph{(multi-)graded seed} is a triple $\big( \wt{\und{x}}, \wt{B}, \und{g} \big)$ where $\big(\wt{\und{x}}, \wt{B} \big)$ is a cluster in $\A$ and $\und{g}$ is a $r \times (n+m)$ integer matrix satisfying 
\begin{equation}
\label{eqn:grading_condition}
\und{g} \wt{B} = \und{0}.
\end{equation}
Such a matrix $\und{g}$ is called a \emph{grading} for $\big(\wt{\und{x}}, \wt{B} \big)$. The $i^\text{th}$ column of $\und{g}$ is the \emph{degree} of $x_i$ and the degree of a homogeneous rational function in the variables of a given graded cluster (possibly with coefficients in frozen variables) is defined in the obvious way.  
The condition \pref{eqn:grading_condition} ensures that exchange relations, and therefore all cluster variables, will be homogeneous. Then, given a cluster variable $z$, we will denote its degree with respect to $\und{g}$ by $\deg_{\und{g}} (z)$. 

\emph{Mutation of a graded seed} in direction $k$ is defined in the same way as in Definition \ref{dnf:matrix_mutation}, with the addition that the $k^\text{th}$ entry of $\und{g}$ is replaced by the degree of the new cluster variable. Alternatively, we may directly mutate $\und{g}$: we replace the variables in Equation \ref{eqn:cluster_mutation} by their degrees, and replace multiplication, division and taking powers of variables by addition, subtraction and multiplication, respectively, of their degrees.
 
The following objects are useful when the only information we need about some cluster variables pertains to their degrees.
A \emph{degree seed} in $\A$ is a pair $\big( \und{g}, \wt{B} \big)$ such that $\big( \wt{\und{x}}, \wt{B} ,\und{g} \big) $ is a graded seed occurring in $\A$. A \emph{graded cluster} or \emph{degree cluster} (a pair of the form $(\wt{\und{x}}, \und{g})$ or the tuple $\und{g}$, respectively) is defined in a similar way.
 
Finally, a \emph{graded cluster algebra} will mean a cluster algebra generated from an initial graded seed.
\end{dfn}

\begin{rmkdfn} \label{rmkdfn:gradings_space_Z_vs_Q}
While we stipulate that the entries of $\und{g}$ are integers, there is no reason in principle why they could not be elements of an arbitrary abelian group, provided $\und{g}$ satisfies \pref{eqn:grading_condition} (see \cite{GP} for more discussion along these lines). We may wish to refer to the \emph{grading space} of $\wt{B}$, which we define as the set of matrices $\und{g}$ over some field satisfying \pref{eqn:grading_condition} (i.e.~the nullspace of $\wt{B}^T$). When we do so, we will think of this vector space as being over $\mathbb{Q}$. Often we will be interested in working with grading matrices $\und{g}$ whose rows form a basis of the grading space (or possibly a subspace of the gradings space). We will call these rows \emph{grading vectors}. In practice, we will scale these basis vectors so that they consist of integer entries (which is always possible if the entries are in $\mathbb{Q}$).
\end{rmkdfn}

\begin{ntn} 
\label{ntn:g_matrix_or_tuple_of_degrees}
We will sometimes use the same notation to mean both the matrix $\und{g}$ and the tuple whose entries are the columns of the matrix (though it should be clear from the context which version we are referring to at any given time). More precisely, let $g_1, \dots, g_{n+m}$ be the columns of $\und{g}$. We will write $\und{g}=(g_1,\dots,g_{n+m})$ and refer to  $g_1, \dots, g_{n+m}$ as the \emph{entries} of $\und{g}$. We do this so $\und{g}$ can be thought of as a tuple containing the degrees of the variables in the corresponding cluster. To fit better with this overloaded notation, we have transposed the matrix $\und{g}$ from that in \cite{G15} so that, if $r=1$, say, $\und{g}$ is a tuple of integer degrees.

We extend Notation \ref{ntn:sd_cl_B_p} to the graded setting in the natural way. Aside from this, if we wish to denote the degree of a variable, or just the degree seed or cluster, we add the prefix $``\deg."$ to the respective notation, and we may simply write $\deg_{\wt{B}} \pth{p}$ to mean $\deg.\var_{\wt{B}} \pth{p}$. Again, $\pth{p} ({\und{x}},\wt{B},\und{g})$ will mean the graded seed obtained on applying the mutation path $\pth{p}$ to $ ({\und{x}},\wt{B},\und{g})$.
\end{ntn} 

\begin{prop}[{\cite[Proposition 3.2]{G15}}]
Let $\Sigma=(\wt{\und{x}}, \wt{B}, \und{g})$ be a graded seed. The graded cluster algebra $\A(\Sigma)$ is a $\Z^r$-graded algebra, and every cluster variable of $\A(\Sigma)$ is homogeneous with respect to this grading.
\end{prop}

\begin{dfn} \label{dfn:standard_grading}
The grading matrix $\und{g}$ is called a \emph{standard} grading if its rows form a basis of the grading space. A seed or degree seed containing $\und{g}$ is called a \emph{standard} graded seed or degree seed.
\end{dfn}

As mentioned in Remark/Definition \ref{rmkdfn:gradings_space_Z_vs_Q}, we will often choose our grading matrix $\und{g}$ such that the rows form a basis of the grading space. We do so since understanding the behaviour of such a grading essentially lets us understand all gradings on the given cluster algebra. The next two results show that this is the case.

\begin{lem}[{\cite[Lemma 3.7]{G15}}]
Any mutation of a standard graded seed or standard degree seed is again standard.
\end{lem}

\begin{lem}[{\cite[Lemma 3.8]{G15}}]
Let $\Sigma=\big( \wt{\und{x}}, \wt{B}, \und{g} \big)$ be a standard graded seed and let $\und{h}$ be any grading for $\big( \wt{\und{x}}, \wt{B} \big)$. There exists an integer matrix $M=M(\und{g}, \und{h})$ such that for any cluster variable $z$ in $\A \big( \wt{\und{x}}, \wt{B}, \und{h} \big)$ we have 
\begin{equation}
\deg_{\und{h}} (z) = M\deg_{\und{g}} (z).
\end{equation}
\end{lem}

\begin{dfn}
Two graded seeds  $\big( \wt{\und{x}}, \wt{A}, \und{g} \big)$ and  $\big( \wt{\und{y}}, \wt{B}, \und{h} \big)$ are essentially equivalent if their underlying seeds are.
Two degree seeds $(\und{g},\wt{A})$ and $(\und{h},\wt{B})$ are essentially equivalent by $\sigma$ if $\wt{B} \esseq \sigma(\wt{A})$ and \pref{eqn:esseq} holds for all $k$ but for degrees rather than variables.
\end{dfn}

Of course, essential equivalence of degree seeds is much weaker than essential equivalence of graded seeds: two graded seeds that are not essentially equivalent may have underlying degree seeds that are. Indeed, we are often interested in finding infinite sequences of non-equivalent graded seeds whose underlying degree seeds are essentially equivalent. This allows us to prove the existence of infinitely many different cluster variables of a particular degree or set of degrees.

Finally, note the simple relationship between the mutation of a degree cluster and its negative.

\begin{lem} \label{lem:mut_of_deg_seed_vs_negative}
Let  $(\und{g}, \wt{B})$ be a degree seed and $\pth{p}$ a mutation path. We have 
\[\deg.\cl_{\wt{B}}\pth{p}(-\und{g}) = - \deg.\cl_{\wt{B}}\pth{p}(\und{g}) \]
\end{lem}

\begin{proof}
For $l(\pth{p})=1$, this is immediate from Equation \ref{eqn:cluster_mutation}. The full result then follows inductively.
\end{proof}

\section{Quivers}

If $\wt{B}$ is an extended exchange matrix whose principal part is skew-symmetric, we can represent $\wt{B}$ with a quiver, a directed weighted graph with no loops or $2$-cycles. Working with a quiver is often easier than working with the corresponding matrix, and it lets us represent a whole seed in one diagram by identifying the cluster variables (or their degrees) with the vertices of the quiver. For background on quivers, see, for example, \cite{SchifflerQR}.

Unless otherwise stated, we will assume any quiver we refer to is connected.

\begin{dfn} \label{dfn:quiver}
Let $B$ be an $n \times n$ skew-symmetric matrix. The quiver $Q(B)$ associated to $B$ has $n$ vertices (labelled $1,\dots, n$) and an arrow of weight $b_{ij}$ from vertex $i$ to vertex $j$ (a negative entry counts as an arrow in the opposite direction and a weight $0$ arrow counts as no arrow).

We can also assign a quiver $Q(\wt{B})$ to an extended exchange matrix $\wt{B}$. This quiver has $n+m$ vertices (the last $m$ of which correspond to frozen variables) and arrows are assigned in the same way as above. The frozen vertices are denoted by drawing a square around them.

A quiver is \emph{acyclic} if it contains no oriented cycles. A rank $3$ quiver is \emph{cyclic} if it is an oriented cycle.

A \emph{subquiver} of $Q$ is a quiver $R$ whose vertices are a subset of those in $Q$, such that if $v_1$ and $v_2$ are vertices in $R$ then the edge from $v_1$ to $v_2$ in $R$ is of the same weight and direction as it is in $Q$. What we define as a subquiver is often called a \emph{full} subquiver in the literature.
\end{dfn}

\begin{dfn}[Quiver mutation]
The \emph{mutation} of a quiver $Q(B)$ in direction $k$ is the quiver $Q'$ corresponding to $B_{[k]}$. Equivalently, we obtain $Q'$ from $Q$ as follows:
\begin{enumerate}[(1)]
\item Reverse each arrow incident to vertex $k$.
\item  For each directed path of the form 
\Image{
\begin{tikzpicture}[node distance=1.6cm, auto, font=\footnotesize]                                          
  \node (L) []{$i$};                                                                                                                     
  \node (C) [right of = L]{$k$};                                                                               
  \node (R) [right of = C]{$j$};   
  \draw[-latex] (L) to node {$b_{ik}$} (C);      
  \draw[-latex] (C) to node {$b_{kj}$} (R);    
\end{tikzpicture},}
add $b_{ik} b_{kj} $ to the weight of the arrow from vertex $i$ to $j$ (which may involve creating an arrow or reversing or deleting a present arrow).
\end{enumerate}
\end{dfn}

In general, we will use a quiver and its corresponding matrix interchangeably, but we make a couple of specific definitions for working with the corresponding objects.

\begin{dfn} \label{dfn:quiver_seed_degree_balanced}
A \emph{(graded) quiver seed} will mean a (graded) seed in which we identify vertex $i$ of the quiver associated to the extended exchange matrix with the cluster variable in position $i$ of the cluster. We represent this in a diagram by placing the cluster variables where the corresponding vertices would be. When we wish to emphasise the vertex numbering, we will place each cluster variable inside parentheses with the corresponding vertex number as a subscript. In terms of quivers, an equivalent condition to \pref{eqn:grading_condition} for the matrix $\und{g}$ is that each non-frozen vertex $k$ is \emph{balanced}. By this we mean that 
\begin{equation} \label{eqn:degree_quiver_balanced_requirement}
 \sum_{i \overset{b_{ik}}{\rightarrow} k} b_{ik} \deg{(i)}  = \sum_{k \overset{b_{ki}}{\rightarrow} i} b_{ki} \deg{(i)} ,
 \end{equation}
where the left sum runs over all weighted arrows into vertex $k$ and the right sum does the same for arrows out. We will sometimes refer to the left (right) hand side of \pref{eqn:degree_quiver_balanced_requirement} as the \emph{weight into (out of) $k$}, which we will denote $w_{\rightarrow}(k)$ ($w_{\leftarrow}(k)$).

A \emph{degree quiver} will mean a degree seed in which we identify vertex $i$ of the associated quiver with the $i^\text{th}$ entry of $\und{g}$. In the corresponding diagram, the $i^\text{th}$ vertex will be replaced by $(d)_i$, where $(d) =g_i$ (though sometimes we may omit the subscript if it is not needed). A degree quiver is \emph{positive} if the degree corresponding to each vertex is non-negative.

We extend the definition of \emph{essentially equivalent} to the above objects in the obvious way: two (graded) quiver seeds are essentially equivalent if their underlying (graded) seeds are, and two degree quivers are essentially equivalent if their underlying degree seeds are. Alternatively, two quivers are essentially equivalent if they are the same up to relabelling vertices and taking the opposite quiver (i.e.~reversing all arrows), and similarly for (graded) quiver seeds and degree quivers.
\end{dfn}

\begin{ex}
\sloppy Let $B = \ssmat{-5}{-4}{3}$ and let $\Sigma=\big((x_1,x_2,x_3),B,(4,3,5) \big)$. 
Then $\Sigma$ is a graded seed. (Strictly, in accordance with Notation \ref{ntn:g_matrix_or_tuple_of_degrees}, the degree seed in $\Sigma$ is the tuple $((4),(3),(5))$ of $1 \times 1$ matrices or the corresponding $1 \times 3$ matrix $(4~3~5)$). The quiver seed corresponding to $\Sigma$ is
\Image{
\begin{tikzpicture}[node distance=1.6cm, auto, font=\footnotesize]                                          
  \node (A2) []{$(x_2)_2$};                                                                                                                     
  \node (A1) [below left of = A2]{$(x_1)_1$};                                                                               
  \node (A3) [below right of = A2]{$(x_3)_3$};   
  \draw[latex-] (A2) to node [swap] {$5$} (A1);      
  \draw[latex-] (A3) to node [swap] {$4$} (A2);    
  \draw[-latex] (A3) to node [] {$3$} (A1);      
\end{tikzpicture}
}
and the degree seed is $Q_1=$
\Image{
\begin{tikzpicture}[node distance=1.6cm, auto, font=\footnotesize]                                          
  \node (A2) []{$(3)_2$};                                                                                                                     
  \node (A1) [below left of = A2]{$(4)_1$};                                                                               
  \node (A3) [below right of = A2]{$(5)_3$};   
  \draw[latex-] (A2) to node [swap] {$5$} (A1);      
  \draw[latex-] (A3) to node [swap] {$4$} (A2);    
  \draw[-latex] (A3) to node [] {$3$} (A1);      
\end{tikzpicture}
}.
For vertex $1$, we have 
\[ \sum_{i \overset{b_{i1}}{\rightarrow} 1} b_{i1} \deg{(i)}  =3.5 =5.3 = \sum_{1 \overset{b_{1i}}{\rightarrow} i} b_{1i} \deg{(i)} ,\]
which shows that this vertex is balanced, as are the other two. The degree seed
$Q_2=$
\Image{
\begin{tikzpicture}[node distance=1.6cm, auto, font=\footnotesize]                                          
  \node (A2) []{$(5)_2$};                                                                                                                     
  \node (A1) [below left of = A2]{$(4)_1$};                                                                               
  \node (A3) [below right of = A2]{$(3)_3$};   
  \draw[latex-] (A2) to node [swap] {$3$} (A1);      
  \draw[latex-] (A3) to node [swap] {$4$} (A2);    
  \draw[-latex] (A3) to node [] {$5$} (A1);      
\end{tikzpicture}
}
is essentially equivalent to $Q_1$: we have $Q_2 \esseq \sigma(Q_1)$, where $\sigma =(2~3)$.
\end{ex}

\section{Denominator vectors}

The \emph{denominator vector} of a cluster variable is a well-known invariant which we will make frequent use of. Denominator vectors were introduced by Fomin and Zelevinsky and we take the definition from \cite{FZCAIV}.

\begin{dfn} \label{dfn:denominator_vector}
Let $\und{x}=(x_1,\dots,x_n)$ be a cluster in a cluster algebra $\A$. By Theorem \ref{thm:laurent_phenomenon} we can express any given cluster variable $z$ uniquely as
\begin{equation}
z = \frac{N(x_1,\dots,x_n)}{x_1^{d_1}, \dots, x_n^{d_n}},
\end{equation}
where $N(x_1,\dots,x_n)$ is a polynomial with coefficients in $\Z \mathcal{P}$ that is not divisible by any $x_i$.	
The vector
 \[ \dv_{ \und{x} } (z) =  \begin{bmatrix} d_1 \\ \vdots \\d_n \end{bmatrix} \]
is called the \emph{denominator vector} of $z$ with respect to $\und{x}$. Its $i^{\text{th}}$ entry will be referred to as $\dv (z)|_i$ or $\dv (z)_i$. If the cluster with respect to which we are writing our denominator vector is clear, we will omit the subscript $\und{x}$ and write $\dv (z)$ (if not specified, we assume $\und{x}$ is the initial cluster from which we are generating our cluster algebra).
\end{dfn}

When working with a denominator vector we may write it as a tuple rather than a column vector.

\begin{ex}
With respect to the cluster $\und{x}=(x_1,\dots,x_n)$, the variable $z=\frac{x_2^a+x_3^c}{x_1} $ has denominator vector $ \dv_{ \und{x} } (z) =  (1,0,0)$.
\end{ex}

\begin{dfn}
Let $\und{x}=(x_1,\dots,x_n)$ be a cluster of rank $n$. Then $\dcl (\und{x})$ will denote the corresponding \emph{denominator cluster}, that is, the $n$-tuple whose $i^{\text{th}}$ entry is $\dv (x_i)$ (with respect to some initial cluster). 
 We will refer to the $i^{\text{th}}$ entry of $\dcl (\und{x})$ by $\dcl (\und{x})^i$ or $\dcl (\und{x})|^i$. 
 We will also use $\dv_B\pth{p}$ and $\dcl_B \pth{p}$ for the denominator vector or cluster of the variable or cluster, respectively, obtained after mutation path $\pth{p}$.
 
 A \emph{denominator quiver} will mean a quiver corresponding to a given seed $(\und{x},B)$ but where we identify the $i^\text{th}$ vertex with $\dcl(\und{x})|^i$. 
 In the corresponding diagram, the $i^\text{th}$ vertex will be replaced by $(\dcl(\und{x})|^i)_i$ (sometimes omitting the subscript if it is clear from context). 
 \end{dfn}
 
\begin{dfn}
\label{dfn:partial_order_on_denominator_vectors}
We define a partial order on the set of denominator vectors in $\Z^n$ by setting $\dv(w) \geq \dv(z)$ if each entry of $\dv(w)$ is greater than or equal to the corresponding entry of $\dv(z)$, and we also say $\dv(w) > \dv(z)$ if $\dv(w) \geq \dv(z)$ and at least one entry of $\dv(w)$ is strictly greater than the corresponding entry of $\dv(z)$.
\end{dfn}

Sometimes we will be interested in finding mutation paths which produce growing denominator vectors. When proving that denominator vectors are growing, it becomes cumbersome to work with the entire cluster of denominator vectors, so we single out only one entry of each denominator vector. This simplifies calculations, and showing that this entry grows is enough to distinguish the variables. This motivates the following definition.

\begin{dfn}
The $j^{\text{th}}$ \emph{denominator slice} of $\dcl (\und{x})$ is the vector whose $i^{\text{th}}$ entry is the  $j^{\text{th}}$ entry of  $\dcl (\und{x})^i$. We will denote this using angle brackets as  
\[ \dslof{j}{\und{x}}=\la \dcl(\und{x})^1_j, \dots, \dcl(\und{x})^n_j \ra= \la \dcl(\und{x})|^1_j, \dots, \dcl(\und{x})|^n_j \ra .\]

Given a seed $(\und{x},B)$, the $j^{\text{th}}$ \emph{denominator quiver} will mean the quiver corresponding to $B$ but where we identify the $i^\text{th}$ vertex with $\dcl(\und{x})|^i_j$.
In the corresponding diagram, the $i^\text{th}$ vertex will be replaced by $\la \dcl(\und{x})|^i_j \ra_i$ (again, sometimes omitting the subscript $i$).
\end{dfn}

\begin{ex}
Consider the seed 
\[ \Sigma = \left(\und{z}= \Big( \frac{x_2^2 + x_3}{x_1}, \frac{x_1^2x_3 + x_2^4 + 2x_2^2x_3 + x_3^2}{x_1^2x_2},  x_3 \Big), B=\ssmat{-2}{-1}{1} \right). \]
The denominator cluster with respect to $\und{x}=(x_1, x_2, x_3)$ is 
\[ \dcl_{ \und{x} } (\und{z}) = \left( ( 1, 0, 0 ),( 2, 1, 0 ), ( 0, 0, -1 ) \right). \]
The denominator quiver is
\Image{
\begin{tikzpicture}[node distance=1.6cm, auto, font=\footnotesize]                                          
  \node (A2) []{$(2,1,0)_2$};                                                                                                                     
  \node (A1) [below left of = A2]{$(1,0,0)_1$};                                                                               
  \node (A3) [below right of = A2]{$(0,0,-1)_3$};   
  \draw[latex-] (A2) to node [swap] {$2$} (A1);      
  \draw[latex-] (A3) to node [swap] {$1$} (A2);    
  \draw[-latex] (A3) to node [] {$1$} (A1);      
\end{tikzpicture}
}.
The first denominator slice of $\und{z}$ is $\la 1,2,0 \ra$ and the first denominator quiver is
\Image{
\begin{tikzpicture}[node distance=1.8cm, auto, font=\footnotesize]                                          
  \node (A2) []{$\la 2 \ra_2$};                                                                                                                     
  \node (A1) [below left of = A2]{$\la 1 \ra_1$};                                                                               
  \node (A3) [below right of = A2]{$\la 0 \ra_3$};   
  \draw[latex-] (A2) to node [swap] {$2$} (A1);      
  \draw[latex-] (A3) to node [swap] {$1$} (A2);    
  \draw[-latex] (A3) to node [] {$1$} (A1);      
\end{tikzpicture}
}.
\end{ex}

The next formula lets us mutate denominator vectors directly.

\begin{lem}[{\cite[Equation 7.7]{FZCAIV}}]
\label{lem:denominator_mutation}
The new cluster variable $x_k'$ obtained after mutation of the seed $(\und{x}=(x_1,\dots,x_n), B)$ in direction $k$ satisfies
\begin{equation}
\label{eqn:denominator_mutation}
\dv(x_k')= -\dv(x_k) +  \max \left( \sum_i [b_{ik}]_+ \dv(x_i), \sum_i [-b_{ik}]_+ \dv(x_i)  \right),
\end{equation}
where $[a]_+ := \begin{cases} a &\text{ if $a \geq 0$,}\\ 0 &\text{ otherwise,} \end{cases}$
and where $\max$ is taken component-wise.
\end{lem}

\biblio

\chapter{Gradings on finite type cluster algebras \label{chap:finite_type_cluster_algebras}} \label{chp:finite_type}

In this chapter we classify gradings on the cluster algebras of type $B_n$ and $C_n$. This work was completed prior to the publication of \cite{G15} and the results are reported in the same paper. The classification for the other finite type cases that admit nontrivial gradings, namely $A_n$, $D_n$ and $E_n$, can be also found in \cite{G15}.

\section{Review of background}

\begin{dfn}[{\cite[Section 1.3]{FZCAII}}]
Let $B=(b_{ij})$ be an $n \times n$ skew-symmetrisable matrix. The \emph{Cartan counterpart of $B$}, $A(B)=(a_{ij})$, is defined by setting
\begin{equation}
a_{ij}=
\begin{cases}
2 &\text{ if $i=j$}, \\
-|b_{ij}| &\text{ otherwise.}
\end{cases}
\end{equation}
\end{dfn}

\begin{thm}[{\cite[Theorem 1.8]{FZCAII}}]
A skew-symmetrisable matrix $B$ gives rise to a finite type cluster algebra  $\mathcal{A}(\und{x},B)$ if and only if $B$ is mutation equivalent to some matrix $B'$ whose Cartan counterpart $A(B')$ is of finite type (i.e.\ $A$ is irreducible and all its principal minors are strictly positive). 
\end{thm}

\begin{dfn}
A skew-symmetrisable matrix $B$ is \emph{bipartite} if $b_{ij}b_{ik} \geq 0$ for all $i,j,k$.
\end{dfn}

\begin{lem}[{\cite[Theorem 1.9]{FZCAII}}] \label{lem:bijection_almost_pos_roots_cluster_vars_finite_type}
Let $B$ be a bipartite skew-symmetrisable matrix and $\mathcal{A}(\und{x},B)$ a finite type cluster algebra. There is a bijection $\alpha \mapsto x[\alpha]$ between the almost positive roots (i.e.\ the positive roots along with the negation of the simple roots) of the semi-simple Lie algebra corresponding to $A(B)$ and the cluster variables of $\mathcal{A}$. Write $\alpha = \lambda_1 \alpha_1 + \dots + \lambda_n \alpha_n$ as a linear combination of the simple roots $\alpha_1, \dots, \alpha_n$. The cluster variable $x[\alpha]$ is then expressed in terms of the initial cluster $\und{x}=(x_1, \dots, x_n)$ as
\begin{equation}
x[\alpha] = \frac{P_{\alpha} (\und{x})}{x^{\alpha}}
\end{equation}
where ${P_{\alpha} (\und{x})}$ is a polynomial over $\Z$ with non-zero constant term and $x^{\alpha}:= \prod_{i=1}^n x_i^{\lambda_i}$.
Under this bijection, $x[-\alpha_i]=x_i$. That is, the negative simple roots correspond to the initial cluster variables.
\end{lem}

Note for types $B_n$ and $C_n$ there are $n^2+n$ almost positive roots.

The following allows us to calculate the degree of a cluster variable in terms of its corresponding root.

\begin{lem}[{\cite[Corollary 4.2]{G15}}]  \label{lem:deg_in_terms_of_simples}
Let  $\mathcal{A}\big(\und{x},B, \und{g}=(g_1,\dots,g_n) \big)$ be a graded cluster algebra of finite type where $B$ is bipartite. Let $\alpha$ be an almost positive root and write $\alpha = \lambda_1 \alpha_1 + \dots + \lambda_n \alpha_n$ as a linear combination of the simple roots.
Then 
\begin{equation}
\label{eqn:deg_in_terms_of_simples}
\deg(x[\alpha])= -\sum_{i=1}^n \lambda_i  g_i.
\end{equation} 
\end{lem}

It turns out that any gradings on cluster algebras of type $B_n$ or $C_n$, or indeed any finite type cluster algebra, must be \emph{balanced}, by which we mean that there is a bijection between the cluster variables of degree $d$ and degree $-d$. For an explanation of this fact, see \cite[Corollary 5.7]{G15} and the remark following. We will verify this for types $B_n$ and $C_n$ below.

\section{Gradings in type $B_n$}

The appropriate bipartite matrix for type $B_n$ is the $n \times n$ matrix is given by
\[
B= \scalemath{0.85}{\begin{pmatrix}
0  & 1     &     &      &      &    &  &  \\
-1  &  0 & -1 &   &   &   &  &  \\
  & 1 & 0  &1  &  &    & &   \\
  &   & -1 & 0  &  &    & &   \\
   & &  & &  &\ddots & 1 \\   
&    &        &         &      &-1  & 0  &-1    \\
&    &        &         &      &   &2   &  0   
\end{pmatrix}} \text{ or } 
B= \scalemath{0.85}{\begin{pmatrix}
0  & 1     &     &      &      &    &  &  \\
-1  &  0 & -1 &   &   &   &  &  \\
  & 1 & 0  &1  &  &    & &   \\
  &   & -1 & 0  &  &    & &   \\
   & &  & &  &\ddots & -1 \\   
&    &        &   &  &1  & 0  &1    \\
&    &        &   &  &   &-2   &  0   
\end{pmatrix}}
\]
corresponding to whether $n$ is odd or even, respectively.
One checks that if $n$ is even, $B$ has full rank and so there is no non-trivial grading. If $n$ is odd, $B$ has rank $n-1$. The grading space is then spanned by a one dimensional vector. There are two cases for the grading vector $\und{g}$:
 \[\und{g}=
 \begin{cases}
 (2,0,-2,0,\dots,0,-2,0,1) &\text{ if } n=4k+1,\\
  (2,0,-2,0,\dots,0,2,0,-1) &\text{ if } n=4k+3.\\
 \end{cases}
 \]
We will assume from now on that $n=4k+3$, but the final results we obtain are the same in both cases.

We take the following description of the root system  for $B_n$ (as found in~\cite{Carter-Book}):
\[ \{ \pm\beta_i \pm \beta_j \mid i,j=1\dots n , i \neq j \} \cup \{\pm \beta_i \mid i= 1 \dots n \}. \]
Here $\beta_i$ is the function which selects the $(i,i)$ entry of an $n \times n$  diagonal matrix, and we have taken
\[\alpha_1=\beta_1-\beta_2 ,\dots, \alpha_{n-1}=\beta_{n-1}-\beta_n, \text{ and } \alpha_n=\beta_n ,\]
where the $\alpha_i$ denote the simple roots.
From this description, the positive roots are of one of the following forms:
\begin{align}
\beta_i-\beta_j =&  \alpha_i + \cdots +\alpha_{j-1}& &(i<j),  \label{eqn:BnForm1}\\
\beta_i+\beta_j=&  \alpha_i + \dots +\alpha_{j-1}+2(\alpha_j+\cdots +\alpha_n) & &(i<j) ,  \label{eqn:BnForm2}\\
\beta_i =&  \alpha_i + \dots +\alpha_{n}& &(i=1,\dots,n)  . \label{eqn:BnForm3}
\end{align}

To find the degrees of the remaining variables in $\mathcal{A}(\und{x},B)$ we now apply Equation \ref{eqn:deg_in_terms_of_simples} to the above lists. This, along with the negative simple roots, tells us the degree of the cluster variable corresponding to each almost positive root. Counting the number of variables of each degree, we summarise the results in the Table \ref{tab:classificationBn} below.
\newline \newline
\begin{minipage}{\textwidth}
\begin{center}
{%
\setlength{\tabcolsep}{3pt} 
\renewcommand{\arraystretch}{1.2} 
    \begin{tabular}{ | c | c | c | c | c | c |}
    \hline
    Degree &Type (\ref{eqn:BnForm1}) & Type (\ref{eqn:BnForm2}) &Type (\ref{eqn:BnForm3}) & negative simples&Total \\ \hline
   -2 &\multicolumn{1}{c}{ $\frac{(n+1)^2}{8}$ } & \multicolumn{1}{c}{ $\frac{(n-1)(n-3)}{8}$} &\multicolumn{1}{c}{$0$} &$\frac{n-3}{4}$ & $\frac{(n+1)(n-1)}{4}$  \\ 
   -1&\multicolumn{1}{c}{0} & \multicolumn{1}{c}{0} &\multicolumn{1}{c}{ $\frac{n-1}{2}$ } &1 & $\frac{n+1}{2}$  \\ 
    0&\multicolumn{1}{c}{$\frac{(n+1)^2}{4}$} & \multicolumn{1}{c}{ $\frac{(n+1)(n-1)}{4}$} &\multicolumn{1}{c}{0} &$\frac{n-1}{2}$ & $\frac{(n+1)(n-1)}{2}$  \\ 
    1&\multicolumn{1}{c}{0} & \multicolumn{1}{c}{0} &\multicolumn{1}{c}{$\frac{n+1}{2}$} &0&$\frac{n+1}{2}$ \\ 
    2&\multicolumn{1}{c}{  $\frac{(n+1)(n-3)}{8}$ } & \multicolumn{1}{c}{ $\frac{(n+1)(n-1)}{8}$} &\multicolumn{1}{c}{0} &$\frac{n+1}{4}$&$\frac{(n+1)(n-1)}{4}$   \\ \hline
    Total&\multicolumn{1}{c}{$\frac{n(n-1)}{2}$} & \multicolumn{1}{c}{$\frac{n(n-1)}{2}$} &\multicolumn{1}{c}{$n$} &$n$&$n^2+n $\\  \hline  
    \end{tabular}
}    
    \captionof{table}{\label{tab:classificationBn}}
\end{center}    
\end{minipage}
\newline \newline
Thus, as expected, this grading is balanced.

\section{Gradings in type $C_n$}

The appropriate bipartite matrix for type $C_n$ is the $n \times n$ matrix is given by
\[
B= \scalemath{0.85}{\begin{pmatrix}
0  & 1     &     &      &      &    &  &  \\
-1  &  0 & -1 &   &   &   &  &  \\
  & 1 & 0  &1  &  &    & &   \\
  &   & -1 & 0  &  &    & &   \\
   & &  & &  &\ddots & 1 \\   
&    &        &         &      &-1  & 0  &-2    \\
&    &        &         &      &   &1   &  0   
\end{pmatrix}} \text{ or } 
B= \scalemath{0.85}{\begin{pmatrix}
0  & 1     &     &      &      &    &  &  \\
-1  &  0 & -1 &   &   &   &  &  \\
  & 1 & 0  &1  &  &    & &   \\
  &   & -1 & 0  &  &    & &   \\
   & &  & &  &\ddots & -1 \\   
&    &        &   &  &1  & 0  &2    \\
&    &        &   &  &   &-1   &  0   
\end{pmatrix}}
\]
corresponding to whether $n$ is odd or even, respectively.
If $n$ is even, $C$ has full rank and so again there is no non-trivial grading. If $n$ is odd, $C$ has rank $n-1$ and the grading space is again spanned by a one dimensional vector. There are two cases for the grading vector $\und{g}$:
 \[\und{g}=
 \begin{cases}
 (1,0,-1,0,\dots,0,-1,0,1) &\text{ if } n=4k+1,\\
  (1,0,-1,0,\dots,0,1,0,-1) &\text{ if } n=4k+3.\\
 \end{cases}
 \]
In both case, the distribution of degrees we obtain is the same.

This time the description of the root system is
\[ \{ \pm\beta_i \pm \beta_j \mid i,j=1\dots n , i \neq j \} \cup \{\pm 2 \beta_i \mid i= 1 \dots n \}. \]
Here we have taken
\[\alpha_1=\beta_1-\beta_2 ,\dots, \alpha_{n-1}=\beta_{n-1}-\beta_n, \text{ and } \alpha_n=2\beta_n .\]
We then find that the positive roots are of one of the following forms:
\begin{align}
2\beta_i =&  2(\alpha_i + \cdots +\alpha_{n-1})+\alpha_n& &(i=1,\dots, n),  \label{eqn:CnForm1}\\
\beta_i-\beta_j=&  \alpha_i + \dots +\alpha_{j-1}& &(i<j) ,  \label{eqn:CnForm2}\\
\beta_1+\beta_j =&   \alpha_i + \dots +\alpha_{j-1}+2(\alpha_j+\cdots +\alpha_{n-1})+\alpha_n& &(i=1,\dots,n)  . \label{eqn:CnForm3}
\end{align}
Applying Equation \ref{eqn:deg_in_terms_of_simples} to the above and counting the number of variables of each degree again, we obtain the following distribution:
there are 
\begin{itemize}
\item $\frac{(n+1)^2}{4}$ variables of degree $1$ and and also of degree $-1$, and
\item $\frac{(n+1)(n-1)}{2}$ variables of degree $0$.
\end{itemize}
Again, this grading is balanced as expected.

\biblio

\chapter{Gradings in the rank $3$ case} \label{chp:3_by_3}

The goal of this chapter is to classify the graded cluster algebras generated by $3 \times 3$ skew-symmetric matrices. This classification is in terms of the cardinality of the set of occurring degrees and how the cluster variables are distributed with respect to these degrees. We do not consider any $2 \times 2$ matrices, since they can only give rise to trivial gradings. For this reason also, we will assume that any matrix we consider corresponds to a connected quiver; cluster algebras generated by quivers that are not connected (i.e.~from direct sums of matrices) are easily understood in terms of the cluster algebras generated by the corresponding connected subquivers. Thus, only matrices corresponding to connected quivers can give rise to nontrivial grading behaviours. So, for example, we will not consider matrices such as $\left(\begin{smallmatrix}0&a&0\\-a& 0&0\\0&0&0\end{smallmatrix} \right)$.

\section{Structure of the problem}

A graded cluster algebra can be classified by the following criteria:
\begin{itemize}
\item The cardinality of the set of degrees that occur as degrees of cluster variables may be finite or infinite.
\item Of the degrees that occur, there may be
\begin{enumerate}[(1)]
\item infinitely many cluster variables of each degree,
\item  finitely many variables in each degree (which in the case of finitely many degrees means the cluster algebra is of finite type), or,
\item  a ``mixed" case where some degrees have finitely many variables while others have infinitely many. 
\end{enumerate}
\end{itemize}

We first note that, for the purposes of classifying the graded cluster algebras above, there is an immediate reduction in the number of cases of generating matrices that need to be considered. By the proof of Proposition \ref{prop:esseq_isom_graded_algebras}, we see that essentially equivalent matrices generate isomorphic graded algebras (we simply relabel the variables in the initial cluster), so we need only consider matrices up to essential equivalence. It is also clear that we need only consider matrices up to mutation equivalence. To see this, let $\A_1=\A\big( (x_1,x_2,x_3),A,(g_1,g_2,g_3) \big)$ and $\A_2=\A\big( (x_1',x_2',x_3'),B,(g'_1,g'_2,g'_3) \big)$ and suppose $B=A_{\pth{p}}$ for some mutation path $\pth{p}$. 
Set $(y_1,y_2,y_3)=\cl_A{\pth{p}}$ in $\A_1$. Then, since the initial seed of a cluster algebra can be chosen arbitrarily, $\A_1 = \A\big( (y_1,y_2,y_3 ), B,( g'_1,g'_2,g'_3) \big)$ (as if we mutate this cluster along $[\und{p}^{-1}]$, we get back to $\big((x_1,x_2,x_3),A,(g_1,g_2,g_3) \big)$). But this is clearly isomorphic to $\A_2$ as a graded algebra.

Given the above, it turns out that there is just one class of $3 \times 3$ matrices that needs to be studied. (Note that this does not mean that we will not need to look in detail at any cluster algebras generated by other classes of matrices in the process, rather that once we have classified those corresponding to the class of matrices above, we are done.)

\begin{dfn}
Let $A$ be a matrix. A column of $A$ will be called \emph{sign-coherent} if all of its non-zero entries are of the same sign. Otherwise the column is of \emph{mixed sign}. If $A$ is skew-symmetric, it is called \emph{acyclic} if its associated quiver is acyclic. If $A$ is also $3 \times 3$, it is called \emph{cyclic} if its associated quiver is cyclic. We call $A$ \emph{mutation-cyclic} if every matrix in its mutation class is cyclic and \emph{mutation-acyclic} if there is a matrix in its mutation class that is acyclic.
\end{dfn}

\begin{rmk}
If $A$ is skew-symmetric and $3 \times 3$, it is easy to see that $A$ is acyclic if and only if it has a sign-coherent column: if $Q(A)$ is acyclic then, up to essential equivalence, $A$ is of the form $\ssmat{a}{b}{c}$ or $\left(\begin{smallmatrix}0&-a&0\\a& 0&-b\\0&b&0\end{smallmatrix} \right)$, where $a,b,c \geq 0$. On the other hand, if $A$ has a sign-coherent column then it must also be of one of these two forms, again up to essential equivalence.
\end{rmk}

\begin{prop}
\label{prop:all_mats_esseq_to_this}
Every $3 \times 3$ skew-symmetric matrix is either essentially equivalent to the matrix
$A=\left(\begin{smallmatrix}0&a&-c\\-a& 0&b\\c&-b&0\end{smallmatrix} \right)$, for some $a,b,c \in \N_0$ and $a \geq b \geq c$, or mutation equivalent to a matrix which is essentially equivalent to $A$.
\end{prop}
\begin{proof}
Let $B$ be any given $3 \times 3$ skew-symmetric matrix and write $B=\ssmat{d}{e}{f}$ for some $d,e,f \in \Z$. It is easy to check directly that $B$ is in one of the following essential equivalence classes:
\begin{align*}
& \left\{ 
   \ssmat{-a}{-b}{c}, \ssmat{-a}{-c}{b}, \ssmat{-b}{-a}{c}, \ssmat{-b}{-c}{a}, \ssmat{-c}{-a}{b}, \ssmat{-c}{-b}{a}
\right\}, \\
&\left\{ 
   \ssmat{a}{b}{c}, \ssmat{b}{a}{c}, \ssmat{c}{-a}{b}, \ssmat{c}{-b}{a}, \ssmat{-a}{c}{b}, \ssmat{-b}{c}{a}
\right\},	 \\
&\left\{ 
   \ssmat{a}{c}{b}, \ssmat{c}{a}{b}, \ssmat{b}{-a}{c}, \ssmat{b}{-c}{a}, \ssmat{-a}{b}{c}, \ssmat{-c}{b}{a}
\right\}	, \\
&\left\{ 
   \ssmat{b}{c}{a}, \ssmat{c}{b}{a}, \ssmat{a}{-b}{c}, \ssmat{a}{-c}{b}, \ssmat{-b}{a}{c}, \ssmat{-c}{a}{b}
\right\}	,
\end{align*}
for some $a \geq b \geq c \geq 0$. (Note that these classes can be indexed by which entries appear in columns of mixed sign. Also, strictly, these sets are not the full equivalence classes, only up to sign; however, as noted previously, the sign of an exchange matrix is of no importance.)  
So it is now enough to consider just the representatives $A=\ssmat{-a}{-b}{c}$, $X=\ssmat{a}{b}{c}$, $Y=\ssmat{a}{c}{b}$ and $Z=\ssmat{b}{c}{a}$, say. Mutating $X$, $Y$ and $Z$ in direction 2 we obtain 
$\mmut{X}{2} = \ssmat{-a}{-b}{(ab+c)}$, $\mmut{Y}{2} = \ssmat{-a}{-c}{(ac+b)}$ and $\mmut{Z}{2} = \ssmat{-c}{-b}{(cb+a)}$, and each of these matrices is clearly essentially equivalent to a matrix with the same form as $A$.
\end{proof}
	
\begin{dfn}
\label{dfn:standard_form}
As we see from the proof of Proposition \ref{prop:all_mats_esseq_to_this}, every $3 \times 3$ skew-symmetric matrix $A$ is essentially equivalent to (at least) one of $\ssmat{-a}{-b}{c}$, $\ssmat{-a}{b}{c}$, $\ssmat{a}{-b}{c}$ or $\ssmat{a}{b}{c}$, for some $a\geq b \geq c \geq 0$. We define the \emph{standard form} of $A$, denoted $\sfr{A}$, to be the first element in the above list that $A$ is essentially equivalent to. We will say $A$ is \emph{in standard form} if $A= \sfr{A}$.
\end{dfn}

\begin{rmk}
\label{rmk:only_one_grading_up_to_sign_3by3}
Let $a,b,c \in \Z$ such that we do not have $a=b=c=0$. Up to sign and integer scaling, the degree vector $(g_1,g_2,g_3)$ is uniquely determined in the graded seed 
\[\left(\cluster{x_1,x_2,x_3}, \ssmat{-a}{-b}{c}, (g_1,g_2,g_3)\right),\]
namely $(g_1,g_2,g_3)=(b,c,a)$. Thus the degree vector is uniquely determined up to sign (and scaling, though the scale must match that of the initial degree vector) in every graded seed of any cluster algebra corresponding to any $3 \times 3$ skew-symmetric matrix.
\end{rmk}

A consequence of this is the following lemma. 

\begin{lem} \label{lem:degree_is_in_matrix}
Let $\pth{p}=[p_n,\dots,p_1]$ be a mutation path (without repetition) and let $z=\var_A \pth{p}$ be a cluster variable in $\A\left((x_1,x_2,x_3), \ssmat{-a}{-b}{c}, (b,c,a)\right)$, where $a,b,c \in \Z$. Then $\deg (z)$ is an entry of $A_{\pth{p}}$:
\begin{equation}
  \deg(z)=\begin{cases}
    (-1)^{l(\pth{p})+1}(A_{\pth{p}})_{32}, & p_n=1,\\
    (-1)^{l(\pth{p})}(A_{\pth{p}})_{31}, & p_n=2,\\
    (-1)^{l(\pth{p})+1}(A_{\pth{p}})_{21}, & p_n=3. 
  \end{cases}
  \label{eqn:degree_is_in_matrix}
\end{equation}
\end{lem}

\begin{proof}
Equation \ref{eqn:degree_is_in_matrix} is easily checked directly if $l(\pth{p})=1$. Now suppose it is true for each path of length $r$. Let $\pth{q}$ be such a path. Assume that $r$ is even (the proof is similar if it is odd).
Write $A_{\pth{q}}=\ssmat{-d}{-e}{f}$ for some integers $d$, $e$ and $f$. Then by \pref{eqn:degree_is_in_matrix} the corresponding degree vector is $(e,f,d)$. 
Mutating this matrix and degree vector in directions $1$, $2$, and $3$, we have
\begin{align*}
 &&A_{[1,\und{q}]} &= \ssmat{d}{(df-e)}{-f} \text{, } &\deg.\cl_{A} [1,\und{q}]&=(df-e,f,d) ,&& \\
 && A_{[2,\und{q}]} &= \ssmat{d}{e}{-(de-f)} \text{, }& \deg.\cl_{A} [2,\und{q}]&=(e,de-f,d), && \\                    
\text{and} && A_{[3,\und{q}]} &= \ssmat{(ef-d)}{e}{-f} \text{, } &\deg.\cl_{A} [3,\und{q}]&=(e,f,ef-d). && 
\end{align*}
Thus \pref{eqn:degree_is_in_matrix} is true for all paths of length $r+1$.
\end{proof}

In particular, Lemma \ref{lem:degree_is_in_matrix} implies the following.

\begin{cor} \label{cor:mutation_infinite_3x3_iff_inf_degrees}
A graded cluster algebra generated by a $3 \times 3$ skew-symmetric matrix $A$ has infinitely many degrees if and only if $A$ is mutation-infinite.
\end{cor}

Finally, we note a condition that lets us infer that a grading is balanced.

\begin{thm}[{\cite[Corollary 5.7]{G15}}]
Let $Q$ be a quiver and $\A(Q)$ the cluster algebra generated by $Q$. If $\A(Q)$ admits a cluster categorification then every grading for $\A(Q)$ is balanced.
\end{thm}

We will not cover background on cluster categorification (see \cite{BMMRT}), but we note the following, which is what is relevant to our situation. 

\begin{lem} \label{lem:acyclic_implies_balanced}
If $Q$ is a finite connected acyclic quiver then $\A(Q)$ admits a cluster categorification.
\end{lem}

This fact follows from work of Palu in \cite{Pal}. All of the quivers we deal with will be finite and connected, so whenever we consider a grading arising from an acyclic quiver (or a quiver mutation equivalent to one), we automatically know it is balanced.

The main result of this chapter is the following.

\begin{thm}
\sloppy A partial classification of graded cluster algebras of the form $\mathcal{A}\left(\cluster{x_1,x_2,x_3}, \left(\begin{smallmatrix}0&a&-c\\-a& 0&b\\c&-b&0\end{smallmatrix} \right), (b,a,c)\right)$, with $a,b,c \in \N_0$ and $a \geq b \geq c$, is given in Table \ref{tab:classification3x3}.
\end{thm}

As noted earlier, the restrictions $a,b,c \in \N_0$ and $a \geq b \geq c$ above do not exclude any graded cluster algebras generated by $3 \times 3$ skew-symmetric matrices as far as questions of classification are concerned, by Proposition \ref{prop:all_mats_esseq_to_this}. 

\begin{center}
{
\setlength{\tabcolsep}{7pt} 
\renewcommand{\arraystretch}{1.5} 
\resizebox{\textwidth}{!}{%
\begin{tabular}{ | l | l | l | }
	\hline
	\multicolumn{3}{ |c| }{\bf{{For mutation-finite matrices (which give rise to finitely many degrees) }}} \\ 
	\hline
	\footnotesize{\bf{Finite type}} & \footnotesize{\bf{Mixed}} &\footnotesize{\bf{Infinitely many variables per degree}} \\
	\hline
	\rule[-2.5ex]{0pt}{7ex} $A_3 =\left(\begin{smallmatrix}0&1&0\\-1& 0&1\\0&-1&0\end{smallmatrix} \right)$
	&
	$\left(\begin{smallmatrix}0&2&-1\\-2& 0&1\\1&-1&0\end{smallmatrix} \right)$
	& 
	$\left(\begin{smallmatrix}0&2&-2\\-2& 0&2\\2&-2&0\end{smallmatrix} \right)$ \\ 
	\hline \hline
	\multicolumn{3}{ |c| }{\bf{ { For mutation-infinite matrices (which give rise to infinitely many degrees) }}}\\ 
	\hline
	\footnotesize{\bf{Finitely many variables per degree}}& \footnotesize{\bf{Mixed}} & \footnotesize{\bf{Infinitely many variables per degree}} \\ 
	\hline                                      
	\rule[-2.5ex]{0pt}{7ex} $\left(\begin{smallmatrix}0&a&-c\\-a& 0&b\\c&-b&0\end{smallmatrix} \right)$, mutation-cyclic with $c>2$
	 &
	 $\emptyset$ 
	 &
	$\left(\begin{smallmatrix}0&a&-c\\-a& 0&b\\c&-b&0\end{smallmatrix} \right)$,  mutation-acyclic \\ 
	\hline
\end{tabular}
}
\captionof{table}{\label{tab:classification3x3}}
}
\end{center}	

The entries of this table cover all possible cases except for one: $\left(\begin{smallmatrix}0&a&-2\\-a& 0&a\\2&-a&0\end{smallmatrix} \right)$ with $a \geq 3$, which we will call the \emph{singular cyclic} case. That all other cases are covered follows from results we will see in Section \ref{sec:mut_cyclic_and_acyclic_matrices}: Lemma \ref{lem:c=2_mutation_acyclic_or_cyclic} shows that if $c=2$ then $\left(\begin{smallmatrix}0&a&-c\\-a& 0&b\\c&-b&0\end{smallmatrix} \right)$ is either mutation acyclic or a singular cyclic matrix, and Lemma \ref{lem:c=1_mutation_acyclic} shows that if $c=1$ then $\left(\begin{smallmatrix}0&a&-c\\-a& 0&b\\c&-b&0\end{smallmatrix} \right)$ is mutation-acyclic. We also need \cite[Corollary 2.3]{ABBS}, which says precisely that the three matrices in the top row are the only $3 \times 3$ mutation-finite matrices.

We conjecture that the singular cyclic case should be placed in the lower right cell of the table:

\begin{cnj} \label{cnj:singular_cyclic_infinitely_many_in_each}
For $a \geq 3$, that is, in the singular cyclic case, the graded cluster algebra $\A \left( \und{x}, \left(\begin{smallmatrix}0&a&-2\\-a& 0&a\\2&-a&0\end{smallmatrix} \right), (a,2,a) \right)$ has infinitely many variables in each occurring degree.
\end{cnj}

In the remainder of this chapter we will prove that Table \ref{tab:classification3x3} is correct. Note that $\left(\begin{smallmatrix}0&1&0\\-1& 0&1\\0&-1&0\end{smallmatrix} \right)$ has been covered in Chapter \ref{chp:finite_type}, and $\ssmat{-2}{-2}{2}$ gives rise to the well-known Markov type cluster algebra in which all (infinitely many) variables have degree $2$. We will address what can be said about the singular cyclic case in Section \ref{sec:singular_cyclic}.

\section{Finitely many degrees: the mixed case}
\label{sec:fin_deg_mixed}

We now let
$A=\left( \begin{smallmatrix}0 & 2&-1\\-2 & 0&1 \\1 & -1&0\end{smallmatrix} \right)$ and $G=(1,1,2)$ and consider the associated cluster algebra
$ \mathcal{A}\left( (x_1,x_2,x_3), \ssmat{-2}{-1}{1}, (1,1,2)\right).$ We will show that this graded cluster algebra is of mixed type 
with finitely many degrees. More specifically, there is exactly one variable of degree $2$ and one of degree $-2$, but infinitely many variables of degree $1$ and $-1$ each; these are the only degrees which occur.

\begin{prop} \label{prop:finite_degrees_mixed_finite_degrees}
The only degrees that occur in the graded cluster algebra 
 $\mathcal{A}\left( (x_1,x_2,x_3), \ssmat{-2}{-1}{1}, (1,1,2)\right)$
  are $\pm 1$ and $\pm 2$.
\end{prop}

\begin{proof}
One way to find all occurring degrees is by attempting to compute the exchange tree for $\mathcal{A}$, but closing each branch whenever we obtain a degree seed that is essentially equivalent, up to the sign of the degree cluster, to one that has already occurred. In this case, by Lemma \ref{lem:mut_of_deg_seed_vs_negative}, any subsequent degrees obtained on the particular branch will have already occurred up to sign. This process terminates and is recorded in the diagram below. A black box indicates we have found a degree seed that is essentially equivalent (up to sign) to some previously occurring one. A circle in a degree cluster indicates the newly obtained degree after mutation. We see that the only degrees that can occur are $\pm1$ and $\pm2$. 
\begingroup
\begin{center}
\begin{tikzpicture}[node distance=2.4cm,auto]
  \node (O) {$\big((1,1,2),A\big)$};
  \node (M1) [below of=O, left of=O, left of = O] {$\big( (\circtx{1},1,2) ,-A\big) \blacksquare$};
  \node (M2) [below of=O] {$\big( (1,\circtx{1},2) ,-A\big) \blacksquare$};
  \node (M3) [below of=O, right of = O, right of = O] {$\left( (1,1,\circtx{$-1$}),\left(\begin{smallmatrix}0 &1&1\\-1& 0&-1\\-1&1&0\end{smallmatrix} \right)\right)$};
  \node (M32) [below=7cm of O, right of = O, right of = O] {$\left( (1,\circtx{-1},-1),\left(\begin{smallmatrix}0 &-1&1\\1& 0&1\\-1&-1&0\end{smallmatrix} \right)\right)  \blacksquare$};
  \node (M31) [below=7cm of O, left of = O, left of = O] {$\left(-(\circtx{1},-1,1),\left(\begin{smallmatrix}0 &-1&-1\\1& 0&-1\\1&1&0\end{smallmatrix} \right)\right)  \blacksquare$};  
  \draw[-] (O) to node [swap] {$1$} (M1);
  \draw[-] (O) to node {$2$} (M2);
  \draw[-] (O) to node {$3$} (M3);
  \draw[-] (M3) to node {$2$} (M32);  
  \draw[-] (M3) to node {$1$} (M31);   
\end{tikzpicture}
\end{center}
\captionof{figure}{}
\label{fig:degtree211}
\endgroup
\end{proof}

Notice that, after negating the degree clusters, the degree seeds obtained after the mutation paths $[1,3]$ and $[2,3]$ are essentially equivalent to the one obtained by mutating the initial degree seed in direction $3$. So we can obtain from either of these new degree seeds the negative of any degree we have encountered so far (thus degree $-2$ must occur).
The diagram above also shows there are matrices in the mutation class of $\ssmat{-2}{-1}{1}$ that are sign-coherent. That is, matrices whose corresponding quivers are acyclic. Therefore, by Lemma \ref{lem:acyclic_implies_balanced}, any grading on the cluster algebra we are considering is balanced.

\begin{prop}
 $\mathcal{A}\left( (x_1,x_2,x_3), \ssmat{-2}{-1}{1} \right)$ has infinitely many cluster variables of degree $1$ and degree $-1$.
\end{prop}
\begin{proof}
We will show that $\var_A[(2,1)_r]$ is distinct for all $r$. Denote the $i^\text{th}$ entry of the cluster $\cl_A[(2,1)^j]$ by $x_i^{(j)}$. These cluster variables are rational functions in $x_1$, $x_2$ and $x_3$ and we will evaluate them at $(1,1,1)$.
We claim that $x_1^{(j)}|_{(1,1,1)}$  and $x_2^{(j)}|_{(1,1,1)}$ are strictly increasing in $j$. To simplify notation we will write $f(x_1,x_2,x_3) > g(x_1,x_2,x_3)$ to mean  $f(x_1,x_2,x_3)|_{(1,1,1)} > g(x_1,x_2,x_3)|_{(1,1,1)}$, for rational functions $f$ and $g$. It is easy to show that $\mmut{A}{(2,1)_r}=\pm A$, so we have the following exchange relations for all $j$:
\begin{align}
x_1^{(j+1)}&=\frac{(x_2^{(j)})^2+x_3}{x_1^{(j)}},\\
x_2^{(j+1)}&=\frac{(x_1^{(j+1)})^2+x_3}{x_2^{(j)}}.
\end{align}
The claim will follow by induction. For the base case, note $\cl_A[2,1]|_{(1,1,1)}=(2,5,1)$ and  $\cl_A[(2,1)^2]|_{(1,1,1)}=(13,34,1)$. For the induction hypothesis assume  
\begin{equation} \label{eqn:indassump1} 
x_2^{(j)}>x_1^{(j)}, \end{equation}
(which implies $ (x_1^{(j+1)})/x_1^{(j)}>1$),
and
\begin{equation}\label{eqn:indassump2} 
x_1^{(j+1)}>x_2^{(j)},
\end{equation}
(which implies $(x_2^{(j+1)})/x_2^{(j)}> 1$).
We wish to show (\ref{eqn:indassump1}) and (\ref{eqn:indassump2}) for $j+1$.
First, using (\ref{eqn:indassump2}), we have
\begin{equation}
\label{eqn:indassump1n}
x_2^{(j+1)}=\frac{ (x_1^{(j+1)})^2 +x_3}{x_2^{(j)}}=x_1^{(j+1)} \frac{x_1^{(j+1)}}{x_2^{(j)}}+\frac{x_3}{x_2^{(j)}}>x_1^{(j+1)}.
\end{equation}
Now, using (\ref{eqn:indassump1n}),
\begin{equation}
x_1^{(j+2)}=\frac{ (x_2^{(j+1)})^2 +x_3}{x_1^{(j+1)}}=x_2^{(j+1)} \frac{x_2^{(j+1)}}{x_1^{(j+1)}}+\frac{x_3}{x_1^{(j+1)}}>x_2^{(j+1)}.
\end{equation}
Thus, we have our result for $j+1$. Now since $x_1^{(j)}|_{(1,1,1)}$ and $x_2^{(j)}|_{(1,1,1)}$ are strictly increasing in $j$, $x_1^{(j)}$ and $x_2^{(j)}$ must be distinct for all $j$. 

Now, since $\deg.\var_A[(2,1)_r]=1$ for all $r$, we have infinitely many variables of degree $1$.

Next, note that $\deg.\var_A[3,(2,1)_r]=-1$ for all $r$. To find infinitely many variables of degree $-1$, let $x_3^{(j)} = \var [3,(2,1)^j]$; we will show $x_3^{(j)}|_{(1,1,1)}$ is strictly increasing in $j$. We have $x_3^{(j)}=\frac{x_1^{(j)}+x_2^{(j)}}{x_3}$, so since we have shown that $x_1^{(j)}|_{(1,1,1)}$ and $x_2^{(j)}|_{(1,1,1)}$ are strictly increasing in $j$, this is also the case for $x_3^{(j)}|_{(1,1,1)}$. Thus there are infinitely many variables of degree $-1$.
\end{proof}

In the proof above, we have shown that cluster variables are distinct by considering their numerators. In future, we will usually instead use the more efficient method of considering their denominators via denominator vectors.

To complete the classification for $\mathcal{A}\left( (x_1,x_2,x_3), \ssmat{-2}{-1}{1}, (1,1,2)\right)$ it remains to prove that there is just one variable each in degrees $2$ and $-2$. For this we make use of Example 7.8 of \cite{FZCAI}. Here, the authors consider the cluster algebra generated by the matrix $\ssmat{-1}{-1}{-1}$ and show that the corresponding exchange graph is given by Figure \ref{fig:exgraph211}. In this graph, the initial cluster at vertex $t_0$ is $(y_1,y_2,y_3)$ and the cluster at a general vertex is given (up to permutation) by the three variables in the three regions adjacent to the vertex. 

\begin{prop} \label{prop:finite_degrees_mixed_one_var_deg_2}
$\mathcal{A}\left( (x_1,x_2,x_3), \ssmat{-2}{-1}{1}, (1,1,2)\right)$ has exactly one variable of degree $2$ and one variable of degree $-2$.
\end{prop}
\begin{proof}
\sloppy Since, up to essential equivalence, the degree seed $\left((1,1,2), \ssmat{-2}{-1}{1} \right)$ is mutation equivalent to $\left((1,-1,1), \ssmat{-1}{-1}{-1} \right)$ (which can be inferred from Figure \ref{fig:degtree211}), any results concerning the behaviour of the grading on $\mathcal{A}\left( (x_1,x_2,x_3), \ssmat{-2}{-1}{1}, (1,1,2)\right)$ may alternatively be obtained by considering the grading on $\mathcal{A} \left( (y_1,y_2,y_3), \ssmat{-1}{-1}{-1}, (1,-1,1) \right)$.
Thus it will now be enough to show that the latter graded cluster algebra has just one variable each in degrees $2$ and $-2$. Referring again to the example in \cite{FZCAI}, we have the following exchange relations for the cluster variables in $\mathcal{A} \left( (y_1,y_2,y_3), \ssmat{-1}{-1}{-1}, (1,-1,1) \right)$ for all $m \in \Z$:

\begin{align}
& w y_{2m} =y_{2m-1} + y_{2m+1},    \label{eqn:exrel211a}\\
& y_{2m-1}y_{2m+3} =y^2_{2m+1}+w,   \label{eqn:exrel211b}\\
& y_m y_{m+3} =y_{m+1} y_{m+2}+1,   \label{eqn:exrel211c}\\
& y_{2m-2} y_{2m+2} = y^2_{2m} + z, \label{eqn:exrel211d}\\
& y_{2m-1} z = y_{2m-2} + y_{2m}.   \label{eqn:exrel211e}
\end{align}
From these equations we may recursively calculate the degrees of all the cluster variables in $\mathcal{A} \left( (y_1,y_2,y_3), \ssmat{-1}{-1}{-1}, (1,-1,1) \right)$. To start, we have
$w= \var[2] = \frac{y_1 + y_3}{y_2}$, which has degree $2$. (We know that $w$ is obtained by mutation of the initial cluster in direction $2$ since the corresponding cluster is obtained by replacing the second element of the initial cluster, $y_2$, by $w$.) Similarly we find that the variable $z$ has degree $-2$. Now we just need to check that these are the only such variables. By \pref{eqn:exrel211b} we have $y_1 y_5 = y_3^2+w$, so $y_5$ has degree $1$. Then by \pref{eqn:exrel211a}, $wy_4=y_3 + y_5$, so $y_4$ has degree $-1$. By \pref{eqn:exrel211d}, $y_2 y_6 = y^2_4 + z$, so $y_6$ also has degree $-1$.  Similarly we find $\deg(y_7)=1$, $\deg(y_8)=-1$ and $\deg(y_9)=1$. This is now enough to infer that  $\deg(y_{2m})=-1$ and $\deg(y_{2m+1})=1$ for all $m$. Therefore $w$ is the only cluster variable of degree $2$ and $z$ the only variable of degree $-2$.
\end{proof}

\begin{cor}
$\mathcal{A}\left( (x_1,x_2,x_3), \ssmat{-2}{-1}{1}, (1,1,2)\right)$ is of mixed type. The occurring degrees are $\pm 1 $ and $\pm 2$. The graded cluster algebra has one variable each in degrees $2$ and $-2$ and infinitely many variables each in degrees $1$ and $-1$. The grading is also balanced.
\end{cor}

\begin{figure}[H] 
\setlength{\unitlength}{2pt}
\begin{center} 
\begin{picture}(140,55)(0,0)
\multiput(0,10)(0,20){3}{\line(1,0){140}}
\multiput(10,30)(40,0){4}{\line(0,1){20}}
\multiput(30,10)(40,0){3}{\line(0,1){20}}

\multiput(10,30)(40,0){4}{\circle*{2}}
\multiput(10,50)(40,0){4}{\circle*{2}}
\multiput(30,30)(40,0){3}{\circle*{2}}
\multiput(30,10)(40,0){3}{\circle*{2}}

\put(10,20){\makebox(0,0){$y_0$}}
\put(50,20){\makebox(0,0){$y_2$}}
\put(90,20){\makebox(0,0){$y_4$}}
\put(130,20){\makebox(0,0){$y_6$}}

\put(30,40){\makebox(0,0){$y_1$}}
\put(70,40){\makebox(0,0){$y_3$}}
\put(110,40){\makebox(0,0){$y_5$}}

\put(53,33){\makebox(0,0){$t_0$}}
\put(70,2){\makebox(0,0){$z$}}
\put(70,58){\makebox(0,0){$w$}}

\end{picture} 
\end{center} 
\caption{\cite[Figure 4.]{FZCAI}} 
\label{fig:exgraph211} 
\end{figure} 

\section{Mutation-cyclic and mutation-acyclic matrices} \label{sec:mut_cyclic_and_acyclic_matrices}

According to Table \ref{tab:classification3x3}, we will need to be able to determine whether a given $3 \times 3$ matrix is mutation-cyclic or not. In this section we give a simple algorithm that allows us to do this. We also prove that (excluding the singular cyclic case) mutation-cyclic matrices give rise to graded cluster algebras with finitely many variables in each degree.

At this stage we should mention \cite{ABBS}, which is concerned with parametrising and testing for $3 \times 3$ mutation-cyclic matrices, and covers some of the same material that we are about to. While the results in this section were not derived from \cite{ABBS}, some can be. In particular, Algorithm \ref{alg:suff_finitely_each_deg} is very similar in spirit to the algorithm in Section 2.5 of \cite{ABBS}, which detects whether an input matrix is mutation-cyclic. However, in general we treat the material from a slightly different perspective. 

Algorithm \ref{alg:suff_finitely_each_deg} is motivated by the observation that, when mutating a matrix $A=\ssmat{-a}{-b}{c}$, $a\geq b \geq c \geq 2$ (and then rearranging to standard form), there is only one chance to obtain a matrix of the same form but with smaller absolute values for the new values corresponding to the $b$ or $c$ entries: if this does not happen when we mutate in direction $3$, it can never happen, and therefore $A$ must be mutation-cyclic. 
For this reason, we use the slightly weaker condition of $bc-a \geq b$ (rather than $bc-a \geq a$, as used in \cite{ABBS}) to detect whether a matrix is mutation-cyclic during our algorithm. This allows the possibility of detecting the property a step earlier but means we may have to perform an additional mutation if we want to obtain a \emph{minimal} mutation-cyclic matrix. 
In addition to detecting mutation-cyclicity, Algorithm \ref{alg:suff_finitely_each_deg} gives an output that is a standard representative for the mutation class of $A$. (The algorithm in \cite{ABBS} gives the same output if $A$ is mutation-cyclic, but does not necessarily otherwise.) This lets us determine whether two matrices are mutation equivalent and tells us where in Table \ref{tab:classification3x3} an input matrix should be placed. In particular, for any input matrix, it tells us the classification of the corresponding graded cluster algebra (modulo the conjecture about the singular cyclic cases). If the matrix is mutation-acyclic, our output also tells us which case (i.e.~which of \pref{case:acyclic1}--\pref{case:acyclic5}---see the next section) its mutation class falls into.

The idea at the heart of both algorithms is very simple: if the input matrix $A$ is not a fork (defined in the next section), it is not mutation-cyclic, and if it is a fork, we repeatedly mutate at the point of return until we get either a ``minimal fork", which means $A$ is mutation-cyclic, or a matrix that is not a fork, which means $A$ is mutation-acyclic.

Finally, if the mutation class of $A$ does not need to be determined, \cite{BBH} gives a condition that allows us to read off whether $A$ is mutation-cyclic or not from its entries.

\begin{thm}[{\cite[Theorem 1.2]{BBH}}] \label{thm:BBH}
Let $A=\ssmat{-a}{-b}{c}$ with $a,b,c \in \N_0$. $A$ is mutation-cyclic if and only if $a,b,c \geq 2$ and $C(a,b,c) \leq 4$, where
\begin{equation}
C(a,b,c):=a^2+b^2+c^2-acb
\end{equation}
is the \emph{Markov constant} of the triple $(a,b,c)$.
\end{thm}

\begin{rmk}
Note the following fact about mutation of a degree seed of the form
\begin{equation}
\left((b,c,a),\left( \begin{smallmatrix}0 & a&-c\\-a & 0&b \\c & -b&0\end{smallmatrix} \right) \right), \label{eqn:degseedformbca}
\end{equation} 
with $a \geq b \geq c \geq 2$.
We have 
\begin{align}\deg.\sd_A[1] &= 
\left((ca-b,c,a), \left(\begin{smallmatrix}0 & -a&c\\a & 0&-(ca-b) \\-c & ca-b&0\end{smallmatrix}\right)
\right) 
\\ 
&\esseq \left((a,c,ca-b), \left(\begin{smallmatrix}0 & ca-b&-c\\-(ca-b) & 0&a \\c & -a &0\end{smallmatrix} \right)
\right), \label{eqn:equivdegcluster1}
\end{align}
where $ca-b \geq a \geq c$, and 
\begin{align}
\deg.\sd_A[2]& = 
\left((b,ab-c,a), \left(\begin{smallmatrix}0 & -a&ab-c\\a & 0&-b \\-(ab-c) & b&0\end{smallmatrix}\right )
\right) \\
&\esseq \left((a,b,ab-c), \left(\begin{smallmatrix}0 & ab-c&-b\\-(ab-c) & 0&a \\b & -a&0\end{smallmatrix}  \right)
\right),
\label{eqn:equivdegcluster2}
\end{align}
where $ab-c \geq a \geq b$. 
So mutation of
$\left((b,c,a),\left( \begin{smallmatrix}0 & a&-c\\-a & 0&b \\c & -b&0\end{smallmatrix} \right) \right)$
in directions $1$ and $2$ yields a degree seed essentially equivalent to one of the same form and a degree at least as large as $a$ (and strictly larger if $c \geq 3$).
However, we have 
\begin{align}
\deg.\sd_A[3] = \left((c,b,bc-a),  -\ssmat{-(bc-a)}{-b}{c}\right),
\label{eqn:equivdegcluster3}
\end{align}
which may not be equivalent to a degree seed of the same form (if $bc-a \leq 1$).
\end{rmk}

\begin{alg}
\label{alg:suff_finitely_each_deg}
Take as input any $3 \times 3$ skew-symmetric matrix $A$. Assume $A$ is in standard form; if it is not, perform the algorithm on $\sfr{A}$ instead. Write
 $A=\ssmat{-a}{-b}{c}$ with $a,b,c \in \Z$ such that $|a| \geq |b| \geq |c| \geq 0$.
\begin{enumerate}[1.]
\item Set the counter $i$ to the value 1. Let $A_1=A$. If $a \geq b \geq c \geq 2$ and $bc-a \geq b$ then say that $A$ has passed the algorithm and go to Step 4. (Note if $c=2$ that $bc-a \geq b$ means $a = b$.) Next, if any of the columns of $A$ are sign-coherent or if $c \leq 2$, say that $A$ has failed and go to Step 4.
\item If there exists a matrix essentially equivalent to $(A_i)_{[3]}$ which is of the form $\ssmat{-d}{-e}{f}$, where $d \geq e \geq f \geq 2$, then let ${A_{i+1}}$ be such a matrix. If not, then say that $A$ has failed and go to Step 4.
Now, if $f \geq 3$ and $ef-d \geq e$, then say that $A$ has passed and go to Step 4. 
If $f=2$, say that $A$ has passed if $d=e$ and failed if  $d>e$, and go to Step 4.
\item Increment the counter $i$ by $1$ and repeat Step 2. (Note that we reach this step if $f  \geq 3$ but $ef-d < e$.)
\item Define the output of the algorithm, which we will denote $m(A)$, as follows. 
If $A$ passed the algorithm and $ef-d \geq d$, $m(A)$ is the matrix $A_i$ computed most recently.  
If $A$ passed and $ef-d < d$, $m(A):=(A_i)_{[3]}$. 
Otherwise, continue computing $A_i$ for higher values of $i$ in the same way as above until we find $k$ such that $A_k$ is acyclic (we will show later that this must eventually happen---see the proofs of Lemma \ref{lem:c=1_mutation_acyclic} and Lemma \ref{lem:c=2_mutation_acyclic_or_cyclic}). Then $m(A): = A_{k}$.
\end{enumerate}
\end{alg}

Let us first check that this algorithm terminates after finitely many steps. As noted, once the algorithm reaches Step 4, we can conclude it will terminate, so we just need to check that $A$ passes or fails in a finite number of iterations. Suppose it does not pass or fail in Step 1 and does not fail in Step 2 when $i=1$. Consider $A_2:=\ssmat{-d}{-e}{f}$. Since $A_1$ did not satisfy $bc-a \geq b$, we must have $d =b$ and either $e <b$ or $f<c$ (i.e.\ mutation in direction $3$ has produced a degree strictly smaller than $b$). So each of the subdiagonal entries of $A_2$ is no larger in absolute value than the corresponding entry of $A_1$, and at least one is strictly smaller. Now $A_2$ will either pass the algorithm or $A_3$ must have some strictly smaller entry than $A_2$ (provided $A_3$ exists---if not, $A$ fails). Continuing in this way, either $A$ must pass after computing (or showing there is no) $A_i$ for some $i$, or we must eventually obtain a matrix in which at least one entry strictly smaller than than 3, at which point we can immediately determine whether $A$ passes or fails. 

\begin{dfn}
We will call the output $m(A)$ of Algorithm \ref{alg:suff_finitely_each_deg} a \emph{minimal mutation representative} of $A$, denoted $\mrep{A}$. 
\end{dfn}

\begin{rmk}
Algorithm \ref{alg:suff_finitely_each_deg} allows us to determine if any two $3 \times 3$ skew-symmetric matrices are in the same mutation class as follows. We will see below that $\mrep{A}$ is unique up to essential equivalence if $A$ passes the algorithm. If it fails the algorithm, there may be up to three minimal representatives, up to equivalence. (However, all will share the same set of entries.) In this case, it easy to work out all the representatives of $\mrep{A}$ (c.f.\ Equations \ref{eqn:mut_special_acyclic_path_along_1}--\ref{eqn:mut_special_acyclic_path_along_321} in the next section). 
\end{rmk}

Algorithm \ref{alg:suff_finitely_each_deg} tests whether $A$ is mutation-cyclic: we will show below that $A$ fails the algorithm if and only if it is mutation-acyclic.

\begin{dfn}
Let $A=\pm \ssmat{-a}{-b}{c}$ be a $3 \times 3$ skew-symmetric matrix. We will call $A$ \emph{$c$-cyclic} if $a \geq b \geq c \geq 1$ and \emph{cyclic-preserving} if it is $c$-cyclic with $bc-a \geq b$. We will also call a matrix $B$ \emph{essentially $c$-cyclic} or \emph{essentially cyclic-preserving} if $\sfr{B}$ is $c$-cyclic or cyclic-preserving, respectively.
\end{dfn}

Note $A$ cannot be cyclic-preserving if $c=1$, and if $c=2$, $A$ is cyclic-preserving if and only if $a=b$.

\begin{lem}
\label{lem:minrep_is_c_cyclic_on_success}
\mbox{}
\begin{enumerate}[(i)]
\item If $\sfr{A}$ passes Algorithm \ref{alg:suff_finitely_each_deg} then we may write $\mrep{A}:=\ssmat{-a}{-b}{c}$ for some $a \geq b \geq c \geq 2$ and we have that $\mrep{A}$ is cyclic-preserving.
\item If $A$ is cyclic-preserving then either $A=\mrep{A}$ or $A_{[3]}=\mrep{A}$ and if it is essentially cyclic-preserving then either $\sfr{A}=\mrep{A}$ or $\sfr{A}_{[3]}=\mrep{A}$. 
\end{enumerate}
\end{lem}

\begin{proof}
This is immediate considering any of the possible forms of $A_n$ that are required in order that $A$ passes the algorithm. 
\end{proof}

\begin{lem}
\label{lem:standard_form_direction_3}
Let $A$ be $c$-cyclic where $c \geq 2$ and let $i \in \{1,2\}$. Let $\sigma$ be such that $A_{[i]} \esseq \sigma(\sfr{(A_{[i]})})$. If $c \geq 3$ or $c=2$ and $a>b$, then $\sigma$ is uniquely determined and $\sigma(i)=3$. If $c=2$ and $a=b$, then $\sigma$ can be chosen such that $\sigma(i)=3$.
\end{lem}

\begin{proof}
This follows by considering \pref{eqn:equivdegcluster1}--\pref{eqn:equivdegcluster3} at the beginning of this subsection.
\end{proof}

\begin{lem}
\label{lem:c_cyclic_preserved_in_dir_1_2}
If $A$ is $c$-cyclic with $c \geq 2$ then $\sfr{(A_{[1]})}$ and $\sfr{(A_{[2]})}$ are $d$-cyclic, where $d \geq c$. If $c>2$ or $a>b$ (i.e.\ if $A$ is not the singular cyclic case) then $d>c$.
\end{lem}

\begin{proof}
This follows from the forms of \pref{eqn:equivdegcluster1}--\pref{eqn:equivdegcluster2}.
\end{proof}

\begin{cor}
\label{cor:c_cyclic_preserved_in_paths_starting_1_2}
Suppose $A$ is $c$-cyclic with $c \geq 2$. Let $\pth{p}$ be a mutation path without repetitions such that $p_1=1$ or $p_1=2$ and suppose $[\und{p}']$ is either a rooted proper subpath of $\pth{p}$ or the empty path. Then $\sfr{(A_{\pth{p}})}$ is $d$-cyclic, where $d \geq c$, and $\sfr{( A_{[\und{p}']} )}$ is $d'$-cyclic with $d \geq d' \geq c$. 

Suppose further that $c>2$ or $a>b$ (where $A=\ssmat{-a}{-b}{c}$).
Then $\sfr{( A_{[\und{p}']} )}$ is $d'$-cyclic with $d > d' \geq c$.
\end{cor}

\begin{proof}
Assume without loss that $p_1=1$. By Lemma \ref{lem:c_cyclic_preserved_in_dir_1_2}, $\sfr{(A_{[1]})}$ is $d$-cyclic with $d \geq c$. Let $\delta$ be such that $A_{[1]}=\delta(\sfr{(A_{[1]})})$ and $\delta(p_1)=3$; such a permutation exists by Lemma \ref{lem:standard_form_direction_3}. This serves as a base case.

Let $[\und{p(r)}]$ denote the rooted subpath of $p$ of length $r$. 
Assume for induction there exists an $n$-cyclic matrix $C$ which is essentially equivalent to  $A_{[\und{p(r)}]}$ and write $A_{[\und{p(r)}]} = \rho(C)$. Assume also that $\rho$ is such that $\rho(p_r)=3$. Now mutating in direction $p_{r+1}$ we obtain, 
\begin{equation}
A_{[\und{p(r+1)}]} = \left(\rho(C) \right)_{[p_{r+1}]} = \rho(C_{[\rho( p_{r+1} )]}),
\end{equation}
using Part (a) of Lemma \ref{lem:esseq} for the second equality.
Now since we know $\rho( p_{r+1} ) \neq 3$ (otherwise $\pth{p}$ must contain repetitions of direction), $\sfr{( C_{[\rho( p_{r+1} )]} )}$ is $m$-cyclic, where $m \geq n$, by Lemma \ref{lem:c_cyclic_preserved_in_dir_1_2}. But $\sfr{( C_{[\rho( p_{r+1} )]} )}$ is essentially equivalent to $C_{[\rho( p_{r+1} )]}$, say $C_{[\rho( p_{r+1} )]}= \sigma \left( \sfr{( C_{[\rho( p_{r+1} )]} )} \right)$. 
So we have 
\begin{equation}
A_{[\und{p(r+1)}]} = \rho \left( \sigma\left( \sfr{( C_{[\rho( p_{r+1} )]} )} \right) \right) = \rho \left( \sigma\left( C' \right) \right),
\end{equation}
where $C':=\sfr{( C_{[\rho( p_{r+1} )]} )}$. 
This shows that $\sfr{(   A_{[\und{p(r+1)}]}  )}$ is $m$-cyclic.
By Lemma \ref{lem:esseq_transitivity}, we may write $A_{[\und{p(r+1)}]} = (\sigma \rho) (C')$.
Finally, since $C'$ is $m$-cyclic, we know by Lemma \ref{lem:standard_form_direction_3} that ($\sigma$ can be chosen such that) $\sigma\left( \rho \left( p_{r+1} \right) \right) = (\sigma \rho)(p_{r+1})=3$. 

For the second part we use Lemma \ref{lem:c_cyclic_preserved_in_dir_1_2} again: if $c>2$, we may replace the inequalities involving $c,d,n$ and $m$ above by strict inequalities. We may do the same also if $a>b$, since, after mutation, this property is preserved for the corresponding entries of the mutated matrix (once re-arranged into standard form).
\end{proof}

A cyclic-preserving matrix is no more than one mutation away from being a ``minimal" mutation-cyclic matrix, in the following sense.

\begin{prop}
\label{prop:c_cyclic_form_preserving}
Suppose $A=\ssmat{-a}{-b}{c}$ is $c$-cyclic and cyclic-preserving. Let $\pth{p}$ be any mutation path without repetitions and suppose $[\und{p}']$ is either a rooted proper subpath of $\pth{p}$ or the empty path. Then $\sfr{(A_{\pth{p}})}$ is $d$-cyclic with $d \geq c$, and $\sfr{( A_{[\und{p}']} )}$ is $d'$-cyclic with $d \geq d' \geq c$. In particular, $A$ is mutation-cyclic.

Suppose further that $l(p) \geq 2$ or $p_1 \neq 3$. If we also have $c>2$ or $a > b$, then $\sfr{( A_{[\und{p}']} )}$ is $d'$-cyclic with $d > d' \geq c$.
\end{prop}

\begin{proof}
If $\pth{p}$ starts with direction $1$ or $2$ then the result holds by Corollary \ref{cor:c_cyclic_preserved_in_paths_starting_1_2}.
Suppose $p_1=3$. We have $A_{[3]} = \ssmat{(bc-a)}{b}{-c}$.
But $bc-a \geq b$ as $A$ is cyclic-preserving, so $A_{[3]}$ is $c$-cyclic. We now need to show that $\sfr{( (A_{[3]})_{[\und{q}']} )}$ is $d$-cyclic, where $[\und{q}']$ is the subpath obtained by deleting the first mutation in $\pth{p}$ (i.e.~its rightmost entry, $p_1$). This is indeed the case, by Corollary \ref{cor:c_cyclic_preserved_in_paths_starting_1_2}, since the first mutation in $[\und{q}']$ is either direction $1$ or $2$.

For the strict inequality, suppose $c>2$ or $a > b$. If $p_1=1$ or $2$ then the result holds by Corollary \ref{cor:c_cyclic_preserved_in_paths_starting_1_2} again. Assume $p_1=3$ and $l(\pth{p}) \geq 2$. We have that the inequality corresponding to $bc-a>b$ holds in $\sfr{(A_{[3]})}$. Thus, to obtain $\sfr{(A_{\pth{p}})}$, we are now mutating a $c$-cyclic matrix of the form $\ssmat{-d}{-a}{c}$ with $d>a$ (or $c>2$) along a path that does not start with $3$. Therefore by Corollary \ref{cor:c_cyclic_preserved_in_paths_starting_1_2}, $\sfr{(A_{\pth{p}})}$ is $d$-cyclic with $d > c$.
\end{proof}

We list a few results which tell us whether certain matrices are mutation-cyclic or mutation-acyclic.

\begin{lem}
\label{lem:c=3_mutation_cyclic}
Let $A=\ssmat{-a}{-b}{c}$ with $a \geq b \geq c \geq 3$ and $bc-a \geq b$. Then $A$ is mutation-cyclic.
\end{lem}

\begin{proof}
Since $A$ is $c$-cyclic and cyclic-preserving, this holds by Proposition \ref{prop:c_cyclic_form_preserving}.
\end{proof}

\begin{lem}
\label{lem:c=1_mutation_acyclic}
Let $A=\ssmat{-a}{-b}{1}$ with $a \geq b \geq 1$. Then $A$ is mutation-acyclic.
\end{lem}

\begin{proof}
We have $A_{[3]} = \ssmat{(b-a)}{b}{-c}$ which is acyclic since $b-a \leq 0 $.
\end{proof}

\begin{lem}
\label{lem:c=2_mutation_acyclic_or_cyclic}
Let $A=\ssmat{-a}{-b}{2}$ with $a \geq b \geq 2$. Then $A$ is mutation-cyclic if and only if $a=b$.
\end{lem}
 
\begin{proof}
If $a=b$ then $A$ is $2$-cyclic and cyclic-preserving. Then the result holds by Proposition \ref{prop:c_cyclic_form_preserving}.

Suppose $a>b$. In this case, we can infer that $A$ is mutation-acyclic using Theorem \ref{thm:BBH}. However, we will give an alternative proof that shows how to directly mutate $A$ to obtain an acyclic matrix (which we wish to be able to do for part 4.~of Algorithm \ref{alg:suff_finitely_each_deg}). We have $A_{[3]}=\ssmat{(2b-a)}{b}{-2}$. Three cases can occur. If $2b-a \leq 0$, then $A_{[3]}$ is acyclic as it has a (weakly) sign-coherent column. If $2b-a=1$, then $A_{[3]} \sim \ssmat{b}{2}{-1}$, which is mutation-acyclic by Lemma \ref{lem:c=1_mutation_acyclic}. If $2b-a \geq 2$ then 
\[\sfr{(A_{[3]})} = \ssmat{-b}{-(2b-a)}{2} = :\ssmat{-a^{(1)}}{-b^{(1)}}{2}=:A^{(1)},\] 
where 
$a^{(1)} > b^{(1)}$. We have strict inequality here since $a \neq b$. This matrix has the same form as $A$, but $a^{(1)}<a$ and $b^{(1)}<b$. Next consider 
$(A^{(1)})_{[3]} = \ssmat{-b^{(1)}}{-\left(2b^{(1)}-a^{(1)}\right)}{2}$.
This matrix is either immediately seen to be mutation-acyclic, in the same way as above, or $2b^{(1)}-a^{(1)} \geq 2$. Then 
\[ \sfr{((A^{(1)})_{[3]})} = \ssmat{b^{(1)}}{\left(2b^{(1)}-a^{(1)}\right)}{-2}=:-\ssmat{-a^{(2)}}{-b^{(2)}}{2},\]
where $a^{(2)}> b^{(2)}$. Continuing this process, we must eventually have $2b^{(n)}-a^{(n)} \leq 1$ for some $n$, at which point we may conclude that $A$ is mutation-acyclic.
\end{proof}

\begin{prop}
\label{prop:A_passes_alg_iff_mutation_cyclic}
An input matrix $A$ passes Algorithm \ref{alg:suff_finitely_each_deg} if and only if it is mutation-cyclic.
\end{prop}

\begin{proof}
Suppose  $A$ passes the algorithm. Then $A$ is essentially mutation equivalent to $\ssmat{-d}{-e}{f}$ where either $d \geq e \geq f \geq 3$ and $ef-d \geq e$, or $f=2$ and $d=c$. In the first case $A$ must then be mutation-cyclic by Lemma \ref{lem:c=3_mutation_cyclic} and in the second it must be mutation-cyclic by Lemma \ref{lem:c=2_mutation_acyclic_or_cyclic}.

Suppose $A$ fails. Then $A$ is essentially mutation equivalent to a matrix which either has some sign-coherent column or is of the form $\ssmat{-d}{-e}{2}$ where $d>e\geq 2$. In the first case $A$ is essentially mutation equivalent to an acyclic matrix and in the second it must be mutation-acyclic by Lemma \ref{lem:c=2_mutation_acyclic_or_cyclic}.
\end{proof}

Now we apply the above to graded cluster algebras generated by mutation-cyclic matrices.

\begin{thm}
Let $A$ be mutation-cyclic. Let $\mrep{A}=\ssmat{-a}{-b}{c}$ and suppose $c >2$ (i.e.~that $\mrep{A}$ is not a singular cyclic matrix). Then the graded cluster algebra arising from $A$ has finitely many variables in each occurring degree. We also have that the grading is not balanced.
\end{thm}

\begin{proof}
Since $\mrep{A}$ is essentially mutation equivalent to $A$, it is enough to consider $\A\big((x_1,x_2,x_3),\mrep{A},(b,c,a)\big)$. Since $A$ is mutation-cyclic, it passes Algorithm \ref{alg:suff_finitely_each_deg} by Proposition \ref{prop:A_passes_alg_iff_mutation_cyclic}. Then $\mrep{A}$ is $c$-cyclic-preserving by Lemma \ref{lem:minrep_is_c_cyclic_on_success}. 
Let $\pth{p}$ be a mutation path of length $n\geq 1$. Then, as $c >2$, the second part of Proposition \ref{prop:A_passes_alg_iff_mutation_cyclic} implies we must have $\deg_{\mrep{A}} \pth{p} > n$, since $\deg_{\mrep{A}} \pth{p}$ is bounded below by the smallest entry of $\mrep{A}_{\pth{p}}$. Thus, there are only finitely many possible mutation paths that can result in any given degree, and so there are only finitely many possible cluster variables of that degree.

The grading is clearly not balanced since all occurring degrees are positive.
\end{proof}

\section{Cluster algebras generated by mutation-infinite acyclic matrices}

Gradings on cluster algebras generated by mutation-infinite $3 \times 3$ skew-symmetric acyclic matrices all have the property that there are infinitely many variables in each occurring degree. A reason why one may expect this behaviour is the following: cluster algebras associated with such acyclic matrices have a special mutation sequence which gives rise to infinitely many seeds that are not equal (or essentially equivalent) to the initial seed, but whose corresponding degree seeds are. Given any degree, there exists a mutation path applied to the initial seed which results in a cluster variable of this degree. But then (accounting for essential equivalence) the same mutation path applied to each of these different seeds gives us infinitely many different ways of obtaining variables of this degree. The purpose of this section is to prove that the variables we obtain in this way are in fact all distinct.

This turns out to be nontrivial, and we rely on a property of the exchange graphs of acyclic cluster algebras in order to prove our result. All such exchange graphs have been identified by Warkentin in Chapter 10 of \cite{War}. The following theorem is the crucial property we need for our result. It appears in \cite{GSV} as Conjecture 3.47, but is proved in our setting, for acyclic matrices, in \cite[Corollary 3]{CalderoKeller}.

\begin{thm} \label{thm:exchange_graph_connected}
Let $\A$ be a cluster algebra arising from a mutation-acyclic matrix. For any cluster variable $x$ in $\A$, the seeds whose clusters contain $x$ form a connected subgraph of the exchange graph. 
\end{thm}

We will show that certain variables are distinct by showing that they define different subgraphs, or \emph{regions}, of the exchange graph.

\begin{dfn}
Let $x$ be a cluster variable in a cluster algebra $\A$. Call the subgraph of $\Gamma(\A)$ formed by the set of seeds whose clusters contain $x$ the \emph{region} defined by $x$, denoted by $R(x)$.
\end{dfn}

By Theorem \ref{thm:exchange_graph_connected}, if $A$ is a (mutation-)acyclic matrix and $\A=\A(\und{x}, A)$, then $R(x)$ is connected for each cluster variable $x$ in $\A$.

It is clear that the regions of $\A$ are in bijection with the cluster variables. We also have that, in acyclic cluster algebras of rank $3$, a region must either be a loop (or polygon) or an infinite line of vertices. We may see this as follows. If a region $R(x)$ is not a loop then let $t_1$ be a vertex of $R(x)$ and let $L$ be the list of vertices encountered so far, initially containing only $t_1$. Exactly two of the vertices connected to $t_1$, say $t_2$ and $t_3$, are also in $R(x)$, since mutation in exactly two out of three possible directions does not replace $x$. Now add $t_2$ and $t_3$ to the list $L$. Then exactly two of the neighbours of $t_2$ and $t_3$ must also be in $R(x)$, and since $R(x)$ is not a loop, these new vertices must not already be in $L$. (To be precise, this is only true if $R(x)$ contains no loops, rather than is not a loop, but it is clear that $R(x)$ must be a loop if it contains one.) Continuing in this way, we see that $R(x)$ is an infinite line.

Let us fix a representative for our acyclic matrices, to which we may restrict our attention. We will assume any such matrix is of the form 
\begin{equation}
A=\ssmat{a}{b}{c},
\end{equation}
 with $a \geq b \geq c \geq 0$. We may assume $|a| \geq |b| \geq |c|$ because of essential equivalence; if $A$ does not satisfy this to begin with, we can always permute its entries to get an essentially equivalent matrix for which this will be the case. Now note the mutation equivalence
 \begin{flalign}
&&A_{[1]}=\ssmat{-a}{b}{-c}, &&  \label{eqn:mut_special_acyclic_path_along_1} \\ 
&& A_{[2,1]}=\ssmat{a}{-b}{-c},\label{eqn:mut_special_acyclic_path_along_21} &&\\
  \text{and }&& A_{[3,2,1]}=A, \label{eqn:mut_special_acyclic_path_along_321} &&
\end{flalign} 
which shows that if any of $a,b,c$ are negative, we can mutate $A$ to make them all positive (or at least all of the same sign, in which case we can then negate the whole matrix). Thus we can also assume that $a,b$ and $c$ are non-negative. This choice of representative is the same as that in \cite{War}.

The exchange graphs we will consider are made up of two main elements: a series of infinite $2$-trees and a ``floor" of varying complexity, to which each $2$-tree is attached by another edge. We wish to be able to formally distinguish between these two components of the graph. The following definition helps with this. It is made in Chapter 2 of \cite{War}.

\begin{dfn}
Let $Q$ be the quiver corresponding to an $n \times n$ skew-symmetric matrix. $Q$, or its corresponding matrix, is called a \emph{fork} if it is not acyclic, $|q_{jk}| \geq 2$ for all $j \neq k$, and there is a vertex $r$ called the \emph{point of return} such that
\begin{enumerate}
\item For all $i \in Q^-(r)$ and $j \in Q^+(r)$ we have $q_{ji}> q_{ir}$ and $q_{ji}>  q_{rj}$, and
\item $Q^-(r)$ and $Q^+(r)$ are acylic,
\end{enumerate}
where $Q^+(r)$ (resp. $Q^-(r)$) is the full subquiver of $Q$ induced by all vertices with an incoming arrow from $r$ (resp. outgoing arrow to $r$).
\end{dfn}

The property of forks that is relevant to us is that mutation of the corresponding matrix in any direction excluding that of the point of return gives rise to another ``larger" fork. In the exchange graphs we are considering, forks are found at the base of each of the infinite $2$-trees (see \cite{War}, Lemma 2.8) that make up one component of a graph.

\begin{dfn}
We will say a vertex of an exchange graph $\Gamma$ is in the \emph{canopy} of $\Gamma$ if the associated matrix is a fork, otherwise we will say it is in the \emph{floor}. Given a vertex $t$ in the canopy of $\Gamma$, there is a unique path of shortest length connecting $t$ to a vertex $t'$ which is in the floor of $\Gamma$. We will call the infinite $2$-tree of which $t$ is a vertex the \emph{branch} of $\Gamma$ associated with $t'$.
\end{dfn}

If $(\und{x},A)$ is a seed where $A$ is a fork with point of return $r$, then we know that part of the exchange graph looks as in the following figure, where the dotted (but not dashed) lines represent an infinite continuation of the $2$-tree.

\begin{figure}[H]
\begin{center}
\begin{tikzpicture}[level distance=8mm, node distance=2cm, auto, font=\footnotesize] 
  \tikzstyle{level 1}=[sibling distance=0mm]
  \tikzstyle{level 2}=[sibling distance=10mm]
  \tikzstyle{level 3}=[sibling distance=8mm]
  \node (O) {}[grow'=up];
  \node (O1) [right of = O,solid node, fill=black, label=below:{}]{}[grow'=up]
    child 
    { node [ solid node, fill=black, label=right:{$(\und{x}, A)$}] {}
      child 
      {  node [ solid node, fill=black, label=right:{}] {}
         child[dotted] {}
         child[dotted] {}
         edge from parent node[left, xshift = 0]{}
      }
      child 
      { node [ solid node, fill=black, label=right:{}] {}
         child[dotted] {}
         child[dotted] {} 
         edge from parent node[right, xshift = 0]{}   
      }
      edge from parent node[left, xshift = 0]{$r$}
    };    
  \node (O2) [right of = O1, label=below:{}]{}[grow'=up];          
\draw[dashed] (O) to node  {} (O1);      
\draw[dashed] (O1) to node {} (O2); 
\end{tikzpicture}
\end{center}
\captionof{figure}{}
\label{fig:fork}
\end{figure}

We will depict such a seed corresponding to a fork along with the associated infinite $2$-tree by a rectangle as in the following.

\begin{figure}[H]
\begin{center}
\begin{tikzpicture}[level distance=8mm, node distance=1cm, auto, font=\footnotesize] 
  \tikzstyle{level 1}=[sibling distance=0mm]
  \tikzstyle{level 2}=[sibling distance=10mm]
  \tikzstyle{level 3}=[sibling distance=8mm]
  \node (O) {}[grow'=up];
  \node (O1) [right of = O,solid node, fill=black, label=below:{}]{}[grow'=up];    
  \node (O2) [right of = O1, label=below:{}]{}[grow'=up];          
  \node (A) [above of = O1, solid node, rectangle, minimum size=5pt, fill=white, label=right:{$(\und{x},A)$}]{}[grow'=up];    
\draw[dashed] (O) to node  {} (O1);      
\draw[dashed] (O1) to node {} (O2); 
\draw[-] (O1) to node {} (A);
\end{tikzpicture}
\end{center}
\captionof{figure}{}
\label{fig:fork_reduced}
\end{figure}

It turns out that there are five cases of acyclic matrices $\ssmat{a}{b}{c}$ which we will  need to consider:
\begin{align}
& a,b,c \geq 2,\label{case:acyclic1} \\ 
& a,b \geq 2, \text{ } c=1, \label{case:acyclic2}\\
& a \geq 2,\text{ } b=c=1, \label{case:acyclic3}\\
& a,b \geq 2,\text{ } c=0,\label{case:acyclic4}\\
& a \geq 2,\text{ } b=1,\text{ } c=0.\label{case:acyclic5}
\end{align}
Any other cases correspond to finite type cluster algebras, matrices of the form $\left(\begin{smallmatrix}0&a&0\\-a& 0&0\\0&0&0\end{smallmatrix} \right)$, or the mixed case $\ssmat{1}{1}{1}$ which was considered in Section \ref{sec:fin_deg_mixed}.

As mentioned above, these acyclic matrices have a special mutation path which produces equal degree seeds.

\begin{lem} \label{lem:acyclic_special_mutation_path}
Let $A$ be one of the matrices in the above list. Then
\[\deg.\sd[(3,2,1)^{n}]\left(\und{x},A\right) = \left(-\deg(\und{x}),A\right),\]
and so
\[\deg.\sd[(3,2,1)^{2n}]\left(\und{x},A\right) = \left(\deg(\und{x}),A\right).\]
\end{lem}

\begin{proof}
Noting that (up to sign and scaling) the grading vector is $(b,-c,a)$ and using Lemma \ref{lem:degree_is_in_matrix}, this follows immediately from direct computation of the case for $n=1$ (which we have already done by Equations \ref{eqn:mut_special_acyclic_path_along_1}--\ref{eqn:mut_special_acyclic_path_along_321}).
\end{proof}

\begin{rmk}
We could alternatively have used the path $[1,2,3]$ in Lemma \ref{lem:acyclic_special_mutation_path} above, which does the same thing in the ``opposite" direction. In fact, starting from the initial seed and mutating along both $[(3,2,1)^{n}]$ and $[(1,2,3)^{n}]$ for all $n$, we traverse precisely the vertices of the exchange graph that comprise the seeds of the floor containing acyclic matrices. (This idea of traversing the exchange graph will be made more precise when we discuss walks on the exchange graphs below.)
\end{rmk}

\begin{rmk}
While we do not make use of it in this section, it is not difficult to show that $\dv_{\und{x}}(\var_A[(3,2,1)_r]) < \dv_{\und{x}}(\var_A [(3,2,1)_{r+1}])$ for all $r \geq 1$ (using the partial order in Definition \ref{dfn:partial_order_on_denominator_vectors}). In other words, the variables in the clusters obtained under the path in Lemma \ref{lem:acyclic_special_mutation_path} ``grow" indefinitely and so are distinct for each $n$.
\end{rmk}

In all cases, dealing with the infinite $2$-trees beginning at a fork can be done in the same way.

\begin{prop} \label{prop:different_trees_different_regions}
Let $A=\ssmat{a}{b}{c}$ be a matrix with entries satisfying one of \pref{case:acyclic1}--\pref{case:acyclic5}. Let $\A(\und{x},A)$ be the corresponding cluster algebra and $\Gamma(\A)$ the exchange graph. Let $z$ and $z'$ be two cluster variables that appear in the seeds of the vertices $t$ and $t'$, respectively, and suppose these vertices belong to different branches in $\Gamma$. Write $z=\var_A\pth{p}(\und{x})$ and $z' =\var_A [\und{p}'](\und{x})$ for some minimal length mutation paths $\pth{p}$ and $[\und{p}']$ and suppose that $A_{\pth{p}}$ and $A_{[\und{p}']}$ are forks (in other words, that $z$ and $z'$ are obtained by mutation paths terminating in the canopy). Then $z$ and $z'$ are not the same cluster variable.
\end{prop}

\begin{proof}
Assume without loss of generality that $t$ occurs in the branch associated with the initial cluster. We just need to show $R(z) \neq R(z')$. Let $\pth{p}$ be a mutation path of minimal length such that $\var_A \pth{p} (\und{x}) = z$. Let $(\und{x}_t,B)=\sd_A \pth{p} (\und{x})$ be (an equivalence class representative of) the seed associated with $t$. Then the seed $\sd_B[p_n](\und{x}_t)$ does not contain $z$. But then since the component of $\Gamma$ containing $z$ is connected by Theorem \ref{thm:exchange_graph_connected}, it must be separate from the component containing $z'$, as (excluding repeated mutation directions) there is only one path from $t$ to the initial cluster. Therefore $z$ and $z'$ cannot be equal.
\end{proof}

Before turning our attention to the vertices in the floor of a graph, we make note of  \cite[Lemma 7.4]{War} (a restatement of \cite[Theorem 7.7]{FZCAI} for quivers) which is crucial.

\begin{thm} \label{thm:exchgraph_cycles}
Let $A$ be a skew-symmetric matrix and let $i \neq j$. Then, starting at the seed $(\und{x}, A)$, alternating mutations in directions $i$ and $j$ produces a cycle in the corresponding exchange graph if and only if either $a_{ij} =0$ or $|a_{ij}| =1$ . The cycle has length 4 in the first case and length 5 in the second.
\end{thm}
 
Our next aim is to prove the following.

\begin{prop} \label{prop:floor_variables_different_regions}
Let $A=\ssmat{a}{b}{c}$ be a matrix with entries satisfying one of \pref{case:acyclic1}--\pref{case:acyclic5}. Let $z$ be a cluster variable in $\A(\und{x},A)$ such that $z = \var_A\pth{p}(\und{x})$ for some mutation path $\pth{p}$. Suppose that the vertex $t \in \Gamma(\A)$ corresponds to the seed $\sd_A\pth{p}(\und{x})$, and that $t$ is in the floor of $\Gamma(\A)$. Then, for $n \geq 1$, $z':= \var_A [\und{p},(3,2,1)^{2n}] (\und{x})$ is not the same cluster variable as $z$.
\end{prop}

We will only need to prove the above for $n=1$; due to the manner in which we prove it for $n=1$, which will use Theorem \ref{thm:exchange_graph_connected}, the result for all higher values of $n$ will automatically be implied.
We will also make use of Theorem 10.3 of \cite{War}, which tells us what the required exchange graphs are. Note since $A_{[(3,2,1)^{2n}]}=A$, the exchange graph ``repeats" every six vertices along the mutation path whose corresponding vertices make up this part of the floor of the graph. (In fact, $A_{[(3,2,1)^{n}]}=A$, so it actually repeats every three vertices, but we think of the repetition as over six vertices as we need six mutations to obtain equal degree seeds.) Therefore we may restrict our attention to just one such segment of the graph. We will call the segment containing the initial seed the \emph{initial segment}. We will label vertices $t_i$, with $t_0$ being the vertex corresponding to the initial cluster. 

Our argument will proceed as follows. Suppose $z= \var_A\pth{p}(\und{x})$ as above. By taking the fixed representative $(\und{x},A)$ for the initial seed at $t_0$ (which we will denote $\sd(t_0)$---and similarly for seed representatives at the other vertices), the mutation path $\pth{p}$ identifies a walk and a terminal vertex in $\Gamma$. We may identify such a walk since, although the exchange graph is not labelled with mutation directions, we will still be able to deduce the mutation direction that results in a seed in the same class as that of the required vertex at every step of the walk. So we can identify $z$ in a fixed seed representative at a particular vertex $t_i$ in $\Gamma$. 
Next we identify the corresponding vertex $t_i'$ for $z'$. We may simply read this off the graph since it corresponds to the mutation path $[\und{p},(3,2,1)^{2}]$, so, as the graph repeats, it is the end point of the same walk on $\Gamma$ as for $z$ but extended by six vertices along from the initial cluster. We then aim to show that this vertex does not belong to $R(z)$. We will do so by showing that any walk in $\Gamma$ from $t_i$ to $t_i'$ must pass through a vertex whose corresponding cluster no longer contains $z$. Since $R(z)$ is connected by Theorem \ref{thm:exchange_graph_connected}, this accomplishes our aim and shows $z \neq z'$. As mentioned, this also implies that the same holds for all larger values of $n$: if we take a larger value of $n$ and consider the new corresponding variable $z'$ and vertex $t_i'$, it is clearly still the case that any walk from $t_i$ to $t_i'$ must pass through a vertex whose cluster no longer contains $z$.

Since we have a fixed seed representative for $t_i$, we can superimpose mutation paths onto $\Gamma$ in the same way as above and show that $z$ must be replaced along any walk resulting in a seed in the same equivalence class as that of $t_i'$. We will denote (part of) such a walk as a directed series of labelled arrows highlighted in blue and ending at a blue square vertex (by which we will mean that this vertex no longer belongs to $R(z)$).
For the initial vertex $t_0$, we will need to show that each of the three variables in the cluster representative must be replaced at some point in any walk to $t_0'$. In other words, that all three regions adjacent to $t_0$ are distinct from all three regions adjacent to $t_0'$. For other vertices, we will only need to do this for one entry of the corresponding cluster. For example, if a cluster representative at the vertex $t_1$ is obtained by one mutation of the initial cluster at $t_0$, then two of its entries will have already been dealt with in $t_0$, and so only the new entry needs to be considered. 

Let us now apply our method each of the cases \pref{case:acyclic1}--\pref{case:acyclic5}.

\begin{case}[\ref{case:acyclic1}]
The exhange graph is partially drawn in Figure \ref{fig:exchgraph_case1}. We will show that $t_0$ and $t_0'$ are in different regions (the proof for $t_1$ to $t_5$ is essentially the same). This is easy to see since a representative of the seed class at $t_3$ is obtainable from the initial seed, $\sd(t_0)$, by the mutation path $[3,2,1]$, and $\cl_A[3,2,1](\und{x})$ no longer contains any of the variables in $\und{x}$. Since any walk from $t_0$ to $t_0'$ passes through $t_3$, we are done. \hfill \qed

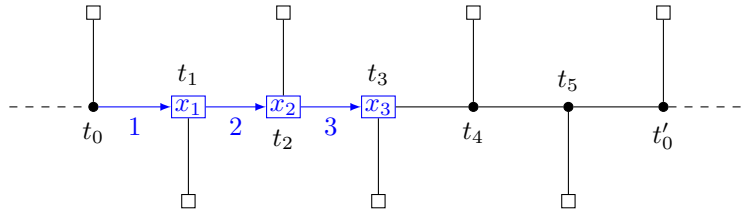
\begin{figure}[H]
\begin{center}
\begin{tikzpicture}[level distance=6mm, node distance=1.25cm, auto, font=\footnotesize] 
  \tikzstyle{level 1}=[sibling distance=0mm]
  \tikzstyle{level 2}=[sibling distance=10mm]
  \tikzstyle{level 3}=[sibling distance=8mm]
  \node (N) [left of=O]{}[grow'=up];
  \node (O) [solid node, label=below:{$t_0$}]{}[grow'=up];
  \node (Oup) [above of = O, solid node, rectangle, minimum size=5pt, fill=white, label=right:{}]{}[grow'=up];   
  \node (O1) [right of = O,solid node, rectangle, blue, fill=white, minimum size=5pt, label=above:{$t_1$}]{$x_1$}[grow'=down];
\node (O1down) [below of = O1, solid node, rectangle, minimum size=5pt, fill=white, label=right:{}]{}[grow'=up];     
  \node (O2) [right of = O1,solid node, rectangle, blue, fill=white, minimum size=5pt, label=below:{$t_2$}]{$x_2$}[grow'=up];  
 \node (O2up) [above of = O2, solid node, rectangle, minimum size=5pt, fill=white, label=right:{}]{}[grow'=up];       
  \node (O3) [right of = O2,solid node, rectangle, blue, fill=white, minimum size=5pt, label=above:{$t_3$}]{$x_3$}[grow'=down]; 
\node (O3down) [below of = O3, solid node, rectangle, minimum size=5pt, fill=white, label=right:{}]{}[grow'=up];    
  \node (O4) [right of = O3,solid node, label=below:{$t_4$}]{}[grow'=up];
\node (O4up) [above of = O4, solid node, rectangle, minimum size=5pt, fill=white, label=right:{}]{}[grow'=up];         
  \node (O5) [right of = O4,solid node, label=above:{$t_5$}]{}[grow'=down]; 
\node (O5down) [below of = O5, solid node, rectangle, minimum size=5pt, fill=white, label=right:{}]{}[grow'=up];     
  \node (O6) [right of = O5,solid node, label=below:{$t_0'$}]{}[grow'=up];                     
\node (O6up) [above of = O6, solid node, rectangle, minimum size=5pt, fill=white, label=right:{}]{}[grow'=up];   
  \node (B) [right of = O6]{}[grow'=up];

\draw[dashed] (N) to node {} (O);     
\draw[-latex, blue] (O) to node [swap] {$1$} (O1);      
\draw[-] (O) to node [swap] {} (Oup);
\draw[-latex, blue] (O1) to node [swap] {$2$} (O2); 
\draw[-] (O1) to node [swap] {} (O1down); 
\draw[-latex, blue] (O2) to node [swap] {$3$} (O3); 
\draw[-] (O2) to node [swap] {} (O2up); 
\draw[-] (O3) to node [swap] {} (O4); 
\draw[-] (O3) to node [swap] {} (O3down);
\draw[-] (O4) to node [swap] {} (O5); 
\draw[-] (O4) to node [swap] {} (O4up); 
\draw[-] (O5) to node [swap] {} (O6); 
\draw[-] (O5) to node [swap] {} (O5down);
\draw[dashed] (O6) to node [swap] {} (B);
\draw[-] (O6) to node [swap] {} (O6up);  
\end{tikzpicture}
\end{center}
\captionof{figure}{Case \ref{case:acyclic1} --- part of the exchange graph.}
\label{fig:exchgraph_case1}
\end{figure}
\end{case}

\begin{case}[\ref{case:acyclic2}]
We will use Figure \ref{fig:exchgraph_case2_t5_var[2,1]} to give a detailed explanation of how we deduce the mutation directions in the walks that are superimposed onto the exchange graph and so show that certain vertices are in separate regions. This approach will also be how we proceed in all subsequent cases. Figure \ref{fig:exchgraph_case2_t5_var[2,1]} corresponds to the variable $\var[2,1]$ in $\cl(t_5)$, the cluster class at vertex $t_5$. The other two variables in this cluster, $\var[1]$ and $x_3$, will have already been dealt with by the time we consider $\var[2,1]$.

We start the case by taking the representative $\left((x_1,x_2,x_3), A=\ssmat{a}{b}{c}\right)$ for the initial seed at vertex $t_0$. 
We need to find the vertex whose seed class contains $\sd_A[2,1](x_1,x_2,x_3)$.
First note the following: by Theorem \ref{thm:exchgraph_cycles}, it is clear that mutating the initial seed in direction 1 will produce a seed in the same class as either $\sd(t_1)$ or $\sd(t_4)$. Since the exchange graph is ``symmetric" about $t_0$, we may simply assume that direction 1 corresponds to $\sd(t_4)$ rather than $\sd(t_1)$. However, once this choice has been fixed, we will need to use it consistently throughout all diagrams in the case. To obtain a seed in the same class as $\sd(t_5)$ we must mutate in direction 2, since (again by Theorem \ref{thm:exchgraph_cycles}) direction 3 keeps us on the pentagon comprising $R(x_2)$. Therefore $t_5$ is the desired vertex. We have now obtained $\var[2,1]$ as the second entry of our cluster, and our seed is in the same class as $\sd(t_5)$. We mark the mutation walk leading up to the variable we are considering in red on the diagram, and we will consider diagrams in an order such that we will have already dealt with the other two entries of our cluster. 

Now we wish to prove that the vertex $t_5'$ is not part of the region $R(\var[2,1])$ by finding ``variable breaks" in the exchange graph, which show that $t_5'$ is not in the connected subgraph containing $\var[2,1]$. Such breaks will occur the next time we mutate in direction 2 along any walk, as this is when the second entry of our cluster representative is replaced. (We do not need to worry about the variable appearing as a different entry before this point, as this would violate that the variables in a cluster form a transcendence base.) 

Clearly any walk from $t_5$ to $t_5'$ must pass through $t_1'$, which we can show is not in $R(\var[2,1])$. To get from $t_5$ to $t_6$ we next mutate in direction 3. 
To see that it is direction 3 that takes us to $t_6$ rather than direction 1, recall that $A_{[1,2,3]}=A$. So the exchange graph must have the same form about the vertex we end up at after $[3,2,1]$ as it does about the initial vertex $t_0$---this is not the case for $t_9$, so direction 3 must indeed take us to $t_6$. (Alternatively, we may argue by contradiction that if this were not the case then the walk corresponding to $[3,2,1,3,2,1]$ would in fact take us to a seed which is a fork, but the matrix $A_{[3,2,1,3,2,1]}=A$ is not a fork.) Next, to get to $t_7$, we mutate in direction 1. Finally we must mutate in direction 2 to arrive at $t_1'$. In the process, we replace the second entry of our cluster representative, $\var[2,1]$. Now that we know $\cl(t_1')$ does not contain $\var[2,1]$ we mark it with a blue square, and we are done.

We check all the other variables in the diagrams below. The regions corresponding to the initial variables are indicated as a visual aid, although we bear in mind that regions are abstract sets rather than regions of the plane in the geometric sense. In some of the figures, we simultaneously superimpose two mutation paths onto the exchange graph. 
Note that a representative of $\cl(t_4)$ is given by $\left(\var[1,3,1,3] = x_3 ,x_2,\var[3,1,3]\right)$, and these three variables are accounted for in Figures \ref{fig:exchgraph_case2_t0_var1}--\ref{fig:exchgraph_case2_t3_var[3,1,3]}. (This is suggested by the diagram, as we can see that the regions adjacent to $t_4$ are $R(x_2)$, $R(x_3)$ and $R(\var[3,1,3])$.) Alternatively, we could have taken another representative of $\cl(t_4)$: $\left(\var[1], x_2, x_ 3\right)$. But our previous choice makes it more obvious that $\var[1]$ (which is $\var[3,1,3]$) has been accounted for, given that we have already considered all variables along the path $[3,1,3]$. In a similar way, the variables of $\cl(t_9)$ are already accounted for elsewhere. Therefore we do not to concern ourselves with the variables in $t_4$ and $t_9$.

\begin{figure}[H]
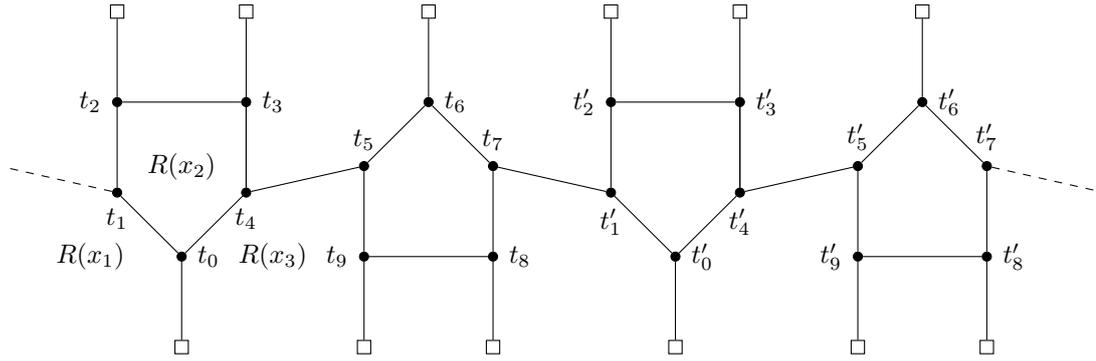

\begin{center}

\end{center}
\captionof{figure}{Case \ref{case:acyclic2} --- two segments of the exchange graph for this matrix.}
\label{fig:exchgraph_case2}
\end{figure}

\begin{figure}[H]
\begin{center}
%
\end{center}
\captionof{figure}{Case \ref{case:acyclic2} --- variable $x_1$ of the initial cluster $\cl(t_0)$.}
\label{fig:exchgraph_case2_t0_var1}
\end{figure}

\begin{figure}[H]
\begin{center}
%
\end{center}
\captionof{figure}{Case \ref{case:acyclic2} --- variable $x_2$ of the initial cluster $\cl(t_0)$.}
\label{fig:exchgraph_case2_t0_var2}
\end{figure}

\begin{figure}[H]
\begin{center}
%
\end{center}
\captionof{figure}{Case \ref{case:acyclic2}. Variable $x_3$ of the initial cluster $\cl(t_0)$.}
\label{fig:exchgraph_case2_t0_var3}
\end{figure}

\begin{figure}[H]
\begin{center}
%
\end{center}
\captionof{figure}{Case \ref{case:acyclic2} --- variable $\var[3]$ of $\cl(t_1)$.}
\label{fig:exchgraph_case2_t1_var[3]}
\end{figure}

\begin{figure}[H]
\begin{center}
%
\end{center}
\captionof{figure}{Case \ref{case:acyclic2} --- variable $\var[1,3]$ of $\cl(t_2)$.}
\label{fig:exchgraph_case2_t2_var[1,3]}
\end{figure}

\begin{figure}[H]
\begin{center}
%
\end{center}
\captionof{figure}{Case \ref{case:acyclic2} --- variable $\var[3,1,3]$ of $\cl(t_3)$.}
\label{fig:exchgraph_case2_t3_var[3,1,3]}
\end{figure}

\begin{figure}[H]
\begin{center}
%
\end{center}
\captionof{figure}{Case \ref{case:acyclic2} --- variable $\var[2,1]$ of $\cl(t_5)$.}
\label{fig:exchgraph_case2_t5_var[2,1]}
\end{figure}

\begin{figure}[H]
\begin{center}
%
\end{center}
\captionof{figure}{Case \ref{case:acyclic2} --- variable $\var[3,2,1]$ of $\cl(t_6)$.}
\label{fig:exchgraph_case2_t6_var[3,2,1]}
\end{figure}

\begin{figure}[H]
\begin{center}
%
\end{center}
\captionof{figure}{Case \ref{case:acyclic2} --- variable $\var[1,3,2,1]$ of $\cl(t_7)$.}
\label{fig:exchgraph_case2_t7_var[1,3,2,1]}
\end{figure}

\begin{figure}[H]
\begin{center}
%
\end{center}
\captionof{figure}{Case \ref{case:acyclic2} --- variable $\var[3,1,3,2,1]$ of $\cl(t_8)$.}
\label{fig:exchgraph_case2_t9_var[3,1,3,2,1]}
\end{figure}

In each case, we see that the variable in question at vertex $t_i$ is in a separate region to those of the variables at $t_i'$.

\hfill \qed
\end{case}

\begin{case}[\ref{case:acyclic3}]
We deal with all variables in the initial cluster at the same time in Figure \ref{fig:exchgraph_case3_t0_var123}. We then consider $\var[1]$ of $\cl(t_4)$, $\var[3]$ of $\cl(t_1)$, $\var[1,3]$ of $\cl(t_2)$ and $\var[3,1]$ of $\cl(t_3)$. The other cases are similar and are not shown.
\begin{figure}[H]
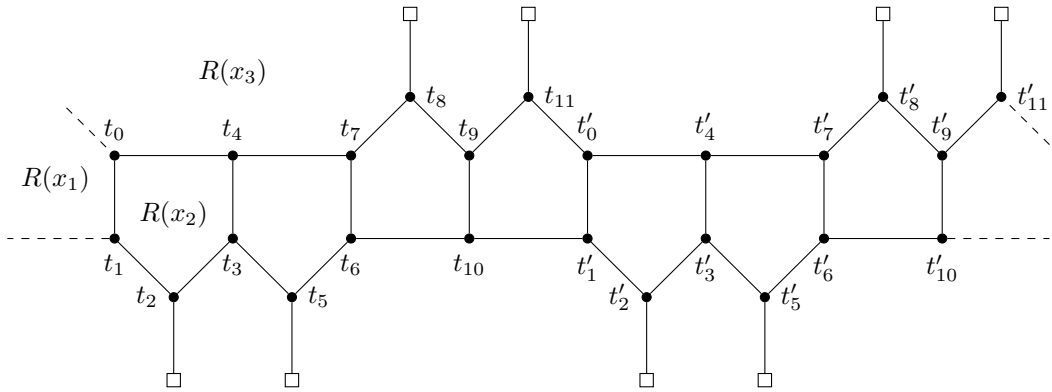

\begin{center}

\end{center}
\captionof{figure}{Case \ref{case:acyclic3} --- two segments of the exchange graph.}
\label{fig:exchgraph_case3}
\end{figure}

\begin{figure}[H]
\begin{center}
%
\end{center}
\captionof{figure}{Case \ref{case:acyclic3} --- variables $x_1$, $x_2$ and $x_3$ of the initial cluster $\cl(t_0)$.}
\label{fig:exchgraph_case3_t0_var123}
\end{figure}

\begin{figure}[H]
\begin{center}
%
\end{center}
\captionof{figure}{Case \ref{case:acyclic3} --- variable $\var[1]$ of $\cl(t_4)$.}
\label{fig:exchgraph_case3_t4_var[1]}
\end{figure}

\begin{figure}[H]
\begin{center}
%
\end{center}
\captionof{figure}{Case \ref{case:acyclic3} --- variable $\var[3]$ of $\cl(t_1)$.}
\label{fig:exchgraph_case3_t1_var[3]}
\end{figure}

\begin{figure}[H]
\begin{center}
%
\end{center}
\captionof{figure}{Case \ref{case:acyclic3} --- variable $\var[1,3]$ of $\cl(t_2)$.}
\label{fig:exchgraph_case3_t2_var[1,3]}
\end{figure}

\begin{figure}[H]
\begin{center}
%
\end{center}
\captionof{figure}{Case \ref{case:acyclic3} --- variable $\var[3,1]$ of $\cl(t_3)$.}
\label{fig:exchgraph_case3_t3_var[3,1]}
\end{figure}

\hfill \qed
\end{case}

\begin{case}[\ref{case:acyclic4}]
In this case we show the relevant diagrams for the initial variables as well as $\var[1]$ of $\cl(t_3)$, $\var[2,1]$ of $\cl(t_4)$,  $\var[3]$ of $\cl(t_1)$ and $\var[1,3]$ of $\cl(t_2)$. Other cases are similar.
\begin{figure}[H]
\begin{center}
%
\end{center}	
\captionof{figure}{Case \ref{case:acyclic4} --- two segments of the exchange graph.}
\label{fig:exchgraph_case4}
\end{figure}

\begin{figure}[H]
\begin{center}
%
\end{center}	
\captionof{figure}{Case \ref{case:acyclic4} --- variables $x_1$, $x_2$ and $x_3$ of the initial cluster $\cl(t_0)$.}
\label{fig:exchgraph_case4_t0_var123}
\end{figure}

\begin{figure}[H]
\begin{center}
%
\end{center}	
\captionof{figure}{Case \ref{case:acyclic4} --- variable $\var[1]$ of $\cl(t_3)$.}
\label{fig:exchgraph_case4_t3_var[1]}
\end{figure}

\begin{figure}[H]
\begin{center}
%
\end{center}	
\captionof{figure}{Case \ref{case:acyclic4} --- variable $\var[2,1]$ of $\cl(t_4)$.}
\label{fig:exchgraph_case4_t4_var[2,1]}
\end{figure}

\begin{figure}[H]
\begin{center}
%
\end{center}	
\captionof{figure}{Case \ref{case:acyclic4} --- variable $\var[1,3]$ of $\cl(t_2)$.}
\label{fig:exchgraph_case4_t2_var[1,3]}
\end{figure}

\hfill \qed
\end{case}

\begin{case}[\ref{case:acyclic5}]
This case has a slightly different property to the others: the vertices $t_4$ and $t_4'$ share a common region, namely, $R(\var [2,1,3])$. (Similarly, so do $t_9$ and $t_9'$.) This means we can not prove that the cluster variables of $t_4$ are distinct from those of $t_4'$. This does not cause any problems, though; while it is not the case that the variable is replaced as we move one segment to the right, this does occur as we move two segments to the right. Showing this is not difficult (see Figure \ref{fig:exchgraph_case5_t1_var[2,1,3]}) and is enough to prove that we have infinitely many variables in all degrees. Indeed, if we are given a variable in the canopy, this shows it is distinct from all the variables in the corresponding cluster two segments to the right (which, as usual, is the same degree cluster). As before, this behaviour propagates inductively, so that the corresponding cluster in every second segment must contain a variable unique to just one pair of adjacent segments, but of the same degree as the others.

We show the diagrams for the initial variables as well as for $\var [3]$ of $\cl(t_1)$, $\var [1,3]$ of $\cl(t_3)$, $\var [2,1,3]$ of $\cl(t_4)$ and $\var [3,2,1,3]$ of $\cl(t_5)$. 

\begin{figure}[H]
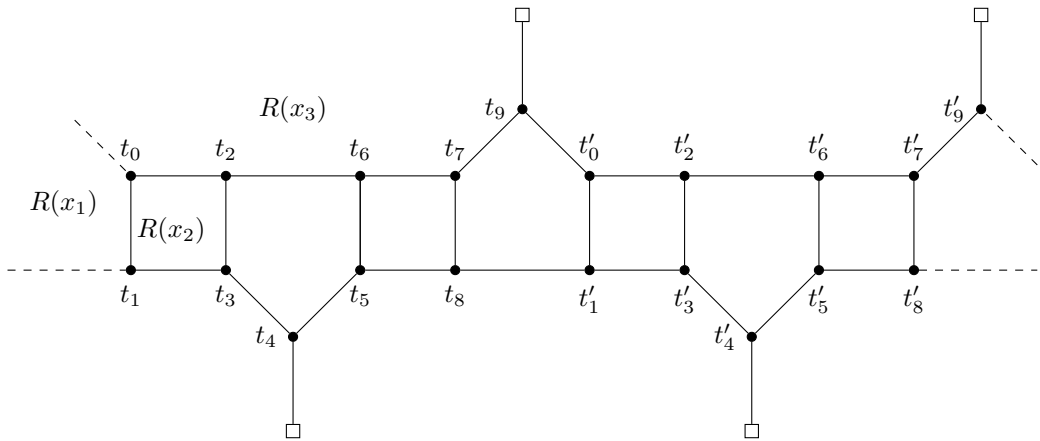

\begin{center}

\end{center}	
\captionof{figure}{Case \ref{case:acyclic5} --- two segments of the exchange graph.}
\label{fig:exchgraph_case5}
\end{figure}

\begin{figure}[H]
\begin{center}
%
\end{center}	
\captionof{figure}{Case \ref{case:acyclic5} --- variables $x_1$, $x_2$ and $x_3$ of the initial cluster $\cl(t_0)$.}
\label{fig:exchgraph_case5_t0_var123}
\end{figure}

\begin{figure}[H]
\begin{center}
%
\end{center}	
\captionof{figure}{Case \ref{case:acyclic5} --- variable $\var[3]$ of $\cl(t_1)$.}
\label{fig:exchgraph_case5_t1_var[3]}
\end{figure}

\begin{figure}[H]
\begin{center}
%
\end{center}	
\captionof{figure}{Case \ref{case:acyclic5} ---variable $\var[1,3]$ of $\cl(t_3)$.}
\label{fig:exchgraph_case5_t1_var[1,3]}
\end{figure}

\begin{figure}[H]
\begin{center}
%
\end{center}	
\captionof{figure}{Case \ref{case:acyclic5} --- variable $\var[2,1,3]$ of $\cl(t_4)$.}
\label{fig:exchgraph_case5_t1_var[2,1,3]}
\end{figure}

\begin{figure}[H]
\begin{center}
%
\end{center}	
\captionof{figure}{Case \ref{case:acyclic5} ---variable $\var[3,2,1,3]$ of $\cl(t_5)$.}
\label{fig:exchgraph_case5_t1_var[3,2,1,3]}
\end{figure}
\hfill \qed
\end{case}

We have now shown what we want for all cases, which establishes Proposition \ref{prop:floor_variables_different_regions}. Putting together Proposition \ref{prop:different_trees_different_regions} and Proposition \ref{prop:floor_variables_different_regions}, we have the following.

\begin{thm} \label{thm:acyclic_infinitely_many_in_each}
Let $A$ be a $3 \times 3$ mutation-acyclic matrix. Then the graded cluster algebra $\A \big( \und{x}, A, \und{g}\big)$ has infinitely many variables in each occurring degree. The grading is also balanced.
\end{thm}

\begin{proof}
Given the above work, we now need only justify that the grading is balanced. Of course, this follows from Lemma \ref{lem:acyclic_implies_balanced} since $A$ is acyclic, but here we are to able see this directly. Let $z$ be a cluster variable in $\A$. Then $z =\var_A [\und{p},(3,2,1)^n]$ for some $n$ and some minimal length mutation path $\pth{p}$ (strictly, we may need to replace the path $[(3,2,1)^n]$ by $[(1,2,3)^n]$). But then $z':=\var_A [\und{p},(3,2,1)^{n+1}]$ satisfies $\deg_{\und{g}}(z')=-\deg_{\und{g}}(z)$. This gives a bijection between the variables of degree $d$ and degree $-d$ for any occurring degree $d$.
\end{proof}

\section{The singular cyclic case} \label{sec:singular_cyclic}

The case of $A=\left(\begin{smallmatrix}0&a&-2\\-a& 0&a\\2&-a&0\end{smallmatrix} \right)$, $a \geq 3$, is special. Mutation in direction $1$ or $3$ gives $-A$, while mutation in direction $2$ gives a fork. Thus $A$ is mutation-cyclic, but unlike the other mutation-cyclic cases we can mutate along a repetition-free path without obtaining a strictly increasing sequence of degrees. Thus the arising graded cluster algebra is not forced to have finitely many variables in each degree (though, like the other mutation-cyclic cases, all degrees will be positive). Indeed, as in the acyclic case, we can obtain a floor of equal degree seeds with different underlying clusters, and so we might expect the graded cluster algebra to behave as one generated by an acyclic matrix. However, as Theorem \ref{thm:exchange_graph_connected} is only conjectured for mutation-cyclic matrices, this can not be proved in the same way as above. Let us now consider what can be said about this case.

For the following lemma, recall the relevant notation regarding denominator vectors from Definition \ref{dfn:denominator_vector}.

\begin{lem}
\label{lem:increasing_floor}
We have $\dv_{\und{x}}(\var_A[(3,1)_k]) < \dv_{\und{x}}(\var_A [(3,1)_{k+1}])$ for all $k \geq 1$. 
\end{lem}

\begin{proof}
We will show
\[
\dcl[(3,1)_k]=
\begin{cases}
\big( (k-1,0,k-2), (0,-1,0), (k,0,k-1) \big) & \text{ if $k$ is even},\\
\big((k,0,k-1),(0,-1,0), (k-1,0,k-2) \big) & \text{ if $k$ is odd}.
\end{cases}
 \] 
 As a base case note that $\dcl[(3,1)_1]=\dcl[1]= \big((1,0,0),(0,-1,0),(0,0,-1) \big)$. Assume the result is true for $k$ and assume without loss of generality that $k$ is even. Then the current degree cluster is $\big( (k-1,0,k-2), (0,-1,0), (k,0,k-1) \big)$ and the current matrix is $A$. We then have
\begin{align*}
\dv [(3,1)_{k+1}]& =\dv[1,(3,1)_k] \\
 &= - ( k-1,0,k-2 ) + \max \big(a(0,-1,0), 2(k,0,k-1)) \big)\\
 &= - ( k-1,0,k-2 ) + (2k,0,2k-2) \\
 &= (k+1,0,k). 
\end{align*}
So the new denominator cluster is $\big((k+1,0,k),(0,-1,0), (k,0,k-1) \big) $ as required.
\end{proof}

Thus, since $A_{[(3,1)_k]}=(-1)^k A$, the path $[(3,1)^n]$ is the special mutation path for $A$ analogous to the one in Lemma \ref{lem:acyclic_special_mutation_path} for an acyclic matrix.

\begin{cor}
In the graded cluster algebra $\A \left( \und{x}, \left(\begin{smallmatrix}0&a&-2\\-a& 0&a\\2&-a&0\end{smallmatrix} \right), (a,2,a) \right)$, there are infinitely many degrees that each contain infinitely many variables.
\end{cor}

\begin{proof}
Let $\pth{p}$ be a mutation path in $\A$ whose corresponding terminal vertex is at least three mutations above the floor of the exchange graph (i.e.~three mutations above a vertex in the equivalence class of $\sd [(3,1)_r] $.) We must have that $\cl_A\pth{p}$ and $\cl_A[\und{p},3,1]$ are distinct. If not, then $\cl_A\pth{p} = \cl_A[\und{p},3,1]$ so that, applying the reverse path $[\und{p}^{-1}]$ to both clusters, $\cl_A [\und{p}^{-1}, \und{p}] = \cl_A[\und{p}^{-1},\und{p},3,1]$. This implies  $\cl_A[3,1]= (x_1,x_2,x_3)$, a contradiction by Lemma \ref{lem:increasing_floor}.
Now by prepending the path $[3,1]$ to $\pth{p}$ infinitely many times, we obtain infinitely many distinct clusters whose corresponding degree clusters are all equal. By the pigeonhole principle, this produces infinitely many variables in at least one of the degrees in $\deg.\cl \pth{p}$. Now say $\pth{p}$ has length $n$ and let $q_1,q_2$ and $q_3$ be mutation directions such that $q_1 \neq p_n$ and $\{q_1,q_2,q_3\}=\{1,2,3\}$. (Applying the mutation path $[q_3,q_2,q_1]$ to the current cluster will then replace each each cluster variable.) Repeating the above process with the path $[q_3,q_2,q_1,\und{p}]$ then yields infinitely many variables of a different degree to the previous one, since degrees strictly increase as we increase the mutation distance from the floor of the exchange graph. We may continue extending the mutation path in a similar way an indefinite number of times and then repeat the same prepending process each time. This produces infinitely many degrees each containing infinitely many variables.
\end{proof}

\begin{rmk}
Given that we may no longer use Theorem \ref{thm:exchange_graph_connected}, the most obvious approach when attempting to prove Conjecture \ref{cnj:singular_cyclic_infinitely_many_in_each} might be that of an inductive argument. We have shown in Lemma \ref{lem:increasing_floor} that the denominator clusters along the floor of the exchange graph grow as we move along the special mutation path, and this would appear to be a suitable base case for the following argument. Suppose $\pth{p}$ is a mutation path resulting in degree $d$, starting from the initial cluster, and that the denominator vector obtained after $\pth{p}$ is smaller than that obtained after $[\und{p},3,1]$ (i.e.~after starting with a larger denominator cluster). We may then expect to be able to prove that if $[\und{p}']$ is a path obtained by extending $\pth{p}$ by one mutation, then $[\und{p}',3,1]$ also gives a larger denominator vector than $[\und{p}]$ does. In practice, showing this is true seems more difficult than might initially be hoped, the difficulty being that when we mutate a larger denominator cluster, we must also subtract a larger denominator vector when computing the new variable.

An alternative approach may involve trying to find an invariant of a cluster variable that is more appropriate to our task. This invariant would need to capture some notion of the size of a cluster variable, while having a mutation formula that circumvents the above problem. Yet another approach may involve considering $\pth{p}$ as a rational function of a denominator cluster (or possibly of some other object capturing relevant information about a cluster) and attempting to show that, given two such objects as inputs, the larger of the two must give a larger output.
\end{rmk}

\biblio

\chapter{Degree subquivers with growth} \label{chp:growing_subquivers}

It is useful to make note of certain small quivers and degree quivers which, if embedded as a subquiver in a larger quiver, will provide a way to find increasing arrows and degrees. This will reduce the problem of finding infinitely many degrees to one of finding such a subquiver in the mutation class of our initial exchange quiver. This reduction will, in particular, be useful in the Chapter \ref{chp:matrix_grassmannian}. Note that all results concerning quivers below are automatically also true of their respective opposite quivers.

\begin{ntn}
Given a quiver $Q$ with vertices $v_1, \dots, v_n$, we will use $w_Q (v_i, v_j)$, to mean the $(i,j)$ entry of the corresponding matrix. That is, the number (or weight) of arrows from vertex $v_i$ to vertex $v_j$, with an arrow in the opposite direction counting as negative.
\end{ntn}

In what follows we will refer to \emph{degree subquivers}. This means a subquiver (in the sense of Definition \ref{dfn:quiver}) of a degree quiver, where we require that the degrees of the vertices of the subquiver must match those of the corresponding vertices of the quiver. It is important to note that a degree subquiver need not be a degree quiver itself.

\section{Subquivers with growing arrows}
\begin{prop}
\label{prop:growing_subquiver}
Consider the quiver
\begin{equation}
Q=
\begin{tikzpicture}[baseline=(current  bounding  box.center),node distance=1.6cm, auto, font=\footnotesize]                                          
  \node (A7) []{$v_2$};                                                                                                                     
  \node (A1) [below left of = A7]{$v_3$};                                                                               
  \node (A14) [below right of = A7]{$v_1$};   
  \draw[-latex] (A7) to node {} (A1);      
  \draw[-latex] (A14) to node [swap] {$2$} (A7);    
\end{tikzpicture}.
\end{equation}
(Recall that no label on an arrow means it has weight $1$.) Let $\pth{p}$ be a repetition-free mutation path such that $p_1 \neq v_1$. Let  $R=Q_{[\und{p},v_1,v_2]}$. Then $|w_R (v_i, v_j)| \geq |w_Q (v_i, v_j)|$ for each $i \neq j$, and this inequality is strict for at least one choice of $i$ and $j$.
Moreover, if $[\und{p}']$ is a proper rooted subpath of $\pth{p}$ and $R':=Q_{[\und{p'},v_1,v_2]}$, then $|w_R (v_i, v_j)| > |w_{R'} (v_i, v_j)|$ for some pair $i \neq j$.
\end{prop}

\begin{proof}
For notational simplicity we will use $i$ and $v_i$ interchangeably when working with mutation directions. Mutation at vertex $v_2$ turns $Q$ into the quiver
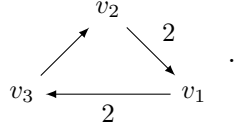
\[ \Image{
\begin{tikzpicture}[node distance=1.6cm, auto, font=\footnotesize]                                          
  \node (A7) []{$v_2$};                                                                                                                     
  \node (A1) [below left of = A7]{$v_3$};                                                                               
  \node (A14) [below right of = A7]{$v_1$};   
  \draw[latex-] (A7) to node {} (A1);      
  \draw[latex-] (A14) to node [swap] {$2$} (A7);    
  \draw[-latex] (A14) to node [] {$2$} (A1);      
\end{tikzpicture}
}.
\]
 Mutation at $v_1$ then turns this into 
\[ \Image{
\begin{tikzpicture}[node distance=1.6cm, auto, font=\footnotesize]                                          
  \node (A7) []{$v_2$};                                                                                                                     
  \node (A1) [below left of = A7]{$v_3$};                                                                               
  \node (A14) [below right of = A7]{$v_1$};   
  \draw[-latex] (A7) to node [swap] {$3$} (A1);      
  \draw[-latex] (A14) to node [swap] {$2$} (A7);    
  \draw[latex-] (A14) to node [] {$2$} (A1);      
\end{tikzpicture}
},
\]
which corresponds to the matrix  $\ssmat{-2}{-3}{2} $. Note that $\ssmat{-2}{-3}{2}  \esseq \sigma \ssmat{-3}{-2}{2} $, where $\sigma$ is such that $\sigma(1)=3$ (i.e. such that if $k \in \{2,3\}$ then $\sigma(k) \in \{1,2\}$). The matrix $\ssmat{-3}{-2}{2}$ is $2$-cyclic, and since $\pth{p}$ starts with mutation of vertex $v_2$ or $v_3$, we may apply Corollary \ref{cor:c_cyclic_preserved_in_paths_starting_1_2} to $ \ssmat{-3}{-2}{2} $ with the path $[\sigma(\und{p})]$. This gives the desired result for an essentially equivalent quiver, and hence for our original quiver.
\end{proof} 

Since we applied Corollary \ref{cor:c_cyclic_preserved_in_paths_starting_1_2} in the proof of Proposition \ref{prop:growing_subquiver}, we immediately obtain a further corollary: 

\begin{cor} \label{cor:quiver_remains_cyclic}
Let $Q$ and $\pth{p}$ be as in Proposition \ref{prop:growing_subquiver}. Then $Q_{[\und{p},v_1,v_2]}$ is cyclic.
\end{cor}

\section{Growing degrees}

Next we give some conditions in which the same quiver allows us to obtain an increasing sequence of degrees.

\begin{rmk}
Let $Q$ be a positive degree quiver and $R$ a degree subquiver. Notice that when carrying out degree seed mutation at a particular vertex $v$ in $R$, we subtract the degree of $v$ and then only add multiples of degrees of adjacent vertices. Since we have assumed $Q$ is positive, the degrees of these adjacent vertices are positive. Any arrows external to $R$ can therefore only add to the degrees that would be obtained under mutation if $R$ was the entire quiver, and if we only mutate at vertices inside $R$, these are the only vertices whose degrees may change. Thus, if we are interested in showing that mutation paths give increasing degrees, we may safely ignore any arrows that are adjacent to any vertex outside $R$.
\end{rmk}

\begin{prop}
\label{prop:growing_degree_subquiver}
Let Q be a positive degree quiver which contains the degree subquiver
\begin{equation}
R=
\begin{tikzpicture}[baseline=(current  bounding  box.center),node distance=1.6cm, auto, font=\footnotesize]                                          
  \node (A7) []{$(d_2)_{v_2}$};                                                                                                                     
  \node (A1) [below left of = A7]{$ (d_3)_{v_3}$};                                                                               
  \node (A14) [below right of = A7]{$(d_1)_{v_1}$};   
  \draw[-latex] (A7) to node {} (A1);      
  \draw[-latex] (A14) to node [swap] {$2$} (A7);    
\end{tikzpicture}.
\end{equation}
Let $\pth{p}$ be a repetition free mutation path such that $p_1 \neq v_1$. Suppose
\begin{align}
 (d_i) \geq 1  \text{ for } i = 1,2,3, \\
  (d_3)> (d_1) \geq (d_2).
\end{align}
Then $\deg [\und{p},v_1,v_2] > \max(d_3,d_2)$, and if $[\und{p}']$ is a proper rooted subpath of $\pth{p}$, then $\deg [\und{p},v_1,v_2]> \deg [\und{p}',v_1,v_2]$.
\end{prop}

\begin{proof}
We will prove the result by induction. First, mutating at vertex $v_2$, we obtain the degree subquiver
\[ \begin{tikzpicture}[baseline=(current  bounding  box.center),node distance=1.6cm, auto, font=\footnotesize]                                          
  \node (A7) []{$(d_2')$};                                                                                                                     
  \node (A1) [below left of = A7]{$(d_3)$};                                                                               
  \node (A14) [below right of = A7]{$ (d_1)$};   
  \draw[latex-] (A7) to node {} (A1);      
  \draw[latex-] (A14) to node [swap] {$2$} (A7);    
  \draw[-latex] (A14) to node [] {$2$} (A1);      
\end{tikzpicture},
\]
where $(d_2') \geq 2 (d_1) - (d_2)>(d_1)$.
Mutating next at $v_1$ we obtain
\[\begin{tikzpicture}[baseline=(current  bounding  box.center),node distance=1.6cm, auto, font=\footnotesize]                                          
  \node (A7) []{$(d_2')$};                                                                                                                     
  \node (A1) [below left of = A7]{$(d_3)$};                                                                               
  \node (A14) [below right of = A7]{$ (d_1')$};   
  \draw[-latex] (A7) to node [swap] {$3$} (A1);      
  \draw[-latex] (A14) to node [swap] {$2$} (A7);    
  \draw[latex-] (A14) to node [] {$2$} (A1);      
\end{tikzpicture}, 
\]
where $(d_1') \geq 2 (d_2') - (d_1)>(d_2') +(d_1)- (d_1)=(d_2')$
and $(d_1') \geq 2(d_3) - (d_1) > (d_3)$.

For the induction step, let $\pth{q}$ be a proper rooted subpath of $\pth{p}$ and suppose the result is true for $\pth{q}$. Up to isomorphism, our degree subquiver is then
\begin{equation}
\label{fig:growing_degree_quiver}
\begin{tikzpicture}[baseline=(current  bounding  box.center),node distance=1.6cm, auto, font=\footnotesize]                                          
  \node (A7) []{$(e_2)$};                                                                                                                     
  \node (A1) [below left of = A7]{$(e_3)$};                                                                               
  \node (A14) [below right of = A7]{$ (e_1)$};   
  \draw[-latex] (A7) to node [swap] {$c$} (A1);      
  \draw[-latex] (A14) to node [swap] {$a$} (A7);    
  \draw[latex-] (A14) to node [] {$b$} (A1);      
\end{tikzpicture},
\end{equation}
where, depending on whether the last mutation of $\pth{q}$ is $v_1$, $v_2$ or $v_3$, we have either
\begin{align}
& (e_1) > (e_2) , (e_3), \label{case:growing_degree_quiver_1}\\ 
& (e_2) > (e_1), (e_3) \text{, or}  \label{case:growing_degree_quiver_2}\\
& (e_3) > (e_1), (e_2),  \label{case:growing_degree_quiver_3}
\end{align}
respectively, and where $a, b, c \geq 2$ by Proposition \ref{prop:growing_subquiver}. Let $\pth{\hat{q}}$ be the rooted subpath of $\pth{p}$ of length $m:=l(\pth{q})+1$. We wish to show that mutation in the direction of $\hat{q}_m$ results in a degree larger than each $(e_i)$, and another degree quiver of the form \pref{fig:growing_degree_quiver} satisfying one of \pref{case:growing_degree_quiver_1} -- \pref{case:growing_degree_quiver_3}.
First we will write down what is obtained upon mutation of \pref{fig:growing_degree_quiver} in each direction.

\begin{align}
\text{Direction $v_1$: } & 
\begin{tikzpicture}[baseline=(current  bounding  box.center),node distance=1.6cm, auto, font=\footnotesize]                                          
  \node (A7) []{$(e_2)$};                                                                                                                     
  \node (A1) [below left of = A7]{$(e_3)$};                                                                               
  \node (A14) [below right of = A7]{$ (e_1')$};   
  \draw[latex-] (A7) to node [swap] {$ab-c$} (A1);      
  \draw[latex-] (A14) to node [swap] {$a$} (A7);    
  \draw[-latex] (A14) to node [] {$b$} (A1);      
\end{tikzpicture}, 
 \text{ where } (e_1') \geq \max(a(e_2),b(e_3))-(e_1). \\
\text{Direction $v_2$: } & 
\begin{tikzpicture}[baseline=(current  bounding  box.center),node distance=1.6cm, auto, font=\footnotesize]                                          
  \node (A7) []{$(e_2')$};                                                                                                                     
  \node (A1) [below left of = A7]{$(e_3)$};                                                                               
  \node (A14) [below right of = A7]{$ (e_1)$};   
  \draw[latex-] (A7) to node [swap] {$c$} (A1);      
  \draw[latex-] (A14) to node [swap] {$a$} (A7);    
  \draw[-latex] (A14) to node [] {$ac-b$} (A1);      
\end{tikzpicture}, 
 \text{ where } (e_2') \geq \max(a(e_1),c(e_3))-(e_2). \\
 \text{Direction $v_3$: } & 
\begin{tikzpicture}[baseline=(current  bounding  box.center),node distance=1.6cm, auto, font=\footnotesize]                                          
  \node (A7) []{$(e_2)$};                                                                                                                     
  \node (A1) [below left of = A7]{$(e_3')$};                                                                               
  \node (A14) [below right of = A7]{$ (e_1)$};   
  \draw[latex-] (A7) to node [swap] {$c$} (A1);      
  \draw[latex-] (A14) to node [swap] {$bc-a$} (A7);    
  \draw[-latex] (A14) to node [] {$b$} (A1);      
\end{tikzpicture}, 
 \text{ where } (e_3') \geq \max(c(e_2),b(e_1))-(e_3).
\end{align}

In the degree subquivers above, all edge weights are again at least two, by Proposition \ref{prop:growing_subquiver}.
Now if $\hat{q}_m$ is direction $v_1$, then the previous mutation direction was $v_2$ or $v_3$, so that, by the induction hypothesis, either $(e_2)$ or $(e_3)$ is strictly greater than the other two degrees. Then $\max(a(e_2),b(e_3))-(e_1) > \max((e_2), (e_3))>(e_1)$, since $a,c \geq 2$, so $(e_1')>\max((e_2),(e_3)) > (e_1)$. Similarly, if  $\hat{q}_m$ is direction $v_2$, we find $(e_2') \geq \max(a(e_1),c(e_3))-(e_2) > \max((e_1),(e_3))>(e_2)$ and if  $\hat{q}_m$ is direction $v_3$ then $(e_3') \geq \max(c(e_2),b(e_1))-(e_3) > \max((e_1),(e_2))>(e_3)$.
\end{proof}

We can relax our assumptions and still obtain infinitely many degrees, though possibly without a strictly increasing sequence, and requiring a specific mutation path.

\begin{prop}
\label{prop:degree_subquiver_tends_infty_relaxed}
Let $Q$ be a positive
degree quiver which contains the degree subquiver
\begin{equation}
R=
\begin{tikzpicture}[baseline=(current  bounding  box.center),node distance=1.6cm, auto, font=\footnotesize]                                          
  \node (A7) []{$(d_2)_{v_2}$};                                                                                                                     
  \node (A1) [below left of = A7]{$ (d_3)_{v_3}$};                                                                               
  \node (A14) [below right of = A7]{$( d_1)_{v_1}$};   
  \draw[-latex] (A7) to node {} (A1);      
  \draw[-latex] (A14) to node [swap] {$2$} (A7);    
\end{tikzpicture},
\end{equation}
where $(d_i) \geq 1  \text{ for } i = 1,2,3$.
Then
\begin{align*}
&\deg [(v_1,v_2)_n] \rightarrow \infty \text{ as } n \rightarrow \infty & \text{ if } (d_1) \geq (d_2),\\
&\deg [(v_1,v_2)_n,v_3,v_1] \rightarrow \infty \text{ as } n \rightarrow \infty & \text{ if } (d_2) > (d_1).
\end{align*}
\end{prop}

\begin{proof}
We will split our proof into three cases:
\begin{enumerate}[(i)]
\item $(d_1)>(d_2)$  \label{item:1}
\item $(d_2)>(d_1) $ \label{item:2}
\item $(d_1)=(d_2)$ \label{item:3}
\end{enumerate}

\textbf{\pref{item:1}.} Mutate $R$ in direction $v_2$ to obtain 
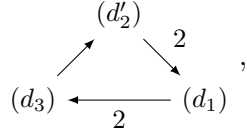
\[ \begin{tikzpicture}[baseline=(current  bounding  box.center),node distance=1.6cm, auto, font=\footnotesize]                                          
  \node (A7) []{$(d_2')$};                                                                                                                     
  \node (A1) [below left of = A7]{$(d_3)$};                                                                               
  \node (A14) [below right of = A7]{$ (d_1)$};   
  \draw[latex-] (A7) to node {} (A1);      
  \draw[latex-] (A14) to node [swap] {$2$} (A7);    
  \draw[-latex] (A14) to node [] {$2$} (A1);      
\end{tikzpicture},
\]
where $(d_2') \geq 2(d_1)-(d_2) > (d_1)$. Next mutate in direction $v_1$ to obtain
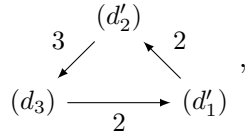
\[ \begin{tikzpicture}[baseline=(current  bounding  box.center),node distance=1.6cm, auto, font=\footnotesize]                                          
  \node (A7) []{$(d_2')$};                                                                                                                     
  \node (A1) [below left of = A7]{$(d_3)$};                                                                               
  \node (A14) [below right of = A7]{$ (d_1')$};   
  \draw[-latex] (A7) to node [swap] {$3$} (A1);      
  \draw[-latex] (A14) to node [swap] {$2$} (A7);    
  \draw[latex-] (A14) to node [] {$2$} (A1);      
\end{tikzpicture},
\]
where $(d_1') \geq 2(d_2')-(d_1) > (d_2')$. The result is now easy to see by induction. Notice in this case that we did not need any assumptions about any edge weights or directions other than $|w_Q (v_1, v_2)|$, which is always $2$ under this mutation path. However, we already know that our quiver will remain cyclic by Corollary \ref{cor:quiver_remains_cyclic} (so that there will always be exactly one incoming and one outgoing weighted arrow adjacent to the vertex we are mutating at), and this means we do not even need to assume that $(d_3) \geq 0$; we can always use the weight $2$ arrow between $v_1$ and $v_2$ to get a lower bound on the degree we will next obtain.

\textbf{\pref{item:2}}. Mutate $R$ in direction $v_1$ to obtain
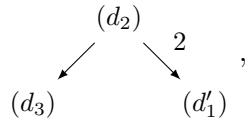
\[ \begin{tikzpicture}[baseline=(current  bounding  box.center),node distance=1.6cm, auto, font=\footnotesize]                                          
  \node (A7) []{$(d_2)$};                                                                                                                     
  \node (A1) [below left of = A7]{$ (d_3)$};                                                                               
  \node (A14) [below right of = A7]{$(d_1')$};   
  \draw[-latex] (A7) to node {} (A1);      
  \draw[latex-] (A14) to node [swap] {$2$} (A7);    
\end{tikzpicture},
\]
where $(d_1') \geq 2(d_2)-(d_1)>(d_2)$. Next mutate in direction $v_3$ to obtain
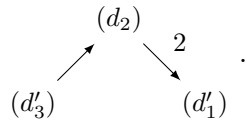
\[ \begin{tikzpicture}[baseline=(current  bounding  box.center),node distance=1.6cm, auto, font=\footnotesize]                                          
  \node (A7) []{$(d_2)$};                                                                                                                     
  \node (A1) [below left of = A7]{$(d_3')$};                                                                               
  \node (A14) [below right of = A7]{$(d_1')$};   
  \draw[latex-] (A7) to node {} (A1);      
  \draw[latex-] (A14) to node [swap] {$2$} (A7);    
\end{tikzpicture}.
\]
This subquiver is $R^{\text{op}}$, but now satisfying case \pref{item:1}, so we are done.

\textbf{\pref{item:3}.} As we alternate mutations between vertices $v_2$ and $v_1$, the weight of the edge between $v_3$ and $v_2$ increases by one each time. Now, suppose after $n$ mutations that we have not obtained a new degree, that is, we have obtained the degree $(d_1)=(d_2)$ again after each mutation. Then the next time we mutate at vertex $v_2$, we obtain the lower bound $m(d_3) - (d_2)$ for the new degree, where $m$ is the edge weight between $v_3$ and $v_2$ in the current degree subquiver. As we continue to mutate, we must either eventually obtain a new degree which is strictly greater than $(d_1)$ (we always obtain one that is at least as great), or otherwise the value of $m$ will eventually grow large enough that  $m(d_3) - (d_2)>(d_1)$, since in this case $(d_3) \geq 1$. Once this occurs, we may proceed in a similar way to the previous cases.
\end{proof}

\begin{cor} \label{cor:sufficient_inf_degrees}
Suppose that $\A$ is a graded cluster algebra generated by the initial quiver $Q$ and that $Q$ contains a subquiver $R$ satisfying the conditions of Proposition \ref{prop:growing_degree_subquiver} or Proposition \ref{prop:degree_subquiver_tends_infty_relaxed}. Then $\A$ contains cluster variables of infinitely many different degrees.
\end{cor}

\begin{rmk}
Of course, if we mutate a degree subquiver satisfying the hypothesis of Proposition \ref{prop:growing_degree_subquiver} or of Proposition \ref{prop:degree_subquiver_tends_infty_relaxed}, we obtain another degree subquiver whose presence is an equivalent condition for showing that a cluster algebra has variables of infinitely many degrees. So Corollary \ref{cor:sufficient_inf_degrees} still holds if we replace $Q$ by a quiver mutation equivalent to it.
\end{rmk}

\begin{rmk}
It may well be possible to further relax the assumptions in Proposition \ref{prop:degree_subquiver_tends_infty_relaxed} and still retain the infinitely many degrees property. In some situations, we may not even need the quiver $Q$ to be positive outside of the subquiver $R$. Having weaker assumptions here would make it more likely that we could read off from a given degree quiver that the associated graded cluster algebra has infinitely many degrees, and so this line of enquiry may be a good candidate for further research.
\end{rmk}

 Finally, we give a sufficient condition for a graded cluster algebra to have all positive degrees, based on the presence of a degree subquiver. As with  Proposition \ref{prop:degree_subquiver_tends_infty_relaxed}, this will be applied in Chapter \ref{chp:matrix_grassmannian}.

\begin{ntn} \label{ntn:flowing_triangles}
Consider the three cyclic subquivers
\begin{equation}
\label{form:stable_triangle1}
\Image{
\begin{tikzpicture}[node distance=1.6cm, auto, font=\footnotesize]                                          
  \node (A1) []{$v_1$};                                                                                                                                                                                              
  \node (A3) [below right of = A1]{$v$};   
  \node (A2) [above right of  = A3]{$v_2$};        
  \draw[-latex] (A1) to node {2} (A2);      
  \draw[-latex] (A3) to node [] {$k+1$} (A1); 
  \draw[ latex-] (A3) to node [swap] {$k$} (A2);   
\end{tikzpicture}},
\end{equation}
\begin{equation}
\label{form:stable_triangle2}
\Image{
\begin{tikzpicture}[node distance=1.6cm, auto, font=\footnotesize]                                          
  \node (A1) []{$v_1$};                                                                                                                                                                                              
  \node (A3) [below right of = A1]{$v$};   
  \node (A2) [above right of  = A3]{$v_2$};        
  \draw[-latex] (A1) to node {2} (A2);      
  \draw[ -latex] (A3) to node [] {$k$} (A1); 
  \draw[latex-] (A3) to node [swap] {$k+1$} (A2);   
\end{tikzpicture}},
\end{equation}
and
\begin{equation}
\label{form:stable_triangle3}
\Image{
\begin{tikzpicture}[node distance=1.6cm, auto, font=\footnotesize]                                          
  \node (A1) []{$v_1$};                                                                                                                                                                                              
  \node (A3) [below right of = A1]{$v$};   
  \node (A2) [above right of  = A3]{$v_2$};        
  \draw[-latex] (A1) to node {2} (A2);      
  \draw[ -latex] (A3) to node [] {} (A1); 
  \draw[ latex-] (A3) to node [] {} (A2);   
\end{tikzpicture}},
\end{equation}
where $k \geq 0$.

Given a quiver $Q$, we will use the notation
\begin{equation}
\label{form:stable_triangle_notation1}
\striL{$v_1$}{$v_2$},
\end{equation}
to mean that for each vertex $v \in Q$ that shares a (non-zero) arrow with either $v_1$ or $v_2$, the triangular subquiver formed by the three vertices is of the form \pref{form:stable_triangle1} or \pref{form:stable_triangle3} (where $k$ is allowed to be different for each $v$). So for each such triangular subquiver, either $W=w+1$, or $W=w=1$.
Similarly, we will use
\begin{equation}
\label{form:stable_triangle_notation2}
\striR{$v_1$}{$v_2$},
\end{equation}
to mean that  for each vertex $v$ that shares an arrow with $v_1$ or $v_2$, the triangular subquiver formed by the three vertices is of the form \pref{form:stable_triangle2} or \pref{form:stable_triangle3}. 

When we refer to \pref{form:stable_triangle_notation1} or \pref{form:stable_triangle_notation2} as being a subquiver of $Q$, we will mean the subquiver containing the vertices $v_1$ and $v_2$ along with all other vertices that share a (non-zero) arrow with $v_1$ or $v_2$ (and by referring to such a subquiver, we will mean that $Q$ has the corresponding property specified above).
\end{ntn} 

Note that mutating \pref{form:stable_triangle1} at $v_1$ yields a subquiver of the form \pref{form:stable_triangle2} with $k$ increased by 1, mutating \pref{form:stable_triangle2} at $v_2$ yields a subquiver of the form \pref{form:stable_triangle1} with $k$ increased by 1
, and mutating \pref{form:stable_triangle3} at either $v_1$ or $v_2$ does not change its form (aside from reversing arrows).
Hence, mutating the subquiver \pref{form:stable_triangle_notation1} at $v_1$ gives a subquiver of the form \pref{form:stable_triangle_notation2}, and mutating the subquiver  \pref{form:stable_triangle_notation2} at $v_2$ gives a subquiver of the form \pref{form:stable_triangle_notation1}. So $[v_2,v_1]$ applied to a subquiver of the form \pref{form:stable_triangle_notation1} is still of the form \pref{form:stable_triangle_notation1}, and similarly for $[v_1,v_2]$ applied to a quiver of the form \pref{form:stable_triangle_notation2}.

\begin{lem}  
\label{lem:linear_increasing_degree_quiver}
Suppose $Q$ is a degree quiver with the subquiver 
\begin{equation}
\label{form:linear_increasing_degree_quiver}
Q'=\striR{$(d_1)_{v_1}$}{$(d_2)_{v_2}$}.
\end{equation}
Then mutation along the path $[(v_1,v_2)_n]$ yields the degree 
\begin{equation}
(d_1-d_2) n + d_1.
\end{equation}
Similarly, mutation of
\begin{equation}
\label{form:linear_increasing_degree_quiver}
Q'=\striL{$(d_1)_{v_1}$}{$(d_2)_{v_2}$},
\end{equation}
along  $[(v_2,v_1)_n]$ yields degree
\begin{equation}
(d_2-d_1) n + d_2.
\end{equation}
\end{lem}

\begin{proof}
It is clear that the result is true for $n=1$ and $n=2$. Suppose the result is true for $n$ and assume without loss $n$ is even (which we may do by the comments in the paragraph below Notation \ref{ntn:flowing_triangles}). Then the degree quiver is of the form
\begin{center}
\begin{tikzpicture}[node distance=1.6cm, auto, font=\footnotesize]                                          
  \node (A1) []{$\Large( [d_1- d_2 ]n+d_1 \Large)$};    
  \node (A3) [below right = 1.8cm of A1]{$*$};                                                                                                                      
  \node (A2) [above right  = 1.8cmof A3]{$\left( [d_1- d_2 ](n-1)+d_1 \right)$};                                                                               
  \draw[-latex] (A1) to node {2} (A2);      
  \draw[dotted, -latex] (A3) to node [] {$w$} (A1); 
  \draw[dotted, latex-] (A3) to node [swap] {$W$} (A2);   
\end{tikzpicture}
\end{center}
and our next mutation is at $v_2$. Carrying out this mutation, we obtain degree 
\[2[  (d_1-d_2)n+d_1 ] -  [(d_1-d_2)(n-1)+d_1] = (d_1-d_2)(n+1)+d_1,\]
 from which the result follows by induction. The result for [$(v_2, v_1)_n$] is immediate from this, as the direction of the weight $2$ arrow does not matter.
\end{proof}

\section{New degree quivers from old}

In the sections above we saw how the presence of a degree subquiver gave us information about cluster algebras arising from larger quivers containing the subquiver. A natural question that arises is whether, in general, understanding properties of rank $n$ graded cluster algebras can be reduced to understanding properties of lower rank cluster algebras. More precisely, we might ask the following:

\begin{qn} \label{qn:reduction_theorem}
Let $R$ be a degree quiver and consider the graded cluster algebra $\A(R)$ it gives rise to. Suppose $Q$ is a degree quiver that contains $R$ as a subquiver. Which properties of the graded cluster algebra $\A(R)$ are inherited by $\A(Q)$? For instance, if $\A(R)$ has infinitely many variables of each occurring degree, must $\A(Q)$ have a corresponding set of degrees each associated with infinitely many variables?

 If, rather than containing $R$ as a subquiver, $Q$ can be constructed from a combination of smaller degree quivers in some way, which properties of these smaller quivers are inherited by $\A(Q)$?
\end{qn}

Progress in this direction would greatly increase our understanding of gradings in general, though establishing an answer to the question does not appear to be easy.
While we don't attempt to tackle this here, we do note an interesting property about degree quivers: that they can often be broken down into sums of smaller degree quivers.

\begin{dfn}
Let $Q$ and $R$ be degree quivers and let $Q_0$ and $R_0$ be the sets of vertices of $Q$ and $R$, respectively. Let $M$ be a subset of $Q_0 \times R_0$, whose elements will be called \emph{matches}, such that for each match $(v,\omega) \in M$, we have that $\deg(v)=\deg(\omega)$ and also that there is no other match involving either $v$ or $\omega$ in $M$. We will call the vertices that do not appear in any match in $M$ \emph{unmatched}. Denote the unmatched vertices of $Q$ and $R$ by $U_Q$ and $U_R$, respectively.

We define $Q+_M R$ to  be the quiver obtained from $Q$ and $R$ as follows: 
\begin{enumerate}[(1)]
\item For each match $(v,\omega) \in M$, create a new vertex $m_{v\omega}$ (in other words, identify vertex $v$ in $Q$ with vertex $\omega$ in $R$). Then define the set of vertices of $Q+_M R$ to be 
\[U_Q \sqcup U_R \sqcup \{m_{v\omega} | (v,\omega) \in M \}. \]
\item For each pair of vertices $m_{v_1 \omega_1}$ and $m_{v_2 \omega_2}$, create an arrow from $m_{v_1 \omega_1}$ to $m_{v_2 \omega_2}$ of weight $w_Q(v_1,v_2) + w_R(\omega_1,\omega_2)$ (with a negatively weighted arrow in one direction meaning a positively weighted arrow in the opposite direction, as usual). For each pair of unmatched vertices $u_1$ and $u_2$, if both come from the same original quiver, set the weight of the arrow from $u_1$ to $u_2$ to its original value, $w(u_1,u_2)$, and if both come from different quivers, do not create an arrow between them.
\end{enumerate}
\end{dfn}

\begin{prop}
Let $Q$ and $R$ be degree quivers. Then $Q+_M R$ is also a degree quiver.
\end{prop}

\begin{proof}
Let $t$ be a vertex in $Q +_M R$. We just need to check that $t$ is balanced, that is, satisfies Equation \ref{eqn:degree_quiver_balanced_requirement}. If $t$ is an unmatched vertex then it is clearly balanced since it was in its original quiver. Suppose $t=m_{v\omega}$ for some match $(v,\omega) \in M$. Then the weight into $m_{v\omega}$ is $w_{\rightarrow}(v)+w_{\rightarrow}(\omega)$ while the weight out of $m_{v\omega}$ is $w_{\leftarrow}(v)+w_{\leftarrow}(\omega)$. Since $w_{\rightarrow}(v) = w_{\leftarrow}(v)$ and $w_{\rightarrow}(\omega) = w_{\leftarrow}(\omega)$, we have $w_{\rightarrow}(m_{v\omega}) =w_{\leftarrow}(m_{v\omega})$, so $t$ is balanced.
\end{proof}

\begin{ex} \label{ex:sum_of_degree_quivers}
Let
\[
Q=
\begin{tikzpicture}[baseline=(current  bounding  box.center),node distance=1.6cm, auto, font=\footnotesize]                                          
  \node (A7) []{$(3)_{v_2}$};                                                                                                                     
  \node (A1) [below left of = A7]{$(5)_{v_1}$};                                                                               
  \node (A14) [below right of = A7]{$ (4)_{v_3}$};   
  \draw[-latex] (A7) to node [swap] {$4$} (A1);      
  \draw[-latex] (A14) to node [swap] {$5$} (A7);    
  \draw[latex-] (A14) to node [] {$3$} (A1);      
\end{tikzpicture}, 
\text{     } R=
\begin{tikzpicture}[baseline=(current  bounding  box.center),node distance=1.6cm, auto, font=\footnotesize]                                          
  \node (BL) []{$(4)_{\omega_1}$};                                                                                                                     
  \node (TL) [above of = BL]{$(3)_{\omega_2}$};                                                                               
  \node (TR) [right of = TL]{$ (4)_{\omega_3}$};   
  \node (BR) [below of = TR]{$ (5)_{\omega_4}$};     
  \draw[-latex] (BL) to node [] { 5} (TL);          
  \draw[-latex] (TL) to node [] { 5} (TR);  
  \draw[-latex] (TR) to node [] {3 } (BR);  
  \draw[-latex] (BR) to node [] {3 } (BL);        
\end{tikzpicture} 
\]
and $M=\{(v_3,\omega_1),(v_2,\omega_2) \}.$
Then 
\[Q+_MR=
\begin{tikzpicture}[baseline=(current  bounding  box.center),node distance=2cm, auto, font=\footnotesize]                                          
  \node (BL) []{$(4)_{m_{v_3 \omega_1}}$};                                                                                                                     
  \node (L) [above left of = BL]{$ (5)_{v_1}$};        
  \node (TL) [above right of = L]{$(3)_{m_{v_2 \omega_2}}$};                                                                               
  \node (TR) [right =1.3cm of TL]{$ (4)_{\omega_3}$};   
  \node (BR) [right = 1.3cm of BL]{$ (5)_{\omega_4}$};     
  \draw[-latex] (BL) to node [swap] {10} (TL);          
  \draw[-latex] (TL) to node [] {5 } (TR);  
  \draw[-latex] (TR) to node [] { 3} (BR);  
  \draw[-latex] (BR) to node [] { 3} (BL);        
  \draw[latex-] (BL) to node [] {3} (L);    
  \draw[latex-] (L) to node [] {4} (TL);      
\end{tikzpicture} \]
is a degree quiver. The vertex $m_{v_2 \omega_2}$ is balanced since 
\[w_{\rightarrow}(m_{v_2 \omega_2}) =10\cdot 4 = 5 \cdot 4+ 4 \cdot 5 = w_{\leftarrow}(m_{v_2 \omega_2}).\] Similarly, so is $m_{v_3 \omega_1}$.
\end{ex}

\begin{dfn}
We will call a degree quiver $S$ \emph{irreducible} if there are no degree quivers $Q$ and $R$ with $Q_0 \subsetneq S_0$ and $R_0 \subsetneq S_0$ such that $S=Q+_M R$ for some $M$, where $Q_0$, $R_0$ and $S_0$ denote the sets of vertices of $Q$, $R$ and $S$, respectively.
\end{dfn}

\begin{ex}
The degree quiver $Q$ in Example \ref{ex:sum_of_degree_quivers} is irreducible as any degree quiver that is an appropriate potential summand of $Q$ must have two vertices, but any degree quiver with two vertices must assign all vertices degree $0$ in order to be balanced. 

On the other hand, the degree quiver $R$ in Example \ref{ex:sum_of_degree_quivers} is not irreducible since
\[
 R=
\begin{tikzpicture}[baseline=(current  bounding  box.center),node distance=1.6cm, auto, font=\footnotesize]                                          
  \node (BL) []{$(4)_{\omega_1}$};                                                                                                                     
  \node (TL) [above of = BL]{$(3)_{\omega_2}$};                                                                               
  \node (TR) [right of = TL]{$ (4)_{\omega_3}$};   
  \node (BR) [below of = TR]{$ (5)_{\omega_4}$};     
  \draw[-latex] (BL) to node [] { 5} (TL);          
  \draw[-latex] (TL) to node [] { 5} (TR);  
  \draw[-latex] (TR) to node [] {3 } (BR);  
  \draw[-latex] (BR) to node [] {3 } (BL);        
\end{tikzpicture} 
=
\begin{tikzpicture}[baseline=(current  bounding  box.center),node distance=1.6cm, auto, font=\footnotesize]                                          
  \node (BL) []{$(4)_{v_3}$};                                                                                                                     
  \node (T) [above of = BL]{$(3)_{v_2}$};                                                                               
  \node (BR) [ right of = BL]{$ (5)_{v_1}$};   
  \draw[-latex] (BL) to node [] {$5$} (T);      
  \draw[-latex] (BR) to node [] {$3$} (BL);    
  \draw[-latex] (T) to node [] {$4$} (BR);      
\end{tikzpicture}
+_M
\begin{tikzpicture}[baseline=(current  bounding  box.center),node distance=1.6cm, auto, font=\footnotesize]                                          
  \node (B) []{$(5)_{v_1}$};                                                                                                                     
  \node (TR) [above of = B]{$(4)_{v_3}$};                                                                               
  \node (TL) [ left of = TR]{$ (3)_{v_2}$};   
  \draw[-latex] (TR) to node [] {$3$} (B);      
  \draw[-latex] (B) to node [] {$4$} (TL);    
  \draw[-latex] (TL) to node [] {$5$} (TR);      
\end{tikzpicture}, 
\]
where $M= \{ (v_1,v_1), (v_2,v_2) \}$. Thus $R= Q+_M Q^{\text{op}}$.
\end{ex}

For degree quivers $S$ that are not irreducible, it might be hoped that if $S= Q+_M R$, then 
\begin{equation} \label{eqn:false_reduction_theorem}
S_{\pth{p}} = Q_{\pth{p}}+_M R_{\pth{p}}, 
\end{equation}
but this turns out not to be the case. However, as we see in the examples below, there are some restricted settings in which information about $\A(S)$ can obviously be obtained from information about $\A(Q)$ and $\A(R)$. It is possible that progress towards Question \ref{qn:reduction_theorem} could involve attempting to expand these restricted settings or attempting to find the correct relationship, similar to Equation \ref{eqn:false_reduction_theorem}, between degree quivers and their summands.

\begin{ex}
Let $Q$ be any degree quiver with a vertex $\omega$ with $\deg(\omega)=2$ and let
\[R=\begin{tikzpicture}[baseline=(current  bounding  box.center),node distance=1.6cm, auto, font=\footnotesize]                                          
  \node (A7) []{$(2)_{v_2}$};                                                                                                                     
  \node (A1) [below left of = A7]{$(2)_{v_1}$};                                                                               
  \node (A14) [below right of = A7]{$ (2)_{v_3}$};   
  \draw[-latex] (A7) to node [swap] {$2$} (A1);      
  \draw[-latex] (A14) to node [swap] {$2$} (A7);    
  \draw[latex-] (A14) to node [] {$2$} (A1);      
\end{tikzpicture}. \]
So $R$ is the quiver of Markov type---recall from Chapter \ref{chp:3_by_3} that the associated graded cluster algebra has infinitely many variables, all of degree $2$. It is easy to show that the mutation path $[(v_2,v_1)_r]$ gives infinitely many variables of degree $2$ in $\A(R)$. Let $M=\{(\omega,v_3)\}$. Then the graded cluster algebra $\A(Q+_M R)$ must also have infinitely many variables of degree $2$. This is since, when mutating along $[(v_2,v_1)_r]$ in $Q+_MR$, the subquiver corresponding to $Q$ cannot change in any way (as we would need to mutate at the vertex $(\omega,v_3)$ in order to affect $Q$), while the variables we obtain are essentially the same as those obtained under the same mutation in $\A(R)$, and of the same degree.
\end{ex}

\begin{ex}
For another example, let ${A=\left(\begin{smallmatrix}0&a&-c\\-a& 0&b\\c&-b&0\end{smallmatrix} \right)}$ with $a \geq b \geq c\geq 3$. Then, as it is mutation-infinite, $A$ gives rise to a cluster algebra with infinitely many degrees by Corollary \ref{cor:mutation_infinite_3x3_iff_inf_degrees}. Let us show that infinitely many degrees can produced by alternately mutating between only two particular directions. Denote $\mmut{A}{(1,2)^k}$ by $A^{(k)}$ and denote the $(i,j)$ entry of $A^{(k)}$ by $A^{(k)}_{ij}$. We will show that the sequence
\begin{equation}
\label{eqn:infdegincreasing}
 |A^{(0)}_{32} |, |A^{(0)}_{31} |, |A^{(1)}_{32} |, |A^{(1)}_{31} |,\dots 
 \end{equation}
is strictly increasing (excluding the first term if $a=b$). Mutating $A$ along the path $[1,2]$ gives $A^{(1)}=\left(\begin{smallmatrix}0&a&c'\\-a& 0&-b'\\-c'&b'&0\end{smallmatrix} \right)$, where $c'=c-ab$ and $b'=b+ac'$. 
Since $a \geq b \geq c \geq 3$, $c'$ is negative. We have $b+c'<0$ since $c'=c-ab=(c-b)-(a-1)b$. Then $b'=b+ac'=b+c'+(a-1)c' < (a-1)c'<c'$. So $b'$ and $c'$ are both negative and $|b'|>|c'|$. Mutating $A^{(1)}$ along $[1,2]$ again we obtain $A^{(2)}=\left(\begin{smallmatrix}0&a&c''\\-a& 0&-b''\\-c''&b''&0\end{smallmatrix} \right)$, where $b''$ and $c''$ satisfy the same respective relations as $b'$ and $c'$. From this, (\ref{eqn:infdegincreasing}) now holds by induction. Thus, by Lemma \ref{lem:degree_is_in_matrix}, we obtain an infinite increasing sequence of degrees by alternately mutating in directions $2$ and $1$.

Suppose $R$ is the quiver associated to $A$ (and label it such that $v_3$ does not correspond to one of the alternating directions). Let $Q$ be any degree quiver with a vertex $\omega$ satisfying $\deg(\omega)=\deg(v_3)$ and let $M=\{(\omega,v_3)\}$. Then $\A(Q+_MR)$ must also have infinitely many degrees. In other words, attaching $R$ to another quiver by vertex $v_3$ endows the arising cluster algebra with infinitely many degrees. For example, we can take
\[R=\begin{tikzpicture}[baseline=(current  bounding  box.center),node distance=1.6cm, auto, font=\footnotesize]                                          
  \node (A7) []{$(3)_{v_2}$};                                                                                                                     
  \node (A1) [below left of = A7]{$(4)_{v_1}$};                                                                               
  \node (A14) [below right of = A7]{$ (5)_{v_3}$};   
  \draw[latex-] (A7) to node [swap] {$5$} (A1);      
  \draw[latex-] (A14) to node [swap] {$4$} (A7);    
  \draw[-latex] (A14) to node [] {$3$} (A1);      
\end{tikzpicture},  \]
for which the two alternating vertices to mutate at are $v_2$ and $v_1$.
\end{ex}

\biblio

\chapter{Gradings on cluster algebra structure for coordinate algebras of matrices and Grassmannians}  \label{chp:matrix_grassmannian}
In this chapter we consider gradings on cluster algebra structures for coordinate algebras of matrices and Grassmannians and prove that, for the infinite type cases, they contain variables of infinitely many different degrees. We also prove that all positive degrees occur in these cluster algebras. 

While we only deal with the classical cluster algebra structures here, any results about gradings on their quantum analogues are also determined by the classical case. This follows from \cite[Theorem 6.1]{BZ} and \cite[Remark 3.6]{GL}.

At various points we use some fairly long mutation paths to prove certain results. Given such paths, it is not difficult to prove these results by hand, but the reader may wonder how such paths were found. So it is worth mentioning that specific code was written to help discover some of these paths, or other patterns that led to them. As noted previously, this code is available with the electronic version of the thesis.

\section{Graded cluster algebra structure for matrices}
\begin{dfn}
Let $M(k,l)$ be the set of $k \times l$ matrices. The coordinate algebra of the $k \times l$ matrices is $\mathcal{O}(M(k,l)) =\C[ x_{i,j}]$, for $1 \leq i \leq k$, $1 \leq j \leq l$, where $x_{ij}$ is the \emph{coordinate function} $M(k,l) \rightarrow \C$ given by $A \mapsto a_{ij}$.
\end{dfn}

$\mathcal{O}(M(k,l))$ has the structure of a cluster algebra, $\A \big(\mathcal{O}(M(k,l)) \big)$, as follows. For row set $I$ and column set $J$ with $|I|=|J|$, let $ {\scriptsize \left[\begin{matrix}J \\ I \end{matrix} \right] }$ denote the corresponding minor of the matrix
\[
\begin{pmatrix}
   x_{11} & x_{12} & \dots & x_{1l} \\
   \vdots & \vdots & \ddots & \vdots \\
   x_{k1} & x_{k2} & \dots & x_{kl}
\end{pmatrix}.
\]

For $1 \leq r \leq k$ and $1 \leq s \leq l$, define the sets 
\begin{eqnarray}
R(r,s)= &\{ k-r+1,k-r+2, \dots, k-r+s\} \cap \{1, \dots, k \},\\
C(r,s) = &\{l-s+1, l-s+2, \dots, l-s+r \} \cap \{1, \dots, l \}.
\end{eqnarray}

For a cluster algebra of this cardinality, it is more efficient to define the initial seed using a quiver. We show the initial quiver for $M(4,4)$ (which we will denote by $Q_{M(4,4)}$) in Figure \ref{fig:init_quiver_M(4,4)} below.  This generalises to $M(k,l)$ by replacing the $4 \times 4$ grid by a $k \times l$ grid of the same form, again with only the bottom and rightmost vertices frozen (see also \cite[p.~715]{GL}). 
The initial cluster is given by defining the entry of the cluster corresponding to the vertex in position $(i,j)$ to be equal to $ {\scriptsize \left[\begin{matrix} C(i,j) \\ R(i,j) \end{matrix} \right] }$. The initial grading vector is then given by assigning the vertex in the $(i,j)$ position degree $\min(i,j)$.

It is easy to check that this definition gives a valid degree quiver. Moreover, since $\mathcal{O}(M(k,l))$ is an $\N$-graded algebra (\cite[p.~716]{GL}) we have that the corresponding cluster algebra is also $\N$-graded (as, similarly, is the cluster algebra corresponding to $\mathcal{O}(Gr(k,k+l))$ which we will introduce in Section \ref{sec:corollaries_for_grassmannians}). For a justification that this initial structure gives rise to a cluster algebra that agrees with the coordinate ring of (quantum) matrices, see \cite[Corollary 12.10]{GLS} or \cite{GL}.

\section{Infinitely many degrees}
We first consider the case for when $k \geq 4$ and $l \geq 4$. Without loss of generality, we may consider the smallest such case: $M(4,4)$. We prove our result by considering a certain subquiver of $Q_{M(4,4)}$. Larger cases will always contain the same subquiver, and it is thus possible to show they also contain infinitely many degrees in exactly the same way.

The initial exchange quiver, $Q_{M(4,4)}$, is given below. (Recall that we denote frozen vertices by placing a square around them.)

\begin{figure}[H]
\begin{center}
\begin{tikzpicture}[node distance=1.8cm, auto, font=\footnotesize]    
  \node (AA) []{1};   
  \node (AB) [right of = AA]{2}; 
  \node (AC) [right of = AB]{3};   
  \node (AD) [solid node, rectangle, fill=white, right of = AC]{16};  
  \node (BA) [below of = AA]{4};     
  \node (BB) [right of = BA]{5}; 
  \node (BC) [right of = BB]{6};   
  \node (BD) [solid node, rectangle, fill=white, right of = BC]{15}; 
  \node (CA) [below of = BA]{7};     
  \node (CB) [right of = CA]{8}; 
  \node (CC) [right of = CB]{9};   
  \node (CD) [solid node, rectangle, fill=white, right of = CC]{14}; 
  \node (DA) [solid node, rectangle, fill=white, below of = CA]{10};     
  \node (DB) [solid node, rectangle, fill=white, right of = DA]{11}; 
  \node (DC) [solid node, rectangle, fill=white, right of = DB]{12};   
  \node (DD) [solid node, rectangle, fill=white, right of = DC]{13}; 
  
  \draw[-latex] (AA) to node {} (AB);             
  \draw[-latex] (AB) to node {} (AC);
  \draw[-latex] (AC) to node {} (AD);
  \draw[-latex] (BA) to node {} (BB);             
  \draw[-latex] (BB) to node {} (BC);
  \draw[-latex] (BC) to node {} (BD);
  \draw[-latex] (CA) to node {} (CB);             
  \draw[-latex] (CB) to node {} (CC);
  \draw[-latex] (CC) to node {} (CD);        

  \draw[-latex] (AA) to node {} (BA);             
  \draw[-latex] (BA) to node {} (CA);             
  \draw[-latex] (CA) to node {} (DA);                        
  \draw[-latex] (AB) to node {} (BB);             
  \draw[-latex] (BB) to node {} (CB);             
  \draw[-latex] (CB) to node {} (DB);              
  \draw[-latex] (AC) to node {} (BC);             
  \draw[-latex] (BC) to node {} (CC);             
  \draw[-latex] (CC) to node {} (DC);                            
  
  \draw[-latex] (BB) to node {} (AA);   
  \draw[-latex] (BC) to node {} (AB);     
  \draw[-latex] (BD) to node {} (AC);   
  \draw[-latex] (CB) to node {} (BA);   
  \draw[-latex] (CC) to node {} (BB);     
  \draw[-latex] (CD) to node {} (BC);    
  \draw[-latex] (DB) to node {} (CA);   
  \draw[-latex] (DC) to node {} (CB);     
  \draw[-latex] (DD) to node {} (CC);         
                               
\end{tikzpicture}
\end{center}	
\captionof{figure}{Initial quiver for $M(4,4)$.}
\label{fig:init_quiver_M(4,4)}
\end{figure}
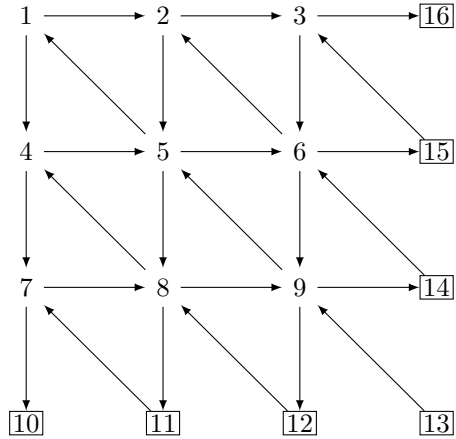

\begin{lem}
\label{lem:MA_4_4_inf_degrees}
The graded cluster algebra structure for $\mathcal{O}(M(4,4))$ has variables in infinitely many degrees.
\end{lem}

\begin{proof} 
We will show that the mutation sequence $[4,8,2,9,5,6,1,5,2]$ transforms $Q$ into a quiver to which we can apply Proposition \ref{prop:degree_subquiver_tends_infty_relaxed}. Note that embedding the same sequence into larger initial quivers (we can embed the subquiver in Figure \ref{fig:subquiver_M(4,4)} at the bottom left of the larger quiver) allows us to immediately generalise any results to $\mathcal{O}(M(k,l))$ for $k,l \geq 4$. (Of course, using the exact same mutation sequence as above would not work in general, since the vertices of the corresponding subquiver are relabelled when embedded in a larger initial quiver, but it is clear how to adjust the mutation path to compensate for this.)

Consider the following subquiver of $Q(M(4,4))$.
\begin{figure}[H]
\begin{center}
\begin{tikzpicture}[node distance=1.6cm, auto, font=\footnotesize]                                          
  \node (A5) []{5};                                                                                                                     
  \node (A4) [left of = A5]{4};                                                                               
  \node (A6) [right of = A5]{6};   
  \node (A2) [above of = A5]{2};     
  \node (A8) [below of = A5]{8};     
  \node (A1) [above of = A4]{1};     
  \node (A9) [right of = A8]{9};     
  \node (A11) [solid node, rectangle, fill=white,below of = A8]{11};       
  \draw[-latex] (A1) to node {} (A2);        
  \draw[-latex] (A4) to node {} (A5);      
  \draw[-latex] (A5) to node {} (A6);      
  \draw[-latex] (A8) to node {} (A9);      
  \draw[-latex] (A1) to node {} (A4);      
  \draw[-latex] (A2) to node {} (A5);      
  \draw[-latex] (A5) to node {} (A8);      
  \draw[-latex] (A8) to node {} (A11);      
  \draw[-latex] (A6) to node {} (A9);      
  \draw[-latex] (A5) to node {} (A1);      
  \draw[-latex] (A9) to node {} (A5);      
  \draw[-latex] (A8) to node {} (A4);        
  \draw[-latex] (A6) to node {} (A2);     
\end{tikzpicture}
\end{center}	
\captionof{figure}{Subquiver of $Q(M(4,4))$.}
\label{fig:subquiver_M(4,4)}
\end{figure}
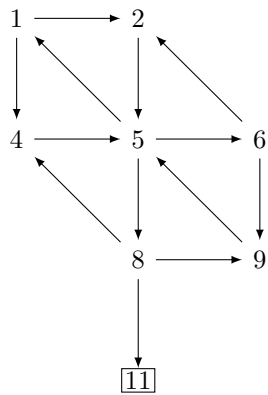
We now perform the sequence of mutations to this subquiver, as shown in the following diagram. After a mutation, we will often no longer be interested in the mutated vertex, and if so we will remove it from our mutated subquiver.
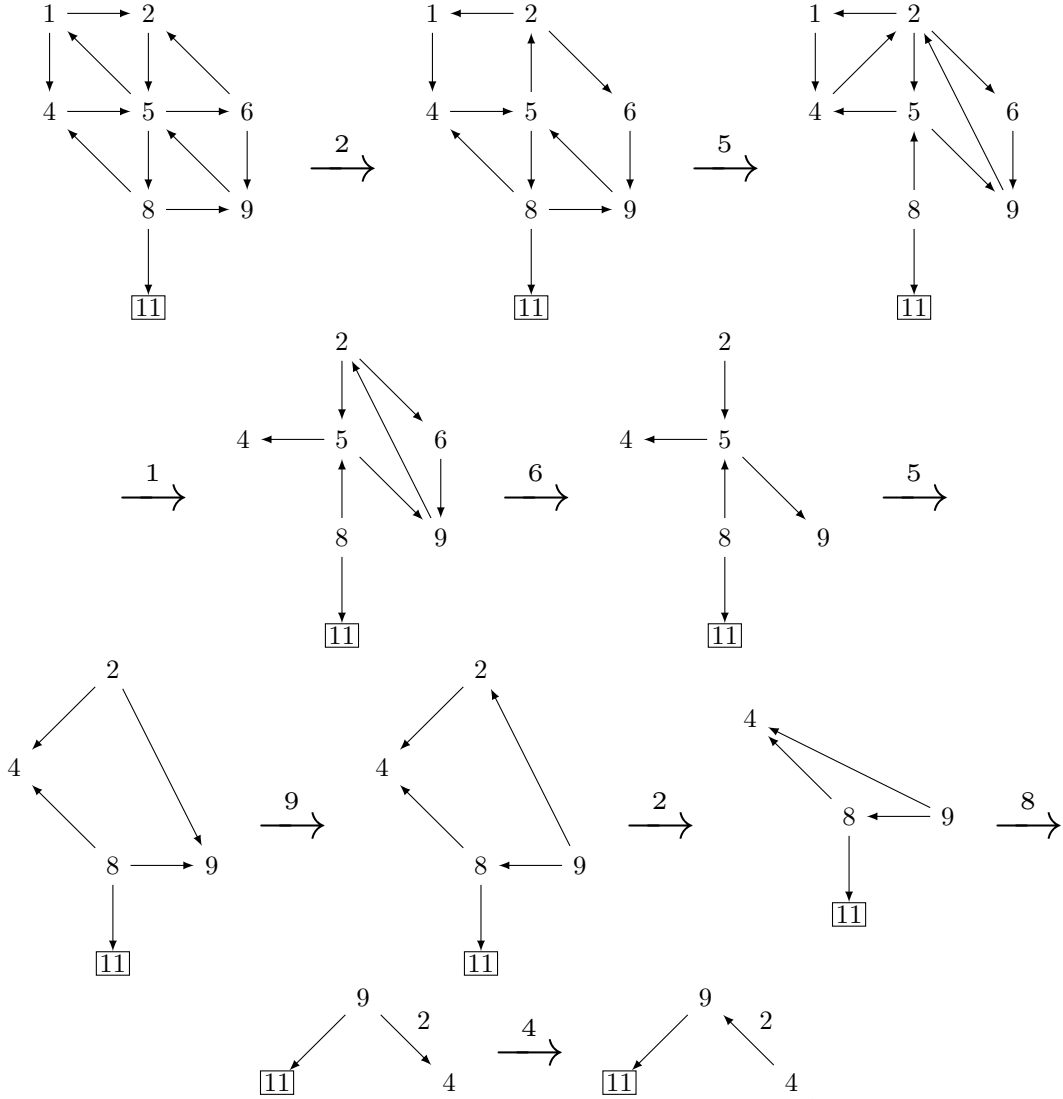
\begin{center}
\Image{
\begin{tikzpicture}[node distance=1.3cm, auto, font=\footnotesize]    
  \node (A5) []{5};                                                                                                                     
  \node (A4) [left of = A5]{4};                                                                               
  \node (A6) [right of = A5]{6};   
  \node (A2) [above of = A5]{2};     
  \node (A8) [below of = A5]{8};     
  \node (A1) [above of = A4]{1};     
  \node (A9) [right of = A8]{9};     
  \node (A11) [solid node, rectangle, fill=white,below of = A8]{11};       
  \draw[-latex] (A1) to node {} (A2);        
  \draw[-latex] (A4) to node {} (A5);      
  \draw[-latex] (A5) to node {} (A6);      
  \draw[-latex] (A8) to node {} (A9);      
  \draw[-latex] (A1) to node {} (A4);      
  \draw[-latex] (A2) to node {} (A5);      
  \draw[-latex] (A5) to node {} (A8);      
  \draw[-latex] (A8) to node {} (A11);      
  \draw[-latex] (A6) to node {} (A9);      
  \draw[-latex] (A5) to node {} (A1);      
  \draw[-latex] (A9) to node {} (A5);      
  \draw[-latex] (A8) to node {} (A4);        
  \draw[-latex] (A6) to node {} (A2);     
\end{tikzpicture}
} \Arrow{2}
\Image{
\begin{tikzpicture}[node distance=1.3cm, auto, font=\footnotesize]    
  \node (A5) []{5};                                                                                                                     
  \node (A4) [left of = A5]{4};                                                                               
  \node (A6) [right of = A5]{6};   
  \node (A2) [above of = A5]{2};     
  \node (A8) [below of = A5]{8};     
  \node (A1) [above of = A4]{1};     
  \node (A9) [right of = A8]{9};     
  \node (A11) [solid node, rectangle, fill=white,below of = A8]{11};       
  \draw[latex-] (A1) to node {} (A2);        
  \draw[-latex] (A4) to node {} (A5);      
  \draw[-latex] (A8) to node {} (A9);      
  \draw[-latex] (A1) to node {} (A4);      
  \draw[latex-] (A2) to node {} (A5);      
  \draw[-latex] (A5) to node {} (A8);      
  \draw[-latex] (A8) to node {} (A11);      
  \draw[-latex] (A6) to node {} (A9);        
  \draw[-latex] (A9) to node {} (A5);      
  \draw[-latex] (A8) to node {} (A4);        
  \draw[latex-] (A6) to node {} (A2); 
\end{tikzpicture}
} \Arrow{5}
\Image{
\begin{tikzpicture}[node distance=1.3cm, auto, font=\footnotesize]    
  \node (A5) []{5};                                                                                                                     
  \node (A4) [left of = A5]{4};                                                                               
  \node (A6) [right of = A5]{6};   
  \node (A2) [above of = A5]{2};     
  \node (A8) [below of = A5]{8};     
  \node (A1) [above of = A4]{1};     
  \node (A9) [right of = A8]{9};     
  \node (A11) [solid node, rectangle, fill=white,below of = A8]{11};       
  \draw[latex-] (A1) to node {} (A2);        
  \draw[latex-] (A4) to node {} (A5);          
  \draw[-latex] (A1) to node {} (A4);      
  \draw[-latex] (A2) to node {} (A5);      
  \draw[latex-] (A5) to node {} (A8);      
  \draw[-latex] (A8) to node {} (A11);      
  \draw[-latex] (A6) to node {} (A9);        
  \draw[latex-] (A9) to node {} (A5);        
  \draw[latex-] (A6) to node {} (A2);  
  \draw[latex-] (A2) to node {} (A4);    
  \draw[latex-] (A2) to node {} (A9);        
\end{tikzpicture}
} \Arrow{1}
\Image{
\begin{tikzpicture}[node distance=1.3cm, auto, font=\footnotesize]    
  \node (A5) []{5};                                                                                                                     
  \node (A4) [left of = A5]{4};                                                                               
  \node (A6) [right of = A5]{6};   
  \node (A2) [above of = A5]{2};     
  \node (A8) [below of = A5]{8};        
  \node (A9) [right of = A8]{9};     
  \node (A11) [solid node, rectangle, fill=white,below of = A8]{11};         
  \draw[latex-] (A4) to node {} (A5);           
  \draw[-latex] (A2) to node {} (A5);      
  \draw[latex-] (A5) to node {} (A8);      
  \draw[-latex] (A8) to node {} (A11);      
  \draw[-latex] (A6) to node {} (A9);        
  \draw[latex-] (A9) to node {} (A5);        
  \draw[latex-] (A6) to node {} (A2);  
  \draw[latex-] (A2) to node {} (A9);  
\end{tikzpicture}
} \Arrow{6}
\Image{
\begin{tikzpicture}[node distance=1.3cm, auto, font=\footnotesize]    
  \node (A5) []{5};                                                                                                                     
  \node (A4) [left of = A5]{4};                                                                               
  \node (A2) [above of = A5]{2};     
  \node (A8) [below of = A5]{8};        
  \node (A9) [right of = A8]{9};     
  \node (A11) [solid node, rectangle, fill=white,below of = A8]{11};         
  \draw[latex-] (A4) to node {} (A5);           
  \draw[-latex] (A2) to node {} (A5);      
  \draw[latex-] (A5) to node {} (A8);      
  \draw[-latex] (A8) to node {} (A11);        
  \draw[latex-] (A9) to node {} (A5);        
\end{tikzpicture}
} \Arrow{5}
\Image{
\begin{tikzpicture}[node distance=1.3cm, auto, font=\footnotesize]    
  \node (A5) []{};                                                                                                                     
  \node (A4) [left of = A5]{4};                                                                               
  \node (A2) [above of = A5]{2};     
  \node (A8) [below of = A5]{8};        
  \node (A9) [right of = A8]{9};     
  \node (A11) [solid node, rectangle, fill=white,below of = A8]{11};               
  \draw[-latex] (A8) to node {} (A11);        
  \draw[-latex] (A2) to node {} (A4);
  \draw[-latex] (A8) to node {} (A4);
  \draw[-latex] (A2) to node {} (A9);    
  \draw[-latex] (A8) to node {} (A9);  
\end{tikzpicture}
}\Arrow{9}
\Image{
\begin{tikzpicture}[node distance=1.3cm, auto, font=\footnotesize]    
  \node (A5) []{};                                                                                                                     
  \node (A4) [left of = A5]{4};                                                                               
  \node (A2) [above of = A5]{2};     
  \node (A8) [below of = A5]{8};        
  \node (A9) [right of = A8]{9};     
  \node (A11) [solid node, rectangle, fill=white,below of = A8]{11};               
  \draw[-latex] (A8) to node {} (A11);        
  \draw[-latex] (A2) to node {} (A4);
  \draw[-latex] (A8) to node {} (A4);
  \draw[latex-] (A2) to node {} (A9);    
  \draw[latex-] (A8) to node {} (A9);    
\end{tikzpicture}
}\Arrow{2}
\Image{
\begin{tikzpicture}[node distance=1.3cm, auto, font=\footnotesize]                                      
  \node (A5) []{};                                                                                                                     
  \node (A4) [left of = A5]{4};                                                                               
  \node (A8) [below of = A5]{8};        
  \node (A9) [right of = A8]{9};     
  \node (A11) [solid node, rectangle, fill=white,below of = A8]{11};               
  \draw[-latex] (A8) to node {} (A11);        
  \draw[-latex] (A8) to node {} (A4);
  \draw[latex-] (A8) to node {} (A9);   
  \draw[latex-] (A4) to node {} (A9);     
\end{tikzpicture}
}\Arrow{8}
\Image{
\begin{tikzpicture}[node distance=1.6cm, auto, font=\footnotesize]                                          
  \node (A9) []{9};                                                                                                                     
  \node (A11) [solid node, rectangle, fill=white,below left of = A9]{11};                                                                               
  \node (A4) [below right of = A9]{4};   
  \draw[-latex] (A9) to node [] {2} (A4);      
  \draw[-latex] (A9) to node [swap] {} (A11);    
\end{tikzpicture}
}\Arrow{4}
\Image{
\begin{tikzpicture}[node distance=1.6cm, auto, font=\footnotesize]                                          
  \node (A9) []{9};                                                                                                                     
  \node (A11) [solid node, rectangle, fill=white,below left of = A9]{11};                                                                               
  \node (A4) [below right of = A9]{4};   
  \draw[latex-] (A9) to node [] {2} (A4);      
  \draw[-latex] (A9) to node [swap] {} (A11);    
\end{tikzpicture}.}
\end{center}	
Direct computation shows that this last subquiver corresponds to the degree subquiver
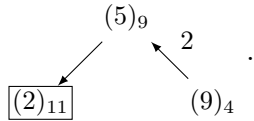
\begin{center}
\Image{
\begin{tikzpicture}[node distance=1.6cm, auto, font=\footnotesize]                                          
  \node (A7) []{$(5)_9$};                                                                                                                     
  \node (A1) [solid node, rectangle, fill=white,below left of = A7]{$(2)_{11}$};           
  \node (A14) [below right of = A7]{$(9)_4$};   
  \draw[latex-] (A1) to node {} (A7);      
  \draw[latex-] (A7) to node [] {2} (A14); 
\end{tikzpicture}}. 
\end{center}
Since $9 \geq 5$, we can now apply Proposition \ref{prop:degree_subquiver_tends_infty_relaxed}, noting that the frozen vertex $11$ is not part of the mutation path stipulated by the proposition.
\end{proof}

\begin{cor}
For $k\geq 4$ and $l \geq 4$, the graded cluster algebra structure for $M(k,l)$ has variables in infinitely many degrees.
\end{cor}

The other infinite type cases that are not covered by the above are $M(3,6)$ and $M(6,3)$. Note that $Q_{M(3,6)} \cong Q_{M(6,3)}$, so we only need to consider $Q_{M(3,6)}$. For this, we can start with the subquiver
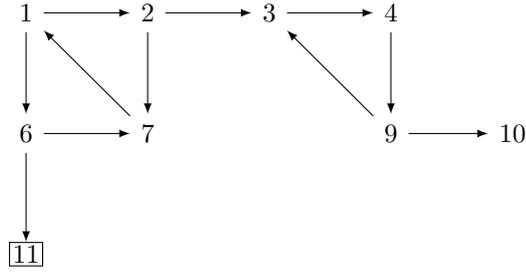
\begin{center}
\Image{
\begin{tikzpicture}[node distance=1.6cm, auto, font=\footnotesize]                                          
  \node (A1) []{1};                                                                                                                     
  \node (A2) [right of = A1]{2};                                                                               
  \node (A3) [right of = A2]{3};   
  \node (A4) [right of = A3]{4};     
  \node (A6) [below of = A1]{6};     
  \node (A7) [right of = A6]{7};     
  \node (A9) [below of = A4]{9};     
  \node (A10) [right of = A9]{10};    
  \node (A11) [solid node, rectangle, fill=white,below of = A6]{11};       
  \draw[-latex] (A1) to node {} (A2);    
  \draw[-latex] (A2) to node {} (A3);  
  \draw[-latex] (A3) to node {} (A4);  
  \draw[-latex] (A6) to node {} (A7);  
  \draw[-latex] (A9) to node {} (A10);  
  \draw[-latex] (A6) to node {} (A11);  
  \draw[-latex] (A7) to node {} (A1);  
  \draw[-latex] (A9) to node {} (A3);
  \draw[-latex] (A1) to node {} (A6);  
  \draw[-latex] (A4) to node {} (A9);  
  \draw[-latex] (A2) to node {} (A7);    
\end{tikzpicture}}
\end{center}
and perform the mutation path $[2,4,3,7,10,6,1,2,7,6,9,3]$ to obtain a quiver containing the subquiver 
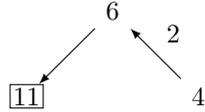
\Image{
\begin{tikzpicture}[node distance=1.6cm, auto, font=\footnotesize]                                          
  \node (A7) []{6};                                                                                                                     
  \node (A1) [solid node, rectangle, fill=white,below left of = A7]{11};                                                                               
  \node (A14) [below right of = A7]{4};   
  \draw[latex-] (A1) to node {} (A7);      
  \draw[latex-] (A7) to node [] {2} (A14); 
\end{tikzpicture}}
(this time we do not write down all the intermediate computations). 
The corresponding degree quiver is  
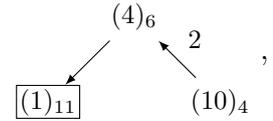
\Image{
\begin{tikzpicture}[node distance=1.6cm, auto, font=\footnotesize]                                          
  \node (A7) []{$(4)_6$};                                                                                                                     
  \node (A1) [solid node, rectangle, fill=white,below left of = A7]{$(1)_{11}$};                                                                               
  \node (A14) [below right of = A7]{$(10)_4$};   
  \draw[latex-] (A1) to node {} (A7);      
  \draw[latex-] (A7) to node [] {2} (A14); 
\end{tikzpicture}},
which satisfies Proposition \ref{prop:degree_subquiver_tends_infty_relaxed}. We again note that vertex $11$ is not part of the stipulated mutation path.

We have now established the following.
\begin{prop}
If the graded cluster algebra associated to $\mathcal{O}(M(k,l))$ is of infinite type, it has cluster variables of infinitely many different degrees.
\end{prop}

\section{Occurring degrees}
A natural question to ask is whether variables of all positive degrees exist in the graded cluster algebra associated to $\mathcal{O}(M(k,l))$. We will prove that this is indeed the case for $\mathcal{O}(M(4,4))$ and $\mathcal{O}(M(3,6))$, which, as before, also proves the result for all larger cases and therefore all infinite type cases.
  
\begin{lem}
\label{lem:M_4_4_all_degrees_occur}
$\A \big( \mathcal{O}(M(4,4)) \big)$ has cluster variables of each degree in $\N$.
\end{lem}

\begin{proof}
We show this by writing down certain mutation sequences that result in degree subquivers of the form of \pref{form:linear_increasing_degree_quiver}. It turns out in this case that these will yield sequences of degrees increasing by $4$ at each mutation.

The initial quiver mutated by the path $[6,4,8,9,1,4,2,6]$ has the subquiver and degree subquiver
\striL{$4$}{$6$} and \striL{$(8)$}{$(12)$}, so we have variables in all degrees of the form $4k$ for $k \geq 2$ by Lemma \ref{lem:linear_increasing_degree_quiver}.

The path $[4,6,8,9,1,4,2,6]$ gives the subquiver and degree subquiver \striR{$4$}{$6$} and \striR{$(10)$}{$(6)$}, so we have variables in all degrees of the form $4k+2$ for $k \geq 1$.

Next, $[1,2,4,8,6,9,1]$ gives the subquiver and degree subquiver \striR{$1$}{$9$} and \striR{$(9)$}{$(5)$}, which leads to variables in all degrees of the form $4k+1$ for $k \geq 1$.

Finally, $[1,9,1,8,6,1,2,6,4,8]$ gives
\striR{$1$}{$9$} and \striR{$(11)$}{$(7)$} which leads to variables in all degrees of the form $4k+3$ for $k \geq 1$.

The initial cluster contains variables of degrees $1,2$ and $3$, and it is easy to find a mutation path that yields a variable of degree $4$. (For example, $\deg [ 6 ]=4$.) This covers all positive degrees.
\end{proof}

\begin{rmk}
In lifting the result of Lemma \ref{lem:M_4_4_all_degrees_occur} to larger cases, we should check that the subquiver (corresponding to $Q_{M(4,4)}$) that we deal with will not interact with the rest of the quiver it is embedded in. In fact it is easy to see that if mutation of $Q_{M(k,l)}$ is restricted to the vertices at $(i,j)$ within the rectangle $1\leq i \leq k$, $1 \leq j \leq l$, then there can be no arrow created between one of these vertices and a vertex in position $(i',j')$ where $i' \geq k+2$ or $j' \geq l+2$. So since vertices we mutate at in $Q_{M(4,4)}$ satisfy $1 \leq i,j \leq 3$ (and since this time, rather than embedding the pattern in the bottom left, we will keep vertices corresponding to the above mutation paths where they are in any larger initial quivers), we do not need to worry about vertices outside the leftmost $4 \times 4$ square in larger cases.
\end{rmk}

\begin{lem}
\label{lem:M_3_6_all_degrees_occur}
$\A \big( \mathcal{O}(M(3,6)) \big)$ has cluster variables of each degree in $\N$.
\end{lem}

\begin{proof}
In this case, we find sequences of degrees increasing by $6$ at each mutation. We might expect that the sequences in this case will increase by a larger amount than in $\mathcal{O}(M(4,4))$ because having so few vertices in which to mutate means achieving a double arrow is not possible in as few mutations as previously, which in turn means degrees have longer to grow along the way. The proof is essentially the same as for $\mathcal{O}(M(4,4))$ so we will summarise it in Table \ref{tab:M_3_6_paths_giving_all_degrees} below.

This leaves degrees $1,2,3,4,6$ and $7$ to consider, but all these degrees are contained in $\cl [ 10, 12, 4, 7, 10, 5, 8, 1, 4, 9, 4 ]$, for example.

\begin{center}
    \setlength{\tabcolsep}{3pt} 
    \renewcommand{\arraystretch}{1.3} 
    \begin{tabular}{ | c  c  c  c |}
    \hline
    Path & Quiver & Degree Quiver & Degrees \\ \hline
   $[9,7,5,2,10,9,5,8,1,4,9]$ & \striL{$7$}{$9$} & \striL{$(12)$}{$(18)$} & $6k$, $k \geq 2$  \\         
   $[10,2,4,7,10,5,8,1,4,9,4]$ & \striR{$10$}{$2$} & \striR{$(19)$}{$(13)$} & $6k+1$, $k \geq 2$  \\ 
   $[7,9,5,2,10,9,5,8,1,4,9]$ & \striR{$7$}{$9$} & \striR{$(14)$}{$(8)$} & $6k+2$, $k \geq 1$  \\     
   $[8,3,5,2,6,9,8,9,3,4,7]$ &\striR{$8$}{$3$} & \striR{$(15)$}{$(9)$} & $6k+3$, $k \geq 1$  \\
   $[10,9,2,7,1,5,9,10,4,8,9]$ &\striR{$10$}{$9$} & \striR{$(16)$}{$(10)$} & $6k+4$, $k \geq 1$  \\    
   $[10,7,9,7,2,9,8,4,9,1,5]$ &\striR{$10$}{$7$} & \striR{$(11)$}{$(5)$} & $6k+5$, $k \geq 0$  \\  \hline  
    \end{tabular}    
    \captionof{table}{\label{tab:M_3_6_paths_giving_all_degrees} Mutation paths leading to all degrees for $M(3,6)$}  
\end{center}

\end{proof}

\begin{cor}
If the graded cluster algebra associated to $\mathcal{O}(M(k,l))$ is of infinite type, it has cluster variables of each degree in $\N$.
\end{cor}

\section{Corollaries for Grassmannians} \label{sec:corollaries_for_grassmannians}
\begin{dfn}
The \emph{Grasmmannian} $Gr(k,k+l)$ is the set of $k$-dimensional subspaces of a $(k+l)$-dimensional vector space $V$ over $\C$. 

After fixing a basis for $V$, such a subspace can be described by a $k \times (k+l)$ matrix whose rows are linearly independent vectors forming a basis of the subspace. Let $I$ be a subset of $\{1,\dots,(k+l)\}$ such that $| I |=k$. The \emph{Pl{\"u}cker coordinate} $x_I$ is the function that maps a $k \times (k+l)$ matrix to its minor indexed by $I$. 

\sloppy The \emph{coordinate ring} $\mathcal{O}(Gr(k,k+l))$ is isomorphic to the subalgebra of  ${\mathcal{O}(M(k,k+l))}$ generated by the Pl{\"u}cker coordinates.
\end{dfn}

The Grasmannian coordinate ring $\mathcal{O}(Gr(k,k+l))$ has the structure of a graded cluster algebra as follows. For the initial exchange quiver, we take the same quiver as for $\mathcal{O}(M(k,l))$, but add a frozen vertex with a single arrow to the vertex at position $(1,1)$. We then assign all vertices degree $1$. By results of \cite{GL}, this quiver gives rise to a cluster algebra for $\mathcal{O}(Gr(k,k+l))$. While an explicit expression for the initial cluster variables is not given in \cite{GL}, one can be obtained by tracing through the construction therein. References \cite{GSV} and \cite{Sc} give more explicit initial clusters but with different quivers. Since we do not need to know the cluster variables explicitly, we do not concern ourselves unduly with this.

Again, it is easy to check that the degree quiver is valid. For further background on how the associated cluster algebra has the desired structure, see \cite{Sc} (for the classical case), \cite{GL} (for the quantum case) or \cite{GSV}.

In terms of proving the existence of infinitely many degrees, the corresponding results for $\mathcal{O}(Gr(k,k+l))$ will immediately follow from the matrix algebra case. Since the grading on $\mathcal{O}(Gr(k,k+l))$ is an $\N$-grading, and since the initial exchange quiver is the same as that for $\mathcal{O}(M(k,l))$ but for an added frozen vertex, the same mutation sequences as in the matrix case will yield quivers which have the same subquivers as above, to which we can again apply Proposition \ref{prop:degree_subquiver_tends_infty_relaxed}. 
Thus we have:

\begin{prop}
If the graded cluster algebra associated to $\mathcal{O}(Gr(k,k+l))$ is of infinite type, it has cluster variables of infinitely many different degrees.
\end{prop}

We also have variables in each degree and prove this in a very similar way as for the matrix case.

\begin{lem}
$\A \big( \mathcal{O}(Gr(4,8)) \big)$ has cluster variables of each degree in $\N$.
\end{lem}

\begin{proof}
In this case we find paths that give degree sequences increasing by $2$ with each mutation.
Applying $[9,1,2,4,8,6,9,1]$ gives the subquiver and degree subquiver
\striR{$9$}{$1$} and \striR{$(6)$}{$(4)$}, so by Lemma \ref{lem:linear_increasing_degree_quiver} we have variables in all degrees of the form $2k$ for $k \geq 2$.
The path $[4,6,8,9,1,4,2,6]$ gives the subquiver and degree subquiver
\striR{$4$}{$6$} and \striR{$(5)$}{$(3)$}, so we have variables in all degrees of the form $2k+1$ for $k \geq 1$.

This only leaves degree $2$, which is easy to find. For example, $\var[6]=2$.
\end{proof}

\begin{lem}
$\A \big( \mathcal{O}(Gr(3,9)) \big)$ has cluster variables of each degree in $\N$.
\end{lem}

\begin{proof}
As in the previous section, the degree sequences for this case increase by a larger amount than for the matrix case. We summarise three  paths which work in Table \ref{tab:Gr_3_9_paths_giving_all_degrees} below.
\begin{center}
    \begin{tabular}{ | c  c  c  c |}
    \hline
    Path & Quiver & Degree Quiver & Degrees \\ \hline
   $[10,7,9,7,2,9,8,4,9,1,5]$ &\striR{$10$}{$7$} & \striR{$(6)$}{$(3)$} & $3k$, $k \geq 1$  \\     
   $[7,4,3,1,4,9,10,6,2,9,3]$ & \striL{$4$}{$7$} & \striL{$(4)$}{$(7)$} & $3k+1$, $k \geq 1$  \\         
   $[10,9,2,7,1,5,9,10,4,8,9]$ &\striR{$10$}{$9$} & \striR{$(8)$}{$(5)$} & $3k+2$, $k \geq 1$  \\    \hline    
    \end{tabular}   
    \captionof{table}{ \label{tab:Gr_3_9_paths_giving_all_degrees} Mutation paths leading to all degrees for $Gr(3,9)$} 
\end{center}    
This leaves degrees $1$ and $2$ which are found in $\cl [7]$, for example.
\end{proof}

\begin{cor}
If the graded cluster algebra associated to $\mathcal{O}(Gr(k,k+l))$ is of infinite type, it has cluster variables of each degree in $\N$.
\end{cor}

\biblio

\chapter{Gradings for finite mutation type quivers not arising from surface triangulations} \label{chp:finite_mutation_type}

Quivers of finite mutation type (i.e.~ones whose associated matrices are mutation-finite) are known to fall into one of two cases: adjacency matrices of triangulations of two-dimensional marked surfaces, and a finite collection of quivers that do not correspond to triangulations of surfaces. More precisely, we have the following:

\begin{thm}[{\cite[Theorem 6.1]{FSTFinite}}] 
A rank $n$ quiver $(n \geq 3)$ of finite mutation type either corresponds to an adjacency matrix of a triangulation of a bordered two-dimensional surface or is mutation equivalent to one of the following quivers: $E_6$, $E_7$, $E_8$, $\wt{E}_6$, $\wt{E}_7$, $\wt{E}_8$, $E_6^{(1,1)}$, $E_7^{(1,1)}$, $E_8^{(1,1)}$, $X_6$ and $X_7$.
\end{thm}

In this chapter we will investigate the latter class.

\section{Initial information about the finite list of quivers}
We will try to determine what can be said about the gradings on the cluster algebras associated to some of the quivers in this finite list; ultimately we will consider $X_7$, $\wt{E}_6$, $\wt{E}_7$ and $\wt{E}_8$ in detail (though we give some additional information about $E_6^{(1,1)}$, $E_7^{(1,1)}$ and $E_8^{(1,1)}$ in Table \ref{tab:fin_mut_quiv_info}). Note that the quivers $E_6$, $E_7$ and $E_8$ give rise to finite type cluster algebras and have been considered in \cite{G15}.  Our approach will be to use computer-aided calculation to find the degree quiver class for each case, which will yield an exhaustive list of all occurring degrees, and then to prove which degrees have infinitely many variables. To do this, we will find mutation paths that give rise to repeating degree quivers (up to sign), but which produce growing denominator vectors. 

We start by writing down initial grading bases for the quivers in the list that we have not considered in a previous chapter.

\begin{lem}
  Initial grading bases for $\wt{E}_6$, $\wt{E}_7$, $\wt{E}_8$, $E_6^{(1,1)}$, $E_7^{(1,1)}$, $E_8^{(1,1)}$, $X_6$ and $X_7$ are given in Table \ref{tab:fin_mut_quiv_init_gradings}.

\begin{center}
\scalebox{0.864}{
    \begin{tabular}{ | c  c |}
		  \hline
		  \emph{Quiver} & \emph{Grading} \\ \hline
		  $X_6=$
		  \begin{tikzpicture}[baseline=(current  bounding  box.center),node distance=1.4cm, auto, font=\footnotesize]                                          
			  \node (A2) []{$2$};                                                                                                                     
			  \node (A1) [below left of = A2]{$1$};                                                                               
			  \node (A3) [below right of = A2]{$3$};   
			  \node (A4) [above right of = A3]{$ 4$};
			  \node (A5) [below right of = A4]{$ 5$}; 
			  \node (A6) [below of = A3]{$ 6$};     
			  \draw[latex-] (A2) to node [swap] {$2$} (A1);      
			  \draw[latex-] (A3) to node [] {$ $} (A2);    
			  \draw[-latex] (A3) to node [] {$ $} (A1);   
			  \draw[-latex] (A3) to node [] {$ $} (A4);   
			  \draw[-latex] (A4) to node [] {$2 $} (A5);
			  \draw[-latex] (A6) to node [] {$ $} (A3);
			  \draw[-latex] (A5) to node [] {$ $} (A3);
		  \end{tikzpicture} & $\{(0,0,0,0,0,0) \}$ \\     
		  $X_7=$
		  \begin{tikzpicture}[baseline=(current  bounding  box.center),node distance=1.4cm, auto, font=\footnotesize]                                          
			  \node (A2) []{$2$};                                                                                                                     
			  \node (A1) [below left of = A2]{$1$};                                                                               
			  \node (A3) [below right of = A2]{$3$};   
			  \node (A4) [above right of = A3]{$ 4$};
			  \node (A5) [below right of = A4]{$ 5$}; 
			  \node (A6) [below right of = A3]{$ 6$};     
			  \node (A7) [below left of = A3]{$7$};     			  
			  \draw[latex-] (A2) to node [swap] {$2$} (A1);      
			  \draw[latex-] (A3) to node [] {$ $} (A2);    
			  \draw[-latex] (A3) to node [] {$ $} (A1);   
			  \draw[-latex] (A3) to node [] {$ $} (A4);   
			  \draw[-latex] (A4) to node [] {$2 $} (A5);
			  \draw[latex-] (A6) to node [] {$ $} (A3);
			  \draw[-latex] (A5) to node [] {$ $} (A3);
			  \draw[-latex] (A6) to node [] {$2 $} (A7);			  
			  \draw[-latex] (A7) to node [] {$ $} (A3);			  
		  \end{tikzpicture}& $\{(1,1,2,1,1,1,1)\} $  \\        \hline
		  $\wt{E}_6=$  
		  \begin{tikzpicture}[baseline=(current  bounding  box.center),node distance=1.1cm, auto, font=\footnotesize]                                          
			  \node (A3) []{$3$};                                                                                                                     			  
			  \node (A2) [left of = A3]{$2$};
			  \node (A1) [left of = A2]{$1$};
			  \node (A6) [above of = A3]{$6$};
			  \node (A7) [above of = A6]{$7$};			  
			  \node (A4) [right of = A3]{$4$};
			  \node (A5) [right of = A4]{$5$};
			  \draw[-latex] (A1) to node [] {} (A2);
			  \draw[-latex] (A2) to node [] {} (A3);
			  \draw[-latex] (A6) to node [] {} (A3);
			  \draw[-latex] (A4) to node [] {} (A3);
			  \draw[-latex] (A7) to node [] {} (A6);
			  \draw[-latex] (A5) to node [] {} (A4);		  			  			  			  			  			  
		  \end{tikzpicture}		  & $\{ (1,0,1,0,1,0,1)\}$  \\     
		  $\wt{E}_7=$  
		  \begin{tikzpicture}[baseline=(current  bounding  box.center),node distance=1.1cm, auto, font=\footnotesize]                                          
			  \node (A3) []{$3$};                                                                                                                     			  
			  \node (A2) [left of = A3]{$2$};
			  \node (A1) [left of = A2]{$1$}; 
			  \node (A4) [right of = A3]{$4$};
			  \node (A5) [right of = A4]{$5$};
			  \node (A6) [right of = A5]{$6$};			  
			  \node (A7) [right of = A6]{$7$};			  
			  \node (A8) [above of = A4]{$8$};			  			  
			  \draw[-latex] (A1) to node [] {} (A2);
			  \draw[-latex] (A2) to node [] {} (A3);
			  \draw[-latex] (A3) to node [] {} (A4);
			  \draw[-latex] (A5) to node [] {} (A4);
			  \draw[-latex] (A6) to node [] {} (A5);		  			  			  			  			  			  
			  \draw[-latex] (A7) to node [] {} (A6);			  
			  \draw[-latex] (A8) to node [] {} (A4);			  
		  \end{tikzpicture}		  & \begin{tabular}{c}  $\{(-1,0,-1,0,1,0,1,0)$, \\ $(-1,0,-1,0,0,0,0,1) \}$ \end{tabular}  \\   		  
		  $\wt{E}_8=$  
		  \begin{tikzpicture}[baseline=(current  bounding  box.center),node distance=1.1cm, auto, font=\footnotesize]                                          
			  \node (A3) []{$3$};                                                                                                                     			  
			  \node (A2) [left of = A3]{$2$};
			  \node (A1) [left of = A2]{$1$}; 
			  \node (A4) [right of = A3]{$4$};
			  \node (A5) [right of = A4]{$5$};
			  \node (A6) [right of = A5]{$6$};			  
			  \node (A7) [right of = A6]{$7$};			  
			  \node (A8) [right of= A7]{$8$};			  			  
			  \node (A9) [above of = A3]{$9$};					  
			  \draw[-latex] (A1) to node [] {} (A2);
			  \draw[-latex] (A2) to node [] {} (A3);
			  \draw[-latex] (A4) to node [] {} (A3);
			  \draw[-latex] (A5) to node [] {} (A4);
			  \draw[-latex] (A6) to node [] {} (A5);		  			  			  			  			  			  
			  \draw[-latex] (A7) to node [] {} (A6);			  
			  \draw[-latex] (A8) to node [] {} (A7);			  
			  \draw[-latex] (A9) to node [] {} (A3);			  			  
		  \end{tikzpicture}		  & $\{(0,0,0,-1,0,-1,0,-1,1)\}$\\   		 \hline
		  $E_6^{(1,1)}=$  
		  \begin{tikzpicture}[baseline=(current  bounding  box.center),node distance=1.1cm, auto, font=\footnotesize]                                          
			  \node (B) [] {};                                                                                                                     			  
			  \node (A2) [left of = B]{$2$};
			  \node (A1) [left of = A2]{$1$}; 
			  \node (A4) [right of = B]{$4$ };
			  \node (A5) [right of = A4]{$8$};
			  \node (A6) [right of = A5]{$6$};			  
			  \node (A7) [right of = A6]{$7$};			  
			  \node (A8) [above of = B]{$3$};			  			  
			  \node (A9) [below of = B]{$5$};				  
			  \draw[-latex] (A1) to node [] {} (A2);
			  \draw[-latex] (A5) to node [] {} (A4);	  			  			  			  			  			  
			  \draw[-latex] (A7) to node [] {} (A6);			    
			  \draw[-latex] (A9) to node [] {$2$} (A8);				  
			  \draw[-latex] (A8) to node [] {} (A2);			  
			  \draw[-latex] (A2) to node [] {} (A9);				  
			  \draw[-latex] (A8) to node [] {} (A4);			  
			  \draw[-latex] (A4) to node [] {} (A9);				  
			  \draw[-latex] (A8) to node [] {} (A6);			  
			  \draw[-latex] (A6) to node [] {} (A9);				  			  
		  \end{tikzpicture}		  & $\{(0,0,0,0,0,0,0,0)\}$\\ 		    	 	  
		  $E_7^{(1,1)}=$  
		  \begin{tikzpicture}[baseline=(current  bounding  box.center),node distance=1.1cm, auto, font=\footnotesize]                                          
			  \node (B) [] {};                                                                                                                     			  
			  \node (A2) [left of = B]{$2$};
			  \node (A1) [left of = A2]{$1$}; 
			  \node (A4) [right of = B]{$4$ };  				  	  
			  \node (A6) [right of = A4]{$6$};
			  \node (A7) [right of = A6]{$7$};			  						  
			  \node (A5) [right of = A7]{$8$};				  
			  \node (A8) [above of = B]{$3$};			  			  
			  \node (A9) [below of = B]{$5$};				  
			  \node (A10) [left of = A1]{$9$};				  
			  \draw[-latex] (A1) to node [] {} (A2);
			  \draw[-latex] (A5) to node [] {} (A7);		  			  			  			  			  			  
			  \draw[-latex] (A7) to node [] {} (A6);			    
			  \draw[-latex] (A9) to node [] {$2$} (A8);				  
			  \draw[-latex] (A8) to node [] {} (A2);			  
			  \draw[-latex] (A2) to node [] {} (A9);				  
			  \draw[-latex] (A8) to node [] {} (A4);			  
			  \draw[-latex] (A4) to node [] {} (A9);				  
			  \draw[-latex] (A8) to node [] {} (A6);			  
			  \draw[-latex] (A6) to node [] {} (A9);				  			  
			  \draw[-latex] (A10) to node [] {} (A1);				  
		  \end{tikzpicture}		  &\begin{tabular}{c}  $\{(0,1,0,-1,0,0,0,0,1)$, \\  $(0,0,1,2,1,0,0,0,0)$, \\ $(0,0,0,-1,0,1,0,1,0) \}$\end{tabular}  \\  
		  $E_8^{(1,1)}=$  
		  \begin{tikzpicture}[baseline=(current  bounding  box.center),node distance=1.1cm, auto, font=\footnotesize]                                          
			  \node (B) [] {};                                                                                                                     			  
			  \node (A2) [left of = B]{$2$};
			  \node (A1) [left of = A2]{$1$}; 
			  \node (A4) [right of = B]{$4$ };  				  	  
			  \node (A6) [right of = A4]{$6$};
			  \node (A7) [right of = A6]{$7$};			  						  
			  \node (A5) [right of = A7]{$8$};				  
			  \node (A8) [above of = B]{$3$};			  			  
			  \node (A9) [below of = B]{$5$};				  
			  \node (A10) [right of = A5]{$9$};				  
			  \node (A11) [right of = A10]{$10$};				  			  
			  \draw[-latex] (A1) to node [] {} (A2);
			  \draw[-latex] (A5) to node [] {} (A7);		  			  			  			  			  			  
			  \draw[-latex] (A7) to node [] {} (A6);			    
			  \draw[-latex] (A9) to node [] {$2$} (A8);				  
			  \draw[-latex] (A8) to node [] {} (A2);			  
			  \draw[-latex] (A2) to node [] {} (A9);				  
			  \draw[-latex] (A8) to node [] {} (A4);			  
			  \draw[-latex] (A4) to node [] {} (A9);				  
			  \draw[-latex] (A8) to node [] {} (A6);			  
			  \draw[-latex] (A6) to node [] {} (A9);				  			  
			  \draw[-latex] (A10) to node [] {} (A5);				  
			  \draw[-latex] (A11) to node [] {} (A10);				  			  
		  \end{tikzpicture}		  &\begin{tabular}{c}  $\{(0,0,1,2,1,0,0,0,0,0)$, \\ $(0,0,0,-1,0,1,0,1,0,1) \}$\end{tabular}  \\ 		  
	      \hline  
    \end{tabular}
    }
    \captionof{table}{\label{tab:fin_mut_quiv_init_gradings}}
\end{center}  
\end{lem}

\begin{proof}
This reduces to finding bases for the kernels of the skew-symmetric matrices corresponding to these quivers.
\end{proof}

So we see that $X_6$ and $E_6^{(1,1)}$ only give rise to the zero grading. Next we wish to determine which degrees occur in the corresponding cluster algebras. We do so with a MAGMA algorithm (available in an accompanying file) which computes the mutation class of degree quivers as well as some other information, which we summarise in Proposition \ref{prop:fin_mut_quiv_info}. Note that mutation and degree quiver class sizes are computed up to essential equivalence. In particular, this number will be smaller than that of the classes computed up to quiver isomorphism (as is done in the Java app of \cite{KJA}) since not all quivers are isomorphic to their opposites whereas a quiver is always essentially equivalent to its opposite.

\begin{prop}
\label{prop:fin_mut_quiv_info}
 We have the following for the quivers $E_6$, $E_7$, $E_8$, $\wt{E}_6$, $\wt{E}_7$, $\wt{E}_8$, $E_6^{(1,1)}$, $E_7^{(1,1)}$, $E_8^{(1,1)}$, $X_6$ and $X_7$.

\begin{center}
\setlength{\tabcolsep}{3pt} 
\renewcommand{\arraystretch}{1.3} 
    \begin{tabular}{ | c  c  c  c |}
		  \hline
		  \emph{Quiver} & \emph{Mutation class} & \emph{Degree quivers class} & \emph{Occurring Degrees}\\
		  \hline
		  $X_6$ & $3$ & $3$ & $ \{ 0 \}$ \\
		  $X_7$ & $2$ & $2$ & $ \{ 1,2\}$ \\		  
		  \hline
		  $\wt{E}_6$ & $74$ & $148$ & $ \{ 0,\pm 1, \pm 2 \}$ \\
		  $\wt{E}_7$ & $571$ & $2297$ & $ \{ \CII{0}{0}, \pm \CII{2}{1}, \pm \CII{1}{1}, \pm \CII{1}{0}, \pm \CII{0}{1} \}$ \\		  		  
		  $\wt{E}_8$ & $7560$ & $7634$ & $ \{ 0,\pm 1, \pm 2 \}$ \\
		  \hline
		  $E_6^{(1,1)}$ & $26$ & $26$ & $ \{ 0 \}$ \\
		  $E_7^{(1,1)}$ & $279$ & $13616$ &
		  \scalemath{0.85}{ \begin{tabular}{c}  $\Big\{ \thvec{0}{0}{0}, \pm \thvecal{1}{-2}{1}, \pm \thvecal{1}{-1}{1},$ \\ $\pm \thvecal{0}{1}{-1}, \pm \thvecal{1}{-1}{0}, \pm \thvecal{1}{0}{1}, \pm \thvecal{1}{0}{-1},$\\ $\pm \thvecal{0}{1}{0}, \pm \thvecal{1}{0}{0}, \pm \thvecal{0}{0}{1} \Big\}$ \end{tabular} } \\
		  $E_8^{(1,1)}$ & $5739$ & $>30000$ & \begin{tabular}{c} $\supseteq \big\{ \tvec{0}{0}, \pm \tvecal{1}{0}, \pm \tvecal{0}{1},$ \\ $ \pm \tvecal{1}{1}, \pm \tvecal{1}{-1},\pm \tvecal{2}{-1}, \pm \tvecal{1}{-2} \big\}$ \end{tabular} \\		  	
		  \hline	  		  
    \end{tabular}
    \captionof{table}{\label{tab:fin_mut_quiv_info}}
\end{center}  
\end{prop}

We are now ready to investigate the question of variables per degree for these quivers. Although we are not able to answer this completely in every case, we may still determine part of the answer.  

\section{Grading for $X_7$}
In the case of $X_7$, whose grading behaviour is qualitatively different from that of $\wt{E}_i$, we are able to determine all the information about the associated graded cluster algebra. To start, note that $(X_7)_{[j]} \esseq X_7$ for all $j \neq 3$.

\begin{lem}
\label{lem:X_7_denom_21}
$\dv_{X_7}[(2,1)_n] = (n,n-1,0,0,0,0,0)$.
\end{lem}

\begin{proof}
It is easy to show
 $(X_7)_{[(2,1)_n]}= 
		  \begin{tikzpicture}[baseline=(current  bounding  box.center),node distance=1.4cm, auto, font=\footnotesize]                                          
			  \node (A2) []{$2$};                                                                                                                     
			  \node (A1) [below left of = A2]{$1$};                                                                               
			  \node (A3) [below right of = A2]{$3$};   
			  \node (A4) [above right of = A3]{$ 4$};
			  \node (A5) [below right of = A4]{$ 5$}; 
			  \node (A6) [below right of = A3]{$ 6$};     
			  \node (A7) [below left of = A3]{$7$};     			  
			  \draw[-latex] (A2) to node [swap] {$2$} (A1);      
			  \draw[-latex] (A3) to node [] {$ $} (A2);    
			  \draw[latex-] (A3) to node [] {$ $} (A1);   
			  \draw[-latex] (A3) to node [] {$ $} (A4);   
			  \draw[-latex] (A4) to node [] {$2 $} (A5);
			  \draw[latex-] (A6) to node [] {$ $} (A3);
			  \draw[-latex] (A5) to node [] {$ $} (A3);
			  \draw[-latex] (A6) to node [] {$2 $} (A7);			  
			  \draw[-latex] (A7) to node [] {$ $} (A3);			  
		  \end{tikzpicture}
$ 
for odd $n$ and $(X_7)_{[(2,1)_n]} = X_7$ for even $n$.  
The claim is true for $n=1,2$, by direct computation. Assume true for $n$ and suppose $n$ is even. Then, since $(X_7)_{[(2,1)_n]} = X_7$, we have  
\begin{align*}
 \dv_{X_7}[(2,1)_{n+1}] ={}& \dv_{X_7}[1,(2,1)^{n/2}] \\
                                 ={}& -(n-1,n-2,0,0,0,0,0)\\
                                 	& +\max(2(n,n-1,0,0,0,0,0),(0,0,-1,0,0,0,0)) \\
                                 ={}& (n+1,n,0,0,0,0,0).
\end{align*} 
(Recall $\max$ is taken componentwise as in Equation \ref{eqn:denominator_mutation}.) Now suppose $n$ is odd. Then, using the quiver we wrote down above for odd $n$, we have 
\begin{align*}
\dv_{X_7}[(2,1)_{n+1}] ={}& \dv_{X_7}[(2,1)^{n/2+1}] \\
                                ={}&-(n-1,n-2,0,0,0,0,0) \\
                                  &+\max(2(n,n-1,0,0,0,0,0),(0,0,-1,0,0,0,0)) \\
                                ={}&(n+1,n,0,0,0,0,0) \\
\end{align*}
again. This proves the result for all $n$.
\end{proof}

\begin{cor}
$\A((1,1,2,1,1,1,1),X_7)$ has infinitely many variables in degree $1$.
\end{cor}

\begin{proof}
This follows since $\deg_{X_7}[(2,1)_n]=1$, which is easy to show.
\end{proof}

\begin{cor}
$\dv_{X_7}[5,4,3,(2,1)_n] = (2n,2(n-1),1,1,1,0,0)$.
\end{cor}

\begin{proof}
Assume $n$ is even (the proof for $n$ odd is very similar). By Lemma \ref{lem:X_7_denom_21}, the denominator cluster $\dcl_{X_7}[(2,1)_n]$ is then
\[   \scalemath{0.915}{ \big( (n-1,n-2,0,0,0,0,0), (n, n-1, 0,0,0,0,0), (0,0,-1,0,0,0,0), \dots, (0,0,0,0,0,0,-1) \big). }  \]
Mutating in direction $3$ we get 
\begin{align*}
\dv _{ }[3,(2,1)_n]=& -\dcl _{ }[(2,1)_n]|^3\\
&+\max \Big(\dcl_{ }[(2,1)_n]|^1+\dcl_{ }[(2,1)_n]|^4+\dcl_{ }[(2,1)_n]|^6,  \\ 
& \phantom{\max .......} \dcl_{ }[(2,1)_n]|^2+\dcl_{ }[(2,1)_n]|^5+\dcl_{ }[(2,1)_n]|^7 \Big)\\
                    =& -\dcl_{ }[(2,1)_n]|^3  +\dcl_{ }[(2,1)_n]|^2 \\
                    =& -(0,0,-1,0,0,0,0)+(n,n-1,0,0,0,0,0)\\
                    =&\phantom{.} (n,n-1,1,0,0,0,0).
\end{align*}
We have
\[(X_7)_{[3,(2,1)_n]}= 
\left( \begin{smallmatrix*}[r]
0 & 1  &1&  0 &-1  &0 &-1\\
-1 & 0 &-1&  1  &0  &1 & 0\\
-1  &1 & 0& -1  &1 &-1 & 1\\
 0 &-1 & 1 & 0  &1 & 0 &-1\\
 1  &0 &-1 &-1  &0 & 1 & 0\\
 0 &-1 & 1 & 0 &-1 & 0 & 1\\
 1 & 0& -1 & 1&  0 &-1&  0
\end{smallmatrix*} \right)\]
(we express this as a matrix as the quiver is not planar). 
So 
\begin{align*}
\dv _{ }[4,3,(2,1)_n] ={}& -\dcl_{ }[3,(2,1)_n]|^4\\
                                &+\max \Big(\dcl_{ }[3,(2,1)_n]|^2+\dcl_{ }[3,(2,1)_n]|^7,  \\ 
                                & \phantom{\max .......} \dcl_{ }[3,(2,1)_n]|^3+\dcl_{ }[3,(2,1)_n]|^5 \Big)\\
                            ={}& -\dcl_{ }[3,(2,1)_n]|^4  +\dcl_{ }[3,(2,1)_n]|^3 \\
                            ={}& -(0,0,0,-1,0,0,0)+(n,n-1,1,0,0,0,0)\\
                            ={}& (n,n-1,1,1,0,0,0).
\end{align*}
Next we have
\[(X_7)_{[4,3,(2,1)_n]}= 
\left( \begin{smallmatrix*}[r]
 0 & 1 & 1  &0 &-1  &0 &-1\\
-1  &0 & 0 &-1 & 1  &1 & 0\\
-1 & 0 & 0 & 1  &1 &-1 & 0\\
 0 & 1 &-1  &0 &-1  &0 & 1\\
 1 &-1 &-1  &1  &0  &1 &-1\\
 0 &-1  &1 & 0 &-1  &0 & 1\\
 1  &0 & 0& -1&  1 &-1&  0\\
\end{smallmatrix*} \right),\]
so
\begin{align*} 
\dv _{ }[5,4,3,(2,1)_n] ={}& -\dcl_{ }[4,3,(2,1)_n]^5\\
&\resizebox{0.74\hsize}{!}{$+\max \Big(\dcl_{ }[4,3,(2,1)_n]^2+\dcl_{ }[4,3,(2,1)_n]^3+\dcl_{ }[4,3,(2,1)_n]^7, $} \\ 
& \resizebox{0.75\hsize}{!}{$\phantom{\max ......} \dcl_{ }[4,3,(2,1)_n]^1+\dcl_{ }[4,3,(2,1)_n]^4+\dcl_{ }[4,3,(2,1)_n]^6 \Big) $}\\
={}&-(0,0,0,0,-1,0,0)\\
&+\max \Big( (n,n-1,0,0,0,0,0)+(n,n-1,1,0,0,0,0), \\ 
& \phantom{\max .......}(n-1,n-2,0,0,0,0,0)+(n,n-1,1,1,0,0,0) \Big) \\
                    ={}& (0,0,0,0,1,0,0)+(2n,2(n-1),1,1,0,0,0) \\
                    ={}& (2n,2(n-1),1,1,1,0,0),
\end{align*}
as required. (We omitted $\dcl_{ }[4,3,(2,1)_n]^6$ and $\dcl_{ }[4,3,(2,1)_n]^7$ in the second $\max$ above as they only have zero or negative entries.)
\end{proof}

\begin{cor}
$\A((1,1,2,1,1,1,1),X_7)$ has infinitely many variables in degree $2$.
\end{cor}

\begin{proof}
This follows since $\deg_{X_7}[5,4,3,(2,1)_n] = \deg_{X_7}[5,4,3]=2$.
\end{proof}

\begin{prop}
 $\A((1,1,2,1,1,1,1),X_7)$ has infinitely many variables of each occurring degree.
\end{prop}

\begin{proof}
\sloppy This is immediate from the previous two Corollaries, since $\A((1,1,2,1,1,1,1),X_7)$ only has variables in degrees $1$ and $2$, as noted in Table \ref{tab:fin_mut_quiv_init_gradings}.
\end{proof}

\section{Gradings for $\wt{E}_6$, $\wt{E}_7$ and $\wt{E}_8$}

Now let us consider $\wt{E}_i$. For these cases, there are certain degrees which we expect to correspond with only one variable, while other degrees have infinitely many variables. We are able to show that the latter is true, and we do so by finding a mutation path which gives us infinitely many examples. Such a path is comparable to the special path $[(3,2,1)^{2n}]$ in Lemma \ref{lem:acyclic_special_mutation_path}. However, while there was only one minimal path of this kind in the cases in Chapter 3, we can find more than one in these examples, and it is not always the case that such a minimal path will always have length equal to the rank. Such paths are also usually much more difficult to find for these cases. 

All of these facts are due to the higher rank: for the rank 3 acyclic quivers, the floor of the exchange graph consists of repeating segments comprised of two-dimensional polytopes along which the special path traverses. In higher ranks, we expect the floor to consist of polytopes whose dimension is close to the rank. Thus there may be many routes between one ``segment" of the floor and another, but due to the size and complexity of the segments, finding paths that move between the initial vertex of each segment is more difficult.

In general, little is known about the structure of the exchange graph of cluster algebras of rank greater than 3. The paths we provide below show that the exchange graph for these example does consist of repeating segments. We conjecture that, similar to the rank 3 mixed case (but in contrast to the rank 3 acyclic cases), these segments are bounded by hyperplanes (rather than being enclosed by a canopy of infinitely many trees whose vertices are forks).

\begin{rmk}
The graded cluster algebras arising from $\wt{E}_6$, $\wt{E}_7$ and $\wt{E}_8$ are all balanced. This is by Lemma \ref{lem:acyclic_implies_balanced}, since the corresponding quivers are acyclic.
\end{rmk}

\subsection{Grading for $\wt{E}_6$}
In $\wt{E}_6$ we will show that degrees $0$ and $\pm 1$ have infinitely many different variables. We conjecture that degrees $\pm 2$ correspond to one variable each.
 
\begin{lem}
\label{lem:R_aff_6_degseed_loop}
Let $\pth{p}=[3,4,6,2,1,7,5]$. Then 
\[ \dsd_{\wt{E}_6}\pth{p^n}= \big((-1)^n (1,0,1,0,1,0,1), \wt{E}_6 \big) .\]
That is, mutating the initial degree seed of  $\wt{E}_6$ along $\pth{p}$ results in the the negative of the initial degree seed again.
\end{lem}

\begin{proof}
Although it is straightforward, we will write out the computation for $n=1$ since we will need to refer to the intermediate quivers involved again later. After this, the result is immediate.
We have
\begin{align*}
\dsd_{\wt{E}_6}[5] &=  
\left( 
( 1, 0, 1, 0, -1, 0, 1 ),
\left( \begin{smallmatrix}
 0&  1&  0&  0&  0&  0&  0\\
-1  &0  &1  &0  &0  &0&  0\\
 0 &-1 & 0& -1&  0& -1&  0\\
 0 & 0 & 1 & 0 & 1 & 0 & 0\\
 0&  0&  0& -1&  0&  0&  0\\
 0 & 0&  1&  0&  0&  0& -1\\
 0  &0 & 0 & 0 & 0 & 1 & 0\\
\end{smallmatrix} \right) 
\right), \\
\dsd_{\wt{E}_6}[7,5] &=  
\left( 
( 1, 0, 1, 0, -1, 0, -1  ),
\left( \begin{smallmatrix}
 0&  1&  0&  0&  0&  0&  0\\
-1  &0  &1  &0  &0  &0&  0\\
 0 &-1 & 0& -1&  0& -1&  0\\
 0 & 0 & 1 & 0 & 1 & 0 & 0\\
 0&  0&  0& -1&  0&  0&  0\\
 0 & 0&  1&  0&  0&  0& 1\\
 0  &0 & 0 & 0 & 0 & -1 & 0\\
\end{smallmatrix} \right) 
\right), \\
\dsd_{\wt{E}_6}[1,7,5] &=  
\left( 
(  -1, 0, 1, 0, -1, 0, -1  ),
\left( \begin{smallmatrix}
 0&  -1&  0&  0&  0&  0&  0\\
1  &0  &1  &0  &0  &0&  0\\
 0 &-1 & 0& -1&  0& -1&  0\\
 0 & 0 & 1 & 0 & 1 & 0 & 0\\
 0&  0&  0& -1&  0&  0&  0\\
 0 & 0&  1&  0&  0&  0& 1\\
 0  &0 & 0 & 0 & 0 & -1 & 0\\
\end{smallmatrix} \right) 
\right), \\
\dsd_{\wt{E}_6}[2,1,7,5] &=  
\left( 
( -1, 0, 1, 0, -1, 0, -1   ),
\left( \begin{smallmatrix}
 0&  1&  0&  0&  0&  0&  0\\
-1  &0  &-1  &0  &0  &0&  0\\
 0 &1 & 0& -1&  0& -1&  0\\
 0 & 0 & 1 & 0 & 1 & 0 & 0\\
 0&  0&  0& -1&  0&  0&  0\\
 0 & 0&  1&  0&  0&  0& 1\\
 0  &0 & 0 & 0 & 0 & -1 & 0\\
\end{smallmatrix} \right) 
\right), \\
\dsd_{\wt{E}_6}[6,2,1,7,5] &=  
\left( 
( -1, 0, 1, 0, -1, 0, -1),
\left( \begin{smallmatrix}
 0&  1&  0&  0&  0&  0&  0\\
-1  &0  &-1  &0  &0  &0&  0\\
 0 &1 & 0& -1&  0& 1&  0\\
 0 & 0 & 1 & 0 & 1 & 0 & 0\\
 0&  0&  0& -1&  0&  0&  0\\
 0 & 0&  -1&  0&  0&  0& -1\\
 0  &0 & 0 & 0 & 0 & 1 & 0\\
\end{smallmatrix} \right) 
\right), \\
\dsd_{\wt{E}_6}[4,6,2,1,7,5] &=  
\left( 
( -1, 0, 1, 0, -1, 0, -1),
\left( \begin{smallmatrix}
 0&  1&  0&  0&  0&  0&  0\\
-1  &0  &-1  &0  &0  &0&  0\\
 0 &1 & 0& 1&  0& 1&  0\\
 0 & 0 & -1 & 0 & -1 & 0 & 0\\
 0&  0&  0& 1&  0&  0&  0\\
 0 & 0&  -1&  0&  0&  0& -1\\
 0  &0 & 0 & 0 & 0 & 1 & 0\\
\end{smallmatrix} \right) 
\right), \\
\dsd_{\wt{E}_6}[3,4,6,2,1,7,5] &=  
\left( 
( -1, 0, -1, 0, -1, 0, -1),
\wt{E}_6
\right),
\end{align*}
as required.
\end{proof}

Note that $l(\pth{p})$ is equal to the rank of $\wt{E}_6$. (We conjecture that this is a minimal length path with the desired property.)

We will soon require the use of floor and ceiling functions. Note the following identities involving these functions, which we will make frequent use of in this chapter. We have
\begin{flalign}
&& \left\lceil {\frac {n}{m}}\right\rceil&=\left\lfloor {\frac {n-1}{m}}\right\rfloor +1,& \label{eqn:floor_to_ceil}\\
{  \text{(or equivalently)} } && \left\lfloor {\frac {n}{m}}\right\rfloor& =\left\lceil {\frac {n+1}{m}}\right\rceil -1,& \label{eqn:ceil_to_floor}
\end{flalign}
and
 \begin{equation}
n=\left\lceil {\frac {n}{m}}\right\rceil +\left\lceil {\frac {n-1}{m}}\right\rceil +\dots +\left\lceil {\frac {n-m+1}{m}}\right\rceil, \label{eqn:n_as_sum_of_ceils}
\end{equation}
for any $n, m \in \N$. (Equation \ref{eqn:n_as_sum_of_ceils} is given in \cite[p.~85]{GKP}.)

\begin{lem}
\label{lem:E_aff_6_increasing_denom_slice}
Let $\pth{p}=[3,4,6,2,1,7,5]$. Then, for $n \geq 1$, 
\begin{equation}
\label{eqn:E_aff_6_denom_slice}
\dslpm{\wt{E}_6}{[(\und{p})^n]}{3}=\left\langle  \floor[\Big]{\frac{n-1}{2}}, n-1, \ceil[\Big]{\frac{3n}{2}}-1, n-1,\floor[\Big]{\frac{n-1}{2}} , n-1, \floor[\Big]{\frac{n-1}{2}} \right\rangle.
\end{equation}
\end{lem}

\begin{proof}
We refer to the matrices in Lemma \ref{lem:R_aff_6_degseed_loop} throughout. As a base case, we have that
$\dslpm{\wt{E}_6}{[\und{p}]}{3}=\left\langle 0,0,1,0,0,0,0 \right\rangle.$
 Assume the claim is true for $n$. Then:

\begin{itemize}
\item$ \dv[5 ,(\und{p})^n]|_3 = -\floor[\Big]{\frac{n-1}{2}}+ \max \left( n-1,0 \right)= n-1-\floor[\Big]{\frac{n-1}{2}} = \ceil[\Big]{\frac{n-1}{2}},\\
\text{so }  \dslpm{\wt{E}_6}{[5,(\und{p})^n]}{3} =\left\langle  \floor[\Big]{\frac{n-1}{2}}, n-1, \ceil[\Big]{\frac{3n}{2}-1}, n-1,\ceil[\Big]{\frac{n-1}{2}} , n-1, \floor[\Big]{\frac{n-1}{2}} \right\rangle.$
\item $\dv[7,5 ,(\und{p})^n]|_3 = -\floor[\Big]{\frac{n-1}{2}}+ \max \left( n-1,0 \right)= n-1-\floor[\Big]{\frac{n-1}{2}} = \ceil[\Big]{\frac{n-1}{2}},\\
\text{so }  \dslp{[7,5,(\und{p})^n]}{3} =\left\langle  \floor[\Big]{\frac{n-1}{2}}, n-1, \ceil[\Big]{\frac{3n}{2}-1}, n-1,\ceil[\Big]{\frac{n-1}{2}} , n-1, \ceil[\Big]{\frac{n-1}{2}} \right\rangle.$
\item $\dv[1,7,5 ,(\und{p})^n]|_3 = -\floor[\Big]{\frac{n-1}{2}}+ \max \left( n-1,0 \right)=  \ceil[\Big]{\frac{n-1}{2}},\\
\text{so }  \dslp{[1,7,5,(\und{p})^n]}{3} =\left\langle  \ceil[\Big]{\frac{n-1}{2}}, n-1, \ceil[\Big]{\frac{3n}{2}-1}, n-1,\ceil[\Big]{\frac{n-1}{2}} , n-1, \ceil[\Big]{\frac{n-1}{2}} \right\rangle .$
\item $\dv[2,1,7,5 ,(\und{p})^n]|_3 = -(n-1)+ \max \left(\ceil[\Big]{\frac{3n}{2}} +\ceil[\Big]{\frac{n-1}{2}} -1 ,0 \right) \\
\phantom{\dv[2,1,7,5 ,(\und{p})^n]|_3} = \ceil[\Big]{\frac{3n}{2}}+ \ceil[\Big]{\frac{n-1}{2}}-(n-1)-1=n,\\
\text{so }  \dslp{[2,1,7,5,(\und{p})^n]}{3} = \left\langle  \ceil[\Big]{\frac{n-1}{2}}, n, \ceil[\Big]{\frac{3n}{2}-1}, n-1,\ceil[\Big]{\frac{n-1}{2}} , n-1, \ceil[\Big]{\frac{n-1}{2}} \right\rangle .$
\item$ \dv[6,2,1,7,5 ,(\und{p})^n]|_3 = -(n-1)+ \max \left(\ceil[\Big]{\frac{3n}{2}} +\ceil[\Big]{\frac{n-1}{2}} -1 ,0 \right) =n,\\
\text{so }  \dslp{[6,2,1,7,5,(\und{p})^n]}{3} =\left\langle  \ceil[\Big]{\frac{n-1}{2}}, n, \ceil[\Big]{\frac{3n}{2}-1}, n-1,\ceil[\Big]{\frac{n-1}{2}} , n, \ceil[\Big]{\frac{n-1}{2}} \right\rangle .$
\item $\dv[4,6,2,1,7,5 ,(\und{p})^n]|_3 = -(n-1)+ \max \left(\ceil[\Big]{\frac{3n}{2}} +\ceil[\Big]{\frac{n-1}{2}} -1 ,0 \right) =n,\\
\text{so }  \dslp{[4,6,2,1,7,5,(\und{p})^n]}{3} =\left\langle  \ceil[\Big]{\frac{n-1}{2}}, n, \ceil[\Big]{\frac{3n}{2}-1}, n,\ceil[\Big]{\frac{n-1}{2}} , n, \ceil[\Big]{\frac{n-1}{2}} \right\rangle.$
\end{itemize}

Finally
\begin{align*}
&\dv[\und{p}^{n+1}]|_3 = -\left(\ceil[\Big]{\frac{3n}{2}}-1 \right)+ \max \left(n+n+n ,0 \right) \\
&                             = 3n-\ceil[\Big]{\frac{3n}{2}}+1 =\ceil[\Big]{\frac{3n}{2}}+\ceil[\Big]{\frac{3n-1}{2}}-\ceil[\Big]{\frac{3n}{2}}+1 =\ceil[\Big]{\frac{3(n+1)}{2}}-1 ,
\end{align*}
where we have used the identity \pref{eqn:n_as_sum_of_ceils} with $m=2$.
So
\begin{align*}
 \dslp{[\und{p}^{n+1}]}{3} &=\left\langle  \ceil[\Big]{\frac{n-1}{2}}, n,\ceil[\Big]{\frac{3(n+1)}{2}}-1, n,\ceil[\Big]{\frac{n-1}{2}} , n, \ceil[\Big]{\frac{n-1}{2}} \right\rangle\\
 &=\left\langle  \floor[\Big]{\frac{n}{2}}, n,\ceil[\Big]{\frac{3(n+1)}{2}}-1, n, \floor[\Big]{\frac{n}{2}} , n,  \floor[\Big]{\frac{n}{2}} \right\rangle.
\end{align*}
Thus the result is also true for $n+1$. Therefore the claim is proved.
\end{proof}

\begin{prop}
$\A((1,0,1,0,1,0,1),\wt{E_6})$ has infinitely many variables of degree $0$ and $\pm 1$.
\end{prop}

\begin{proof}
By the proof of Lemma \ref{lem:R_aff_6_degseed_loop} we see that repeating the path $\pth{p}$ gives infinitely many variables of degree $m$, where $m \in \{ 0, \pm 1\}$. By Lemma \ref{lem:E_aff_6_increasing_denom_slice} we have that infinitely many of these paths must result in different variables as their denominator vectors are different.
\end{proof}

\begin{cnj}
$\A \big((1,0,1,0,1,0,1),\wt{E}_6 \big)$ has only one variable in degree $2$ and one variable in degree $-2$.
\end{cnj}

This statement may be compared with the analogous result for $\ssmat{-2}{-1}{1}$ in Proposition \ref{prop:finite_degrees_mixed_one_var_deg_2}. There, it was not difficult to establish the structure of the exchange graph, and its form was simple enough that we could deduce our result using the relevant set of recurrence relations. In this case, the higher rank makes it much more difficult to attempt the same method, but we still expect the result will hold. Computer aided calculation gives some further confidence in our conjecture: after computing thousands of different mutation paths that result in degree $2$ and $-2$, only one variable was found in each degree. Similar considerations also apply to Conjecture \ref{cnj:E_7} and Conjecture \ref{cnj:E_8}, which we will make later in the chapter.

\subsection{Grading for $\wt{E}_7$}
 
Let us now consider $\wt{E}_7$. This has a 2-dimensional grading under which we will show that the degrees $\CII{0}{0}$, $\pm \CII{1}{1}$, $\pm \CII{1}{0}$ and $\pm \CII{0}{1}$ all contain infinitely many variables. We conjecture that the degrees $ \pm \CII{2}{1}$ correspond to only one variable each.

\begin{lem}
\label{lem:E_aff_7_degseed_loop}
Let $\pth{p}=[7,1,6,2,5,8,3,4]$. Then \[\dsd_{\wt{E}_7}\pth{p^n}= \big((-1)^n \left(\CII{-1}{-1},\CII{0}{0},\CII{-1}{-1},\CII{0}{0},\CII{1}{0},\CII{0}{0}, \CII{1}{0}, \CII{0}{1} \right), \wt{E}_7 \big).\] That is, mutating the initial degree seed of $\wt{E}_7$ along $\pth{p}$ results in the negation of the initial degree seed.
\end{lem}

\begin{proof}
See Proof \ref{apdx:E_aff_7_degseed_loop} of the appendix for the proof, including the list of matrices obtained (which we will use again below).
\end{proof}

\begin{lem}
\label{lem:E_aff_7_increasing_denom_slice}
Let $\pth{p}=[7,1,6,2,5,8,3,4]$. Then 
\begin{equation}
\label{eqn:E_aff_7_denom_slice}
\dslpm{\wt{E}_7}{[(\und{p})^n]}{4}=\left\langle  \ceil[\Big]{\frac{n}{3}}, \ceil[\Big]{\frac{2n}{3}}, n,\ceil[\Big]{\frac{4n}{3}}-1 , n, \ceil[\Big]{\frac{2n}{3}}, \ceil[\Big]{\frac{n}{3}}, \ceil[\Big]{\frac{2n-1}{3}} \right\rangle.
\end{equation}
\end{lem}

\begin{proof}
The result is clearly true for $n=0$ which provides a base case. Assume true for $n$. We refer to the matrices in Proof \ref{apdx:E_aff_7_degseed_loop} throughout. This time we will not write down the new denominator slice after each mutation.

We have the following:
\begin{align*}
\dv[4 ,(\und{p})^n]|_4 &= -\left( \ceil[\Big]{\frac{4n}{3}} - 1 \right)+ \max \left\{ n+n+\ceil[\Big]{\frac{2n-1}{3}},0 \right\} \\
&=  2n + \ceil[\Big]{\frac{2n-1}{3}} - \ceil[\Big]{\frac{4n}{3}}+1 \\
&=n + \ceil[\Big]{\frac{2n-1}{3}} - \ceil[\Big]{\frac{n}{3}}+1 .\\
\dv[3,4 ,(\und{p})^n]|_4 &= -n + \max \left\{\ceil[\Big]{\frac{2n}{3}} + \left(n + \ceil[\Big]{\frac{n+1}{3}} \right),0 \right\} \\
&=  \ceil[\Big]{\frac{2n}{3}} +  \ceil[\Big]{\frac{n+1}{3}}.\\
\dv[8,3,4 ,(\und{p})^n]|_4 &= \ceil[\Big]{\frac{2n-1}{3}} + \max \left\{ n + \ceil[\Big]{\frac{n+1}{3}} ,0 \right\} \\
 &= n \ceil[\Big]{\frac{n+1}{3}} -1  \ceil[\Big]{\frac{2n-1}{3}}.\\
\dv[5,8,3,4 ,(\und{p})^n]|_4 &= -n + \max \left\{ \left(n + \ceil[\Big]{\frac{n+1}{3}} \right) + \ceil[\Big]{\frac{2n}{3}},0 \right\}\\
 &= \ceil[\Big]{\frac{n+1}{3}}+ \ceil[\Big]{\frac{2n}{3}}. \\
\dv[2,5,8,3,4 ,(\und{p})^n]|_4 &=  -\ceil[\Big]{\frac{2n}{3}} + \max \left\{  \ceil[\Big]{\frac{n}{3}} + \left( \ceil[\Big]{\frac{2n}{3}} + \ceil[\Big]{\frac{n+1}{3}} \right) ,0 \right\}\\
& = \ceil[\Big]{\frac{n+1}{3}}+ \ceil[\Big]{\frac{n}{3}}. 
\end{align*}
\begin{align*}
\dv[6,2,5,8,3,4 ,(\und{p})^n]|_4 &=  -\ceil[\Big]{\frac{2n}{3}} + \max \left\{  \left( \ceil[\Big]{\frac{n+1}{3}} + \ceil[\Big]{\frac{2n}{3}} \right) +\ceil[\Big]{\frac{n}{3}} ,0 \right\}\\
 &= \ceil[\Big]{\frac{n+1}{3}}+ \ceil[\Big]{\frac{n}{3}}. \\
\dv[1,6,2,5,8,3,4 ,(\und{p})^n]|_4 &=  -\ceil[\Big]{\frac{n}{3}} + \max \left\{  \left( \ceil[\Big]{\frac{n+1}{3}} + \ceil[\Big]{\frac{n}{3}} \right),0 \right\} = \ceil[\Big]{\frac{n+1}{3}}. \\
\dv[\und{p}^{n+1}]|_4 &=  -\ceil[\Big]{\frac{n}{3}} + \max \left\{  \left( \ceil[\Big]{\frac{n+1}{3}} + \ceil[\Big]{\frac{n}{3}} \right),0 \right\} = \ceil[\Big]{\frac{n+1}{3}}. 
\end{align*}

So the new denominator slice, $\dslpm{\wt{E}_7}{[(\und{p})^n]}{4}$, is

\begin{align*}
\Bigg\langle & \ceil[\Big]{\frac{n+1}{3}}, \ceil[\Big]{\frac{n+1}{3}}+\ceil[\Big]{\frac{n}{3}}, \ceil[\Big]{\frac{2n}{3}}+ \ceil[\Big]{\frac{n+1}{3}}, n + \ceil[\Big]{\frac{2n-1}{3}} - \ceil[\Big]{\frac{n}{3}}+1 ,\\
&  \ceil[\Big]{\frac{n+1}{3}}+\ceil[\Big]{\frac{2n}{3}}, \ceil[\Big]{\frac{n+1}{3}}+\ceil[\Big]{\frac{n}{3}}, \ceil[\Big]{\frac{n+1}{3}}, n+\ceil[\Big]{\frac{n+1}{3}}-\ceil[\Big]{\frac{2n-1}{3}} \Bigg\rangle.
\end{align*}

For the induction step, we need to check that

\begin{enumerate}[(i)]
\item $ \ceil[\Big]{\frac{n+1}{3}}+\ceil[\Big]{\frac{n}{3}} =  \ceil[\Big]{\frac{2(n+1)}{3}}$,
\item $ \ceil[\Big]{\frac{2n}{3}}+ \ceil[\Big]{\frac{n+1}{3}} = n+1$,
\item $n + \ceil[\Big]{\frac{2n-1}{3}} - \ceil[\Big]{\frac{n}{3}}+1 = \ceil[\Big]{\frac{4(n+1)}{3}}-1 $, and
\item $n + \ceil[\Big]{\frac{n+1}{3}} - \ceil[\Big]{\frac{2n-1}{3}} = \ceil[\Big]{\frac{2(n+1)-1}{3}} $.
\end{enumerate}

This is straightforward. We write $n=3k+i$, for some $k \in \N$ and some $i \in \{0,1,2\}$ and split into cases based on $n \bmod 3$. We show the calculation for (i) and assert that the other cases are easily shown in a similar way.
\begin{align*}
&\text{For } n=3k \text{: } \\
&  \ceil[\Big]{\frac{n+1}{3}} +  \ceil[\Big]{\frac{n}{3}} = \ceil[\Big]{\frac{3k+1}{3}} +  \ceil[\Big]{\frac{3k}{3}} = k +   \ceil[\Big]{\frac{1}{3}} +k = 2k+1 \\
& \phantom{ \ceil[\Big]{\frac{n+1}{3}} +  \ceil[\Big]{\frac{n}{3}}} = 2k+  \ceil[\Big]{\frac{2}{3}} = \ceil[\Big]{\frac{6k+2}{3}} =   \ceil[\Big]{\frac{2(n+1)}{3}}. \\ 
&\text{For } n=3k+1 \text{: } \\
&\ceil[\Big]{\frac{n+1}{3}} +  \ceil[\Big]{\frac{n}{3}} = \ceil[\Big]{\frac{3k+2}{3}} +  \ceil[\Big]{\frac{3k+1}{3}} = 2k+2 = \ceil[\Big]{\frac{6k+4}{3}} =   \ceil[\Big]{\frac{2(n+1)}{3}}. \\
&\text{For } n=3k+2 \text{: } \\
&\ceil[\Big]{\frac{n+1}{3}} +  \ceil[\Big]{\frac{n}{3}} = \ceil[\Big]{\frac{3k+3}{3}} +  \ceil[\Big]{\frac{3k+2}{3}} = 2k+2 = \ceil[\Big]{\frac{6k+6}{3}} =   \ceil[\Big]{\frac{2(n+1)}{3}}. 
\end{align*}
Thus we obtain the induction step which gives our result.
\end{proof}

\begin{prop}
$\A\left(  \left( \CII{-1}{-1},\CII{0}{0},\CII{-1}{-1},\CII{0}{0},\CII{1}{0},\CII{0}{0}, \CII{1}{0}, \CII{0}{1} \right),\wt{E}_7  \right)$ has infinitely many variables of degree $\CII{0}{0}$, $\pm \CII{1}{1}$, $\pm \CII{1}{0}$ and $\pm \CII{0}{1}$.
\end{prop}

\begin{proof}
This follows by combining Lemma \ref{lem:E_aff_7_degseed_loop}  and Lemma \ref{lem:E_aff_7_increasing_denom_slice}
\end{proof}

\begin{cnj} \label{cnj:E_7}
$\A\left(  \left( \CII{-1}{-1},\CII{0}{0},\CII{-1}{-1},\CII{0}{0},\CII{1}{0},\CII{0}{0}, \CII{1}{0}, \CII{0}{1} \right),\wt{E_7}  \right)$ has only one variable in degree $ \CII{2}{1}$ and one variable in degree $ \CII{-2}{-1}$.
\end{cnj}

\subsection{Grading for $\wt{E}_8$}

In $\wt{E}_8$, we again expect degrees $\pm 2$ correspond to one variable each, while degrees $0$ and $\pm 1$ have infinitely many different variables.
For this case it is more difficult to write down a simple formula for the entries of the third denominator slice. The entries after $[(\und{p})^n]$ depend on $n \bmod 5$.

\begin{lem} \label{lem:E_aff_8_degseed_loop}
Let $\pth{p}=[3,4,2,6,9,1,6,5,6,7,8]$. Then \[\dsd_{\wt{E}_8}[(\und{p})^n] = \big((-1)^n (0,0,0,-1,0,-1,0,-1,1), \wt{E}_8 \big). \]
\end{lem}

\begin{proof}
See Proof \ref{apdx:E_aff_8_degseed_loop} of the appendix, where again we include the list of matrices obtained.
\end{proof}

\begin{lem} \label{lem:E_aff_8_increasing_denom_slice}
Let $\pth{p}=[3,4,2,6,9,1,6,5,6,7,8]$. The entries of $\dslpm{\wt{E}_8}{[(\und{p})^n]}{3}$ are as follows.
\begin{align}
  \dcl^1_3 [(\und{p})^n] &= \begin{cases} \label{eqn:E_aff_8_entry_1}
    2 \ceil[\Big]{\frac{n}{5}}-2, & n=5k+i, (i=1,2)\\
    2 \ceil[\Big]{\frac{n}{5}}-1, & n=5k+j, (j=0,3,4).
  \end{cases}
\\
  \dcl^2_3 [(\und{p})^n] &= \begin{cases}  \label{eqn:E_aff_8_entry_2}
   4 \ceil[\Big]{\frac{n}{5}}-4, & n=5k+1\\
   4 \ceil[\Big]{\frac{n}{5}}-3, & n=5k+2\\
   4 \ceil[\Big]{\frac{n}{5}}-2, & n=5k+i, (i=3,4)\\
   4 \ceil[\Big]{\frac{n}{5}}-1, & n=5k       .  
  \end{cases}
\\
  \dcl^3_3 [(\und{p})^n] &= \begin{cases}  \label{eqn:E_aff_8_entry_3}
   6 \ceil[\Big]{\frac{n}{5}}-5, & n=5k+1\\
   6 \ceil[\Big]{\frac{n}{5}}-4, & n=5k+2\\
   6 \ceil[\Big]{\frac{n}{5}}-3, & n=5k+3\\
   6 \ceil[\Big]{\frac{n}{5}}-2, & n=5k+4\\      
   6 \ceil[\Big]{\frac{n}{5}}-1, & n=5k     .    
  \end{cases}  
\end{align}  
\begin{align}
  \dcl^4_3 &= n-1  \label{eqn:E_aff_8_entry_4}
\\
  \dcl^5_3 [(\und{p})^n] &= \begin{cases}  \label{eqn:E_aff_8_entry_5}
   4 \ceil[\Big]{\frac{n}{5}}-4, & n=5k+i, (i=1,2)\\
   4 \ceil[\Big]{\frac{n}{5}}-3, & n=5k+3\\
   4 \ceil[\Big]{\frac{n}{5}}-2, & n=5k+4\\
   4 \ceil[\Big]{\frac{n}{5}}-1, & n=5k  .       
  \end{cases}
\\  
  \dcl^6_3 [(\und{p})^n] &= \begin{cases}  \label{eqn:E_aff_8_entry_6}
   3 \ceil[\Big]{\frac{n}{5}}-3, & n=5k+i, (i=1,2,3)\\
   3 \ceil[\Big]{\frac{n}{5}}-2, & n=5k+4\\
   3 \ceil[\Big]{\frac{n}{5}}-1, & n=5k.
  \end{cases}  
\\
  \dcl^7_3 [(\und{p})^n] &= \begin{cases}  \label{eqn:E_aff_8_entry_7}
    2 \ceil[\Big]{\frac{n}{5}}-2, & n=5k+i, (i=1,2,3,4)\\
    2 \ceil[\Big]{\frac{n}{5}}-1, & n=5k.
  \end{cases}  
\\
  \dcl^8_3[\und{p}^{n}] &=  \ceil[\Big]{\frac{n}{5}}-1.  \label{eqn:E_aff_8_entry_8}
\\  
  \dcl^9_3 [(\und{p})^n] &= \begin{cases}  \label{eqn:E_aff_8_entry_9}
   3 \ceil[\Big]{\frac{n}{5}}-3, & n=5k+1\\
   3 \ceil[\Big]{\frac{n}{5}}-2, & n=5k+i , (i=2,3)\\
   3 \ceil[\Big]{\frac{n}{5}}-1, & n=5k+j , (j=0,4).
  \end{cases}    
\end{align}
\end{lem}

\begin{proof}
The result is clearly true for $n=0$ which provides a base case. Assume true for $n$ when $n=5k$. (Other cases can be proved in a similar way.) We refer to the matrices in Proof \ref{apdx:E_aff_8_degseed_loop}. Note, since  $n=5k$, $\ceil[\Big]{\frac{n+1}{5}} = \ceil[\Big]{\frac{n}{5}}+1$ and  $\ceil[\Big]{\frac{n}{5}}=\frac{n}{5}$.

We compute the denominator entries along the mutation path.
\begin{align*}
\dv^3[8 ,(\und{p})^n] &= -\left( \ceil[\Big]{\frac{n}{5}} - 1 \right)+ \max \left(2\ceil[\Big]{\frac{n}{5}}-1,0 \right)\\
                             &= \ceil[\Big]{\frac{n}{5}} .\\
\dv^3[7,8 ,(\und{p})^n] &= -\left(2 \ceil[\Big]{\frac{n}{5}} -1 \right)+ \max \left( \left(3\ceil[\Big]{\frac{n}{5}}-1 \right) + \ceil[\Big]{\frac{n}{5}},0 \right)\\
                                &= 2\ceil[\Big]{\frac{n}{5}}.\\
\dv^3[6,7,8 ,(\und{p})^n] &= -\left(3 \ceil[\Big]{\frac{n}{5}} -1 \right)+ \max \left( \left(4\ceil[\Big]{\frac{n}{5}}-1 \right) + 2\ceil[\Big]{\frac{n}{5}},0 \right)\\
                                  &= 3\ceil[\Big]{\frac{n}{5}}.\\
\dv^3[5,6,7,8 ,(\und{p})^n] &= -\left(4 \ceil[\Big]{\frac{n}{5}} -1 \right)+ \max \left( (n-1) + 3\ceil[\Big]{\frac{n}{5}},0 \right)\\
                                  &= n-\ceil[\Big]{\frac{n}{5}}=n-\frac{n}{5}=4\ceil[\Big]{\frac{n}{5}}.\\
\dv^3[6,5,6,7,8 ,(\und{p})^n] &= -3 \ceil[\Big]{\frac{n}{5}} + \max \left( \left(n- \ceil[\Big]{\frac{n}{5}} \right), 2\ceil[\Big]{\frac{n}{5}} \right)\\
                                  &= n-4\ceil[\Big]{\frac{n}{5}} =\ceil[\Big]{\frac{n}{5}}.
\end{align*}                                  
\begin{align*}                                  
\dv^3[1,6,5,6,7,8 ,(\und{p})^n] &= -\left(2 \ceil[\Big]{\frac{n}{5}} -1 \right)+ \max \left( \left(4\ceil[\Big]{\frac{n}{5}}-1 \right) ,0 \right)\\
                                          &= 2\ceil[\Big]{\frac{n}{5}}.\\
\dv^3[9,1,6,5,6,7,8 ,(\und{p})^n] &= -\left(3 \ceil[\Big]{\frac{n}{5}} -1 \right)+ \max \left( \left(6\ceil[\Big]{\frac{n}{5}}-1 \right) ,0 \right)\\
                                             &= 3\ceil[\Big]{\frac{n}{5}}.\\
\dv^3[6,9,1,6,5,6,7,8 ,(\und{p})^n] &= \resizebox{0.6\hsize}{!}{$-\left(n- 4\ceil[\Big]{\frac{n}{5}} -1 \right)+ \max \left( \left(n-\ceil[\Big]{\frac{n}{5}}-1 \right) ,2\ceil[\Big]{\frac{n}{5}} \right) $}\\
                                                &= 4\ceil[\Big]{\frac{n}{5}}.\\
\dv^3[2,6,9,1,6,5,6,7,8 ,(\und{p})^n] &= -\left(4 \ceil[\Big]{\frac{n}{5}} -1 \right)+ \max \left(  2\ceil[\Big]{\frac{n}{5}}+\left(6\ceil[\Big]{\frac{n}{5}}-1 \right),0 \right)\\
                                                   &= 4\ceil[\Big]{\frac{n}{5}}.\\
\dv^3[4,2,6,9,1,6,5,6,7,8 ,(\und{p})^n] &= -(n-1)+ \max \left(\left(6\ceil[\Big]{\frac{n}{5}}-1 \right) +  \left(n-\ceil[\Big]{\frac{n}{5}} \right),0 \right)\\
                                                     &= 5\ceil[\Big]{\frac{n}{5}}.
\end{align*}
Finally,
\begin{align*}
\dv^3[\und{p}^{n+1}] &= -\left(6\ceil[\Big]{\frac{n}{5}}-1 \right) + \max \left(4\ceil[\Big]{\frac{n}{5}}  +5\ceil[\Big]{\frac{n}{5}}  +3\ceil[\Big]{\frac{n}{5}}  ,0 \right)\\
                                   &= 6\ceil[\Big]{\frac{n}{5}}+1.
\end{align*}

So the new denominator slice is 
\[\dslpm{\wt{E}_8}{[\und{p}^{n+1}]}{3} =\ceil[\Big]{\frac{n}{5}} \langle 2,4,6,5,4,3,2,1,3 \rangle + \langle 0,0,1,0,0,0,0,0,0 \rangle .\] 
Considering the first entry, we have $2\ceil[\Big]{\frac{n}{5}} = 2\left(\ceil[\Big]{\frac{n}{5}}+1 \right)-2 = 2\ceil[\Big]{\frac{n+1}{5}} - 2$, so \pref{eqn:E_aff_8_entry_1} is true for $n+1$ (which is of the form $5k+1$). 
Similarly, so are \pref{eqn:E_aff_8_entry_2}--\pref{eqn:E_aff_8_entry_9}.
So, when $n=5k$, the result is true for $n+1$.
\end{proof}

\begin{prop}
 $\A\big((0,0,0,-1,0,-1,0,-1,1), \wt{E}_8 \big)$ has infinitely many variables of degree $0$ and $\pm 1$.
\end{prop}

\begin{proof}
By the proof of Lemma \ref{lem:E_aff_8_degseed_loop} we see that repeating the path $\pth{p}$ gives infinitely many variables of degree $m$, where $m \in \{ 0, \pm 1\}$. By Lemma \ref{lem:E_aff_8_increasing_denom_slice} we have that infinitely many of these paths must result in different variables since the denominator vectors obtained grow indefinitely.
\end{proof}

\begin{cnj} \label{cnj:E_8}
$\A \big((0,0,0,-1,0,-1,0,-1,1), \wt{E}_8 \big)$ has only one variable in degree $2$ and one variable in degree $-2$.
\end{cnj}

\biblio

\chapter{Gradings for quivers arising from surface triangulations} \label{chp:surface_type}

\section{Graded cluster algebras associated to surfaces} \label{sec:surface_theory}

There exists a class of cluster algebras associated to oriented bordered surfaces with marked points. In \cite{FST} the authors describe the process by which a cluster algebra arises from such a surface. As is explained in \cite{Muller}, these cluster algebras may be given a grading by assigning each cluster variable a degree which is the sum of valued marked points (though this grading does not agree precisely with the definition of grading we have assumed in this thesis). In this chapter we will review the theory of cluster algebras arising from surfaces, following \cite{FST}, and adapt a definition from \cite{Muller} to our setting in order to define a grading (in our sense) on such a cluster algebra arising from the properties of its associated surface. We will then apply this theory to study the graded cluster algebra structure of some classes of examples.

The basic idea when associating a surface with a cluster algebra is that the surface will be triangulated by arcs, each of which represents a cluster variable, and the configuration of the triangulation then also encodes an exchange quiver so that the triangulated surface can be considered as a seed. This is made more precise with the following definitions.

\begin{dfn}[\emph{Marked surfaces}]
A \emph{bordered surface with marked points} is a pair $(\bf{S},\bf{M})$ where $\bf{S}$ is a connected oriented $2$-dimensional Riemann surface with boundary and $\bf{M}$ is a finite, non-empty set of $\emph{marked points}$ in the closure of $\bf{S}$ such that each connected component of the boundary of $\bf{S}$ has at least one marked point. 

However, for technical reasons, the following particular cases are excluded from the definition: a sphere with one or two punctures, an unpunctured or once-punctured monogon, an unpunctured digon, an unpunctured triangle and a sphere with three punctures. This excludes cases that are not able to be \emph{triangulated}---which we will define shortly.

A marked point which is not on a boundary component is called a \emph{puncture}. (However, when we consider gradings arising from surfaces later, we will not allow surfaces with punctures.)
\end{dfn}

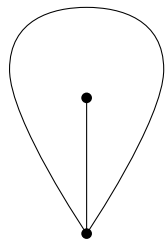
\begin{figure}[h]
\begin{center}
\Image{ \begin{tikzpicture}
  \def\rotation{90} 
  \def\radius{3cm}
  \def\innerradius{ {(4/5)*\radius} }  
  \coordinate (O) at (0,0);
  \path (O) ++(\rotation+25:\innerradius) coordinate (A1);
  \path (O) ++(\rotation+0:\radius) coordinate (A2);
  \path (O) ++(\rotation+335:\innerradius) coordinate (A3); 
  \path (O) ++(\rotation+0:{(3/5)*\radius}) coordinate (M);   
  \fill[black] (O) circle[radius=2pt] ++(\rotation+90:0.8em) node {\scriptsize $ $};
  \fill[black] (M) circle[radius=2pt] ++(\rotation+90:0.8em) node {\scriptsize $ $};  
  \fill[black] (A2)  ++(\rotation+180:1.5em) node {\scriptsize $ $};
     \draw [black, xshift=4cm] plot [smooth, tension=.8] coordinates { (O) (A1) (A2) (A3) (O)};  
     \draw [black, xshift=4cm] plot [smooth, tension=.8] coordinates { (O) (M)};  
\end{tikzpicture}}
\end{center}
\caption{The self-folded triangle.} 
\label{fig:self_folded_triangle} 
\end{figure}

\begin{dfn}[\emph{Arcs and ideal triangulations}]
An arc in $(\bf{S},\bf{M})$ is a curve (up to isotopy) in $\bf{S}$ that does not intersect itself (except, possibly, for its endpoints) whose endpoints (and only its endpoints) are marked points in $\bf{M}$. An arc is not allowed to cut out an unpunctured monogon or unpunctured digon, so that an arc cannot be contractible to a marked point or to the boundary of $\bf{S}$. We use $\bf{A}^\circ (\bf{S},\bf{M})$ to denote the set of all arcs in $(\bf{S},\bf{M})$ (which is typically infinite---see \cite[Proposition 2.3]{FST}).

Two arcs are called \emph{compatible} if they have representatives in their respective isotopy classes that do not intersect in the interior of $\bf{S}$. An \emph{ideal triangulation} of  $(\bf{S},\bf{M})$ is a maximal collection of distinct, pairwise compatible arcs, and these arcs cut $\bf{S}$ into \emph{ideal triangles}. Thus, each side of an ideal triangle is an arc or a segment of a boundary component between two marked points. Deviating from \cite{FST} slightly, we will use the term \emph{boundary arc} to refer to such a segment, even though it is not a genuine arc. The three sides of an ideal triangle are not required to be distinct; such a triangle (shown in Figure \ref{fig:self_folded_triangle}) is called \emph{self-folded} (although, as above, these will not be allowed when we consider gradings). Furthermore, two triangles can share more than one side.
\end{dfn}

The following defines the analogue of seed mutation for triangulated surfaces.

\begin{dfn}[\emph{Arc complex and flips}]
The \emph{arc complex} is the simplicial complex on the ground set $\bf{A}^\circ(\bf{S},\bf{M})$ whose simplices are collections of distinct mutually compatible arcs and whose maximal simplices are ideal triangulations. We denote the arc complex by $\Delta^\circ(\bf{S},\bf{M})$ and its dual graph by $\bf{E}^\circ(\bf{S},\bf{M})$. So the vertices of $\bf{E}^\circ(\bf{S},\bf{M})$ are in bijection with the ideal triangulations of $(\bf{S},\bf{M})$.

A \emph{flip} is a transformation of an ideal triangulation of $(\bf{S},\bf{M})$ that removes a particular arc and replaces it with a (unique) new arc that forms a new ideal triangulation of  $(\bf{S},\bf{M})$ together with the remaining arcs. Thus, the edges of $\bf{E}^\circ(\bf{S},\bf{M})$ above correspond to flips.
\end{dfn}

It is clear that a flip of a triangulated surface is involutive, though not every arc can be flipped: there is no suitable replacement arc for the self-folded edge (interior to the loop) of the self-folded triangle mentioned above. This is in contrast with seed mutation, in which every (non-frozen) entry of the seed may be mutated. In \cite[Section 7]{FST}, the authors resolve this issue by introducing \emph{tagged arcs} and their flips. However, for the class of examples from which our gradings arise, every arc will be flippable (indeed, as we will not allow punctures, the self-folded triangle can not arise in our examples).

\begin{prop}
\label{prop:dual_arc_complex_connected}
Any two ideal triangulations are related by a sequence of flips. That is, $\bf{E}^\circ(\bf{S},\bf{M})$ is connected.
\end{prop}

\begin{proof}
This follows from results of \cite{Har}, \cite{Hat} and \cite{Mos}.
\end{proof}

Compare the following to Theorem \ref{thm:exchgraph_cycles}.

\begin{thm}[{\cite[Theorem 1.1]{Che}}]
The fundamental group of $\bf{E}^\circ(\bf{S},\bf{M})$ is generated by cycles of length 4 and 5, pinned down to a base point.
\end{thm}

Still following \cite{FST}, we now explain how to associate an exchange matrix to a triangulated surface.

\begin{dfn}[{\cite[Definition 4.1]{FST}}]
Let $T$ be an ideal triangulation of $(\bf{S},\bf{M})$. We define the \emph{signed adjacency matrix} $B(T)$ as follows. We label the rows and columns of $B(T)$ by the arcs in $T$ (strictly, by $l(1), \dots, l(n)$ for some labelling $l$ as in Definition \ref{dfn:labelled_set}). For an arc labelled $i$, the map $\pi_T(i)$ is defined as follows: if there is a self-folded ideal triangle in $T$ in which $i$ is the self-folded edge, then $\pi_T(i)$ is the arc that is the remaining edge of the triangle. Otherwise, $\pi_T(i) = i$. Then, for each ideal triangle $\Delta$ in $T$ which is not self-folded, define the $n \times n$ matrix $B^\Delta$ by 
\begin{equation}
b_{ij}^\Delta=
\begin{cases}
 1      &\text{if $\Delta$ has sides labeled $\pi_T(i)$ and $\pi_T(j)$}\\
         & \text{ with $\pi_T(i)$ followed by $\pi_T(j)$ in the clockwise order}, \\
 -1     &\text{if the same is true but in the counter-clockwise order,}\\
 0      &\text{otherwise.}
\end{cases}
\end{equation}
Then we set $ B(T) = \sum_\Delta B^\Delta$, where $\Delta$ ranges over all ideal triangles in $T$ that are not self-folded. 
\end{dfn}

We have that $B(T)$ is skew-symmetric and all its entries are in $\{0, \pm 1, \pm 2 \}$. 

\begin{rmk}
\label{rmk:quiver_of_triangulated_surface}
As an alternative to using a matrix, we may draw the corresponding quiver on a triangulated surface by picking a point on each arc and adding an arrow between arcs $\alpha$ and $\alpha'$ each time $\alpha'$ appears after $\alpha$ in the clockwise direction in a triangle (which may happen more than once for a given pair of arcs, leading to a double arrow or cancellation of arrows).
\end{rmk}

\begin{rmk}
Suppose $T$ and $T'$ are two triangulations of the same marked surface. It is not true that $ T \neq T' \implies B(T) \neq B(T')$. Indeed, unless $B(T)$ gives rise to a finite type cluster algebra, there will be infinitely many different triangulations $T'$ such that $B(T')= B(T)$. Note also that two different marked
 surfaces may have the same signed adjacency matrix (for example, see \break {\cite[Example 4.3]{FST}}).
\end{rmk}

\begin{ex}
Consider the triangulated marked surface
\begin{center} $(\bf{S},\bf{M})=$
\Image{ \begin{tikzpicture}
  \coordinate (O) at (0,0);  
  \def\radius{2.5cm}
  \def\midradius{ {(3/5)*\radius} }  
  \def\outermidradius{ {(3.8/5)*\radius} }      
  \def\innermidradius{ {(2.3/5)*\radius} }      
  \def\innerradius{  {(0.4/3)*\radius} }  
  \centerarc[](O)(0:360:\radius)    
  \path (O) ++(180+270:\radius) coordinate (B);
  \path (O) ++(180+90:{(0.85/3)*\radius} ) coordinate (GC);  
  \path (O) ++(180+170:{(0.75/3)*\radius} ) coordinate (GL);  
  \path (O) ++(180+10:{(0.75/3)*\radius} ) coordinate (GR);      
  \path (O) ++(180+90:\radius) coordinate (OU);  
  \path (O) ++(180+150:\radius) coordinate (OL1); \path (O) ++(180+150:\outermidradius) coordinate (GOL1);   
  \path (O) ++(180+210:\radius) coordinate (OL2); \path (O) ++(180+210:\outermidradius) coordinate (GOL2);   
  \path (O) ++(180+30:\radius) coordinate (OR1); \path (O) ++(180+30:\outermidradius) coordinate (GOR1);   
  \path (O) ++(150:\outermidradius) coordinate (GIR2);     
  \path (O) ++(180+330:\radius) coordinate (OR2); \path (O) ++(180+330:\outermidradius) coordinate (GOR2);       
  \path (O) ++(180+90:\innermidradius) coordinate (GIR1);                 
  \path (O) ++(180+210:\innermidradius) coordinate (GIR2);                 
  \path (O) ++(180+330:\innermidradius) coordinate (GIR3);                   
   \fill[black] (OU) circle[radius=2pt] ++(180+90:1em) node {\scriptsize $ $};
   \fill[black] (OL1) circle[radius=2pt] ++(310:1em) node {\scriptsize $ $};   
   \fill[black] (OL2) circle[radius=2pt] ++(30:1em) node {\scriptsize $ $};
   \fill[black] (OR1) circle[radius=2pt] ++(200:1em) node {\scriptsize $ $};   
   \fill[black] (OR2) circle[radius=2pt] ++(160:1em) node {\scriptsize $ $};
   \fill[black] (B) circle[radius=2pt] ++(180+270:1em) node {\scriptsize $ $};   
   \fill[black] (GIR1)  ++(270:1em) node {\small $\alpha_1$};              
   \fill[black] (GIR2)  ++(0:1em) node {\small $\alpha_2$};      
   \fill[black] (GIR3)  ++(180:1em) node {\small $\alpha_3$};            
   \draw [black, xshift=4cm] plot [smooth, tension=0.5] coordinates { (B) (OL1)};            
   \draw [black, xshift=4cm] plot [smooth, tension=0.5] coordinates { (B) (OR1)};     
   \draw [black, xshift=4cm] plot [smooth, tension=0.5] coordinates { (OL1) (OR1)};        
\end{tikzpicture}}.
\end{center}
This is the six-gon. The ideal triangulation $T$ consists of the three non-boundary arcs $\alpha_1$, $\alpha_2$ and $\alpha_3$ and we have
$B(T) = \ssmat{1}{1}{-1}.$ Flipping the arc $\alpha_1$, we obtain the new ideal triangulation $T'=\{ \alpha_1',\alpha_2,\alpha_3 \}$, as in
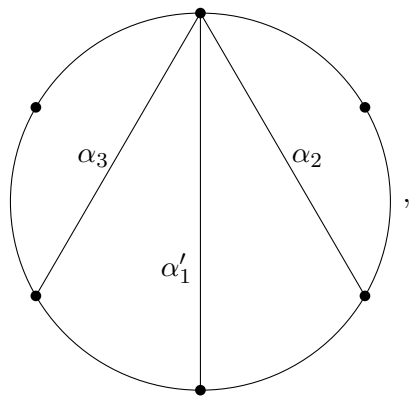
\begin{center}
\Image{ \begin{tikzpicture}
  \coordinate (O) at (0,0);  
  \def\radius{2.5cm}
  \def\midradius{ {(3/5)*\radius} }  
  \def\outermidradius{ {(3.8/5)*\radius} }      
  \def\innermidradius{ {(2.3/5)*\radius} }      
  \def\innerradius{  {(0.4/3)*\radius} }  
  \centerarc[](O)(0:360:\radius)    
  \path (O) ++(180+270:\radius) coordinate (B);
  \path (O) ++(180+90:{(0.85/3)*\radius} ) coordinate (GC);  
  \path (O) ++(180+170:{(0.75/3)*\radius} ) coordinate (GL);  
  \path (O) ++(180+10:{(0.75/3)*\radius} ) coordinate (GR);      
  \path (O) ++(180+90:\radius) coordinate (OU);  
  \path (O) ++(180+150:\radius) coordinate (OL1); \path (O) ++(180+150:\outermidradius) coordinate (GOL1);   
  \path (O) ++(180+210:\radius) coordinate (OL2); \path (O) ++(180+210:\outermidradius) coordinate (GOL2);   
  \path (O) ++(180+30:\radius) coordinate (OR1); \path (O) ++(180+30:\outermidradius) coordinate (GOR1);   
  \path (O) ++(150:\outermidradius) coordinate (GIR2);     
  \path (O) ++(180+330:\radius) coordinate (OR2); \path (O) ++(180+330:\outermidradius) coordinate (GOR2);       
  \path (O) ++(180+90:\innermidradius) coordinate (GIR1);                 
  \path (O) ++(180+210:\innermidradius) coordinate (GIR2);                 
  \path (O) ++(180+330:\innermidradius) coordinate (GIR3);                   
   \fill[black] (OU) circle[radius=2pt] ++(180+90:1em) node {\scriptsize $ $};
   \fill[black] (OL1) circle[radius=2pt] ++(310:1em) node {\scriptsize $ $};   
   \fill[black] (OL2) circle[radius=2pt] ++(30:1em) node {\scriptsize $ $};
   \fill[black] (OR1) circle[radius=2pt] ++(200:1em) node {\scriptsize $ $};   
   \fill[black] (OR2) circle[radius=2pt] ++(160:1em) node {\scriptsize $ $};
   \fill[black] (B) circle[radius=2pt] ++(180+270:1em) node {\scriptsize $ $};   
   \fill[black] (GIR1)  ++(140:1em) node {\small $\alpha_1'$};              
   \fill[black] (GIR2)  ++(0:1em) node {\small $\alpha_2$};      
   \fill[black] (GIR3)  ++(180:1em) node {\small $\alpha_3$};            
   \draw [black, xshift=4cm] plot [smooth, tension=0.5] coordinates { (B) (OL1)};            
   \draw [black, xshift=4cm] plot [smooth, tension=0.5] coordinates { (B) (OR1)};     
   \draw [black, xshift=4cm] plot [smooth, tension=0.5] coordinates { (B) (OU)};        
\end{tikzpicture}},
\end{center}
and we have $B(T')=\left( \begin{smallmatrix}0 & 1&-1\\-1 & 0&0 \\1 & 0&0\end{smallmatrix} \right)$. The corresponding quiver is an orientation of a Dynkin diagram of type $A_3$, and indeed the cluster algebra associated to this marked surface is of type $A_3$. (In general, the cluster algebra associated to the $n$-gon with no punctures is of type $A_n$---see \cite[Table 1]{FST}.)
\end{ex}

Notice in the example above that $B(T')=B(T)_{[1]}$, that is, the matrix corresponding to flipping the arc $\alpha_1$ in $T$ is the same as the one obtained by mutating $B(T)$ in direction $1$. This happens in general, as we see from the next result, which relates flips of arcs to matrix mutation.

\begin{prop}[{\cite[Proposition 4.8]{FST}}]
\label{prop:flip_vs_matrix_mutation}
Let $T$ be an ideal triangulation and suppose the ideal triangulation $T'$ is obtained from $T$ by flipping an arc labelled $k$. Then $B(T') = B(T)_{[k]}$. 
\end{prop}

The following is a corollary of Proposition \ref{prop:dual_arc_complex_connected} and Proposition \ref{prop:flip_vs_matrix_mutation}.

\begin{prop}[{\cite[Proposition 4.10]{FST}}]
Let $T$ be a triangulation of $(\bf{S},\bf{M})$. The mutation equivalence class of $B(T)$ is independent of $T$ and depends only on $(\bf{S},\bf{M})$. 
\end{prop}

\begin{dfn}
We will write $\A (\bf{S}, \bf{M})$ to unambiguously refer to the cluster algebra arising from a particular triangulated surface. 
\end{dfn}

We have not yet precisely established the relationship between the arcs in the class of ideal triangulations of a given surface and the cluster variables in the corresponding cluster algebra, that is, the relationship between the arc complex $\Delta^\circ(\bf{S},\bf{M})$ of a triangulated surface and the cluster complex of the cluster algebra arising from that surface. The next theorem, which is the main result of \cite{FST}, addresses this. The version presented here is in fact only a special case of that theorem (which allows for the possibility of self-folded triangles by using tagged triangulations), but it is all we will need as we will ultimately be excluding surfaces with punctures. 

\begin{thm}[{\cite[Theorem 7.11]{FST}}]
\label{thm:FST_main_theorem}
Let $(\bf{S},\bf{M})$ be a bordered surface with marked points without any punctures and let $\A$ be a cluster algebra generated by $B(T)$, where $T$ is a triangulation of $(\bf{S},\bf{M})$. Then the arc complex $\Delta^\circ(\bf{S},\bf{M})$ is isomorphic to the cluster complex of $\A$, and $\bf{E}^\circ(\bf{S},\bf{M})$ is isomorphic to the exchange graph of $\A$.
\end{thm} 

Apart from one or two particular cases, this result can be extended to allow for surfaces with punctures. We will exclude such surfaces in our study as they do not work with the set-up we are going to introduce that allows a grading to arise from a triangulation. 

We are now ready to introduce the mechanism by which such a grading arises. The idea is to assign values to the marked points of a surface and then set the degree of an arc (and thus its corresponding cluster variable) to the sum of its endpoints. First note the following, which lets us read off the dimension of the grading from $(\bf{S},\bf{M})$.

\begin{thm}[{\cite[Theorem 14.3]{FST}}]
\label{thm:FST_even_boundaries_gives_rank}
Let $T$ be an ideal triangulation of $(\bf{S},\bf{M})$. Then the corank of $B(T)$ is the number of punctures in $(\bf{S},\bf{M})$ plus the number of boundary components having an even number of marked points.
\end{thm}

We will refer to boundary components with an even number of marked points as \emph{even} boundary components and ones with an odd number of points as \emph{odd} components. Thus, for our class of examples, the dimension of a grading is given by the number of even boundary components.

The set $\mathsf{E}$ in Definition \ref{dfn:surface_grading_space_E} below is where we make an alteration to the definition in \cite[Section 3.5]{Muller}. This set will essentially give us our grading space, but since we have changed the definition from \cite{Muller}, we need to prove this and any subsequent results for our altered setting. We will do so here in a combinatorial way without needing the background theory in \cite{Muller}. 

\begin{dfn}
Let $\alpha$  be an arc in an ideal triangulation of $(\bf{S},\bf{M})$ and consider the marked points (or common marked point) $m$ and $m'$ attached to each end of $\alpha$. We will call $m$ and $m'$ the \emph{endpoint(s)} of $\alpha$.
We will call functions ${s,t : \bf{A}^\circ (\bf{S},\bf{M}) \rightarrow \bf{M}}$ such that $\{ s(\alpha), t(\alpha) \} = \{ m,m'\}$ for all $\alpha \in \bf{A}^\circ (\bf{S},\bf{M})$ \emph{endpoint maps} of $ \bf{A}^\circ (\bf{S},\bf{M})$. Given $ \bf{A}^\circ (\bf{S},\bf{M})$, unless otherwise stated we will tacitly fix two endpoint maps and refer to them as $s$ and $t$.
\end{dfn}

\begin{dfn}
\label{dfn:surface_grading_space_E}
Let $(\bf{S},\bf{M})$ be a marked surface without punctures with an ideal triangulation $T$. Then
\[
\mathsf{E}:= \left\{  f: {\bf M} \rightarrow \mathbb{Q}  \mid f(s(\beta))+f(t(\beta))=0 \text{ for all boundary arcs } \beta \in T \right\}.
\]
We will call a function $f \in \mathsf{E}$ a \emph{valuation function} and $(\bf{S},\bf{M})$ along with $f$ (or a tuple of valuation functions $(f_1, \dots, f_n)$)  a \emph{valued marked surface}. Although we allow $f$ to take arbitrary rational values (so that $\mathsf{E}$ is a genuine vector space), in practice we will choose values such that the degree of any arc, which we will define presently, is an integer (c.f.~Remark/Definition \ref{rmkdfn:gradings_space_Z_vs_Q}).

For an arc $\alpha$ in a valued marked surface with valuation function $f$, we define the function $\deg_f: \bf{A}^\circ (\bf{S},\bf{M}) \rightarrow \Z$ by $\deg_f(\alpha) :=  f(s(\alpha))+f(t(\alpha))$. We will call $\deg_f(\alpha)$ the \emph{degree} of $\alpha$.
\end{dfn}

We will now turn to proving that a valuation function gives rise to a grading on the cluster algebra associated to a marked surface. For a seed to be graded, the exchange relations for mutation in each direction need to be homogeneous. We will show that the corresponding exchange relations will be homogeneous for any arc in an ideal triangulation of a valued marked surface. First, we need to introduce some new combinatorial objects that let us translate the notion of balanced vertices of a quiver to arcs of a triangulation.

\begin{dfn}
\label{dfn:configuration_of_alpha}
Let $\alpha$ be a (non-boundary) arc in a marked surface $(\bf{S},\bf{M})$ without punctures with triangulation $T$. Then $\alpha$ is an edge of two ideal triangles. We define the \emph{configuration} of $\alpha$, $\cfg{\alpha}$ to be the subset of $T$ consisting of these two triangles. A \emph{standard configuration} is one in which the two triangles share only one edge and whose arcs' endpoints make up four distinct points. (In other words, a standard configuration is a square with a diagonal, up to homotopy.) 
A \emph{starred configuration} is one that looks like a standard configuration except that certain marked points, those \emph{starred} by the symbol $*$ or $\star$, are identified with any others starred by the same symbol. We do not allow two or more adjacent boundary arcs to have all their endpoints starred by the same symbol.

If $(\bf{S},\bf{M})$ has a valuation function $f$, we will say $\cfg{\alpha}$ is \emph{balanced} if 
\begin{equation} \label{eqn:balanced_configuration}
 \sum_{\alpha_\text{cw}} \deg_f(\alpha_\text{cw})  = \sum_{\alpha_\text{ac}} \deg_f(\alpha_\text{ac}),
\end{equation} 
 where $\alpha_\text{cw}$ runs over all (non-boundary) arcs that follow $\alpha$ in the clockwise direction inside one of the triangles (adding two summands if this occurs in both triangles) and similarly for $\alpha_\text{ac}$ in the counter-clockwise direction. (Empty sums are defined to be 0.) This may be extended in the obvious way if we have a tuple of valuation functions. Equivalently, when the quiver is superimposed on the configuration, $\cfg{\alpha}$ is balanced if the corresponding vertex of the quiver is balanced---c.f.\ Remark \ref{rmk:quiver_of_triangulated_surface}.

When representing $\cfg{\alpha}$ as a diagram, we will use dotted lines to mean that an edge is a boundary arc. There may be other arcs contained in the region of $\bf{S}$ that falls inside the outer boundary formed by the two triangles; in this case, we may write $T^*$ to indicate that some arbitrary set of arcs is present (although they are not part of the configuration) to make the diagram more clear and avoid depicting something that appears not to be part of an ideal triangulation. 
\end{dfn}

Note that the outer boundary formed by the configuration of an arc may consist of between one and four arcs.

\begin{ex}
\label{ex:configuration}
Consider the valued marked surface
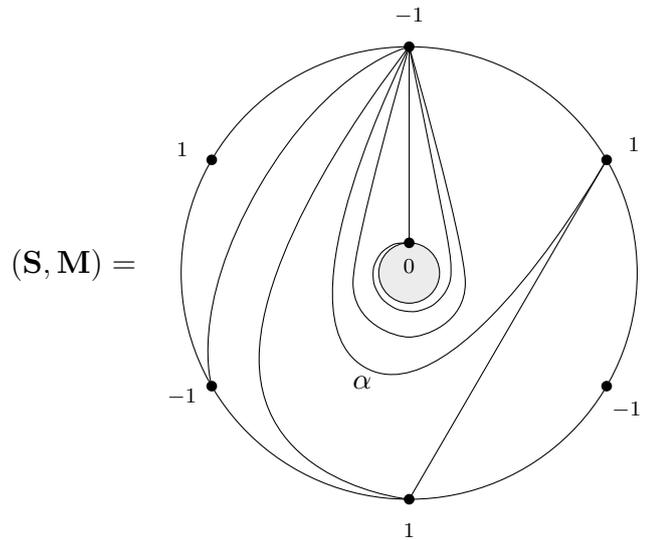
\begin{center} $(\bf{S},\bf{M})=$
\Image{ \begin{tikzpicture}
  \coordinate (O) at (0,0);  
  \def\radius{3cm}
  \def\midradius{ {(3/5)*\radius} }  
  \def\outermidradius{ {(3.8/5)*\radius} }      
  \def\innermidradius{ {(2.3/5)*\radius} }      
  \def\innerradius{  {(0.4/3)*\radius} }  
  \centerarc[](O)(0:360:\radius)  
  \shadedraw[inner color=gray!15,outer color=gray!15, draw=white] (O) circle[radius=\innerradius];    
  \centerarc[](O)(0:360:\innerradius)
  \path (O) ++(180+270:\radius) coordinate (B);
  \path (O) ++(180+90:{(0.85/3)*\radius} ) coordinate (GC);  
  \path (O) ++(180+170:{(0.75/3)*\radius} ) coordinate (GL);  
  \path (O) ++(180+10:{(0.75/3)*\radius} ) coordinate (GR);      
  \path (O) ++(180+90:\radius) coordinate (OU);  
  \path (O) ++(180+150:\radius) coordinate (OL1); \path (O) ++(180+150:\outermidradius) coordinate (GOL1);   
  \path (O) ++(180+210:\radius) coordinate (OL2); \path (O) ++(180+210:\outermidradius) coordinate (GOL2);   
  \path (O) ++(180+30:\radius) coordinate (OR1); \path (O) ++(180+30:\outermidradius) coordinate (GOR1);   
  \path (O) ++(150:\outermidradius) coordinate (GIR2);     
  \path (O) ++(180+330:\radius) coordinate (OR2); \path (O) ++(180+330:\outermidradius) coordinate (GOR2);       
  \path (O) ++(180+270:\innerradius) coordinate (IB); 
  \path (O) ++(180+30:\innerradius) coordinate (IR1);  \path (O) ++(180+65:\innermidradius) coordinate (GIR1);  
  \path (O) ++(50: {(0.75/3)*\radius}   ) coordinate (E0);  
  \path (O) ++(360: {1.1*(0.5/3)*\radius}   ) coordinate (E1);      
 \path (O) ++({333}: {1.1*(0.49/3)*\radius}) coordinate (E2);      
  \path (O) ++({306}:  {1.1*(0.48/3)*\radius}) coordinate (E3);      
  \path (O) ++({279}:  {1.1*(0.47/3)*\radius}) coordinate (E4);      
  \path (O) ++({252}:  {1.1*(0.46/3)*\radius}) coordinate (E5);      
  \path (O) ++({225}:  {1.1*(0.45/3)*\radius}) coordinate (E6);      
  \path (O) ++({198}:  {1.1*(0.44/3)*\radius}) coordinate (E7);      
  \path (O) ++({171}:  {1.1*(0.43/3)*\radius}) coordinate (E8);      
  \path (O) ++({144}:  {1.08*(0.42/3)*\radius}) coordinate (E9);      
  \path (O) ++({117}:  {1.05*(0.41/3)*\radius}) coordinate (E10);                        
   \fill[black] (OU) circle[radius=2pt] ++(180+90:1em) node {\scriptsize $ 1$};
   \fill[black] (OL1) circle[radius=2pt] ++(310:1em) node {\scriptsize $ -1$};   
   \fill[black] (OL2) circle[radius=2pt] ++(30:1em) node {\scriptsize $ 1$};
   \fill[black] (OR1) circle[radius=2pt] ++(200:1em) node {\scriptsize $ -1$};   
   \fill[black] (OR2) circle[radius=2pt] ++(160:1em) node {\scriptsize $1 $};
   \fill[black] (IB) circle[radius=2pt] ++(270:0.75em) node {\scriptsize $0$};      
   \fill[black] (B) circle[radius=2pt] ++(180+270:1em) node {\scriptsize $-1 $};   
   \fill[black] (GIR1)  ++(180+80:0.5em) node {\small $\alpha$};              
   \draw [black, xshift=4cm] plot [smooth, tension=0.5] coordinates { (B) (GL) (GC) (GR) (B)}; 
   \draw [black, xshift=4cm] plot [smooth, tension=0.5] coordinates { (B) (E0) (E1) (E2) (E3) (E4) (E5) (E6) (E7) (E8) (E9) (E10) (IB)}; 
   \draw [black, xshift=4cm] plot [smooth, tension=1] coordinates { (IB) (B) };  
   \draw [black, xshift=4cm] plot [smooth, tension=1.7] coordinates { (OL2) (GIR1) (B) };        
   \draw [black, xshift=4cm] plot [smooth, tension=1] coordinates { (OL2) (OU) };      
   \draw [black, xshift=4cm] plot [smooth, tension=1.3] coordinates { (OR1) (GIR2) (B) };    
   \draw [black, xshift=4cm] plot [smooth, tension=1] coordinates { (OU) (GOR1) (B) };                 
\end{tikzpicture}}.
\end{center}
The shaded area represents a region that is excluded from the surface---a convention we will adopt throughout. So $(\bf{S},\bf{M})$ is an annulus (a class of surfaces we will consider in the next section) with six points on the outer boundary and one point on the inner boundary. The configuration of $\alpha$ is
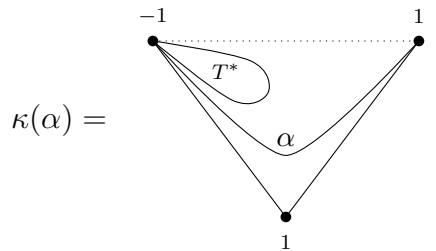
\begin{center}
$\cfg{\alpha}=$
\Image{ \begin{tikzpicture}
  \def\length{3.5cm}
  \coordinate (L) at (0,0);
  \coordinate (R) at (\length,0);
  \coordinate (B) at ( {(0.5)*\length}, {-(2/3)*\length} );
  \coordinate (Bup) at ( {0.5*\length}, {-(1.3/3)*\length} );  
  \coordinate (A1) at ( {(1.1/3)*\length}, {-(1.5/20)*\length} );  
  \coordinate (A2) at ( {(1.3/3)*\length}, {-(3.9/20)*\length} );
  \coordinate (A3) at  ( {(0.9/3)*\length}, {-(4.5/20)*\length} );
  \fill[black] (L) circle[radius=2pt] ++(90:0.8em) node {\scriptsize $ -1$};
  \fill[black] (R) circle[radius=2pt] ++(90:0.8em) node {\scriptsize $1 $};
  \fill[black] (B) circle[radius=2pt] ++(270:0.8em) node {\scriptsize $1 $};     
  \fill[black] (Bup) ++(90:0.5em) node {\small $\alpha $};       
  \fill[black] (A2)  ++(153:1.5em) node {\scriptsize $ T^*$};
     \draw [black, xshift=4cm] plot [smooth, tension=0] coordinates { (L) (B) (R)};  
     \draw [ dotted, black, xshift=4cm] plot [smooth, tension=0] coordinates { (L)  (R)};  
     \draw [black, xshift=4cm] plot [smooth, tension=0.6] coordinates { (L) (A1) (A2) (A3) (L)};  
     \draw [black, xshift=4cm] plot [smooth, tension=0.4] coordinates { (L) (Bup) (R) }; 
\end{tikzpicture}},
\end{center}
which is balanced since $\sum\limits_{\alpha_\text{cw}} \deg_f(\alpha_\text{cw})  = (-1-1) +(1+1)$ while $ \sum\limits_{\alpha_\text{ac}} \deg_f(\alpha_\text{ac}) = (-1+1)$.
\end{ex}

\begin{dfn}
In the setting of Definition \ref{dfn:configuration_of_alpha}, suppose $\cfg{\alpha}$ is such that the two triangles share exactly one edge. Then if  $\cfg{\alpha}$  is not a standard configuration we may use it to form a new starred configuration as follows. For each loop in $\cfg{\alpha}$ that encloses a region that is not the interior of one of the two triangles (i.e.\ a region marked $T^*$), we ``cut" the loop open at its basepoint. More precisely, each loop of the form
\Image{ \begin{tikzpicture}
  \def\rotation{90} 
  \def\radius{1.5cm}
  \def\innerradius{ {(4/5)*\radius} }  
  \coordinate (O) at (0,0);
  \path (O) ++(\rotation+15:\innerradius) coordinate (A1);
  \path (O) ++(\rotation+0:\radius) coordinate (A2);
  \path (O) ++(\rotation+345:\innerradius) coordinate (A3); 
  \fill[black] (O) circle[radius=2pt] ++(\rotation+90:0.8em) node {\scriptsize $ $};
  \fill[black] (A2)  ++(\rotation+180:1.5em) node {\scriptsize $ T^*$};
     \draw [black, xshift=4cm] plot [smooth, tension=0.6] coordinates { (O) (A1) (A2) (A3) (O)};  
\end{tikzpicture}}
 is replaced by 
 \Image{ \begin{tikzpicture}
  \def\rotation{0} 
  \def\radius{1.5cm}
  \coordinate (O) at (0,0);
  \path (O) ++(\rotation:\radius) coordinate (E);
  \fill[black] (O)  ++(\rotation+90:0.8em) node {\scriptsize $ $};
  \node at (O) {$*$};
  \fill[black] (E)  ++(\rotation+90:0.8em) node {\scriptsize $ $};
  \fill[black] (E)  ++(\rotation+270:1em) node {\scriptsize $ $};  
  \node at (E) {$*$};
  \draw [black, xshift=4cm,shorten <=0.12cm,shorten >=0.2cm] plot [smooth, tension=0.6] coordinates { (O) (E)}; 
\end{tikzpicture}}
and any other arcs that were attached to the base point of the loop are now attached to one of the starred points, corresponding to the side of the arc they were originally on. It is important to maintain the distinction between loops with different base points in the starred configuration we obtain, so we will star the corresponding pairs of points with different symbols when cutting loops based at different points. This is a valid starred configuration: the only way we could obtain two or more adjacent boundary arcs with endpoints starred with the same symbol is if the original configuration had multiple boundary loops with the same basepoint, but there is no such configuration. We will call the starred configuration obtained by this process the \emph{opened configuration} of $\alpha$.
\end{dfn}

\begin{ex} \label{ex:opened_configurations}
The configuration of $\alpha$ in Example \ref{ex:configuration} is not a standard configuration. The opened configuration is
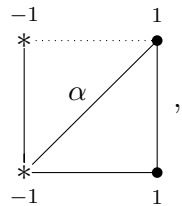
\begin{center}
\Image{ \begin{tikzpicture}
  \def\length{1.75cm}
  \coordinate (TL) at (0,0);
  \coordinate (TR) at (\length,0);
  \coordinate (BL) at ( {0}, {-\length} );
  \coordinate (BR) at ( {\length}, {-\length} );  
  \fill[black] (TL)  ++(90:0.8em) node {\scriptsize $ -1$};
  \node at (TL) {$*$};  
  \fill[black] (TR) circle[radius=2pt] ++(90:0.8em) node {\scriptsize $1 $};
  \fill[black] (BL) ++(270:0.8em) node {\scriptsize $-1 $}; 
  \node [inner sep=0pt] at (BL) {$*$};      
  \fill[black] (BR) circle[radius=2pt] ++(270:0.8em) node {\scriptsize $1 $};       
  \fill[black] (TL) ++(315:2.4em) node {\small $\alpha$};		     
     \draw [dotted, black, xshift=4cm, shorten <=0.15cm,shorten >=0cm] plot [smooth, tension=0] coordinates {(TL) (TR)}; 
     \draw [black, xshift=4cm] plot [smooth, tension=0.5] coordinates { (TR) (BR)};  
     \draw [black, xshift=4cm,shorten <=0.1cm,shorten >=0cm] plot [smooth, tension=0] coordinates {(BL) (BR)}; 
     \draw [black, xshift=4cm, shorten <=0.1cm,shorten >=0.25cm] plot [smooth, tension=0.5] coordinates { (TL) (BL) };
     \draw [black, xshift=4cm,shorten <=0.15cm,shorten >=0cm] plot [smooth, tension=0.5] coordinates { (BL) (TR) };    
\end{tikzpicture}},
\end{center}
which is again balanced.

For a second example, let $(\bf{S},\bf{M})$ be the marked surface
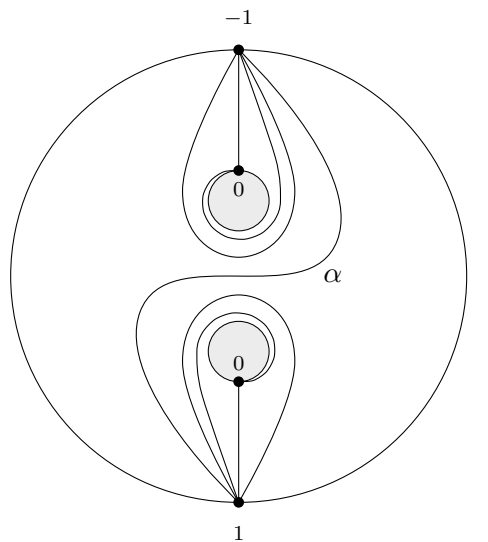
\begin{center}
\Image{ \begin{tikzpicture}
  \coordinate (O) at (0,0);  
  \def\radius{3cm}
  \def\midradius{ {(3/5)*\radius} }  
  \def\outermidradius{ {(3.8/5)*\radius} }      
  \def\innermidradius{ {(2.3/5)*\radius} }      
  \def\innerradius{  {(0.4/3)*\radius} }  
  \def\IBradius{  {(1/3)*\radius} }  
  \centerarc[](O)(0:360:\radius)  
   \path (O) ++(90:\IBradius) coordinate (upperB);  
   \path (O) ++(270:\IBradius) coordinate (lowerB);     
  \shadedraw[inner color=gray!15,outer color=gray!15, draw=white] (upperB) circle[radius=\innerradius];    
  \shadedraw[inner color=gray!15,outer color=gray!15, draw=white] (lowerB) circle[radius=\innerradius];      
  \centerarc[](upperB)(0:360:\innerradius)
  \centerarc[](lowerB)(0:360:\innerradius)  
  \path (O) ++(180+270:\radius) coordinate (B);
  \path (O) ++(200:\innermidradius) coordinate (alphaGL);  
  \path (O) ++(20:\innermidradius) coordinate (alphaGR);    
  \path (upperB) ++(180+90:{(0.75/3)*\radius} ) coordinate (GC);  
  \path (upperB) ++(180+190:{(0.75/3)*\radius} ) coordinate (GL);  
  \path (upperB) ++(180-10:{(0.75/3)*\radius} ) coordinate (GR);     
  \path (lowerB) ++(90:{(0.75/3)*\radius} ) coordinate (LGC);  
  \path (lowerB) ++(190:{(0.75/3)*\radius} ) coordinate (LGL);  
  \path (lowerB) ++(-10:{(0.75/3)*\radius} ) coordinate (LGR);     
  \path (O) ++(180+90:\radius) coordinate (OU);    
  \path (upperB) ++(180+270:\innerradius) coordinate (IB); 
  \path (lowerB) ++(270:\innerradius) coordinate (LIB); 
  \path (upperB) ++(50: {(0.75/3)*\radius}   ) coordinate (E0);  
  \path (upperB) ++(360: {1.1*(0.5/3)*\radius}   ) coordinate (E1);      
  \path (upperB) ++({333}: {1.1*(0.49/3)*\radius}) coordinate (E2);      
  \path (upperB) ++({306}:  {1.1*(0.48/3)*\radius}) coordinate (E3);      
  \path (upperB) ++({279}:  {1.1*(0.47/3)*\radius}) coordinate (E4);      
  \path (upperB) ++({252}:  {1.1*(0.46/3)*\radius}) coordinate (E5);      
  \path (upperB) ++({225}:  {1.1*(0.45/3)*\radius}) coordinate (E6);      
  \path (upperB) ++({198}:  {1.1*(0.44/3)*\radius}) coordinate (E7);      
  \path (upperB) ++({171}:  {1.1*(0.43/3)*\radius}) coordinate (E8);      
  \path (upperB) ++({144}:  {1.08*(0.42/3)*\radius}) coordinate (E9);      
  \path (upperB) ++({117}:  {1.05*(0.41/3)*\radius}) coordinate (E10);      
  \path (lowerB) ++(180+50: {(0.75/3)*\radius}   ) coordinate (F0);  
  \path (lowerB) ++(180+360: {1.1*(0.5/3)*\radius}   ) coordinate (F1);      
  \path (lowerB) ++({180+333}: {1.1*(0.49/3)*\radius}) coordinate (F2);      
  \path (lowerB) ++({180+306}:  {1.1*(0.48/3)*\radius}) coordinate (F3);      
  \path (lowerB) ++({180+279}:  {1.1*(0.47/3)*\radius}) coordinate (F4);      
  \path (lowerB) ++({180+252}:  {1.1*(0.46/3)*\radius}) coordinate (F5);      
  \path (lowerB) ++({180+225}:  {1.1*(0.45/3)*\radius}) coordinate (F6);      
  \path (lowerB) ++({180+198}:  {1.1*(0.44/3)*\radius}) coordinate (F7);      
  \path (lowerB) ++({180+171}:  {1.1*(0.43/3)*\radius}) coordinate (F8);      
  \path (lowerB) ++({180+144}:  {1.08*(0.42/3)*\radius}) coordinate (F9);      
  \path (lowerB) ++({180+117}:  {1.05*(0.41/3)*\radius}) coordinate (F10);                          
   \fill[black] (OU) circle[radius=2pt] ++(180+90:1em) node {\scriptsize $ 1$};  
   \fill[black] (B) circle[radius=2pt] ++(180+270:1em) node {\scriptsize $-1 $};   
   \fill[black] (IB) circle[radius=2pt] ++(270:0.6em) node {\scriptsize $0 $};      
   \fill[black] (LIB) circle[radius=2pt] ++(180+270:0.6em) node {\scriptsize $0 $};     
   \fill[black] (O)  ++(0:3em) node {\small $\alpha$};              
   \draw [black, xshift=4cm] plot [smooth, tension=0.8] coordinates { (B) (GL) (GC) (GR) (B)}; 
   \draw [black, xshift=4cm] plot [smooth, tension=0.8] coordinates { (OU) (LGL) (LGC) (LGR) (OU)}; 
   \draw [black, xshift=4cm] plot [smooth, tension=0.5] coordinates { (B) (E0) (E1) (E2) (E3) (E4) (E5) (E6) (E7) (E8) (E9) (E10) (IB)}; 
   \draw [black, xshift=4cm] plot [smooth, tension=0.5] coordinates { (OU) (F0) (F1) (F2) (F3) (F4) (F5) (F6) (F7) (F8) (F9) (F10) (LIB)}; 
   \draw [black, xshift=4cm] plot [smooth, tension=1] coordinates { (IB) (B) };  
   \draw [black, xshift=4cm] plot [smooth, tension=1] coordinates { (LIB) (OU) };  
   \draw [black, xshift=4cm] plot [smooth, tension=1] coordinates { (OU) (alphaGL) (alphaGR) (B)};  
\end{tikzpicture}} \phantom{.}.
\end{center}
Then 
\begin{center} $\cfg{\alpha}=$
\Image
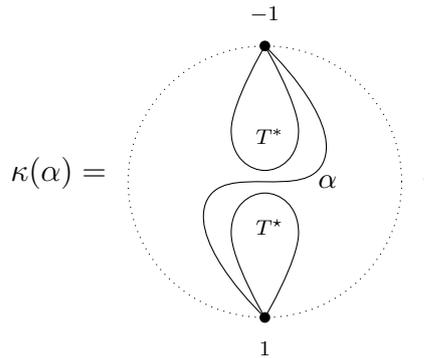
{ \begin{tikzpicture}
  \coordinate (O) at (0,0);  
  \def\radius{1.8cm}
  \def\midradius{ {(3/5)*\radius} }  
  \def\outermidradius{ {(3.8/5)*\radius} }      
  \def\innermidradius{ {(2.3/5)*\radius} }      
  \def\innerradius{  {(0.4/3)*\radius} }  
  \def\IBradius{  {(1/3)*\radius} }  
  \centerarc[dotted](O)(0:360:\radius)  
   \path (O) ++(90:\IBradius) coordinate (upperB);  
   \path (O) ++(270:\IBradius) coordinate (lowerB);     
  \path (O) ++(180+270:\radius) coordinate (B);
  \path (O) ++(200:\innermidradius) coordinate (alphaGL);  
  \path (O) ++(20:\innermidradius) coordinate (alphaGR);    
  \path (upperB) ++(180+90:{(0.75/3)*\radius} ) coordinate (GC);  
  \path (upperB) ++(180+190:{(0.75/3)*\radius} ) coordinate (GL);  
  \path (upperB) ++(180-10:{(0.75/3)*\radius} ) coordinate (GR);     
  \path (lowerB) ++(90:{(0.75/3)*\radius} ) coordinate (LGC);  
  \path (lowerB) ++(190:{(0.75/3)*\radius} ) coordinate (LGL);  
  \path (lowerB) ++(-10:{(0.75/3)*\radius} ) coordinate (LGR);     
  \path (O) ++(180+90:\radius) coordinate (OU);    
  \path (upperB) ++(180+270:\innerradius) coordinate (IB); 
  \path (lowerB) ++(270:\innerradius) coordinate (LIB); 
   \fill[black] (OU) circle[radius=2pt] ++(180+90:1em) node {\scriptsize $ 1$};  
   \fill[black] (B) circle[radius=2pt] ++(180+270:1em) node {\scriptsize $-1 $};   
   \fill[black] (IB)  ++(285:0.6em) node {\scriptsize $T^* $};      
   \fill[black] (LIB) ++(75:0.6em) node {\scriptsize $T^{\star} $};     
   \fill[black] (O)  ++(0:2em) node {\small $\alpha$};              
   \draw [black, xshift=4cm] plot [smooth, tension=0.8] coordinates { (B) (GL) (GC) (GR) (B)}; 
   \draw [black, xshift=4cm] plot [smooth, tension=0.8] coordinates { (OU) (LGL) (LGC) (LGR) (OU)}; 
   \draw [black, xshift=4cm] plot [smooth, tension=1] coordinates { (OU) (alphaGL) (alphaGR) (B)};  
\end{tikzpicture}} \phantom{.},
\end{center}
which is balanced as $\sum\limits_{\alpha_\text{cw}} \deg_f(\alpha_\text{cw})  = 0$ while $ \sum\limits_{\alpha_\text{ac}} \deg_f(\alpha_\text{ac}) = (-1-1)+(1+1)$.
We obtain the opened configuration in two steps. First we open the loop whose endpoints are on the upper marked point of value $-1$ (i.e.~the loop bounding $T^*$), which gives
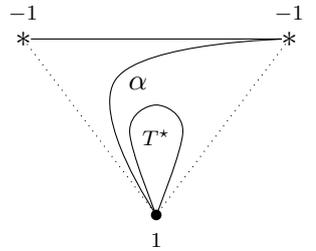
\begin{center}
\Image{ \begin{tikzpicture}
  \def\length{3.5cm}
  \coordinate (L) at (0,0);
  \coordinate (R) at (\length,0);
  \coordinate (B) at ( {(0.5)*\length}, {-(2/3)*\length} );
  \coordinate (Bup) at ( {0.5*\length}, {-(1.3/3)*\length} );  
  \coordinate (A1) at ( {(1.2/3)*\length}, {-(7/20)*\length} );  
  \coordinate (A2) at ( {(1.5/3)*\length}, {-(5/20)*\length} );  
  \coordinate (A3) at  ( {(1.8/3)*\length}, {-(7/20)*\length} );
  \coordinate (G1) at ( {0.35*\length}, {-(3/20)*\length} );    
  \fill[black] (L) ++(90:0.8em) node {\scriptsize $ -1$};
  \node [inner sep=0pt] at (L) {$*$};        
  \fill[black] (R)  ++(90:0.8em) node {\scriptsize $-1 $};
  \node [inner sep=0pt] at (R) {$*$};   
  \fill[black] (B) circle[radius=2pt] ++(270:0.8em) node {\scriptsize $1 $};     
  \fill[black] (Bup) ++(90:0.5em) node {\scriptsize  $ T^{\star}$};       
  \fill[black] (A2) ++(130:.9em) node {\small $\alpha$};
     \draw [ dotted,black, xshift=4cm] plot [smooth, tension=0] coordinates { (L) (B) (R)};  
     \draw [ black, xshift=4cm,shorten <=0.1cm,shorten >=0.15cm] plot [smooth, tension=0] coordinates { (L)  (R)};  
     \draw [ black, xshift=4cm] plot [smooth, tension=0.6] coordinates { (B) (A1) (A2) (A3) (B)};  
     \draw [black, xshift=4cm,shorten >=0.1cm] plot [smooth, tension=0.75] coordinates { (B) (G1) (R)}; 
\end{tikzpicture}}.
\end{center}
Then we open the remaining loop, which gives us the opened configuration
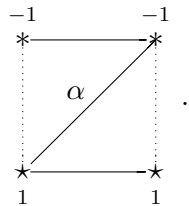
\begin{center}
\Image{ \begin{tikzpicture}
  \def\length{1.75cm}
  \coordinate (TL) at (0,0);
  \coordinate (TR) at (\length,0);
  \coordinate (BL) at ( {0}, {-\length} );
  \coordinate (BR) at ( {\length}, {-\length} );  
  \fill[black] (TL)  ++(90:0.8em) node {\scriptsize $ -1$};
  \node at (TL) {$*$};  
  \fill[black] (TR) ++(90:0.8em) node {\scriptsize $-1 $};
  \node at (TR) {$*$};    
  \fill[black] (BL) ++(270:0.8em) node {\scriptsize $1 $}; 
  \node [inner sep=0pt] at (BL) {$\star$};      
  \fill[black] (BR) ++(270:0.8em) node {\scriptsize $1 $};      
  \node at (BR) {$\star$};    
  \fill[black] (TL) ++(315:2.4em) node {\small $\alpha$};		      
     \draw [ black, xshift=4cm, shorten <=0.1cm,shorten >=0.22cm] plot [smooth, tension=0] coordinates {(TL) (TR)}; 
     \draw [dotted,black, xshift=4cm,shorten <=0.1cm,shorten >=0.25cm] plot [smooth, tension=0.5] coordinates { (TR) (BR)};  
     \draw [black, xshift=4cm,shorten <=0.1cm,shorten >=.25cm] plot [smooth, tension=0] coordinates {(BL) (BR)}; 
     \draw [dotted,black, xshift=4cm, shorten <=0.1cm,shorten >=0.25cm] plot [smooth, tension=0.5] coordinates { (TL) (BL) };
     \draw [black, xshift=4cm,shorten <=0.15cm,shorten >=0.15cm] plot [smooth, tension=0.5] coordinates { (BL) (TR) };    
\end{tikzpicture}}.
\end{center}
This is again balanced. 
\end{ex}

In the examples above, a balanced configuration gave rise to a balanced open configuration. This is not a coincidence, as we will see shortly.

\begin{rmk}
It is possible to obtain any configuration by starting with a starred configuration and ``gluing" along matching symbols such that we create internal loops. For example, for the latter opened configuration in Example  \ref{ex:opened_configurations}, if we reverse the two steps we carried out to obtain it, we get back the original configuration of $\alpha$. The other possible configurations are obtained by starring a standard configuration in all the possible ways (including starring more than two points with the same symbols) and then gluing.
\end{rmk}

Below are some simple results that allow us to restrict our attention to a small number of configurations in a systematic way.

\begin{lem}
\label{lem:opened_config_balanced_iff_balanced}
Let $\alpha$ be an arc in an ideal triangulation of a valued marked surface such that the triangles of $\cfg{\alpha}$ share just one edge.
Then the opened configuration of $\alpha$ is balanced if and only if $\cfg{\alpha}$ is balanced. 
\end{lem}
\begin{proof}
This is clear by inspection.
\end{proof}

\begin{lem}
\label{lem:standard_config_balanced_iff_stared_balanced}
Suppose a given standard configuration is balanced for all valuation functions. Then so is any starred configuration obtained from it by replacing some marked points with starred points (that is, by starring some of the points).
\end{lem}

\begin{proof}
Replacing points of the standard configuration by starred points simply restricts the possible values of these points (points starred with the same symbol must now share a common value), but the configuration is already balanced for all values endowed by a valuation function, and in particular for these restricted values.
\end{proof}

\begin{lem}
Let $\alpha$ be an arc in some ideal triangulation $T$ of a valued marked surface and let $\alpha'$ be the unique new arc in the triangulation $T'$ obtained by flipping $\alpha$ in $T$. We have that $\cfg{\alpha}$ is balanced if and only if $\cfg{\alpha'}$ is balanced.
\end{lem}

\begin{proof}
The arcs that follow $\alpha'$ in the clockwise direction are exactly those that follow $\alpha$ in the counter-clockwise direction, and the arcs that follow $\alpha'$ in the counter-clockwise direction are those that follow $\alpha$ clockwise. Thus exchanging $\alpha$ for $\alpha'$ simply has the effect of swapping the left hand side of Equation \ref{eqn:balanced_configuration} with the right hand side, and vice versa.
\end{proof}

\begin{prop}
\label{prop:configurations_are_balanced}
Let $\alpha$ be an arc in an ideal triangulation of a valued marked surface. Then $\cfg{\alpha}$ is balanced.
\end{prop}

\begin{proof}
In any configuration, the two triangles share either one or two edges. Thus there are two cases to consider.

If the triangles share two edges, there are only two possibilities for  $\cfg{\alpha}$:
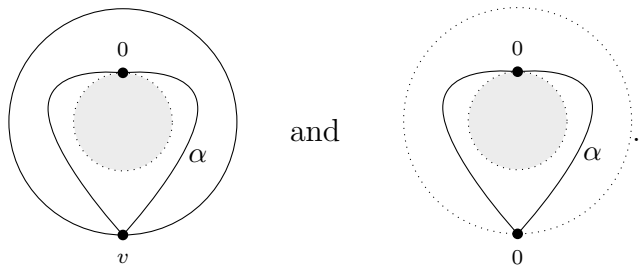
\begin{center}
\Image{ \begin{tikzpicture}
  \coordinate (O) at (0,0);  
  \def\radius{1.5cm}    
  \def\innermidradius{ {(2/3)*\radius} }      
  \def\innerradius{  {(1.3/3)*\radius} }  
  \centerarc[](O)(0:360:\radius)  
  \shadedraw[inner color=gray!15,outer color=gray!15, draw=white] (O) circle[radius=\innerradius];    
  \centerarc[dotted](O)(0:360:\innerradius)
  \path (O) ++(270:\radius) coordinate (B);
  \path (O) ++(170:{\innermidradius} ) coordinate (GL);  
  \path (O) ++(10:{\innermidradius} ) coordinate (GR);      
  \path (O) ++(270:\innerradius) coordinate (IB);
  \path (O) ++(90:\innerradius) coordinate (IU);
   \fill[black] (IU) circle[radius=2pt] ++(90:0.75em) node {\scriptsize $0$};      
   \fill[black] (B) circle[radius=2pt] ++(270:0.75em) node {\scriptsize $v$};   
   \fill[black] (B)  ++(47:3.5em) node {\small $\alpha$};              
   \draw [black, xshift=4cm] plot [smooth, tension=1] coordinates { (B) (GL) (IU)};                         
   \draw [black, xshift=4cm] plot [smooth, tension=1] coordinates { (B) (GR) (IU)};        
\end{tikzpicture}}
\ \ \ \ and \ \ \ \
\Image{ \begin{tikzpicture}
  \coordinate (O) at (0,0);  
  \def\radius{1.5cm}    
  \def\innermidradius{ {(2/3)*\radius} }      
  \def\innerradius{  {(1.3/3)*\radius} }  
  \centerarc[dotted](O)(0:360:\radius)  
  \shadedraw[inner color=gray!15,outer color=gray!15, draw=white] (O) circle[radius=\innerradius];    
  \centerarc[dotted](O)(0:360:\innerradius)
  \path (O) ++(270:\radius) coordinate (B);
  \path (O) ++(170:{\innermidradius} ) coordinate (GL);  
  \path (O) ++(10:{\innermidradius} ) coordinate (GR);      
  \path (O) ++(270:\innerradius) coordinate (IB);
  \path (O) ++(90:\innerradius) coordinate (IU);
   \fill[black] (IU) circle[radius=2pt] ++(90:0.75em) node {\scriptsize $0$};      
   \fill[black] (B) circle[radius=2pt] ++(270:0.75em) node {\scriptsize $0$};   
   \fill[black] (B)  ++(47:3.5em) node {\small $\alpha$};              
   \draw [black, xshift=4cm] plot [smooth, tension=1] coordinates { (B) (GL) (IU)};                         
   \draw [black, xshift=4cm] plot [smooth, tension=1] coordinates { (B) (GR) (IU)};        
\end{tikzpicture}}.
\end{center}
Note the double arrow (between $\alpha$ and the other interior arc) which would be present in the superimposed quiver. Both of these configurations are balanced: in the first we have 
\[ \sum\limits_{\alpha_\text{cw}} \deg_f(\alpha_\text{cw})  = 2(v+0) = (v+v) = \sum\limits_{\alpha_\text{ac}} \deg_f(\alpha_\text{ac}),\]
 while in the second both sums are automatically zero (as both marked points must have zero value under any valuation function, being on boundaries with an odd number of points).
 
Suppose the two triangles share one edge. Then by Lemma \ref{lem:opened_config_balanced_iff_balanced} and Lemma \ref{lem:standard_config_balanced_iff_stared_balanced} we may assume without loss of generality that they form a standard configuration. We therefore just need to list all these configurations and check they are balanced for an arbitrary valuation function. We write down and check this list in Table \ref{tab:configurations_balanced}.
{
\setlength{\tabcolsep}{3pt} 
\renewcommand{\arraystretch}{1.3} %
\begin{table}[H]
\centering
\begin{tabular}{ | l  l |  l l | }
\hline
\multicolumn{4}{|c|}{ Configuration \ \ \ \ \  $ \scalemath{0.85}{\begin{matrix}  \text{LHS: }\sum\limits_{\alpha_\text{cw}} \deg_f(\alpha_\text{cw})  \\  \text{RHS: }  \sum\limits_{\alpha_\text{ac}} \deg_f(\alpha_\text{ac}) \end{matrix} }$ }
\\ \hline
		 \Image{ \begin{tikzpicture}
		  \def\length{1.7cm}
		  \coordinate (TL) at (0,0);
		  \coordinate (TR) at (\length,0);
		  \coordinate (BL) at ( {0}, {-\length} );
		  \coordinate (BR) at ( {\length}, {-\length} );  
		  \fill[black] (TL) circle[radius=2pt] ++(90:0.8em) node {\scriptsize $a_1$};
		  \fill[black] (TR) circle[radius=2pt] ++(90:0.8em) node {\scriptsize $a_4$};
		  \fill[black] (BL) circle[radius=2pt] ++(270:0.8em) node {\scriptsize $a_2$};     
		  \fill[black] (BR) circle[radius=2pt] ++(270:0.8em) node {\scriptsize $a_3$}; 
		  \fill[black] (TL) ++(315:2.4em) node {\small $\alpha$};		          
	      \draw [black] plot [smooth, tension=0] coordinates {(TL) (TR)}; 
	      \draw [black] plot [smooth, tension=0.5] coordinates { (TR) (BR)};  
	      \draw [black] plot [smooth, tension=0] coordinates {(BL) (BR)}; 
	      \draw [black] plot [smooth, tension=0.5] coordinates { (TL) (BL) };
	      \draw [black] plot [smooth, tension=0.5] coordinates { (BL) (TR) };    
		\end{tikzpicture}}  &  $ \scalemath{0.75}{\begin{matrix*}[l] \text{LHS: } (a_1+a_2)+(a_3+a_4) \\ \text{RHS: }  (a_1+a_4)+(a_2+a_3)  \end{matrix*} }$  &
		 \Image{ \begin{tikzpicture}
		  \def\length{1.7cm}
		  \coordinate (TL) at (0,0);
		  \coordinate (TR) at (\length,0);
		  \coordinate (BL) at ( {0}, {-\length} );
		  \coordinate (BR) at ( {\length}, {-\length} );  
		  \fill[black] (TL) circle[radius=2pt] ++(90:0.8em) node {\scriptsize $a_1$};
		  \fill[black] (TR) circle[radius=2pt] ++(90:0.8em) node {\scriptsize $-a_1$};
		  \fill[black] (BL) circle[radius=2pt] ++(270:0.8em) node {\scriptsize $a_2$};     
		  \fill[black] (BR) circle[radius=2pt] ++(270:0.8em) node {\scriptsize $a_3$}; 
		  \fill[black] (TL) ++(315:2.4em) node {\small $\alpha$};       
	      \draw [dotted, black] plot [smooth, tension=0] coordinates {(TL) (TR)}; 
	      \draw [black] plot [smooth, tension=0.5] coordinates { (TR) (BR)};  
	      \draw [black] plot [smooth, tension=0] coordinates {(BL) (BR)}; 
	      \draw [black] plot [smooth, tension=0.5] coordinates { (TL) (BL) };
	      \draw [black] plot [smooth, tension=0.5] coordinates { (BL) (TR) };    
		\end{tikzpicture}}  &  $ \scalemath{0.75}{\begin{matrix*}[l] \text{LHS: } (a_1+a_2)+(a_3-a_1) \\ \text{RHS: }  (a_2+a_3) \end{matrix*} }$  
\\ \hline
		 \Image{ \begin{tikzpicture}
		  \def\length{1.7cm}
		  \coordinate (TL) at (0,0);
		  \coordinate (TR) at (\length,0);
		  \coordinate (BL) at ( {0}, {-\length} );
		  \coordinate (BR) at ( {\length}, {-\length} );  
		  \fill[black] (TL) circle[radius=2pt] ++(90:0.8em) node {\scriptsize $a_1$};
		  \fill[black] (TR) circle[radius=2pt] ++(90:0.8em) node {\scriptsize $-a_1$};
		  \fill[black] (BL) circle[radius=2pt] ++(270:0.8em) node {\scriptsize $a_2$};     
		  \fill[black] (BR) circle[radius=2pt] ++(270:0.8em) node {\scriptsize $-a_2$};     
		  \fill[black] (TL) ++(315:2.4em) node {\small $\alpha$};		      
	      \draw [dotted, black] plot [smooth, tension=0] coordinates {(TL) (TR)}; 
	      \draw [black] plot [smooth, tension=0.5] coordinates { (TR) (BR)};  
	      \draw [dotted,black] plot [smooth, tension=0] coordinates {(BL) (BR)}; 
	      \draw [black] plot [smooth, tension=0.5] coordinates { (TL) (BL) };
	      \draw [black] plot [smooth, tension=0.5] coordinates { (BL) (TR) };    
		\end{tikzpicture}}  &  $ \scalemath{0.75}{\begin{matrix*}[l] \text{LHS: } (a_1+a_2)+(-a_1-a_2) \\ \text{RHS: }  0 \end{matrix*} }$  &
		 \Image{ \begin{tikzpicture}
		  \def\length{1.7cm}
		  \coordinate (TL) at (0,0);
		  \coordinate (TR) at (\length,0);
		  \coordinate (BL) at ( {0}, {-\length} );
		  \coordinate (BR) at ( {\length}, {-\length} );  
		  \fill[black] (TL) circle[radius=2pt] ++(90:0.8em) node {\scriptsize $a_1$};
		  \fill[black] (TR) circle[radius=2pt] ++(90:0.8em) node {\scriptsize $-a_1$};
		  \fill[black] (BL) circle[radius=2pt] ++(270:0.8em) node {\scriptsize $-a_1$};     
		  \fill[black] (BR) circle[radius=2pt] ++(270:0.8em) node {\scriptsize $a_2$};     
		  \fill[black] (TL) ++(315:2.4em) node {\small $\alpha$};		      
	      \draw [dotted, black] plot [smooth, tension=0] coordinates {(TL) (TR)}; 
	      \draw [black] plot [smooth, tension=0.5] coordinates { (TR) (BR)};  
	      \draw [black] plot [smooth, tension=0] coordinates {(BL) (BR)}; 
	      \draw [dotted, black] plot [smooth, tension=0.5] coordinates { (TL) (BL) };
	      \draw [black] plot [smooth, tension=0.5] coordinates { (BL) (TR) };    
		\end{tikzpicture}}  &  $ \scalemath{0.75}{\begin{matrix*}[l] \text{LHS: } (a_2-a_1) \\ \text{RHS: }  (a_2-a_1) \end{matrix*} }$  
\\ \hline
		 \Image{ \begin{tikzpicture}
		  \def\length{1.7cm}
		  \coordinate (TL) at (0,0);
		  \coordinate (TR) at (\length,0);
		  \coordinate (BL) at ( {0}, {-\length} );
		  \coordinate (BR) at ( {\length}, {-\length} );  
		  \fill[black] (TL) circle[radius=2pt] ++(90:0.8em) node {\scriptsize $a_1$};
		  \fill[black] (TR) circle[radius=2pt] ++(90:0.8em) node {\scriptsize $-a_1$};
		  \fill[black] (BL) circle[radius=2pt] ++(270:0.8em) node {\scriptsize $-a_1$};     
		  \fill[black] (BR) circle[radius=2pt] ++(270:0.8em) node {\scriptsize $a_1$};     
		  \fill[black] (TL) ++(315:2.4em) node {\small $\alpha$};		      
	      \draw [dotted, black] plot [smooth, tension=0] coordinates {(TL) (TR)}; 
	      \draw [black] plot [smooth, tension=0.5] coordinates { (TR) (BR)};  
	      \draw [dotted,black] plot [smooth, tension=0] coordinates {(BL) (BR)}; 
	      \draw [dotted, black] plot [smooth, tension=0.5] coordinates { (TL) (BL) };
	      \draw [black] plot [smooth, tension=0.5] coordinates { (BL) (TR) };    
		\end{tikzpicture}}  &  $ \scalemath{0.75}{\begin{matrix*}[l] \text{LHS: } (a_1-a_1) \\ \text{RHS: }  0 \end{matrix*} }$  &
		 \Image{ \begin{tikzpicture}
		  \def\length{1.7cm}
		  \coordinate (TL) at (0,0);
		  \coordinate (TR) at (\length,0);
		  \coordinate (BL) at ( {0}, {-\length} );
		  \coordinate (BR) at ( {\length}, {-\length} );  
		  \fill[black] (TL) circle[radius=2pt] ++(90:0.8em) node {\scriptsize $a_1$};
		  \fill[black] (TR) circle[radius=2pt] ++(90:0.8em) node {\scriptsize $-a_1$};
		  \fill[black] (BL) circle[radius=2pt] ++(270:0.8em) node {\scriptsize $-a_1$};     
		  \fill[black] (BR) circle[radius=2pt] ++(270:0.8em) node {\scriptsize $a_1$};     
		  \fill[black] (TL) ++(315:2.4em) node {\small $\alpha$};		      
	      \draw [dotted, black] plot [smooth, tension=0] coordinates {(TL) (TR)}; 
	      \draw [dotted, black] plot [smooth, tension=0.5] coordinates { (TR) (BR)};  
	      \draw [dotted, black] plot [smooth, tension=0] coordinates {(BL) (BR)}; 
	      \draw [dotted, black] plot [smooth, tension=0.5] coordinates { (TL) (BL) };
	      \draw [black] plot [smooth, tension=0.5] coordinates { (BL) (TR) };    
		\end{tikzpicture}}  &  $ \scalemath{0.75}{\begin{matrix*}[l] \text{LHS: } 0 \\ \text{RHS: } 0 \end{matrix*} }$  
\\ \hline
\end{tabular}
\caption{ the standard configurations.}
\label{tab:configurations_balanced}
\end{table}
}
\end{proof}

\begin{prop}
\label{prop:valued_surfce_gives_graded_algebra}
Let $(\bf{S}, \bf{M})$ be a valued marked surface with valuation function $f$ and initial ideal triangulation $T$. Say the arcs of $T$ are $\alpha_1, \dots, \alpha_n$.  \sloppy Then $\A (\bf{S}, \bf{M})$ is the graded cluster algebra generated by the initial degree seed $\left(  (\deg_f(\alpha_1), \dots, \deg_f(\alpha_n) ) ,B(T) \right)$.
\end{prop}

\begin{proof}
Since arcs are in bijection with cluster variables by \ref{thm:FST_main_theorem}, we just need to justify that the seed above is a degree seed. By  Proposition \ref{prop:configurations_are_balanced}, $\cfg{\alpha_i}$ is balanced for each $i$, but this is exactly the condition required in Equation \ref{eqn:grading_condition}.
\end{proof}

\begin{cor}
Given a tuple of valuation functions, Proposition $\ref{prop:valued_surfce_gives_graded_algebra}$ extends to multi-gradings in the obvious way.
\end{cor}

We have shown that $\deg_f$ takes a triangulation and gives a valid degree cluster for any valuation function $f$. Now we will show that, given a degree cluster, there is a valuation function $f$ such that $\deg_f$ gives this degree cluster.

Suppose $(\bf{S},\bf{M})$ is a marked surface with ideal triangulation $T$, where $T$ has arcs $\alpha_1 \dots, \alpha_n$, and that $B=B(T)$. Assume $\ker(B^T)$ is $k$-dimensional, that is we have a $k$-dimensional grading. Then by Theorem \ref{thm:FST_even_boundaries_gives_rank} we know there are $k$ boundary components of $(\bf{S},\bf{M})$ that have an even number of marked points.

Then $\mathsf{E} = \mathsf{E}(\bf{S},\bf{M})$ is also $k$-dimensional as follows. Denote the even boundary components of $(\bf{S},\bf{M})$ by $\Sigma_1, \dots, \Sigma_k$. For the component $\Sigma_i$, denote the set of its marked points by $\{ m_1^i, \dots m^i_{r_i}\}$ and assume these are labelled such that $m^i_j$ and $m^i_{j+1}$ are adjacent for all $j$ (working $\bmod~r_i$ so that $m_1^i$ and $m^i_{r_i}$ are adjacent). A basis of $\mathsf{E}$ is given by $(f_1, \dots, f_k)$ where 
\begin{equation}
f_i=\delta_{m^i_1} - \delta_{m^i_2} + \delta_{m^i_3 } - \cdots + (-1)^{r_i + 1} \delta_{m^i_{r_i}}.
\end{equation}
(Here, we are defining $\delta_*$ by $\delta_{m}(p)={\begin{cases} 1 & {\text{if }} m=p, \\ 0 & {\text{otherwise,}}\end{cases}}$ for marked points $m$ and $p$.)

\begin{dfn}
For the triangulation $T$, define $\theta_T: \mathsf{E} \rightarrow \ker(B(T)^T)$ by
\begin{equation}
\theta_T(f) = \left( \deg_f(\alpha_1), \dots, \deg_f(\alpha_n) \right).
\end{equation}
\end{dfn}

By Proposition \ref{prop:valued_surfce_gives_graded_algebra}, we know that the image of $f$ is indeed in $\ker(B^T)$, and it is easy to show that $\theta_T$ is linear. We wish to show it is an isomorphism. The crucial step is establishing that it is enough to show the result for any single triangulation in the mutation class.

\begin{lem}
\label{lem:theta_inj_T0_implies_inj_T}
For any mutation direction $j$, if $\theta_T$ is injective then so is $\theta_{T'}$, where $T' = \mu_j (T)$.
\end{lem}

\begin{proof}
Suppose $\theta_T$ is injective. Then $f \neq 0  \implies \theta_T(f) \neq 0$. Thus for any $f \neq 0$ there is some $\alpha_i \in T$ such that $\deg_f(\alpha_i) \neq 0$, that is, such that $f(s(\alpha_i)) + f(t(\alpha_i)) \neq 0$. Fix $f \neq 0 $ and say  $\deg_f(\alpha_i) \neq 0$.  Consider $T':= \mu_j(T)$. If $i \neq j$ or if there is another arc with non-zero degree then clearly we still have $\deg_f(\alpha) \neq 0$ for some $\alpha$ in $T'$. So assume $i=j$ and that $\alpha_i$ is the only arc in $T$ such that  $\deg_f(\alpha_i) \neq 0$.
Consider the configuration of $\alpha_i$ in $T$: we have
\begin{equation}
\cfg{\alpha_i} =
		 \Image{ \begin{tikzpicture}
		  \def\length{1.7cm}
		  \coordinate (TL) at (0,0);
		  \coordinate (TR) at (\length,0);
		  \coordinate (BL) at ( {0}, {-\length} );
		  \coordinate (BR) at ( {\length}, {-\length} );  
		  \fill[black] (TL) circle[radius=2pt] ++(90:0.8em) node {\scriptsize $a$};
		  \fill[black] (TR) circle[radius=2pt] ++(90:0.8em) node {\scriptsize $b$};
		  \fill[black] (BL) circle[radius=2pt] ++(270:0.8em) node {\scriptsize $d$};     
		  \fill[black] (BR) circle[radius=2pt] ++(270:0.8em) node {\scriptsize $c$}; 
		  \fill[black] (TL) ++(315:2.2em) node {\small $\alpha_i$};		          
	      \draw [black] plot [smooth, tension=0] coordinates {(TL) (TR)}; 
	      \draw [black] plot [smooth, tension=0.5] coordinates { (TR) (BR)};  
	      \draw [black] plot [smooth, tension=0] coordinates {(BL) (BR)}; 
	      \draw [black] plot [smooth, tension=0.5] coordinates { (TL) (BL) };
	      \draw [black] plot [smooth, tension=0.5] coordinates { (BL) (TR) };    
		\end{tikzpicture}} ,
\end{equation}
where $a,b,c,d$ are the values of the corresponding marked points. (We do not assume these marked points are distinct and, apart from $\alpha_i$, the arcs may be boundary arcs.) But since $\alpha_i$ is the only arc in $T$ with non-zero degree, we must have $a= -b$, $a=-d$, $d=-c$ and $b=-c$. So the configuration is
\begin{equation}
\cfg{\alpha_i} =
		 \Image{ \begin{tikzpicture}
		  \def\length{1.7cm}
		  \coordinate (TL) at (0,0);
		  \coordinate (TR) at (\length,0);
		  \coordinate (BL) at ( {0}, {-\length} );
		  \coordinate (BR) at ( {\length}, {-\length} );  
		  \fill[black] (TL) circle[radius=2pt] ++(90:0.8em) node {\scriptsize $a$};
		  \fill[black] (TR) circle[radius=2pt] ++(90:0.8em) node {\scriptsize $-a$};
		  \fill[black] (BL) circle[radius=2pt] ++(270:0.8em) node {\scriptsize $-a$};     
		  \fill[black] (BR) circle[radius=2pt] ++(270:0.8em) node {\scriptsize $a$}; 
		  \fill[black] (TL) ++(315:2.2em) node {\small $\alpha_i$};		          
	      \draw [black] plot [smooth, tension=0] coordinates {(TL) (TR)}; 
	      \draw [black] plot [smooth, tension=0.5] coordinates { (TR) (BR)};  
	      \draw [black] plot [smooth, tension=0] coordinates {(BL) (BR)}; 
	      \draw [black] plot [smooth, tension=0.5] coordinates { (TL) (BL) };
	      \draw [black] plot [smooth, tension=0.5] coordinates { (BL) (TR) };    
		\end{tikzpicture}},
\end{equation}
and we must have $\deg(\alpha_i) = -2a \neq 0$. Now the configuration of $\alpha_i'$ in $\mu_j(T')$ is 
\begin{equation}
\cfg{\alpha_i'} =
		 \Image{ \begin{tikzpicture}
		  \def\length{1.7cm}
		  \coordinate (TL) at (0,0);
		  \coordinate (TR) at (\length,0);
		  \coordinate (BL) at ( {0}, {-\length} );
		  \coordinate (BR) at ( {\length}, {-\length} );  
		  \fill[black] (TL) circle[radius=2pt] ++(90:0.8em) node {\scriptsize $a$};
		  \fill[black] (TR) circle[radius=2pt] ++(90:0.8em) node {\scriptsize $-a$};
		  \fill[black] (BL) circle[radius=2pt] ++(270:0.8em) node {\scriptsize $-a$};     
		  \fill[black] (BR) circle[radius=2pt] ++(270:0.8em) node {\scriptsize $a$}; 
		  \fill[black] (TL) ++(298:2.8em) node {\small $\alpha_i'$};		          
	      \draw [black] plot [smooth, tension=0] coordinates {(TL) (TR)}; 
	      \draw [black] plot [smooth, tension=0.5] coordinates { (TR) (BR)};  
	      \draw [black] plot [smooth, tension=0] coordinates {(BL) (BR)}; 
	      \draw [black] plot [smooth, tension=0.5] coordinates { (TL) (BL) };
	      \draw [black] plot [smooth, tension=0.5] coordinates { (BR) (TL) };    
		\end{tikzpicture}},
\end{equation}
which means $\deg(\alpha_i') = 2a \neq 0$. So $\theta_{T'} (f) \neq 0$. Thus $\theta_{T'}$ is also injective.
\end{proof}

This immediately implies the step we wished to establish:

\begin{cor}
\label{cor:inj_on_any_triang_enough}
Let $T_0$ be a fixed ideal triangulation. If $\theta_{T_0}$ is injective then so is $\theta_{T}$ for any ideal triangulation $T$.
\end{cor}

We can now deal uniformly with all but a special class of cases.

\begin{lem}
\label{lem:theta_f_0_implies_f_0_excluding_exaclty_2_pts}
Suppose $(\bf{S},\bf{M})$ is a marked surface. There exists a triangulation $T$ such that the following is true: if a valuation function $f$ is such that $\theta_T(f)=0$, then $f$ is zero on every boundary component that does not have exactly two marked points.
\end{lem}

\begin{proof}
If a boundary component $\Sigma$ has an odd number of marked points, then $f$ must always be zero on $\Sigma$. Let $T$ be a triangulation satisfying the following: for each even boundary component $\Sigma_i$ with more than two marked points, there is an arc $\beta_i$ between the marked points $m^i_1$ and $m^i_3$ (which must both exist as $\Sigma_i$ has at least $4$ marked points). We should justify that $\beta_1,\dots,\beta_k$ are compatible; this is clear since each arc's endpoints are confined to a unique boundary component. It is also easy to see there does indeed exist a triangulation containing these arcs: if $T'=\{\beta_1,\dots,\beta_k \}$ isn't an ideal triangulation already, then by definition there is another arc $\beta \notin T'$ compatible with $T'$.

Let $f \in \mathsf{E}$. Note that since $m^i_1$ and $m^i_3$ are on an even boundary, and by how we have labelled the marked points, we have $f(m^i_1) = f(m^i_3)$. Also, for each $i$ we have $\deg_f(\beta_i) = f(s(\beta_i)) + f(t(\beta_i)) =f(m^i_1) + f(m^i_3)$. Now suppose $\theta_T(f) = 0$. Then in particular $f(m^i_1) + f(m^i_3) = 0$, so we now have both $f(m^i_1) =- f(m^i_3)$ and $f(m^i_1) = f(m^i_3)$ and thus  $f(m^i_1) = 0$. But this implies all the other marked points $m$ on $\Sigma_i$ must also satisfy $f(m)=0$. Thus $f$ is $0$ on each even boundary component of $(\bf{S},\bf{M})$ not containing exactly two marked points.
\end{proof}

Given the above, all that is left to address are boundary components with exactly two marked points, in the surfaces that have them.

\begin{thm}
\label{thm:theta_bijection}
Let $(\bf{S},\bf{M})$ be a marked surface. Then $\theta_T$ is an isomorphism for any ideal triangulation $T$.
\end{thm}

\begin{proof}
We will show $\theta_T$ is injective, and therefore bijective, for any triangulation $T$. If $(\bf{S},\bf{M})$ has no boundary components with exactly two marked points, then, by Corollary \ref{cor:inj_on_any_triang_enough}, we can assume $T$ is a triangulation satisfying Lemma \ref{lem:theta_f_0_implies_f_0_excluding_exaclty_2_pts}. Then we have $\theta_T(f)=0$ implies $f=0$ and we are done. So assume $(\bf{S},\bf{M})$ has boundary components with exactly two marked points. 

Suppose the valuation function $f$ is such that $\theta_T(f)=0$.  We will show that there exists a triangulation containing an arc on $\Sigma$ that goes between marked points of the same value, which implies these points have zero value. This triangulation will be an extension of $\{\beta_1,\dots,\beta_k \}$ from before. Indeed, by Corollary \ref{cor:inj_on_any_triang_enough}, we may assume $T$ is such that for each even boundary component $\Sigma_i$ there is an arc $\beta_i$ between the marked points $m^i_1$ and $m^i_3$. Then, by Lemma  \ref{lem:theta_f_0_implies_f_0_excluding_exaclty_2_pts}, $f$ is zero on any even boundary components with more than two points. As always, $f$ is automatically zero on any odd boundary components.

Now, suppose $\Sigma$ is a boundary component with exactly two marked points. Then $\Sigma$ must fall into one of three cases:
\begin{enumerate}[(i)]
\item $\Sigma$ is the only boundary component (i.e.\ $(\bf{S},\bf{M})$ is the digon), 
\item $\Sigma$ is contained inside the region of $(\bf{S},\bf{M})$ bounded by another component, or,
\item  $\Sigma$ encloses an interior boundary component.
\end{enumerate}
Let $m$ be a fixed marked point on $\Sigma$. Case (i) is trivial as there can be no non-boundary arcs, so the arc complex is empty.  In case (ii), we use Corollary \ref{cor:inj_on_any_triang_enough} again: assume $T$ contains an arc $\alpha$ that is a loop from $m$ to itself encircling $\Sigma$ itself. Such a triangulation exists since $\alpha$ is clearly compatible with the arcs $\beta_1,\dots,\beta_k$, which we have previously assumed are in $T$.
Then $0=f(s(\alpha)) + f(t(\alpha))=f(m)+f(m)$, so $f(m)=0$. Thus, $f$ is zero on $\Sigma$. Case (iii) is topologically equivalent to (ii), so we are done.

This shows we can assume that $T$ is such that $f$ must be zero on every boundary component, that is, $f=0$. Thus, $\theta_T$ is injective and hence bijective. Therefore $\theta_T$ is an isomorphism.
\end{proof}

\begin{rmk}
A potentially worthwhile line of future research is investigating whether the grading we have defined can be extended to surfaces with punctures. A possible approach is the following. By \cite[Lemma 2.13]{FST} such a surface has a triangulation without self folded triangles, and we can choose a valuation function for such a triangulation---problems only occur when the mutated triangulation produces self-folded triangles. If we pass to a tagged triangulation so that self-folded triangles are replaced with notched arcs, a suitable way to extend the grading may be to set the degree of a notched arc to the value of the point at the plain end minus the value of the point at the notched end.
\end{rmk}

\section{Annulus with $n+m$ marked points, $m$ odd} \label{section:Ann_n_m_m_odd}

The annulus $\wt{A}(n,m)$ (where $n,m \geq 1$) has $n$ marked points on the outer boundary component and $m$ points on the inner boundary component. For the grading on the corresponding cluster algebra, there are two different scenarios depending on whether both, or just one, of $n$ and $m$ are even (we assume at least one is even, otherwise no nontrivial gradings can occur by Theorem \ref{thm:FST_even_boundaries_gives_rank}). In the former case, we will assume $n \geq m$, and in the latter, that $n$ is even. We may do so by the following. 

\begin{rmk}
\label{rmk:annulus_triangulation_is_cylinder_triangulation}
The graded cluster algebra associated to the annulus $\wt{A}(n,m)$ is the same as that for $\wt{A}(m,n)$. To see this, simply note that any initial triangulation of  $\wt{A}(n,m)$ is the same as a triangulation of an open cylinder with $n$ marked points on the top circle and $m$ on the base circle. Turning the cylinder upside down then gives an initial triangulation for $\wt{A}(m,n)$, and clearly this does not change the corresponding initial quiver or grading.
\end{rmk}

In this section we consider the case for when only $n$ is even, and we will show that the grading is of mixed type. For our initial set-up we take the valued triangulation
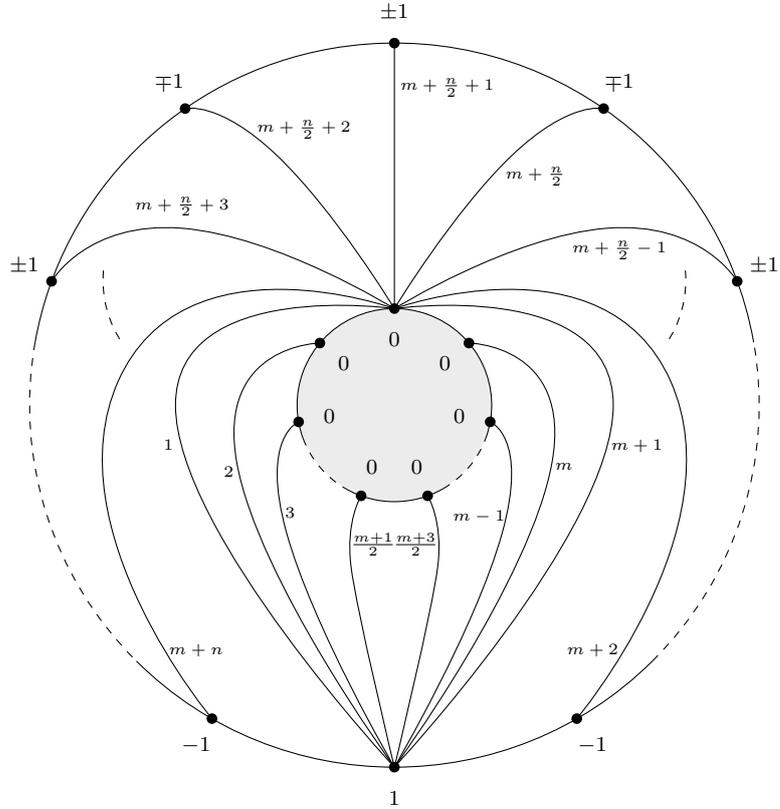
\begin{center}
{ 
\begin{tikzpicture}
  \coordinate (O) at (0,0);
  
  \def\radius{4.8cm}
  \def\midradius{ {(3/5)*\radius} }  
  \def\outermidradius{ {(3.8/5)*\radius} }      
  \def\innermidradius{ {(2.3/5)*\radius} }      
  \def\innerradius{  {(0.8/3)*\radius} }
  
  \centerarc[](O)(10:170:\radius)
  \centerarc[dashed](O)(170:225:\radius)
  \centerarc[](O)(225:315:\radius)
  \centerarc[dashed](O)(315:370:\radius)
  
  \shadedraw[inner color=gray!15,outer color=gray!15, draw=white] (O) circle[radius=\innerradius];  
  
  \centerarc[](O)(-20:200:\innerradius)
  \centerarc[dashed](O)(200:240:\innerradius)
  \centerarc[](O)(240:300:\innerradius)
  \centerarc[dashed](O)(300:340:\innerradius)
  
  \path (O) ++(155:{(4.4/5)*\radius}) coordinate (LA);
  \draw[dashed] (LA) arc (175:215:1.5);
  \path (O) ++(180-155:{(4.4/5)*\radius}) coordinate (RA);
  \draw[dashed] (RA) arc (180-175:180-215:1.5);  

  \path (O) ++(270:\radius) coordinate (B);
  \path (O) ++(90:\radius) coordinate (OU);  
  \path (O) ++(125:\radius) coordinate (OL1); \path (O) ++(115:\outermidradius) coordinate (GOL1);   
  \path (O) ++(160:\radius) coordinate (OL2); \path (O) ++(140:\outermidradius) coordinate (GOL2);   
  \path (O) ++(240:\radius) coordinate (OL3); \path (O) ++(175:\outermidradius) coordinate (GOL3);
  \path (O) ++(180-125:\radius) coordinate (OR1); \path (O) ++(180-115:\outermidradius) coordinate (GOR1);   
  \path (O) ++(180-160:\radius) coordinate (OR2); \path (O) ++(180-140:\outermidradius) coordinate (GOR2);   
  \path (O) ++(180-240:\radius) coordinate (OR3); \path (O) ++(180-175:\outermidradius) coordinate (GOR3);      
  \path (O) ++(90:\innerradius) coordinate (IU);  \path (O) ++(180:\midradius) coordinate (GIUL); \path (O) ++(0:\midradius) coordinate (GIUR);
  \path (O) ++(140:\innerradius) coordinate (IL1);  \path (O) ++(200:\innermidradius) coordinate (GIL1);  
  \path (O) ++(190:\innerradius) coordinate (IL2);  \path (O) ++(230:\innermidradius) coordinate (GIL2);   
  \path (O) ++(250:\innerradius) coordinate (IL3);  \path (O) ++(255:\innermidradius) coordinate (GIL3);  
  \path (O) ++(180-140:\innerradius) coordinate (IR1);  \path (O) ++(180-200:\innermidradius) coordinate (GIR1);  
  \path (O) ++(180-190:\innerradius) coordinate (IR2);  \path (O) ++(180-230:\innermidradius) coordinate (GIR2);   
  \path (O) ++(180-250:\innerradius) coordinate (IR3);  \path (O) ++(180-255:\innermidradius) coordinate (GIR3);    

   \fill[black] (OU) circle[radius=2pt] ++(90:1em) node {\scriptsize $\pm 1$};
   \fill[black] (OL1) circle[radius=2pt] ++(120:1em) node {\scriptsize $\mp 1$};   
   \fill[black] (OL2) circle[radius=2pt] ++(150:1em) node {\scriptsize $\pm 1$};
   \fill[black] (OL3) circle[radius=2pt] ++(240:1em) node {\scriptsize $-1$};    
   \fill[black] (OR1) circle[radius=2pt] ++(180-120:1em) node {\scriptsize $\mp 1$};   
   \fill[black] (OR2) circle[radius=2pt] ++(180-150:1em) node {\scriptsize $\pm 1$};
   \fill[black] (OR3) circle[radius=2pt] ++(180-240:1em) node {\scriptsize $-1$};    
   \fill[black] (IU) circle[radius=2pt] ++(270:1em) node {\scriptsize $0$};
   \fill[black] (IL1) circle[radius=2pt] ++(320:1em) node {\scriptsize $0$};   
   \fill[black] (IL2) circle[radius=2pt] ++(10:1em) node {\scriptsize $0$};               
   \fill[black] (IL3) circle[radius=2pt] ++(70:1em) node {\scriptsize $0$};       
   \fill[black] (IR1) circle[radius=2pt] ++(180-320:1em) node {\scriptsize $0$};   
   \fill[black] (IR2) circle[radius=2pt] ++(180-10:1em) node {\scriptsize $0$};               
   \fill[black] (IR3) circle[radius=2pt] ++(180-70:1em) node {\scriptsize $0$};       
   \fill[black] (B) circle[radius=2pt] ++(270:1em) node {\scriptsize $1$};   
   \draw (GOR2) ++(300:0.8em) node[font=\tiny] {$m+\frac{n}{2}-1$};   
   \draw (GOR1) ++(320:1em) node[font=\tiny] {$m+\frac{n}{2}$};   
   \draw (OU) ++(320:2.2em) node[font=\tiny] {$m+\frac{n}{2}+1$};      
   \draw (GOL1) ++(45:1.2em) node[font=\tiny] {$m+\frac{n}{2}+2$};    
   \draw (GOL2) ++(90:0.7em) node[font=\tiny] {$m+\frac{n}{2}+3$};    
   \draw (GOL3) ++(286:9em) node[font=\tiny] {$m+n$};      
   \draw (GIUL) ++(260:1.3em) node[font=\tiny]  {$1$};
   \draw (GIL1) ++(230:0.4em) node[font=\tiny]  {$2$};
   \draw (GIL2) ++(80:0.65em) node[font=\tiny]  {$3$};  
   \draw (GIL3) ++(45:1em) node[font=\tiny]  {$\frac{m+1}{2}$};  
   \draw (GIR3) ++(135:1em) node[font=\tiny]  {$\frac{m+3}{2}$};   
   \draw (GIR2) ++(145:0.9em) node[font=\tiny]  {$m-1$};  
   \draw (GIR1) ++(340:0.4em) node[font=\tiny]  {$m$};  
   \draw (GIUR) ++(300:1.5em) node[font=\tiny]  {$m+1$};  
   \draw (GOR3) ++(254:9em) node[font=\tiny] {$m+2$};        
     
   \draw [black, xshift=4cm] plot [smooth, tension=1] coordinates { (IU) (OU)};               
   \draw [black, xshift=4cm] plot [smooth, tension=1] coordinates { (IU) (GIUL) (B)}; 
   \draw [black, xshift=4cm] plot [smooth, tension=1] coordinates { (IU) (GIUR) (B)};    
   \draw [black, xshift=4cm] plot [smooth, tension=1] coordinates { (IL1) (GIL1) (B)};        
   \draw [black, xshift=4cm] plot [smooth, tension=1] coordinates { (IL2) (GIL2) (B)};   
   \draw [black, xshift=4cm] plot [smooth, tension=0.6] coordinates { (IL3) (GIL3) (B)};   
   \draw [black, xshift=4cm] plot [smooth, tension=1] coordinates { (IR1) (GIR1) (B)};        
   \draw [black, xshift=4cm] plot [smooth, tension=1] coordinates { (IR2) (GIR2) (B)};          
   \draw [black, xshift=4cm] plot [smooth, tension=0.6] coordinates { (IR3) (GIR3) (B)};    
   \draw [black, xshift=4cm] plot [smooth, tension=1] coordinates { (IU) (GOL1) (OL1)};      
   \draw [black, xshift=4cm] plot [smooth, tension=1] coordinates { (IU) (GOL2) (OL2)};         
   \draw [black, xshift=4cm] plot [smooth, tension=1.3] coordinates { (IU) (GOL3) (OL3)};  
   \draw [black, xshift=4cm] plot [smooth, tension=1] coordinates { (IU) (GOR1) (OR1)};      
   \draw [black, xshift=4cm] plot [smooth, tension=1] coordinates { (IU) (GOR2) (OR2)};         
   \draw [black, xshift=4cm] plot [smooth, tension=1.3] coordinates { (IU) (GOR3) (OR3)};       
\end{tikzpicture}
}
\end{center}
which corresponds to the initial graded quiver
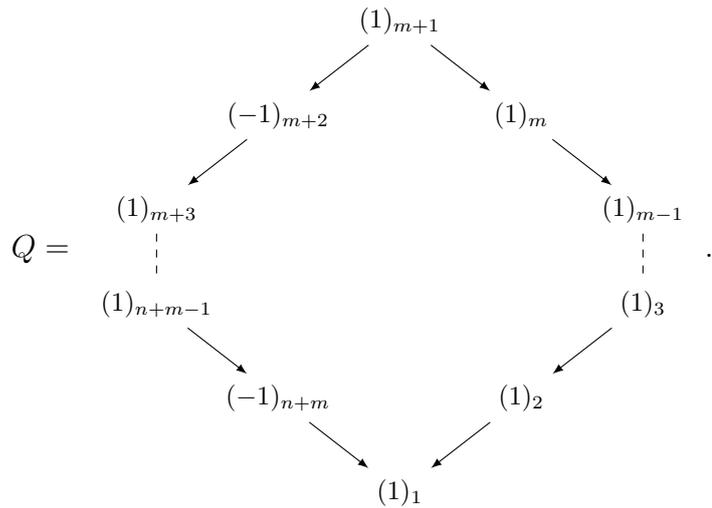
\begin{center}
$Q=$
\Image{
		\begin{tikzpicture}[align=center,node distance=1.25cm and 1.6cm, auto, font=\footnotesize, on grid]    
		  \node (Top) []{$(1)_{m+1}$};                                        			  
		  \node (TopBL) [below left = of Top]{$(-1)_{m+2}$};                                        			  
		  \node (TopBR) [below right = of Top]{$(1)_{m}$};                                        			  		  		  
		  \node (TopBLBL) [below left = of TopBL]{$(1)_{m+3}$};                                        			  		  		  		  
		  \node (TopBRBR) [below right = of TopBR]{$(1)_{m-1}$};
		  \node (LowerL) [below = of TopBLBL]{$(1)_{n+m-1 }$};		                                          			  		  		  		  		  
		  \node (LowerR) [below = of TopBRBR]{$(1)_{3} $};		  
		  \node (LowerLBR) [below right = of LowerL]{$(-1)_{n+m}$};		  		  
		  \node (LowerRBL) [below left = of LowerR]{$(1)_{2}$};			  
		  \node (Bottom) [below right = of LowerLBR]{$(1)_{1}$};		  
		  \draw[-latex] (Top) to node {} (TopBL);      
		  \draw[-latex] (Top) to node {} (TopBR); 		  
		  \draw[-latex] (TopBL) to node {} (TopBLBL);     
		  \draw[-latex] (TopBR) to node {} (TopBRBR);   		  		  
		  \draw[dashed, -] (TopBLBL) to node {} (LowerL);
		  \draw[dashed, -] (TopBRBR) to node {} (LowerR);		  
		  \draw[-latex] (LowerL) to node {} (LowerLBR);		 		  
		  \draw[-latex] (LowerR) to node {} (LowerRBL);		 		  		  
		  \draw[-latex] (LowerRBL) to node {} (Bottom);		 		  		  		  
		  \draw[-latex] (LowerLBR) to node {} (Bottom);		 		  		  		  		  
		\end{tikzpicture} 
}. 
\end{center}

As a tuple, the degree cluster is
$g_{n,m}:=(\overbrace{1, \dots, 1}^{\text{$m$ entries}}, \overbrace{1,-1,\dots,1,-1}^{\text{$n$ entries}} )$.

\begin{lem} \label{lem:Ann_n_m_odd_inf_vars_degs1}
$\A\big(g_{n,m},Q \big)$ has infinitely many variables in degrees $1$ and $-1$.
\end{lem}

\begin{proof}
To see that there are infinitely many variables in degree $-1$, say, fix any marked point $p$ of value $-1$ (which will be on the outer boundary) and any point $q$ on the inner boundary (which will have value $0$). For $r \geq 1$, let $\alpha(r)$ be the arc that starts at $p$, winds $r$ times around the inner boundary in a clockwise fashion (spiralling inwards as it does so), and ends at $q$. For example, in $\wt{A}(6,1)$, if we let $p$ be the point at the top of the outer boundary and $q$ the only point on the inner boundary, then $\alpha(1)$ is the arc as drawn below:
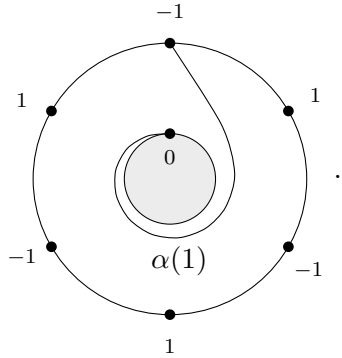
\begin{center}
\Image{ \begin{tikzpicture}
  \coordinate (O) at (0,0);  
  \def\radius{1.8cm}
  \def\midradius{ {(3/5)*\radius} }  
  \def\outermidradius{ {(3.8/5)*\radius} }      
  \def\innermidradius{ {(2.3/5)*\radius} }      
  \def\innerradius{  {(1/3)*\radius} }  
  \centerarc[](O)(0:360:\radius)  
  \shadedraw[inner color=gray!15,outer color=gray!15, draw=white] (O) circle[radius=\innerradius];    
  \centerarc[](O)(0:360:\innerradius)
  \path (O) ++(180+270:\radius) coordinate (B);
  \path (O) ++(180+90:{(0.85/3)*\radius} ) coordinate (GC);  
  \path (O) ++(180+170:{(0.75/3)*\radius} ) coordinate (GL);  
  \path (O) ++(180+10:{(0.75/3)*\radius} ) coordinate (GR);      
  \path (O) ++(180+90:\radius) coordinate (OU);  
  \path (O) ++(180+150:\radius) coordinate (OL1); \path (O) ++(180+150:\outermidradius) coordinate (GOL1);   
  \path (O) ++(180+210:\radius) coordinate (OL2); \path (O) ++(180+210:\outermidradius) coordinate (GOL2);   
  \path (O) ++(180+30:\radius) coordinate (OR1); \path (O) ++(180+30:\outermidradius) coordinate (GOR1);   
  \path (O) ++(150:\outermidradius) coordinate (GIR2);     
  \path (O) ++(180+330:\radius) coordinate (OR2); \path (O) ++(180+330:\outermidradius) coordinate (GOR2);       
  \path (O) ++(180+270:\innerradius) coordinate (IB); 
  \path (O) ++(180+30:\innerradius) coordinate (IR1);  \path (O) ++(180+65:\innermidradius) coordinate (GIR1);  
  \path (O) ++(47: {(1.65/3)*\radius}   ) coordinate (E0);  
  \path (O) ++(5: {2.85*(0.5/3)*\radius}   ) coordinate (E1);      
 \path (O) ++({333}: {2.8*(0.49/3)*\radius}) coordinate (E2);      
  \path (O) ++({306}:  {2.8*(0.48/3)*\radius}) coordinate (E3);      
  \path (O) ++({279}:  {2.8*(0.47/3)*\radius}) coordinate (E4);      
  \path (O) ++({252}:  {2.8*(0.46/3)*\radius}) coordinate (E5);      
  \path (O) ++({225}:  {2.8*(0.45/3)*\radius}) coordinate (E6);      
  \path (O) ++({198}:  {2.8*(0.44/3)*\radius}) coordinate (E7);      
  \path (O) ++({171}:  {2.8*(0.43/3)*\radius}) coordinate (E8);      
  \path (O) ++({144}:  {2.74*(0.42/3)*\radius}) coordinate (E9);      
  \path (O) ++({117}:  {2.6*(0.41/3)*\radius}) coordinate (E10);                        
   \fill[black] (OU) circle[radius=2pt] ++(180+90:1em) node {\scriptsize $ 1$};
   \fill[black] (OL1) circle[radius=2pt] ++(310:1em) node {\scriptsize $ -1$};   
   \fill[black] (OL2) circle[radius=2pt] ++(30:1em) node {\scriptsize $ 1$};
   \fill[black] (OR1) circle[radius=2pt] ++(200:1em) node {\scriptsize $ -1$};   
   \fill[black] (OR2) circle[radius=2pt] ++(160:1em) node {\scriptsize $1 $};
   \fill[black] (IB) circle[radius=2pt] ++(270:0.75em) node {\scriptsize $0$};      
   \fill[black] (B) circle[radius=2pt] ++(180+270:1em) node {\scriptsize $-1 $};   
   \fill[black] (GIR1)  ++(180+145:1.4em) node {\small $\alpha(1)$};              
   \draw [black, xshift=4cm] plot [smooth, tension=0.5] coordinates { (B) (E0) (E1) (E2) (E3) (E4) (E5) (E6) (E7) (E8) (E9) (E10) (IB)}; 
\end{tikzpicture}}.
\end{center}
Then, if $r_1 \neq r_1$, $\alpha(r_1)$ and $\alpha(r_2)$ are not homotopic, so the corresponding cluster variables are distinct by Theorem \ref{thm:FST_main_theorem}. But $\alpha(r)$ has degree $-1$ for any $r$, so $\{\alpha(r) \mid r \in \N \}$ is a set of infinitely many variables of degree $-1$. Degree $1$ may be dealt with in a similar way.
\end{proof}

\begin{lem} \label{lem:Ann_n_m_odd_finite_vars_degs02}
$\A\big(g_{n,m},Q \big)$ has only finitely many variables in degrees $0$ and $\pm 2$.
\end{lem}

\begin{proof}
While we will not write down exact values for upper bounds on the number of variables that can occur in these degrees, it is easy to see that these numbers must be finite. The only way to obtain infinitely many non-homotopic arcs in the annulus with $n+m$ marked points is to produce a sequence of arcs that loop increasingly many times around the inner boundary of the annulus (such as $\alpha(r)$ above). But any such arc must start on one boundary and end on the other and therefore must have degree $\pm 1$, since all points on the inner boundary have value 0 and all points on the outer boundary have value $\pm 1$. Therefore, there are only finitely many arcs of degree $0$ and $\pm 2$.
\end{proof}

\begin{cor}
 $\A\big(g_{n,m},Q \big)$ has infinitely many variables in degrees $\pm 1$ and finitely many in degrees $0$ and $\pm 2$, and these are the only degrees that occur. Thus the graded cluster algebra associated to a triangulation of the annulus with $n+m$ marked points, with $m$ odd, is of mixed type.
\end{cor}

\begin{rmk}
\label{rmk:Ann_n_m_parity_m_doesn't_matter}
It is possible to prove Lemma \ref{lem:Ann_n_m_odd_inf_vars_degs1} above using standard methods involving growing sequences of denominator vectors (though this takes considerably more effort than we needed above), but Lemma \ref{lem:Ann_n_m_odd_finite_vars_degs02}, which involves showing certain degrees have only finitely many variables, would be very difficult to prove without the theory in Section \ref{sec:surface_theory} (as is often the case when trying to establish there are finitely many variables in a given degree)
\end{rmk}

\begin{rmk}
None of the above relies on $m$ being odd. We could thus extend any of these results to the case when $m$ is even (provided we extend the grading $g_{n,m}$ appropriately), although doing so would not alone be enough to show what happens for the graded cluster algebra when $m$ is even, assuming we take a standard grading (recall Definition \ref{dfn:standard_grading}), which will be a $\Z^2$-grading. However, we will make use of this fact in the following section.
\end{rmk}

\section{Annulus with $n+m$ marked points, $m$ even}
Now let us assume both $n$ and $m$ are even. As there are now two boundary components with an even number of marked points, we get a 2-dimensional grading. We take the initial valued triangulation
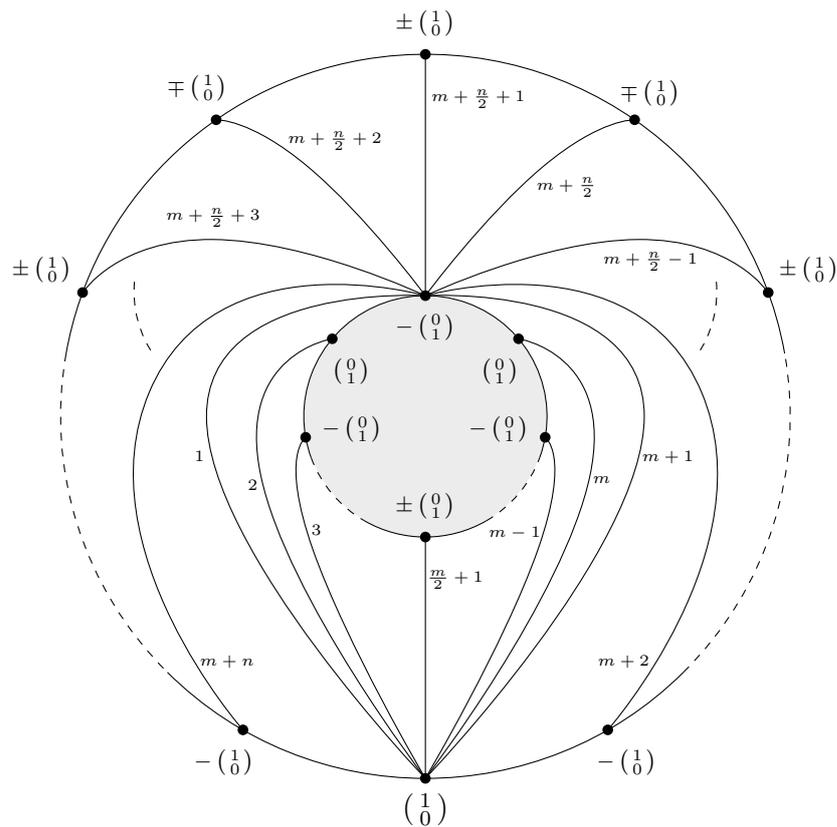
\begin{center}
{ 
\begin{tikzpicture}
  \coordinate (O) at (0,0);
  
  \def\radius{4.8cm}
  \def\midradius{ {(3/5)*\radius} }  
  \def\outermidradius{ {(3.8/5)*\radius} }      
  \def\innermidradius{ {(2.4/5)*\radius} }      
  \def\innerradius{  {(1/3)*\radius} }
  
  \centerarc[](O)(10:170:\radius)
  \centerarc[dashed](O)(170:225:\radius)
  \centerarc[](O)(225:315:\radius)
  \centerarc[dashed](O)(315:370:\radius)
  
  \shadedraw[inner color=gray!15,outer color=gray!15, draw=white] (O) circle[radius=\innerradius];  
  
  \centerarc[](O)(-20:200:\innerradius)
  \centerarc[dashed](O)(200:240:\innerradius)
  \centerarc[](O)(240:300:\innerradius)
  \centerarc[dashed](O)(300:340:\innerradius)
  
  \path (O) ++(155:{(4.4/5)*\radius}) coordinate (LA);
  \draw[dashed] (LA) arc (175:215:1.5);
  \path (O) ++(180-155:{(4.4/5)*\radius}) coordinate (RA);
  \draw[dashed] (RA) arc (180-175:180-215:1.5);  

  \path (O) ++(270:\radius) coordinate (B);
  \path (O) ++(270:\innerradius) coordinate (IB); \path (O) ++(270:\innermidradius) coordinate (GIB);
  \path (O) ++(90:\radius) coordinate (OU);  
  \path (O) ++(125:\radius) coordinate (OL1); \path (O) ++(115:\outermidradius) coordinate (GOL1);   
  \path (O) ++(160:\radius) coordinate (OL2); \path (O) ++(140:\outermidradius) coordinate (GOL2);   
  \path (O) ++(240:\radius) coordinate (OL3); \path (O) ++(175:\outermidradius) coordinate (GOL3);
  \path (O) ++(180-125:\radius) coordinate (OR1); \path (O) ++(180-115:\outermidradius) coordinate (GOR1);   
  \path (O) ++(180-160:\radius) coordinate (OR2); \path (O) ++(180-140:\outermidradius) coordinate (GOR2);   
  \path (O) ++(180-240:\radius) coordinate (OR3); \path (O) ++(180-175:\outermidradius) coordinate (GOR3);      
  \path (O) ++(90:\innerradius) coordinate (IU);  \path (O) ++(180:\midradius) coordinate (GIUL); \path (O) ++(0:\midradius) coordinate (GIUR);
  \path (O) ++(140:\innerradius) coordinate (IL1);  \path (O) ++(200:\innermidradius) coordinate (GIL1);  
  \path (O) ++(190:\innerradius) coordinate (IL2);  \path (O) ++(230:\innermidradius) coordinate (GIL2);   
  \path (O) ++(180-140:\innerradius) coordinate (IR1);  \path (O) ++(180-200:\innermidradius) coordinate (GIR1);  
  \path (O) ++(180-190:\innerradius) coordinate (IR2);  \path (O) ++(180-230:\innermidradius) coordinate (GIR2);   
   
   \fill[black] (OU) circle[radius=2pt] ++(90:1em) node {\scriptsize $\pm \tvec{1}{0}$};
   \fill[black] (OL1) circle[radius=2pt] ++(120:1.2em) node {\scriptsize $\mp \tvec{1}{0}$};   
   \fill[black] (OL2) circle[radius=2pt] ++(150:1.5em) node {\scriptsize $\pm \tvec{1}{0}$};
   \fill[black] (OL3) circle[radius=2pt] ++(240:1.2em) node {\scriptsize $-\tvec{1}{0}$};    
   \fill[black] (OR1) circle[radius=2pt] ++(180-120:1em) node {\scriptsize $\mp \tvec{1}{0}$};   
   \fill[black] (OR2) circle[radius=2pt] ++(180-150:1.5em) node {\scriptsize $\pm \tvec{1}{0}$};
   \fill[black] (OR3) circle[radius=2pt] ++(180-240:1.2em) node {\scriptsize $-\tvec{1}{0}$};    
   \fill[black] (IU) circle[radius=2pt] ++(270:1em) node {\scriptsize $-\tvec{0}{1}$};
   \fill[black] (IL1) circle[radius=2pt] ++(300:1.2em) node {\scriptsize $\tvec{0}{1}$};   
   \fill[black] (IL2) circle[radius=2pt] ++(10:1.5em) node {\scriptsize $-\tvec{0}{1}$};                  
   \fill[black] (IR1) circle[radius=2pt] ++(180-300:1.2em) node {\scriptsize $\tvec{0}{1}$};   
   \fill[black] (IR2) circle[radius=2pt] ++(180-10:1.5em) node {\scriptsize $-\tvec{0}{1}$};                      
   \fill[black] (B) circle[radius=2pt] ++(270:1em) node {$\tvec{1}{0}$};   
   \fill[black] (IB) circle[radius=2pt] ++(90:1em) node {\scriptsize $\pm \tvec{0}{1}$};      
   \draw (GOR2) ++(300:0.8em) node[font=\tiny] {$m+\frac{n}{2}-1$};   
   \draw (GOR1) ++(320:1em) node[font=\tiny] {$m+\frac{n}{2}$};   
   \draw (OU) ++(320:2.2em) node[font=\tiny] {$m+\frac{n}{2}+1$};      
   \draw (GOL1) ++(45:1.2em) node[font=\tiny] {$m+\frac{n}{2}+2$};    
   \draw (GOL2) ++(90:0.7em) node[font=\tiny] {$m+\frac{n}{2}+3$};    
   \draw (GOL3) ++(286:9em) node[font=\tiny] {$m+n$};      
   \draw (GIUL) ++(260:1.3em) node[font=\tiny]  {$1$};
   \draw (GIL1) ++(230:0.4em) node[font=\tiny]  {$2$};
   \draw (GIL2) ++(80:0.65em) node[font=\tiny]  {$3$};  
   \draw (GIR2) ++(145:0.9em) node[font=\tiny]  {$m-1$};  
   \draw (GIR1) ++(340:0.4em) node[font=\tiny]  {$m$};  
   \draw (GIUR) ++(300:1.5em) node[font=\tiny]  {$m+1$};  
   \draw (GOR3) ++(254:9em) node[font=\tiny] {$m+2$};
   \draw (GIB) ++(20:1em) node[font=\tiny] {$\frac{m}{2}+1$};             
     
   \draw [black, xshift=4cm] plot [smooth, tension=1] coordinates { (IU) (OU)};               
   \draw [black, xshift=4cm] plot [smooth, tension=1] coordinates { (IU) (GIUL) (B)}; 
   \draw [black, xshift=4cm] plot [smooth, tension=1] coordinates { (IU) (GIUR) (B)};    
   \draw [black, xshift=4cm] plot [smooth, tension=1] coordinates { (IL1) (GIL1) (B)};        
   \draw [black, xshift=4cm] plot [smooth, tension=1] coordinates { (IL2) (GIL2) (B)};     
   \draw [black, xshift=4cm] plot [smooth, tension=1] coordinates { (IR1) (GIR1) (B)};        
   \draw [black, xshift=4cm] plot [smooth, tension=1] coordinates { (IR2) (GIR2) (B)};           
   \draw [black, xshift=4cm] plot [smooth, tension=1] coordinates { (IU) (GOL1) (OL1)};      
   \draw [black, xshift=4cm] plot [smooth, tension=1] coordinates { (IU) (GOL2) (OL2)};         
   \draw [black, xshift=4cm] plot [smooth, tension=1.3] coordinates { (IU) (GOL3) (OL3)};  
   \draw [black, xshift=4cm] plot [smooth, tension=1] coordinates { (IU) (GOR1) (OR1)};      
   \draw [black, xshift=4cm] plot [smooth, tension=1] coordinates { (IU) (GOR2) (OR2)};         
   \draw [black, xshift=4cm] plot [smooth, tension=1.3] coordinates { (IU) (GOR3) (OR3)};       
   \draw [black, xshift=4cm] plot [smooth, tension=1] coordinates { (IB) (B)};    
\end{tikzpicture}
}
\end{center}
 which corresponds to the initial graded quiver
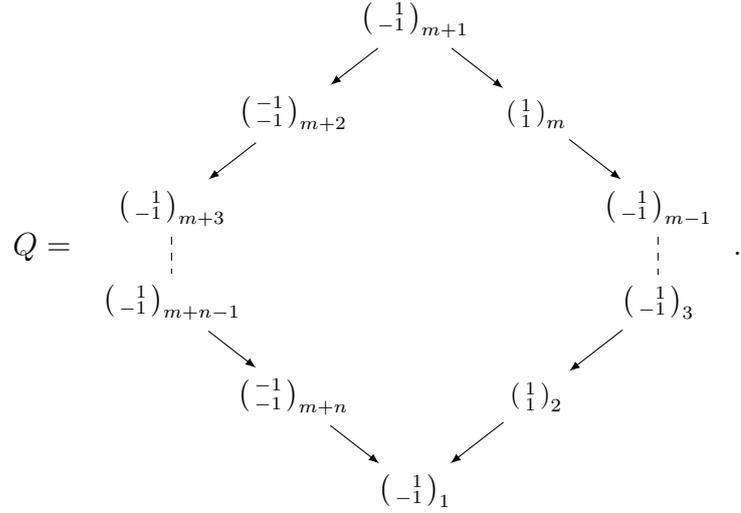
\begin{center}
$Q=$
\Image{
		\begin{tikzpicture}[align=center,node distance=1.25cm and 1.6cm, auto, font=\footnotesize, on grid]    
		  \node (Top) []{$\tvecal{1}{-1}_{m+1}$};                                        			  
		  \node (TopBL) [below left = of Top]{$\tvecal{-1}{-1}_{m+2}$};                                        			  
		  \node (TopBR) [below right = of Top]{$\tvecal{1}{1}_{m}$};                                        			  		  		  
		  \node (TopBLBL) [below left = of TopBL]{$\tvecal{1}{-1}_{m+3}$};                                        			  		  		  		  
		  \node (TopBRBR) [below right = of TopBR]{$\tvecal{1}{-1}_{m-1}$};
		  \node (LowerL) [below = of TopBLBL]{$\tvecal{1}{-1}_{m+n-1 }$};		                                          			  		  		  		  		  
		  \node (LowerR) [below = of TopBRBR]{$\tvecal{1}{-1}_{3} $};		  
		  \node (LowerLBR) [below right = of LowerL]{$\tvecal{-1}{-1}_{m+n}$};		  		  
		  \node (LowerRBL) [below left = of LowerR]{$\tvecal{1}{1}_{2}$};			  
		  \node (Bottom) [below right = of LowerLBR]{$\tvecal{1}{-1}_{1}$};		  
		  \draw[-latex] (Top) to node {} (TopBL);      
		  \draw[-latex] (Top) to node {} (TopBR); 		  
		  \draw[-latex] (TopBL) to node {} (TopBLBL);     
		  \draw[-latex] (TopBR) to node {} (TopBRBR);   		  		  
		  \draw[dashed, -] (TopBLBL) to node {} (LowerL);
		  \draw[dashed, -] (TopBRBR) to node {} (LowerR);		  
		  \draw[-latex] (LowerL) to node {} (LowerLBR);		 		  
		  \draw[-latex] (LowerR) to node {} (LowerRBL);		 		  		  
		  \draw[-latex] (LowerRBL) to node {} (Bottom);		 		  		  		  
		  \draw[-latex] (LowerLBR) to node {} (Bottom);		 		  		  		  		  
		\end{tikzpicture} 
}. 
\end{center}
As a tuple, the degree cluster is
\[ h_{n,m}:= \big( \overbrace{\tvecal{1}{-1}, \tvecal{1}{1},\dots, \tvecal{1}{-1}, \tvecal{1}{1}}^{\text{$m$ entries}},
 \overbrace{\tvecal{1}{-1},\tvecal{-1}{-1},\dots,\tvecal{1}{-1},\tvecal{-1}{-1}}^{\text{$n$ entries}} \big), \] 
and a basis for the grading space is given by 
\[ 
\{g_{n,m}:=( \overbrace{1, \dots, 1}^{\text{$m$ entries}}, \overbrace{1,-1,\dots,1,-1}^{\text{$n$ entries}}  ) ,
f_{n,m}:=( \overbrace{-1,1, \dots, -1,1}^{\text{$m$ entries}}, \overbrace{-1,-1,\dots,-1}^{\text{$n$ entries}}  )  \}.
\]
Note $g_{n,m}$ is the same grading vector as the one in the case for odd $n$.

At first glance, the addition of a dimension to the grading may appear to have a potential effect on the classification of the graded cluster algebra. However, it is not difficult to show that this will not be the case. In fact, the work we have done in Section \ref{section:Ann_n_m_m_odd} is already almost enough to determine the classification for the present case. 

\begin{prop}
Assume $n$ and $m$ are both even. In the graded cluster algebra $\A\big(h_{n,m},Q \big)$, the cardinality of the the set of cluster variables of degree $\tvec{d_1}{d_2}$ is determined by $d_1$: it is equal to the cardinality of the set of cluster variables of degree $d_1$ in $\A\big(g_{n,m}, Q \big)$. This cardinality is the same as for the case when $m$ is odd.  More specifically, $\A\big(h_{n,m},Q \big)$ is of mixed type with infinitely many variables in degrees $\pm \tvec{1}{-1}$ and $\pm \tvec{1}{1}$ and finitely many variables in degrees $\tvec{0}{0}$, $\pm \tvec{0}{2}$ and $\pm \tvec{2}{0}$, and these are all the occurring degrees.
\end{prop}

\begin{proof}
Write $\A = \A\big(h_{n,m},Q \big)$ and  $\A' =\A\big(g_{n,m} , Q \big)$. First note that, as in the case when $m$ is odd, $\A'$ is of mixed type with infinitely many variables in degrees $\pm 1$ and finitely many in degree $\pm 2$ and $0$. (As noted in Remark \ref{rmk:Ann_n_m_parity_m_doesn't_matter}, this may be proved in exactly the same way as for when $m$ is odd.) 

To justify that the degrees listed in the proposition are the only ones that can occur is easy: from the initial valued triangulation we may read off that these are the only possible combinations obtained when summing two valued marked points, to which any degree must correspond.

\emph{A priori}, at least one of $ \tvecal{1}{-1}$ and $\tvec{1}{1}$ must be associated with infinitely many variables, and similarly for $\tvecal{-1}{1}$ and $\tvec{-1}{-1}$. 
This is since, as noted in Lemma \ref{lem:Ann_n_m_odd_inf_vars_degs1}, there are infinitely many non-homotopic arcs of degrees $1$ and $-1$ in $\A'$, and so these same arcs give infinitely many variables in degrees $\tvec{1}{*}$ and $\tvecal{-1}{*}$ in $\A$. It is easy to see that we can in fact obtain infinitely many non-homotopic arcs both of degree $\tvecal{1}{-1}$ and degree $\tvecal{1}{1}$. Again, we do this by fixing appropriate points (of value $\tvec{1}{0}$ on the outer boundary and either $\tvec{0}{-1}$ or $\tvec{0}{1}$ as appropriate on the inner boundary) and taking arcs that wind around the inner boundary component an increasing number of times. Similarly, we have infinitely many variables in degrees $\tvecal{-1}{1}$ and $\tvecal{-1}{-1}$.

On the other hand, there are only finitely many possible arcs that can result in degrees $\pm 2$ or $0$ in $\A'$, and therefore this is the case also for any degree of the form $\tvec{0}{d}$ or $\pm \tvec{2}{d}$ in $\A$, since these same arcs are the only ways to obtain the corresponding first entries. 
\end{proof}

This shows what happens when we consider a standard grading. For a non-standard grading, however, it is possible to get a graded structure with a different behaviour. Consider the valuation on $\wt{A}(n,m)$ given as follows:

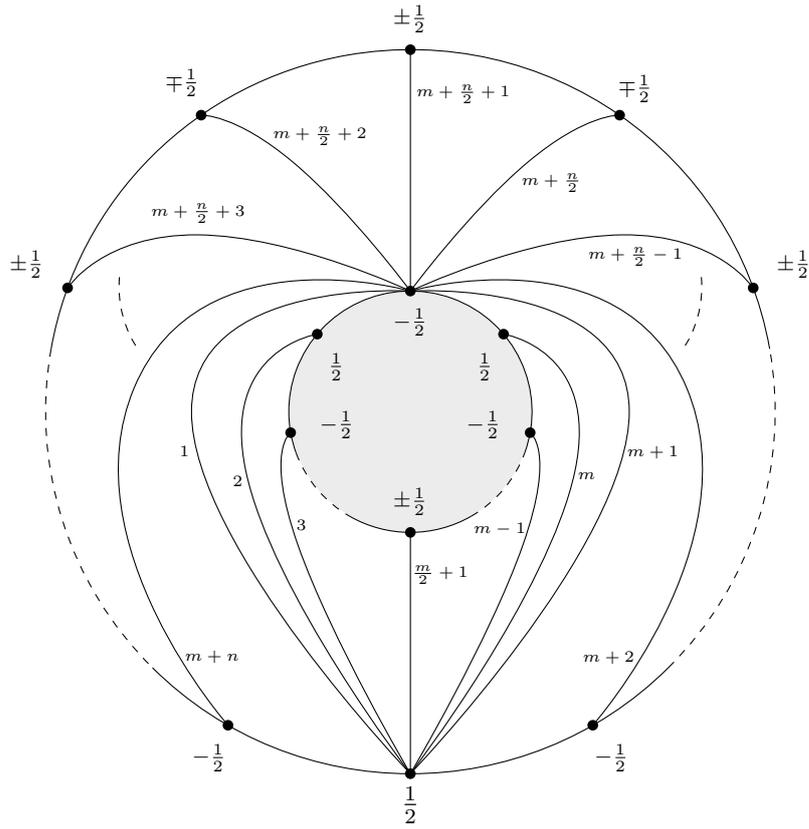
\begin{center}
{ 
\begin{tikzpicture}
  \coordinate (O) at (0,0);
  
  \def\radius{4.8cm}
  \def\midradius{ {(3/5)*\radius} }  
  \def\outermidradius{ {(3.8/5)*\radius} }      
  \def\innermidradius{ {(2.4/5)*\radius} }      
  \def\innerradius{  {(1/3)*\radius} }
  
  \centerarc[](O)(10:170:\radius)
  \centerarc[dashed](O)(170:225:\radius)
  \centerarc[](O)(225:315:\radius)
  \centerarc[dashed](O)(315:370:\radius)
  
  \shadedraw[inner color=gray!15,outer color=gray!15, draw=white] (O) circle[radius=\innerradius];  
  
  \centerarc[](O)(-20:200:\innerradius)
  \centerarc[dashed](O)(200:240:\innerradius)
  \centerarc[](O)(240:300:\innerradius)
  \centerarc[dashed](O)(300:340:\innerradius)
  
  \path (O) ++(155:{(4.4/5)*\radius}) coordinate (LA);
  \draw[dashed] (LA) arc (175:215:1.5);
  \path (O) ++(180-155:{(4.4/5)*\radius}) coordinate (RA);
  \draw[dashed] (RA) arc (180-175:180-215:1.5);  

  \path (O) ++(270:\radius) coordinate (B);
  \path (O) ++(270:\innerradius) coordinate (IB); \path (O) ++(270:\innermidradius) coordinate (GIB);
  \path (O) ++(90:\radius) coordinate (OU);  
  \path (O) ++(125:\radius) coordinate (OL1); \path (O) ++(115:\outermidradius) coordinate (GOL1);   
  \path (O) ++(160:\radius) coordinate (OL2); \path (O) ++(140:\outermidradius) coordinate (GOL2);   
  \path (O) ++(240:\radius) coordinate (OL3); \path (O) ++(175:\outermidradius) coordinate (GOL3);
  \path (O) ++(180-125:\radius) coordinate (OR1); \path (O) ++(180-115:\outermidradius) coordinate (GOR1);   
  \path (O) ++(180-160:\radius) coordinate (OR2); \path (O) ++(180-140:\outermidradius) coordinate (GOR2);   
  \path (O) ++(180-240:\radius) coordinate (OR3); \path (O) ++(180-175:\outermidradius) coordinate (GOR3);      
  \path (O) ++(90:\innerradius) coordinate (IU);  \path (O) ++(180:\midradius) coordinate (GIUL); \path (O) ++(0:\midradius) coordinate (GIUR);
  \path (O) ++(140:\innerradius) coordinate (IL1);  \path (O) ++(200:\innermidradius) coordinate (GIL1);  
  \path (O) ++(190:\innerradius) coordinate (IL2);  \path (O) ++(230:\innermidradius) coordinate (GIL2);   
  \path (O) ++(180-140:\innerradius) coordinate (IR1);  \path (O) ++(180-200:\innermidradius) coordinate (GIR1);  
  \path (O) ++(180-190:\innerradius) coordinate (IR2);  \path (O) ++(180-230:\innermidradius) coordinate (GIR2);   
   
   \fill[black] (OU) circle[radius=2pt] ++(90:1em) node {\scriptsize $\pm  \frac{1}{2}$};
   \fill[black] (OL1) circle[radius=2pt] ++(120:1.2em) node {\scriptsize $\mp\frac{1}{2}$};   
   \fill[black] (OL2) circle[radius=2pt] ++(150:1.5em) node {\scriptsize $\pm \frac{1}{2}$};
   \fill[black] (OL3) circle[radius=2pt] ++(240:1.2em) node {\scriptsize $-\frac{1}{2}$};    
   \fill[black] (OR1) circle[radius=2pt] ++(180-120:1em) node {\scriptsize $\mp \frac{1}{2}$};   
   \fill[black] (OR2) circle[radius=2pt] ++(180-150:1.5em) node {\scriptsize $\pm \frac{1}{2}$};
   \fill[black] (OR3) circle[radius=2pt] ++(180-240:1.2em) node {\scriptsize $-\frac{1}{2}$};    
   \fill[black] (IU) circle[radius=2pt] ++(270:1em) node {\scriptsize $-\frac{1}{2}$};
   \fill[black] (IL1) circle[radius=2pt] ++(300:1.2em) node {\scriptsize $\frac{1}{2}$};   
   \fill[black] (IL2) circle[radius=2pt] ++(10:1.5em) node {\scriptsize $-\frac{1}{2}$};                  
   \fill[black] (IR1) circle[radius=2pt] ++(180-300:1.2em) node {\scriptsize $\frac{1}{2}$};   
   \fill[black] (IR2) circle[radius=2pt] ++(180-10:1.5em) node {\scriptsize $-\frac{1}{2}$};                      
   \fill[black] (B) circle[radius=2pt] ++(270:1em) node {$\frac{1}{2}$};   
   \fill[black] (IB) circle[radius=2pt] ++(90:1em) node {\scriptsize $\pm \frac{1}{2}$};      
   \draw (GOR2) ++(300:0.8em) node[font=\tiny] {$m+\frac{n}{2}-1$};   
   \draw (GOR1) ++(320:1em) node[font=\tiny] {$m+\frac{n}{2}$};   
   \draw (OU) ++(320:2.2em) node[font=\tiny] {$m+\frac{n}{2}+1$};      
   \draw (GOL1) ++(45:1.2em) node[font=\tiny] {$m+\frac{n}{2}+2$};    
   \draw (GOL2) ++(90:0.7em) node[font=\tiny] {$m+\frac{n}{2}+3$};    
   \draw (GOL3) ++(286:9em) node[font=\tiny] {$m+n$};      
   \draw (GIUL) ++(260:1.3em) node[font=\tiny]  {$1$};
   \draw (GIL1) ++(230:0.4em) node[font=\tiny]  {$2$};
   \draw (GIL2) ++(80:0.65em) node[font=\tiny]  {$3$};  
   \draw (GIR2) ++(145:0.9em) node[font=\tiny]  {$m-1$};  
   \draw (GIR1) ++(340:0.4em) node[font=\tiny]  {$m$};  
   \draw (GIUR) ++(300:1.5em) node[font=\tiny]  {$m+1$};  
   \draw (GOR3) ++(254:9em) node[font=\tiny] {$m+2$};
   \draw (GIB) ++(20:1em) node[font=\tiny] {$\frac{m}{2}+1$};             
     
   \draw [black, xshift=4cm] plot [smooth, tension=1] coordinates { (IU) (OU)};               
   \draw [black, xshift=4cm] plot [smooth, tension=1] coordinates { (IU) (GIUL) (B)}; 
   \draw [black, xshift=4cm] plot [smooth, tension=1] coordinates { (IU) (GIUR) (B)};    
   \draw [black, xshift=4cm] plot [smooth, tension=1] coordinates { (IL1) (GIL1) (B)};        
   \draw [black, xshift=4cm] plot [smooth, tension=1] coordinates { (IL2) (GIL2) (B)};     
   \draw [black, xshift=4cm] plot [smooth, tension=1] coordinates { (IR1) (GIR1) (B)};        
   \draw [black, xshift=4cm] plot [smooth, tension=1] coordinates { (IR2) (GIR2) (B)};           
   \draw [black, xshift=4cm] plot [smooth, tension=1] coordinates { (IU) (GOL1) (OL1)};      
   \draw [black, xshift=4cm] plot [smooth, tension=1] coordinates { (IU) (GOL2) (OL2)};         
   \draw [black, xshift=4cm] plot [smooth, tension=1.3] coordinates { (IU) (GOL3) (OL3)};  
   \draw [black, xshift=4cm] plot [smooth, tension=1] coordinates { (IU) (GOR1) (OR1)};      
   \draw [black, xshift=4cm] plot [smooth, tension=1] coordinates { (IU) (GOR2) (OR2)};         
   \draw [black, xshift=4cm] plot [smooth, tension=1.3] coordinates { (IU) (GOR3) (OR3)};       
   \draw [black, xshift=4cm] plot [smooth, tension=1] coordinates { (IB) (B)};    
\end{tikzpicture}
}
\end{center}
 which corresponds to the initial graded quiver
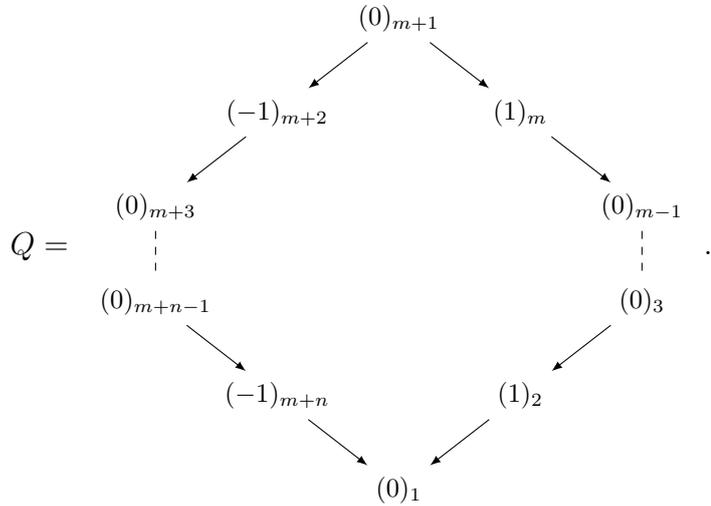
\begin{center}
$Q=$
\Image{
		\begin{tikzpicture}[align=center,node distance=1.25cm and 1.6cm, auto, font=\footnotesize, on grid]    
		  \node (Top) []{$(0)_{m+1}$};                                        			  
		  \node (TopBL) [below left = of Top]{$(-1)_{m+2}$};                                        			  
		  \node (TopBR) [below right = of Top]{$(1)_m$};                                        			  		  		  
		  \node (TopBLBL) [below left = of TopBL]{$(0)_{m+3}$};                                        			  		  		  		  
		  \node (TopBRBR) [below right = of TopBR]{$(0)_{m-1}$};
		  \node (LowerL) [below = of TopBLBL]{$(0)_{m+n-1 }$};		                                          			  		  		  		  		  
		  \node (LowerR) [below = of TopBRBR]{$(0)_{3} $};		  
		  \node (LowerLBR) [below right = of LowerL]{$(-1)_{m+n}$};		  		  
		  \node (LowerRBL) [below left = of LowerR]{$(1)_{2}$};			  
		  \node (Bottom) [below right = of LowerLBR]{$(0)_{1}$};		  
		  \draw[-latex] (Top) to node {} (TopBL);      
		  \draw[-latex] (Top) to node {} (TopBR); 		  
		  \draw[-latex] (TopBL) to node {} (TopBLBL);     
		  \draw[-latex] (TopBR) to node {} (TopBRBR);   		  		  
		  \draw[dashed, -] (TopBLBL) to node {} (LowerL);
		  \draw[dashed, -] (TopBRBR) to node {} (LowerR);		  
		  \draw[-latex] (LowerL) to node {} (LowerLBR);		 		  
		  \draw[-latex] (LowerR) to node {} (LowerRBL);		 		  		  
		  \draw[-latex] (LowerRBL) to node {} (Bottom);		 		  		  		  
		  \draw[-latex] (LowerLBR) to node {} (Bottom);		 		  		  		  		  
		\end{tikzpicture} 
}. 
\end{center}
As a tuple, the degree cluster is thus
$l_{n,m}:=(\overbrace{0,1, \dots,0,1}^{\text{$m$ entries}}, \overbrace{0,-1,\dots,0,-1}^{\text{$n$ entries}} )$.

In this case, the following is easy to see.

\begin{prop}
For $n$ and $m$ even, $\A\big(l_{n,m},Q \big)$ has infinitely many variables in degrees $0$ and $\pm 1$, and these are the only occurring degrees.
\end{prop}

\begin{proof}
That the above are the only possible degrees is clear.

By similar arguments to those used in previous cases we have considered, we can find infinitely many non-homotopic arcs of each of these degrees.
\end{proof}

\begin{rmk}
The theory we have developed in this chapter is likely to be able to successfully attack many other examples of cluster algebras arising from marked surfaces. Some of these examples could include a generalisation of the annulus that has $k$ inner boundaries rather than one, or the torus with a disc or multiple discs removed.
\end{rmk}

\biblio

\appendix
\chapter{Appendix}

\begin{apdxprf}[Lemma \ref{lem:E_aff_7_degseed_loop}] \label{apdx:E_aff_7_degseed_loop}
We will write out the computation for $n=1$, after which the result is immediate.
First we write down the matrices obtained along the mutation path $[7,1,6,2,5,8,3,4]$ from Lemma \ref{lem:E_aff_7_degseed_loop}. We have
\begin{align*}
 (\wt{E}_7)_{[4]} &=
 {
 \scalemath{0.9}{ \left( \begin{xsmallmatrix}{0.35em}
 0&1&0&0&0&0&0&0\\
-1&0&1&0&0&0&0&0\\
 0& -1&0& -1&0&0&0&0\\
 0&0&1&0&1&0&0&1\\
 0&0&0 &-1&0 &-1&0&0\\
 0&0&0&0&1&0 &-1&0\\
 0&0&0&0&0&1&0&0\\
 0&0&0& -1&0&0&0&0\\
\end{xsmallmatrix} \right) }
}, &
 (\wt{E}_7)_{[3,4]} &=
 {
 \scalemath{0.9}{ \left( \begin{smallmatrix} 
0&1&0&0&0&0&0&0\\
-1&0&-1&0&0&0&0&0\\
 0&1&0&1&0&0&0&0\\
 0&0&-1&0&1&0&0&1\\
 0&0&0&-1&0&-1&0&0\\
 0&0&0&0&1&0&-1&0\\
 0&0&0&0&0&1&0&0\\
 0&0&0&-1&0&0&0&0\\
\end{smallmatrix} \right) }
},
\\
 (\wt{E}_7)_{[8,3,4]} &=
 {
 \scalemath{0.9}{ \left( \begin{smallmatrix} 
 0& 1& 0& 0& 0& 0& 0& 0\\
-1& 0&-1& 0& 0& 0& 0& 0\\
 0& 1& 0& 1& 0& 0& 0& 0\\
 0& 0&-1& 0& 1& 0& 0&-1\\
 0& 0& 0&-1& 0&-1& 0& 0\\
 0& 0& 0& 0& 1& 0&-1& 0\\
 0& 0& 0& 0& 0& 1& 0& 0\\
 0& 0& 0& 1& 0& 0& 0& 0\\
\end{smallmatrix} \right) }
}, &
 (\wt{E}_7)_{[5,8,3,4]} &=
 {
 \scalemath{0.9}{ \left( \begin{smallmatrix} 
 0& 1& 0& 0& 0& 0& 0& 0\\
-1& 0&-1& 0& 0& 0& 0& 0\\
 0& 1& 0& 1& 0& 0& 0& 0\\
 0& 0&-1& 0&-1& 0& 0&-1\\
 0& 0& 0& 1& 0& 1& 0& 0\\
 0& 0& 0& 0&-1& 0&-1& 0\\
 0& 0& 0& 0& 0& 1& 0& 0\\
 0& 0& 0& 1& 0& 0& 0& 0\\
\end{smallmatrix} \right) }
},
\\
 (\wt{E}_7)_{[2,5,8,3,4]} &=
 {
 \scalemath{0.9}{ \left( \begin{xsmallmatrix}{0.35em}
 0&-1& 0& 0& 0& 0& 0& 0\\
 1& 0& 1& 0& 0& 0& 0& 0\\
 0&-1& 0& 1& 0& 0& 0& 0\\
 0& 0&-1& 0&-1& 0& 0&-1\\
 0& 0& 0& 1& 0& 1& 0& 0\\
 0& 0& 0& 0&-1& 0&-1& 0\\
 0& 0& 0& 0& 0& 1& 0& 0\\
 0& 0& 0& 1& 0& 0& 0& 0\\
\end{xsmallmatrix} \right) }
}, &
 (\wt{E}_7)_{[6,2,5,8,3,4]} &=
 {
 \scalemath{0.9}{ \left( \begin{smallmatrix} 
 0&-1& 0& 0& 0& 0& 0& 0\\
 1& 0& 1& 0& 0& 0& 0& 0\\
 0&-1& 0& 1& 0& 0& 0& 0\\
 0& 0&-1& 0&-1& 0& 0&-1\\
 0& 0& 0& 1& 0&-1& 0& 0\\
 0& 0& 0& 0& 1& 0& 1& 0\\
 0& 0& 0& 0& 0&-1& 0& 0\\
 0& 0& 0& 1& 0& 0& 0& 0\\
\end{smallmatrix} \right) }
},
\\
 (\wt{E}_7)_{[1,6,2,5,8,3,4]} &=
 {
 \scalemath{0.9}{ \left( \begin{smallmatrix} 
  0& 1& 0& 0& 0& 0& 0& 0\\
 -1& 0& 1& 0& 0& 0& 0& 0\\
  0&-1& 0& 1& 0& 0& 0& 0\\
  0& 0&-1& 0&-1& 0& 0&-1\\
  0& 0& 0& 1& 0&-1& 0& 0\\
  0& 0& 0& 0& 1& 0& 1& 0\\
  0& 0& 0& 0& 0&-1& 0& 0\\
  0& 0& 0& 1& 0& 0& 0& 0\\
\end{smallmatrix} \right) }
}, &
 (\wt{E}_7)_{[7,1,6,2,5,8,3,4]} &=
 {
 \scalemath{0.9}{ \left( \begin{xsmallmatrix}{0.1em}
  0 &1 &0 &0 &0 &0 &0 &0\\
 -1 &0 &1 &0 &0 &0 &0 &0\\
  0&-1 &0 &1 &0 &0 &0 &0\\
  0 &0&-1 &0&-1 &0 &0&-1\\
  0 &0 &0 &1 &0&-1 &0 &0\\
  0 &0 &0 &0 &1 &0&-1 &0\\
  0 &0 &0 &0 &0 &1 &0 &0\\
  0 &0 &0 &1 &0 &0 &0 &0\\
\end{xsmallmatrix} \right) }
}.
\end{align*}
Thus $(\wt{E}_7)_{[7,1,6,2,5,8,3,4]} =\wt{E}_7$, as claimed.

Then, for the degree seeds, we have
\allowdisplaybreaks
\begin{align*}
\deg.\cl_{\wt{E}_7}[4] &=  
\left( \CII{-1}{-1},\CII{0}{0},\CII{-1}{-1},\CII{0}{0},\CII{1}{0},\CII{0}{0}, \CII{1}{0}, \CII{0}{1} \right) ,\\
\deg.\cl_{\wt{E}_7}[3,4] &=  
\left( \CII{-1}{-1},\CII{0}{0},\CII{1}{1},\CII{0}{0},\CII{1}{0},\CII{0}{0}, \CII{1}{0}, \CII{0}{1} \right) ,\\
\deg.\cl_{\wt{E}_7}[8,3,4] &=  
\left( \CII{-1}{-1},\CII{0}{0},\CII{1}{1},\CII{0}{0},\CII{1}{0},\CII{0}{0}, \CII{1}{0}, \CII{0}{-1} \right) ,\\
\deg.\cl_{\wt{E}_7}[5,8,3,4] &=  
\left( \CII{-1}{-1},\CII{0}{0},\CII{1}{1},\CII{0}{0},\CII{-1}{0},\CII{0}{0}, \CII{1}{0}, \CII{0}{-1} \right) ,\\
\deg.\cl_{\wt{E}_7}[2,5,8,3,4] &=  
\left( \CII{-1}{-1},\CII{0}{0},\CII{1}{1},\CII{0}{0},\CII{-1}{0},\CII{0}{0}, \CII{1}{0}, \CII{0}{-1} \right) ,\\
\deg.\cl_{\wt{E}_7}[6,2,5,8,3,4] &=  
\left( \CII{-1}{-1},\CII{0}{0},\CII{1}{1},\CII{0}{0},\CII{-1}{0},\CII{0}{0}, \CII{1}{0}, \CII{0}{-1} \right) ,\\
\deg.\cl_{\wt{E}_7}[1,6,2,5,8,3,4] &=  
\left( \CII{1}{1},\CII{0}{0},\CII{1}{1},\CII{0}{0},\CII{-1}{0},\CII{0}{0}, \CII{1}{0}, \CII{0}{-1} \right) ,\\
\deg.\cl_{\wt{E}_7}[7,1,6,2,5,8,3,4] &=  
\left( \CII{1}{1},\CII{0}{0},\CII{1}{1},\CII{0}{0},\CII{-1}{0},\CII{0}{0}, \CII{-1}{0}, \CII{0}{-1} \right) .\\
\end{align*}
Thus $\deg.\cl_{\wt{E}_7}[7,1,6,2,5,8,3,4]$ is the negation of the initial cluster, as claimed. \hfill \qed
\end{apdxprf}

\begin{apdxprf}[Lemma \ref{lem:E_aff_8_degseed_loop}] \label{apdx:E_aff_8_degseed_loop}
We write down the matrices obtained along the mutation path $[3, 4, 2, 6, 9, 1, 6, 5, 6, 7, 8]$ from Lemma \ref{lem:E_aff_7_degseed_loop}. We have
\begin{align*}
\resizebox{0.9cm}{!}{$ (\wt{E}_7)_{[8]} $} &=
 {
 \scalemath{0.8}{ \left( \begin{xsmallmatrix}{.2em}
  0 &1 &0 &0 &0 &0 &0 &0 &0\\
 -1 &0 &1 &0 &0 &0 &0 &0 &0\\
  0&-1 &0&-1 &0 &0 &0 &0&-1\\
  0 &0 &1 &0&-1 &0 &0 &0 &0\\
  0 &0 &0 &1 &0&-1 &0 &0 &0\\
  0 &0 &0 &0 &1 &0&-1 &0 &0\\
  0 &0 &0 &0 &0 &1 &0 &1 &0\\
  0 &0 &0 &0 &0 &0&-1 &0 &0\\
  0 &0 &1 &0 &0 &0 &0 &0 &0\\
\end{xsmallmatrix} \right) }
}, &
\resizebox{1cm}{!}{$ (\wt{E}_7)_{[7,8]} $}&=
 {
 \scalemath{0.8}{ \left( \begin{xsmallmatrix}{.2em}
  0 &1 &0 &0 &0 &0 &0 &0 &0\\
-1 &0 &1 &0 &0 &0 &0 &0 &0\\
 0 &-1 &0&-1 &0 &0 &0 &0&-1\\
 0 &0 &1 &0&-1 &0 &0 &0 &0\\
 0 &0 &0 &1 &0&-1 &0 &0 &0\\
 0 &0 &0 &0 &1 &0 &1 &0 &0\\
 0 &0 &0 &0 &0&-1 &0&-1 &0\\
 0 &0 &0 &0 &0 &0 &1 &0 &0\\
 0 &0 &1 &0 &0 &0 &0 &0 &0\\
\end{xsmallmatrix} \right) }
},
\\
\resizebox{1.3cm}{!}{$ (\wt{E}_7)_{[6,7,8]} $}&=
 {
 \scalemath{0.8}{ \left( \begin{xsmallmatrix}{.62em}
 0 &1 &0 &0 &0 &0 &0 &0 &0\\
-1 &0 &1 &0 &0 &0 &0 &0 &0\\
 0& 1 &0& 1 &0 &0 &0 &0& 1\\
 0 &0 &1 &0& 1 &0 &0 &0 &0\\
 0 &0 &0 &1 &0 &1 &0 &0 &0\\
 0 &0 &0 &0& 1 &0& 1 &0 &0\\
 0 &0 &0 &0 &0 &1 &0& 1 &0\\
 0 &0 &0 &0 &0 &0 &1 &0 &0\\
 0 &0 &1 &0 &0 &0 &0 &0 &0\\
\end{xsmallmatrix} \right) }
}, &
\resizebox{1.4cm}{!}{$ (\wt{E}_7)_{[5,6,7,8]} $}&=
 {
 \scalemath{0.8}{ \left( \begin{xsmallmatrix}{.2em}
 0 &1 &0 &0 &0 &0 &0 &0 &0\\
-1 &0 &1 &0 &0 &0 &0 &0 &0\\
 0&-1 &0&-1 &0 &0 &0 &0&-1\\
 0 &0 &1 &0 &1 &0 &0 &0 &0\\
 0 &0 &0&-1 &0&-1 &0 &0 &0\\
 0 &0 &0 &0 &1 &0&-1 &0 &0\\
 0 &0 &0 &0 &0 &1 &0&-1 &0\\
 0 &0 &0 &0 &0 &0 &1 &0 &0\\
 0 &0 &1 &0 &0 &0 &0 &0 &0\\
\end{xsmallmatrix} \right) }
},
\\
 \resizebox{1.7cm}{!}{$ (\wt{E}_7)_{[6,5,6,7,8]}$} &=
 {
 \scalemath{0.8}{ \left( \begin{xsmallmatrix}{.14em} 
 0& 1& 0& 0& 0& 0& 0& 0& 0\\
-1& 0& 1& 0& 0& 0& 0& 0& 0\\
 0&-1& 0&-1& 0& 0& 0& 0&-1\\
 0& 0& 1& 0& 1& 0& 0& 0& 0\\
 0& 0& 0&-1& 0& 1&-1& 0& 0\\
 0& 0& 0& 0&-1& 0& 1& 0& 0\\
 0& 0& 0& 0& 1&-1& 0&-1& 0\\
 0& 0& 0& 0& 0& 0& 1& 0& 0\\
 0& 0& 1& 0& 0& 0& 0& 0& 0\\
\end{xsmallmatrix} \right) }
},&
 \resizebox{1.8cm}{!}{$ (\wt{E}_7)_{[1,6,5,6,7,8]}$} &=
 {
 \scalemath{0.8}{ \left( \begin{xsmallmatrix}{.2em} 
 0&-1& 0& 0& 0& 0& 0& 0& 0\\
 1& 0& 1& 0& 0& 0& 0& 0& 0\\
 0&-1& 0&-1& 0& 0& 0& 0&-1\\
 0& 0& 1& 0& 1& 0& 0& 0& 0\\
 0& 0& 0&-1& 0& 1&-1& 0& 0\\
 0& 0& 0& 0&-1& 0& 1& 0& 0\\
 0& 0& 0& 0& 1&-1& 0&-1& 0\\
 0& 0& 0& 0& 0& 0& 1& 0& 0\\
 0& 0& 1& 0& 0& 0& 0& 0& 0\\
\end{xsmallmatrix} \right) }
},
\\
 \resizebox{2.1cm}{!}{$ (\wt{E}_7)_{[9,1,6,5,6,7,8]} $}&=
 {
 \scalemath{0.8}{ \left( \begin{xsmallmatrix}{.2em} 
 0&-1 &0 &0 &0 &0 &0 &0 &0\\
 1 &0 &1 &0 &0 &0 &0 &0 &0\\
 0&-1 &0&-1 &0 &0 &0 &0 &1\\
 0 &0 &1 &0 &1 &0 &0 &0 &0\\
 0 &0 &0&-1 &0 &1&-1 &0 &0\\
 0 &0 &0 &0&-1 &0 &1 &0 &0\\
 0 &0 &0 &0 &1&-1 &0&-1 &0\\
 0 &0 &0 &0 &0 &0 &1 &0 &0\\
 0 &0&-1 &0 &0 &0 &0 &0 &0\\
\end{xsmallmatrix} \right) }
}, &
  \resizebox{2.2cm}{!}{$(\wt{E}_7)_{[6,9,1,6,5,6,7,8]} $}&=
 {
 \scalemath{0.8}{ \left(  \begin{xsmallmatrix}{.28em} 
 0&-1& 0& 0& 0& 0& 0& 0& 0\\
 1& 0& 1& 0& 0& 0& 0& 0& 0\\
 0&-1& 0&-1& 0& 0& 0& 0& 1\\
 0& 0& 1& 0& 1& 0& 0& 0& 0\\
 0& 0& 0&-1& 0&-1& 0& 0& 0\\
 0& 0& 0& 0& 1& 0&-1& 0& 0\\
 0& 0& 0& 0& 0& 1& 0&-1& 0\\
 0& 0& 0& 0& 0& 0& 1& 0& 0\\
 0& 0&-1& 0& 0& 0& 0& 0& 0\\
\end{xsmallmatrix} \right) }
},
\\
 \resizebox{2.6cm}{!}{$ (\wt{E}_7)_{[2,6,9,1,6,5,6,7,8]} $}&=
 {
 \scalemath{0.8}{ \left( \begin{xsmallmatrix}{.27em} 
 0& 1& 0& 0& 0& 0& 0& 0& 0\\
-1& 0&-1& 0& 0& 0& 0& 0& 0\\
 0& 1& 0&-1& 0& 0& 0& 0& 1\\
 0& 0& 1& 0& 1& 0& 0& 0& 0\\
 0& 0& 0&-1& 0&-1& 0& 0& 0\\
 0& 0& 0& 0& 1& 0&-1& 0& 0\\
 0& 0& 0& 0& 0& 1& 0&-1& 0\\
 0& 0& 0& 0& 0& 0& 1& 0& 0\\
 0& 0&-1& 0& 0& 0& 0& 0& 0\\
\end{xsmallmatrix} \right) }
}, &
  \resizebox{2.75cm}{!}{$(\wt{E}_7)_{[4,2,6,9,1,6,5,6,7,8]} $}&=
 {
 \scalemath{0.8}{ \left( \begin{xsmallmatrix}{.3em} 
 0& 1& 0& 0& 0& 0& 0& 0& 0\\
-1& 0&-1& 0& 0& 0& 0& 0& 0\\
 0& 1& 0& 1& 0& 0& 0& 0& 1\\
 0& 0&-1& 0&-1& 0& 0& 0& 0\\
 0& 0& 0& 1& 0&-1& 0& 0& 0\\
 0& 0& 0& 0& 1& 0&-1& 0& 0\\
 0& 0& 0& 0& 0& 1& 0&-1& 0\\
 0& 0& 0& 0& 0& 0& 1& 0& 0\\
 0& 0&-1& 0& 0& 0& 0& 0& 0\\
\end{xsmallmatrix} \right) }
},
\\
 \resizebox{2.95cm}{!}{$(\wt{E}_7)_{[3,4,2,6,9,1,6,5,6,7,8]} $} &=
 {
 \scalemath{0.8}{ \left( \begin{xsmallmatrix}{.15em} 
 0& 1& 0& 0& 0& 0& 0& 0& 0\\
-1& 0& 1& 0& 0& 0& 0& 0& 0\\
 0&-1& 0&-1& 0& 0& 0& 0&-1\\
 0& 0& 1& 0&-1& 0& 0& 0& 0\\
 0& 0& 0& 1& 0&-1& 0& 0& 0\\
 0& 0& 0& 0& 1& 0&-1& 0& 0\\
 0& 0& 0& 0& 0& 1& 0&-1& 0\\
 0& 0& 0& 0& 0& 0& 1& 0& 0\\
 0& 0& 1& 0& 0& 0& 0& 0& 0\\
\end{xsmallmatrix} \right) }
}.
\end{align*}
Thus $(\wt{E}_8)_{[3, 4, 2, 6, 9, 1, 6, 5, 6, 7, 8]} =\wt{E}_8$, as claimed.
\end{apdxprf}

Then, for the degree seeds, we have
\allowdisplaybreaks
\begin{align*}
\deg.\cl_{\wt{E}_7}[8] &=  
\left( 0, 0, 0, -1, 0, -1, 0, 1, 1\right) ,\\
\deg.\cl_{\wt{E}_7}[7,8] &=  
\left(0, 0, 0, -1, 0, -1, 0, 1, 1  \right) ,\\
\deg.\cl_{\wt{E}_7}[6,7,8] &=  
\left(0, 0, 0, -1, 0, 1, 0, 1, 1  \right) ,\\
\deg.\cl_{\wt{E}_7}[5,6,7,8] &=  
\left( 0, 0, 0, -1, 0, 1, 0, 1, 1\right) ,\\
\deg.\cl_{\wt{E}_7}[6,5,6,7,8] &=  
\left(0, 0, 0, -1, 0, -1, 0, 1, 1 \right) ,\\
\deg.\cl_{\wt{E}_7}[1,6,5,6,7,8] &=  
\left( 0, 0, 0, -1, 0, -1, 0, 1, 1\right) ,\\
\deg.\cl_{\wt{E}_7}[9,1,6,5,6,7,8] &=  
\left( 0, 0, 0, -1, 0, -1, 0, 1, -1\right) ,\\
\deg.\cl_{\wt{E}_7}[6,9,1,6,5,6,7,8] &=  
\left(  0, 0, 0, -1, 0, 1, 0, 1, -1\right) ,\\
\deg.\cl_{\wt{E}_7}[2,6,9,1,6,5,6,7,8] &=  
\left( 0, 0, 0, -1, 0, 1, 0, 1, -1\right) ,\\
\deg.\cl_{\wt{E}_7}[4,2,6,9,1,6,5,6,7,8] &=  
\left( 0, 0, 0, 1, 0, 1, 0, 1, -1 \right) ,\\
\deg.\cl_{\wt{E}_7}[3,4,2,6,9,1,6,5,6,7,8] &=  
\left(  0, 0, 0, 1, 0, 1, 0, 1, -1 \right).
\end{align*}
Thus $\deg.\cl_{\wt{E}_7}[3, 4, 2, 6, 9, 1, 6, 5, 6, 7, 8]$ is the negation of the initial cluster, as claimed. \hfill \qed

\biblio

\bibliography{caref}\label{references}

\begin{thebibliography}{}

\bibitem[\protect\astroncite{Assem et~al.}{2008}]{ABBS}
Assem, I., Blais, M., Br\"ustle, T., and Samson, A. (2008).
\newblock Mutation classes of skew-symmetric {$3\times 3$}-matrices.
\newblock {\em Comm. Algebra}, 36(4):1209--1220.

\bibitem[\protect\astroncite{Beineke et~al.}{2011}]{BBH}
Beineke, A., Br\"ustle, T., and Hille, L. (2011).
\newblock Cluster-cyclic quivers with three vertices and the {M}arkov equation.
\newblock {\em Algebr. Represent. Theory}, 14(1):97--112.
\newblock With an appendix by Otto Kerner.

\bibitem[\protect\astroncite{Berenstein et~al.}{2005}]{FZBCAIII}
Berenstein, A., Fomin, S., and Zelevinsky, A. (2005).
\newblock Cluster algebras. {III}. {U}pper bounds and double {B}ruhat cells.
\newblock {\em Duke Math. J.}, 126(1):1--52.

\bibitem[\protect\astroncite{Berenstein and Zelevinsky}{2005}]{BZ}
Berenstein, A. and Zelevinsky, A. (2005).
\newblock Quantum cluster algebras.
\newblock {\em Adv. Math.}, 195(2):405--455.

\bibitem[\protect\astroncite{Buan et~al.}{2006}]{BMMRT}
Buan, A.~B., Marsh, R., Reineke, M., Reiten, I., and Todorov, G. (2006).
\newblock Tilting theory and cluster combinatorics.
\newblock {\em Adv. Math.}, 204(2):572--618.

\bibitem[\protect\astroncite{Caldero and Keller}{2006}]{CalderoKeller}
Caldero, P. and Keller, B. (2006).
\newblock From triangulated categories to cluster algebras. {II}.
\newblock {\em Ann. Sci. \'Ecole Norm. Sup. (4)}, 39(6):983--1009.

\bibitem[\protect\astroncite{Carter}{2005}]{Carter-Book}
Carter, R.~W. (2005).
\newblock {\em Lie algebras of finite and affine type}, volume~96 of {\em
  Cambridge Studies in Advanced Mathematics}.
\newblock Cambridge University Press, Cambridge.

\bibitem[\protect\astroncite{Chekhov and Penner}{2003}]{Che}
Chekhov, L.~O. and Penner, R.~C. (2003).
\newblock Introduction to {T}hurston's quantum theory.
\newblock {\em Uspekhi Mat. Nauk}, 58(6(354)):93--138.

\bibitem[\protect\astroncite{Felikson et~al.}{2012}]{FSTFinite}
Felikson, A., Shapiro, M., and Tumarkin, P. (2012).
\newblock Skew-symmetric cluster algebras of finite mutation type.
\newblock {\em J. Eur. Math. Soc. (JEMS)}, 14(4):1135--1180.

\bibitem[\protect\astroncite{Fomin et~al.}{2008}]{FST}
Fomin, S., Shapiro, M., and Thurston, D. (2008).
\newblock Cluster algebras and triangulated surfaces. {I}. {C}luster complexes.
\newblock {\em Acta Math.}, 201(1):83--146.

\bibitem[\protect\astroncite{Fomin and Zelevinsky}{2002}]{FZCAI}
Fomin, S. and Zelevinsky, A. (2002).
\newblock Cluster algebras. {I}. {F}oundations.
\newblock {\em J. Amer. Math. Soc.}, 15(2):497--529.

\bibitem[\protect\astroncite{Fomin and Zelevinsky}{2003}]{FZCAII}
Fomin, S. and Zelevinsky, A. (2003).
\newblock Cluster algebras. {II}. {F}inite type classification.
\newblock {\em Invent. Math.}, 154(1):63--121.

\bibitem[\protect\astroncite{Fomin and Zelevinsky}{2007}]{FZCAIV}
Fomin, S. and Zelevinsky, A. (2007).
\newblock Cluster algebras. {IV}. {C}oefficients.
\newblock {\em Compos. Math.}, 143(1):112--164.

\bibitem[\protect\astroncite{Gei\ss et~al.}{2013}]{GLS}
Gei\ss, C., Leclerc, B., and Schr\"oer, J. (2013).
\newblock Cluster structures on quantum coordinate rings.
\newblock {\em Selecta Math. (N.S.)}, 19(2):337--397.

\bibitem[\protect\astroncite{Gekhtman et~al.}{2010}]{GSV}
Gekhtman, M., Shapiro, M., and Vainshtein, A. (2010).
\newblock {\em Cluster algebras and {P}oisson geometry}, volume 167 of {\em
  Mathematical Surveys and Monographs}.
\newblock American Mathematical Society, Providence, RI.

\bibitem[\protect\astroncite{Grabowski}{2015}]{G15}
Grabowski, J.~E. (2015).
\newblock Graded cluster algebras.
\newblock {\em J. Algebraic Combin.}, 42(4):1111--1134.

\bibitem[\protect\astroncite{Grabowski and Launois}{2014}]{GL}
Grabowski, J.~E. and Launois, S. (2014).
\newblock Graded quantum cluster algebras and an application to quantum
  {G}rassmannians.
\newblock {\em Proc. Lond. Math. Soc. (3)}, 109(3):697--732.

\bibitem[\protect\astroncite{{Grabowski} and {Pressland}}{2016}]{GP}
{Grabowski}, J.~E. and {Pressland}, M. (2016).
\newblock {Graded Frobenius cluster categories}.
\newblock {\em ArXiv e-prints}.

\bibitem[\protect\astroncite{Graham et~al.}{1994}]{GKP}
Graham, R.~L., Knuth, D.~E., and Patashnik, O. (1994).
\newblock {\em Concrete mathematics}.
\newblock Addison-Wesley Publishing Company, Reading, MA, second edition.
\newblock A foundation for computer science.

\bibitem[\protect\astroncite{Harer}{1986}]{Har}
Harer, J.~L. (1986).
\newblock The virtual cohomological dimension of the mapping class group of an
  orientable surface.
\newblock {\em Invent. Math.}, 84(1):157--176.

\bibitem[\protect\astroncite{Hatcher}{1991}]{Hat}
Hatcher, A. (1991).
\newblock On triangulations of surfaces.
\newblock {\em Topology Appl.}, 40(2):189--194.

\bibitem[\protect\astroncite{Keller}{}]{KJA}
Keller, B.
\newblock Quiver mutation in {J}ava.
\newblock \url{www.math.jussieu.fr/~keller/quivermutation}.

\bibitem[\protect\astroncite{Mosher}{1988}]{Mos}
Mosher, L. (1988).
\newblock Tiling the projective foliation space of a punctured surface.
\newblock {\em Trans. Amer. Math. Soc.}, 306(1):1--70.

\bibitem[\protect\astroncite{Muller}{2016}]{Muller}
Muller, G. (2016).
\newblock Skein and cluster algebras of marked surfaces.
\newblock {\em Quantum Topol.}, 7(3):435--503.

\bibitem[\protect\astroncite{Palu}{2008}]{Pal}
Palu, Y. (2008).
\newblock Cluster characters for 2-{C}alabi-{Y}au triangulated categories.
\newblock {\em Ann. Inst. Fourier (Grenoble)}, 58(6):2221--2248.

\bibitem[\protect\astroncite{Schiffler}{2014}]{SchifflerQR}
Schiffler, R. (2014).
\newblock {\em Quiver representations}.
\newblock CMS Books in Mathematics/Ouvrages de Math\'ematiques de la SMC.
  Springer, Cham.

\bibitem[\protect\astroncite{Scott}{2006}]{Sc}
Scott, J.~S. (2006).
\newblock Grassmannians and cluster algebras.
\newblock {\em Proceedings of the London Mathematical Society}, 92(2):345--380.

\bibitem[\protect\astroncite{Warkentin}{2014}]{War}
Warkentin, M. (2014).
\newblock {\em Exchange Graphs via Quiver Mutation}.
\newblock PhD thesis, Universit{\"a}t Bonn.

\end{thebibliography}


\begin{thebibliography}{1}

\bibitem{G15}
Jan~E. Grabowski.
\newblock Graded cluster algebras.
\newblock {\em J. Algebraic Combin.}, 42(4):1111--1134, 2015.

\bibitem{BMMRT}
Aslak~Bakke Buan, Robert Marsh, Markus Reineke, Idun Reiten, and Gordana
  Todorov.
\newblock Tilting theory and cluster combinatorics.
\newblock {\em Adv. Math.}, 204(2):572--618, 2006.

\bibitem{Pal}
Yann Palu.
\newblock Cluster characters for 2-{C}alabi-{Y}au triangulated categories.
\newblock {\em Ann. Inst. Fourier (Grenoble)}, 58(6):2221--2248, 2008.

\bibitem{FZCAI}
Sergey Fomin and Andrei Zelevinsky.
\newblock Cluster algebras. {I}. {F}oundations.
\newblock {\em J. Amer. Math. Soc.}, 15(2):497--529, 2002.

\bibitem{ABBS}
Ibrahim Assem, Martin Blais, Thomas Br\"ustle, and Audrey Samson.
\newblock Mutation classes of skew-symmetric {$3\times 3$}-matrices.
\newblock {\em Comm. Algebra}, 36(4):1209--1220, 2008.

\bibitem{BBH}
Andre Beineke, Thomas Br\"ustle, and Lutz Hille.
\newblock Cluster-cyclic quivers with three vertices and the {M}arkov equation.
\newblock {\em Algebr. Represent. Theory}, 14(1):97--112, 2011.
\newblock With an appendix by Otto Kerner.

\bibitem{War}
Matthias Warkentin.
\newblock {\em Exchange Graphs via Quiver Mutation}.
\newblock PhD thesis, Universit{\"a}t Bonn, 2014.

\bibitem{GSV}
Michael Gekhtman, Michael Shapiro, and Alek Vainshtein.
\newblock {\em Cluster algebras and {P}oisson geometry}, volume 167 of {\em
  Mathematical Surveys and Monographs}.
\newblock American Mathematical Society, Providence, RI, 2010.

\bibitem{CalderoKeller}
Philippe Caldero and Bernhard Keller.
\newblock From triangulated categories to cluster algebras. {II}.
\newblock {\em Ann. Sci. \'Ecole Norm. Sup. (4)}, 39(6):983--1009, 2006.

\end{thebibliography}


\begin{thebibliography}{}

\end{thebibliography}


\begin{thebibliography}{1}

\bibitem{FSTFinite}
Anna Felikson, Michael Shapiro, and Pavel Tumarkin.
\newblock Skew-symmetric cluster algebras of finite mutation type.
\newblock {\em J. Eur. Math. Soc. (JEMS)}, 14(4):1135--1180, 2012.

\bibitem{G15}
Jan~E. Grabowski.
\newblock Graded cluster algebras.
\newblock {\em J. Algebraic Combin.}, 42(4):1111--1134, 2015.

\bibitem{KJA}
Bernhard Keller.
\newblock Quiver mutation in {J}ava.
\newblock \url{www.math.jussieu.fr/~keller/quivermutation}.

\bibitem{GKP}
Ronald~L. Graham, Donald~E. Knuth, and Oren Patashnik.
\newblock {\em Concrete mathematics}.
\newblock Addison-Wesley Publishing Company, Reading, MA, second edition, 1994.
\newblock A foundation for computer science.

\end{thebibliography}


\begin{thebibliography}{1}

\bibitem{G15}
Jan~E Grabowski.
\newblock Graded cluster algebras.
\newblock {\em arXiv preprint arXiv:1309.6170}, 2015.

\bibitem{FZCAII}
Sergey Fomin and Andrei Zelevinsky.
\newblock Cluster algebras {II}: Finite type classification.
\newblock {\em Inventiones {M}athematicae}, 154(1):63--121, 2003.

\bibitem{Carter-Book}
Roger~William Carter.
\newblock {\em Lie algebras of finite and affine type}.
\newblock Number~96. Cambridge University Press, 2005.

\end{thebibliography}


\begin{thebibliography}{10}

\bibitem{FZCAI}
Sergey Fomin and Andrei Zelevinsky.
\newblock Cluster algebras. {I}. {F}oundations.
\newblock {\em J. Amer. Math. Soc.}, 15(2):497--529, 2002.

\bibitem{FZCAII}
Sergey Fomin and Andrei Zelevinsky.
\newblock Cluster algebras. {II}. {F}inite type classification.
\newblock {\em Invent. Math.}, 154(1):63--121, 2003.

\bibitem{FZBCAIII}
Arkady Berenstein, Sergey Fomin, and Andrei Zelevinsky.
\newblock Cluster algebras. {III}. {U}pper bounds and double {B}ruhat cells.
\newblock {\em Duke Math. J.}, 126(1):1--52, 2005.

\bibitem{FZCAIV}
Sergey Fomin and Andrei Zelevinsky.
\newblock Cluster algebras. {IV}. {C}oefficients.
\newblock {\em Compos. Math.}, 143(1):112--164, 2007.

\bibitem{GLS}
C.~Gei\ss, B.~Leclerc, and J.~Schr\"oer.
\newblock Cluster structures on quantum coordinate rings.
\newblock {\em Selecta Math. (N.S.)}, 19(2):337--397, 2013.

\bibitem{GL}
Jan~E. Grabowski and St\'ephane Launois.
\newblock Graded quantum cluster algebras and an application to quantum
  {G}rassmannians.
\newblock {\em Proc. Lond. Math. Soc. (3)}, 109(3):697--732, 2014.

\bibitem{G15}
Jan~E. Grabowski.
\newblock Graded cluster algebras.
\newblock {\em J. Algebraic Combin.}, 42(4):1111--1134, 2015.

\bibitem{BZ}
Arkady Berenstein and Andrei Zelevinsky.
\newblock Quantum cluster algebras.
\newblock {\em Adv. Math.}, 195(2):405--455, 2005.

\bibitem{GSV}
Michael Gekhtman, Michael Shapiro, and Alek Vainshtein.
\newblock {\em Cluster algebras and {P}oisson geometry}, volume 167 of {\em
  Mathematical Surveys and Monographs}.
\newblock American Mathematical Society, Providence, RI, 2010.

\bibitem{FST}
Sergey Fomin, Michael Shapiro, and Dylan Thurston.
\newblock Cluster algebras and triangulated surfaces. {I}. {C}luster complexes.
\newblock {\em Acta Math.}, 201(1):83--146, 2008.

\bibitem{ABBS}
Ibrahim Assem, Martin Blais, Thomas Br\"ustle, and Audrey Samson.
\newblock Mutation classes of skew-symmetric {$3\times 3$}-matrices.
\newblock {\em Comm. Algebra}, 36(4):1209--1220, 2008.

\bibitem{CalderoKeller}
Philippe Caldero and Bernhard Keller.
\newblock From triangulated categories to cluster algebras. {II}.
\newblock {\em Ann. Sci. \'Ecole Norm. Sup. (4)}, 39(6):983--1009, 2006.

\bibitem{War}
Matthias Warkentin.
\newblock {\em Exchange Graphs via Quiver Mutation}.
\newblock PhD thesis, Universit{\"a}t Bonn, 2014.

\bibitem{FSTFinite}
Anna Felikson, Michael Shapiro, and Pavel Tumarkin.
\newblock Skew-symmetric cluster algebras of finite mutation type.
\newblock {\em J. Eur. Math. Soc. (JEMS)}, 14(4):1135--1180, 2012.

\bibitem{Muller}
Greg Muller.
\newblock Skein and cluster algebras of marked surfaces.
\newblock {\em Quantum Topol.}, 7(3):435--503, 2016.

\end{thebibliography}


\begin{thebibliography}{1}

\bibitem{BZ}
Arkady Berenstein and Andrei Zelevinsky.
\newblock Quantum cluster algebras.
\newblock {\em Adv. Math.}, 195(2):405--455, 2005.

\bibitem{GL}
Jan~E. Grabowski and St\'ephane Launois.
\newblock Graded quantum cluster algebras and an application to quantum
  {G}rassmannians.
\newblock {\em Proc. Lond. Math. Soc. (3)}, 109(3):697--732, 2014.

\bibitem{GLS}
C.~Gei{\ss}, B.~Leclerc, and J.~Schr{\"o}er.
\newblock Cluster structures on quantum coordinate rings.
\newblock {\em Selecta Mathematica}, 19(2):337–397, 2012.

\bibitem{GSV}
Michael Gekhtman, Michael Shapiro, and Alek Vainshtein.
\newblock {\em Cluster algebras and Poisson geometry}, volume 167 of {\em
  Mathematical {S}urveys and {M}onographs}.
\newblock American Mathematical Soc., 2010.

\bibitem{Sc}
Joshua~S. Scott.
\newblock Grassmannians and cluster algebras.
\newblock {\em Proceedings of the London Mathematical Society}, 92(2):345--380,
  2006.

\end{thebibliography}


\begin{thebibliography}{1}

\bibitem{GSV}
Michael Gekhtman, Michael Shapiro, and Alek Vainshtein.
\newblock {\em Cluster algebras and {P}oisson geometry}, volume 167 of {\em
  Mathematical Surveys and Monographs}.
\newblock American Mathematical Society, Providence, RI, 2010.

\bibitem{FZCAI}
Sergey Fomin and Andrei Zelevinsky.
\newblock Cluster algebras. {I}. {F}oundations.
\newblock {\em J. Amer. Math. Soc.}, 15(2):497--529, 2002.

\bibitem{FZCAII}
Sergey Fomin and Andrei Zelevinsky.
\newblock Cluster algebras. {II}. {F}inite type classification.
\newblock {\em Invent. Math.}, 154(1):63--121, 2003.

\bibitem{FZBCAIII}
Arkady Berenstein, Sergey Fomin, and Andrei Zelevinsky.
\newblock Cluster algebras. {III}. {U}pper bounds and double {B}ruhat cells.
\newblock {\em Duke Math. J.}, 126(1):1--52, 2005.

\bibitem{FZCAIV}
Sergey Fomin and Andrei Zelevinsky.
\newblock Cluster algebras. {IV}. {C}oefficients.
\newblock {\em Compos. Math.}, 143(1):112--164, 2007.

\bibitem{G15}
Jan~E. Grabowski.
\newblock Graded cluster algebras.
\newblock {\em J. Algebraic Combin.}, 42(4):1111--1134, 2015.

\bibitem{GP}
J.~E. {Grabowski} and M.~{Pressland}.
\newblock {Graded Frobenius cluster categories}.
\newblock {\em ArXiv e-prints}, September 2016.

\bibitem{SchifflerQR}
Ralf Schiffler.
\newblock {\em Quiver representations}.
\newblock CMS Books in Mathematics/Ouvrages de Math\'ematiques de la SMC.
  Springer, Cham, 2014.

\end{thebibliography}


\begin{thebibliography}{1}

\bibitem{FST}
Sergey Fomin, Michael Shapiro, and Dylan Thurston.
\newblock Cluster algebras and triangulated surfaces. {I}. {C}luster complexes.
\newblock {\em Acta Math.}, 201(1):83--146, 2008.

\bibitem{Muller}
Greg Muller.
\newblock Skein and cluster algebras of marked surfaces.
\newblock {\em Quantum Topol.}, 7(3):435--503, 2016.

\bibitem{Har}
John~L. Harer.
\newblock The virtual cohomological dimension of the mapping class group of an
  orientable surface.
\newblock {\em Invent. Math.}, 84(1):157--176, 1986.

\bibitem{Hat}
Allen Hatcher.
\newblock On triangulations of surfaces.
\newblock {\em Topology Appl.}, 40(2):189--194, 1991.

\bibitem{Mos}
Lee Mosher.
\newblock Tiling the projective foliation space of a punctured surface.
\newblock {\em Trans. Amer. Math. Soc.}, 306(1):1--70, 1988.

\bibitem{Che}
L.~O. Chekhov and R.~Ch. Penner.
\newblock Introduction to {T}hurston's quantum theory.
\newblock {\em Uspekhi Mat. Nauk}, 58(6(354)):93--138, 2003.

\end{thebibliography}


\begin{thebibliography}{1}

\bibitem{G15}
Jan~E. Grabowski.
\newblock Graded cluster algebras.
\newblock {\em J. Algebraic Combin.}, 42(4):1111--1134, 2015.

\bibitem{FSTFinite}
Anna Felikson, Michael Shapiro, and Pavel Tumarkin.
\newblock Skew-symmetric cluster algebras of finite mutation type.
\newblock {\em J. Eur. Math. Soc. (JEMS)}, 14(4):1135--1180, 2012.

\bibitem{FST}
Sergey Fomin, Michael Shapiro, and Dylan Thurston.
\newblock Cluster algebras and triangulated surfaces. {I}. {C}luster complexes.
\newblock {\em Acta Math.}, 201(1):83--146, 2008.

\bibitem{Muller}
Greg Muller.
\newblock Skein and cluster algebras of marked surfaces.
\newblock {\em Quantum Topol.}, 7(3):435--503, 2016.

\end{thebibliography}
\bibliographystyle{apa}
\end{document}